\documentclass[a4paper, 12pt]{article} 

\usepackage{amsmath} 
\usepackage{theorem}
\usepackage{amssymb}
\usepackage{amsfonts}
\usepackage{latexsym}
\usepackage{amscd}
\usepackage{diagramb}
\input GrCalc4.sty
\usepackage{graphicx}
\usepackage{enumerate}
\usepackage{cancel}
\usepackage{pxfonts}  
\usepackage{float}

\usepackage[hyperindex,breaklinks]{hyperref}

\usepackage{stmaryrd}

\usepackage{xspace,colortbl}
\usepackage{color}

\usepackage{arydshln}

\DeclareMathAlphabet{\mathpzc}{OT1}{pzc}{m}{it}

    \newcommand{\axiom}[1]{\textbf{\hypertarget{#1}{(#1)}}}
    \newcommand{\axiomref}[1]{\hyperlink{#1}{(#1)}}

\theoremstyle{plain}
\theoremheaderfont{\normalfont\bfseries}

\newtheorem{thm}{Theorem}[section]

\newtheorem{lma}[thm]{Lemma}
\newtheorem{cor}[thm]{Corollary}
\newtheorem{defn}[thm]{Definition}

\theorembodyfont{\rmfamily}

\newtheorem{rem}[thm]{Remark}
\newtheorem{prop}[thm]{Proposition}

\newcommand{\qed}{\hfill\quad\fbox{\rule[0mm]{0,0cm}{0,0mm}}  \par\bigskip}

\definecolor{rojo}{rgb}{1,0,0}

\newcommand{\x}{\mbox{-}}
\newcommand{\s}{\hspace{0,06cm}}
\newcommand{\w}{\hspace{-0,06cm}}

\newcommand{\comp}{\circ}
\newcommand{\iso}{\cong}
\newcommand{\ot}{\otimes}
\newcommand{\crta}{\overline}

\newcommand{\Aa}{{\mathbb A}}
\newcommand{\Bb}{{\mathbb B}}
\newcommand{\Cc}{{\mathbb C}}
\newcommand{\Dd}{{\mathbb D}}
\newcommand{\Ee}{{\mathbb E}}
\newcommand{\LL}{{\mathbb L}}
\newcommand{\Rr}{{\mathbb R}}
\newcommand{\Zz}{{\mathbb Z}}
\newcommand{\Kk}{{\mathbb K}}
\newcommand{\HH}{{\mathcal H}}
\newcommand{\Z}{{\mathcal Z}}
\newcommand{\F}{{\mathcal F}}
\newcommand{\G}{{\mathcal G}}
\newcommand{\Ps}{\operatorname{\mathbb P s}}
\newcommand{\A}{{\mathcal A}}
\newcommand{\B}{{\mathcal B}}
\newcommand{\C}{{\mathcal C}}
\newcommand{\V}{{\mathcal V}}
\newcommand{\Del}{\boxtimes}

\newcommand{\I}{{\mathcal I}}
\newcommand{\M}{{\mathcal M}}

\def\ul{\underline}

\newcommand{\Kl}{\operatorname {Kl}}
\newcommand{\Id}{\operatorname {Id}}
\newcommand{\id}{\operatorname {id}}
\newcommand{\Cat}{\operatorname {Cat}}
\newcommand{\Fun}{\operatorname {Fun}}
\newcommand{\Lax}{\operatorname{\mathbb L ax}}
\newcommand{\Pseudo}{\operatorname{\mathbb P seudo}}
\newcommand{\Binoidal}{\operatorname{Binoidal}}
\newcommand{\noidal}{\operatorname{noidal}}
\newcommand{\Epsilon}{\varepsilon}
\newcommand{\Dist}{\operatorname{Dist}}
\def\Dd{{\mathbb D}}
\newcommand{\rtr}{\triangleright}
\newcommand{\ltr}{\triangleleft}
\newcommand{\PsNat}{\operatorname{PsNat}}

\usepackage{yfonts}
\usepackage{upgreek}

\newcommand{\lelabel}[1]{\label{le:#1}}
\newcommand{\leref}[1]{Lemma~\ref{le:#1}}
\newcommand{\eqlabel}[1]{\label{eq:#1}}
\newcommand{\equref}[1]{(\ref{eq:#1})}
\newcommand{\thlabel}[1]{\label{th:#1}}
\newcommand{\thref}[1]{Theorem~\ref{th:#1}}
\newcommand{\delabel}[1]{\label{de:#1}}
\newcommand{\deref}[1]{Definition~\ref{de:#1}}
\newcommand{\prlabel}[1]{\label{pr:#1}}
\newcommand{\prref}[1]{Proposition~\ref{pr:#1}}
\newcommand{\colabel}[1]{\label{co:#1}}
\newcommand{\coref}[1]{Corollary~\ref{co:#1}}
\newcommand{\rmlabel}[1]{\label{rm:#1}}
\newcommand{\rmref}[1]{Remark~\ref{rm:#1}}
\newcommand{\selabel}[1]{\label{se:#1}}
\newcommand{\seref}[1]{Section~\ref{se:#1}}
\newcommand{\sslabel}[1]{\label{ss:#1}}
\newcommand{\ssref}[1]{Subsection~\ref{ss:#1}}
\newcommand{\ssslabel}[1]{\label{ss:#1}}
\newcommand{\sssref}[1]{Subsection~\ref{ss:#1}}

\newcommand*{\threefrac}[3]{%
  \begin{array}{@{\,}c@{\,}}%
    #1\\
    \hline
    #2\\
    \hline
    #3%
  \end{array}%
}

\newcommand*{\fourfrac}[4]{%
 \begin{array}{@{\,}c@{\,}c@{\,}}%
    #1\\
    \hline
    #2\\
    \hline
    #3\\
    \hline
		#4%
  \end{array}%
	}

\begin{document}

\title{Premonoidal and Kleisli double categories } 

\author{Bojana Femi\'c \vspace{6pt} \\
{\small Mathematical Institite of  \vspace{-2pt}}\\
{\small Serbian Academy of Sciences and Arts } \vspace{-2pt}\\
{\small Kneza Mihaila 36,} \vspace{-2pt}\\
{\small 11 000 Belgrade, Serbia} \vspace{-2pt}\\
{\small femicenelsur@gmail.com}}

\date{}
\maketitle

\begin{abstract}
We give a double categorical version of the recently introduced notion of premonoidal bicategories. 
We introduce a funny product and a funny type of multicategory on double categories granting them 
a closed funny monoidal structure. We investigate relations between 
various funny type of structures 
and premonoidal double categories. 
We prove that a premonoidal double category $\Dd$ is purely central if and only if its binoidal structure is given by a pseudodouble quasi-functor (a multimap for a Gray type of multicategory) if and only if it admits a monoidal structure. 
For such $\Dd$ we introduce pure center and show that the monoidal structure on $\Dd$ extends to it. We also discuss 
one-sided and general centers. 

Exploiting the companion-lifting properties of vertical structures in a double category into their horizontal counterparts, 
we prove a series of further results simplifying proofs for the corresponding bicategorical findings. 
We introduce vertical strengths on vertical double monads and horizontal strengths on horizontal double monads 
and prove that the former induce the latter. 
We show that vertical strengths induce actions of the induced horizontally monoidal double category on the corresponding Kleisli double category of the induced horizontal double monad. We prove that there is a 1-1 correspondence 
between horizontal strengths and extensions of the canonical action of the double category on itself. Finally, we show that 
for a bistrong vertical double monad the corresponding Kleisli double category is premonoidal. 
\end{abstract}



\maketitle

\tableofcontents

\section{Introduction}

The idea for considering premonoidal categories comes from Computer Science. In 1991 Eugenio Moggi  proposed the use of monads 
for modeling {\em notions of computations} in the Kleisli category of the monad. 
In particular, he showed that the notion of a strong monad $(T,t)$ on a Cartesian (but also any monoidal) category $\C$ 
captures the essence of the structure of effectful programs (\cite{Moggi, Moggi1}). 
Concretely, an effectful program can be modelled by a Kleisli arrow $f: A\to T(A')$ in $\C$ and substitution into another
program $g:B\ot A'\to A''$ can be modelled by the strength $t$ and the Kleisli extension operation: $B\ot A\stackrel{1\ot f}{\to}B\ot T(A')
\stackrel{t}{\to} T(B\ot A') \stackrel{T(g)}{\to} T(A'')$. 
By the framework proposed by Moggi a substantial class of semantic models for various effectful languages can be covered \cite{Moggi1}. 
Though the fact that in denoting a program a monad is involved has as a consequence that the program is not interpreted in the category $\C$ itself, but in the Kleisli category $\C_T$. 

In 1996 John Power and Edmund Robinson reformulated the theory by introducing the concept of a {\em premonoidal category} in \cite{PR}. It  
encodes Moggi's model, namely the structure of $\C_T$ that is necessary for the semantic interpretation. In this way 
the structure of effectful languages is directly reflected without any mention of the monad 
and effectful programs are interpreted simply as morphisms in a premonoidal category. 
That premonoidal categories generalize Moggi's approach with monads is illustrated by the fact that in a symmetric monoidal category $\C$ any strength for a monad $T$ induces a premonoidal category structure on the Kleisli category $\C_T$. However, there are premonoidal categories 
that do not arise from monads, \cite{P,S,T}. This justifies the interest of studying the notion for itself.

Premonoidal categories enable to model sequential composition of programs. This way in the analysis of side effects a program from $A$ to $B$ 
can be modeled by a function $[A]\ot S \to [B]\ot S$, where $S$ is a set of states and $[A],[B]$ are types of the variables. 
If another program is modeled by a function $[A']\ot S \to [B']\ot S$, one obtains two functions $[A]\ot [A']\ot S \to [B]\ot [B']\ot S$, 
so that both can model the composite program, depending on the order of execution of the programs. This essentially amounts to noncommutativity of the following diagram 
\begin{equation} \eqlabel{effects}
\scalebox{0.86}{
\bfig
 \putmorphism(-170,400)(1,0)[A\ot B`A\ot B' `A\ot g]{620}1a
 \putmorphism(450,400)(1,0)[\phantom{A\ot B}`A'\ot B' `f\ot B']{680}1a
 \putmorphism(-150,50)(1,0)[A\ot B`A'\ot B `f\ot B]{600}1a
 \putmorphism(450,50)(1,0)[\phantom{A\ot B}`A'\ot B'. `A'\ot g]{680}1a

\putmorphism(-180,400)(0,-1)[\phantom{Y_2}``=]{350}1r
\putmorphism(1100,400)(0,-1)[\phantom{Y_2}``=]{350}1r
\put(380,210){\fbox{$\not=$}}
\efig}
\end{equation}
To encode the fact that those two functions are different, one needs to avoid bifunctoriality of the tensor product, thus one is led to go beyond a framework of a monoidal category (of denotations of types and programs). This is how the notion of a premonoidal category arose. Summing up, expressions in effectful programs cannot generally be re-ordered even though the data flow allows for it. In this sense 
the composition of morphisms in a premonoidal category should be understood as encoding control flow, as argued in \cite{HF1}. 

\medskip

The framework of premonoidal categories proved very useful for studying existing models with effects and constructing new ones, which all 
can be described by categories. In \cite{HF,HF1} a bicategorical setting for premonoidal structures has been proposed under the notion of 
{\em premonoidal bicategories}. In the cited articles various well-known categorical results concerning Kleisli categories, graded monads,   premonoidality, Freyd categories, just to mention some, have been proved to extend to the bicategorical setting. 
This bicategorical framework proves appropriate for models having more structure, which normally supply more detailed or purposeful  information than the usual categorical models. Plenty of such models appeared in recent years, and they are all bicategorical in nature. Their examples include various kinds of game semantics, linear logic based on profunctors, and models describing the $\beta\eta$-rewrites of the 
simply-typed $\lambda$-calculus. For more details and references on such examples see \cite{HF1}. Thus Paquet-Saville's model extends the 
one of Moggi-Power-Robinson both theoretically and in practice, in the latter case providing a practical tool towards new semantic models 
for programming languages. 

\medskip

When the findings of \cite{HF} were presented at the annual Category Theory 2023 conference the author of the present article has recognized in it the main feature of ``premonoidality'', which is so to say the noncommutativity of \equref{effects}, from Gray's construction of what today is called Gray tensor product (for the reference see \cite{Gray}). More precisely, this noncommutativity appears as a part of the data defining a {\em quasi-functor of two  variables} from \cite{Gray} in the context of 2-categories, which underlies the construction of the inner-hom object and closedness for the closed Gray monoidal structure (for the category of 2-categories). Such quasi-functors have also been studied recently in the context of double categories, see \cite{Fem:Fil, Fem:Bif}. 
(As a matter of fact, they present {\em multimaps} of particular {\em multicategories} that in \cite{Fem:Gray1} we call {\em Gray multicategories}. Nevertheless, 
since, on one hand, we really want to exploit the structure and properties of individual quasi-functors, and 
on the other hand, we introduce also other kinds of multimaps and multicategories in this article, we will stick to the term {\em quasi-functor} throughout.) A deeper study of premonoidal bicategories impelled various questions of curiosity which ended up with the research whose results the reader finds in front him/herself. 

Closely related to the notion of a premonoidal bicategory $\B$ is the notion of {\em central cells} in $\B$, giving rise to a further notion of the {\em center bicategory} $\Z(\B)$. A natural question, extending an existing result in ordinary categories, of whether the center bicategory is monoidal, resulted in a negative answer in \cite{HF} and the necessity to introduce the so-called {\em bicategory of pure maps} $\C_p(\B)$, which is then proved to be monoidal. 
(The findings of the preprint \cite{HF} were in the meantime split into published papers \cite{HF1} and \cite{HF2}, except from the study of the bicategory of pure maps that remained only in \cite{HF}.) 
On the other hand, the authors introduced {\em strengths} for pseudomonads on monoidal bicategories and proved in particular three important results, generalizing famous analogous results from categories. Firstly, the fact that strengths on pseudomonads $T$ on monoidal bicategories $\B$ induce actions of $\B$ on the Kleisli bicategory $\B_T$, secondly, that there is a 1-1 
correspondence between strengths and extensions of the canonical action of $\B$ on itself,  
and thirdly, that for a bistrong monad $T$ the Kleisli bicategory $\B_T$ is premonoidal. In the presented notions and henceforth results merely the formulation of definitions and then much more the carrying out of the proofs involves a substantial set of data, axioms and technicalities of computations, which make the creating process and then also the 
follow-up reading a demanding mission. This is so due to the nature of the subject of monoidal bicategories. Here is where we found the challenge: because of the importance and beauty of the results, it is almost so that they deserve an easier approach, and on the other hand, luckily, precisely an easier approach is provided by working from the point of view of {\em double categories}, which are a treasure by themselves. 

The latter notion is gaining increasingly more interest in the last decade. Moreover, in Shulman's wonderful paper \cite{Shul} it is shown how it is enough to check much less conditions on a double category in order to draw 
desired conclusions that would hold in the underlying horizontal bicategory of the double category in play. Now, it is true, though, that Shulman's \cite[Theorem 5.1]{Shul} is formulated in terms of monoidality of the double category and its underlying bicategory, but in his paper the bases are set for drawing much more results following the same philosophy. We extract the essence for this, the essence of Shulman's machinery behind the proof of his Theorem 5.1, in our \prref{essence} and the following-up \prref{omega*}. Let us spell out what happens. 

For a less instructed reader let us first say that double categories have both horizontal and vertical 1-cells, and henceforth the 2-cells (often called {\em squares}) have both horizontal and vertical composition. 
For this reason many (bi)categorical notions have both their horizontal and vertical double categorical mates. 
The baseline is that one can consider horizontal and vertical natural transformations, horizontal and vertical modifications between the latter two, respectively, but also, so to say, ``mixed'' modifications acting between two horizontal and two vertical transformations. 

In short, the ``hidden treasure'' behind \cite[Theorem 5.1]{Shul} that we make explicit in this paper is the following. If one works with ``vertical structures'' in a double category $\Dd$, that is, with (invertible) vertical strict transformations {\em such that their vertical 1-cell components have companions}, and with identity vertical modifications acting between them, then 
(equivalence) horizontal transformations and horizontal modifications between them are induced, so that the latter satisfy any sensible axiom that can be written down using their compositions. Taking the underlying bicategorical structure from the latter, one obtains, in particular, bicategorical (equivalence) transformations and modifications 
between them that obey any desired axiom. 

\medskip

Armed with this potent tool, we engaged in lifting to the setting of double categories the above-spelled notions and results from \cite{HF2}-\cite{HF}. 
By doing so, not only do we provide simpler proofs for (some of) the corresponding bicategorical results. Rather, by proving double categorical 
versions of results from \cite{HF2}-\cite{HF} we open the way to allow for new applications, and we broaden the mathematical framework so to encompass an even larger class of semantic models for effectful languages than those that have been recognized so far. The particular in double categories with respect to bicategories is their internal category nature, which provides a whole new dimension for additional arrows (the vertical, 
{\em i.e.} strict direction). This has proved at numerous occasions in the literature to enable a more complete, even self-sufficient environment for expressing mathematical concepts.  

The contents of the present article can be divided in two parts. The first one is comprised of the first nine sections, concerning the introduction of the notion of premonoidal and center double categories, double funny functors 
and the corresponding multicategories, as well as pseudodouble quasi-functors (of three and more variables), and their interrelations and properties. 
The second part is concentrated in the last, tenth section and deals with our results on Kleisli double categories. 
After a preliminary part on double categories, we start the first part by defining binoidal structures and central cells. 
We study the relationship between centrality of 1- and 2-cells and companions. 
We define a premonoidal double category so that the structural transformations are taken to be vertical, in order to be able to apply Shulman's recipe to draw the consequences that will encompass the bicategorical results from \cite{HF1}. 
In the premonoidal setting one has three associativity constraints $\alpha$ - one for each of the three variables. 
In \seref{24 axioms} we study 24 axioms which permit to describe the interplay between those three associativity constraints. We study also the consequences of having these 24 axioms. 

In \seref{funny} we introduce two types of double funny products between double categories. To do so, we introduce two types of 
{\em double funny functors}, which will give rise to two distinct versions of multicategories that we call {\em funny multicategories}. They will be closed and representable, so that we obtain two kinds of funny type of monoidal structures on the category $Dbl$ of double categories: a {\em purely funny} and a {\em mixed funny} monoidal structure. They generalize the closed monoidal category of 2-categories from \cite[Section 2]{BG}. A monoid in the purely funny monoidal category structure on $Dbl$ is a {\em strict premonoidal double category}. On the other  hand, we extend the mixed funny 
monoidal structure on $Dbl$ to a funny monoidal 2-category structure on a 2-category $Dbl_2$ of double categories, we construct a pseudomonid in it and show that it is a premonoidal double category $\Dd$ in which the binoidal structure is strict (meaning that it is given by two strict double functors $A\ltimes-, -\rtimes B: \Dd\to\Dd$). 

In \seref{cubical1} we explore the anticipated relationship between premonoidality and quasi-functors. Here we rely on our Proposition 3.3 and Corollary 3.5 from \cite{Fem:Bif} and conclude in the present \thref{lr central} that a double category 
$\Dd$ has a pseudodouble quasi-functor $H:\Dd\times\Dd\to\Dd$ if and only if $\Dd$ is binoidal and purely central.  
For such a double category we introduce a {\em pure center pseudodouble category}, which is then determined by such quasi-functor and 
entails a functorial choice of centrality structures. We also explore one-sided centers in relation to pseudodouble 
quasi-functors and a (general) center double category induced by two partially related pseudodouble quasi-functors, 
see \ssref{not diff}. At this point we indeed study the mentioned center double categories for {\em binoidal} double categories and defer their study for {\em premonoidal} double categories to \seref{center premon}.

We dedicate \seref{cubical2} to prove that premonoidal double categories whose binoidal structure stems from a pseudodouble quasi-functor allow for a monoidal structure. 
From our \cite[Proposition 5.6]{Fem:Bif} applied to the present setting we know that 
there is a double category equivalence $\F\colon q\x\Ps_{hop}^{st}(\Dd\times\Dd,\Dd) \to \Ps_{hop}(\Dd\times\Dd,\Dd)$ 
between the double categories of pseudodouble quasi-functors (for which certain 2-cells $(u, U)$ are trivial) and the 
double category of pseudodouble functors. We then introduce {\em pseudodouble quasi-functors with three variables} 
and vertical strict transformations between them, 
generalizing to double categories the corresponding 2-categorical gadget from \cite[Section 7]{GPS}. 
We make use of these notions to obtain a three (and four) variable analogue of our mentioned double equivalence result. 
Indeed, we use a category equivalence of the underlying vertical categories, \thref{H equiv ot}. The latter, more importantly, helps us prove 
\thref{premon-mon}, stating the following. A double category $\Dd$ equipped with a quasi-functor $H$ is purely central premonoidal with a binoidal structure given via $H$ if and only if there is a monoidal double category structure $(\Dd, \ot)$ 
(in the sense of \cite[Definition 2.9]{Shul}), where $H$ and $\ot$ are in 1-1 correspondence via the above double equivalence functor $\F$. We draw a conclusion about an analogous result for 2-categories, 
on one hand, and a consequence on the underlying horizontal bicategory $\HH(\Dd)$, on the other hand. 
We also prove that in a premonoidal double category whose binoidal structure comes from a pseudodouble quasi-functor the 
24 axioms studied in \seref{24 axioms} are satisfied. Along the way, in analogy to pseudodouble quasi-functors of $n$-variables 
we also introduce $n$-noidal purely and one-sided central structures as $n$-variable extensions of the respective central  binoidal structures, and their corresponding categories. 
In \thref{gen-1} we extend \thref{premon-mon} by proving that there is a premonoidal and purely central structure on $\Dd$ 
(with different sources of purely central binoidal structures) if and only if there is a monoidal structure on $\Dd$. 


In \seref{center premon} we study centers of premonoidal double categories: 
a new requirement with respect to centers of binoidal double categories is that 1v-cell components of premonoidal constraints should live in the centers. Accordingly, we can only consider pure center $\Zz_p(\Dd)$ and (general) center double category 
$\Zz(\Dd)$. Moreover, we also require interrelatedness of the three associativity constraints via the 24 axioms. It turns out 
that they are automatically fulfilled in purely central double categories, {\em i.e.} in $\Zz_p(\Dd)$. We prove that there are pseudodouble functors $A\ltimes-, -\rtimes B:\Zz_\bullet(\Dd)\to\Zz_\bullet(\Dd)$ for $\Zz_\bullet(\Dd)$ being $\Zz_p(\Dd)$ or $\Zz(\Dd)$. In \seref{mon pure cen} we prove that monoidal structure on a purely central premonoidal double category $\Dd$ extends to a monoidal structure on its pure center $\Zz_p^{st}(\Dd)$ (a pure center in which the structure 2-cells $(u, U)$ of the underlying pseudodouble quasi-functor are trivial). Under the additional assumption that the associativity and unity constraints of $\Dd$ are companion-liftable (meaning that their 1v-cell components have companions), the underlying horizontal bicategory $\HH(\Dd)$ is premonoidal and admits a monoidal structure. Also, there is a monoidal bicategory 
$\HH(\Zz_p^{st}(\Dd)_{hm})$ granting the bicategory of pure maps $C_p(\HH(\Dd))$ of $\HH(\Dd)$ a monoidal bicategory structure 
via a pseudofunctor $\HH(\Zz_p^{st}(\Dd)_{hm})\to C_p(\HH(\Dd))$, \thref{gen-H(B)}, and it has a natural structure of a Freyd bicategory, \coref{Freyd}. For a premonoidal 2-category $\B$ that comes from a premonoidal double category $\Dd$ satisfying the above conditions, the fact proved in \cite[Theorem 3]{HF} that the 2-category of pure maps $\C_p(\B)$ is monoidal can be thus seen as a consequence of our theorem.

What follows is the second part concentrated in \seref{Kleisli} in which we study Kleisli double categories. 
We start by recalling the notions of horizontal and vertical double monads on a double category $\Dd$, and of a horizontal Kleisli double category of a horizontal double monad from \cite{GGV}. An important condition on vertical transformations used in \cite{GGV} 
is that they be ``special'', 
a condition that is fulfilled if they are companion-liftable and invertible as vertical transformations (in the sense of 
\deref{vlt}). We also recall the result from \cite{GGV} saying that a vertical double monad $T$ whose structure vertical transformations are companion-liftable and invertible induces a horizontal double monad $\hat T$. (This is another example of the lifting of a vertical to a horizontal structure in the style of \cite{Shul}.) We then introduce in \ssref{strengths} the notion of a {\em horizontally monoidal} double category, 
so that a (vertically, 
{\em i.e.} in the sense of \cite[Definition 2.9]{Shul}) monoidal double category $(\Dd, \ot, I, \alpha, \lambda, \rho)$ for which $\alpha, \lambda, \rho$ are companion-liftable, induces a horizontally monoidal double category $(\Dd, \ot, I, \hat\alpha, \hat\lambda, \hat\rho)$. 
Furthermore, 
we introduce, on one hand, the notion of a vertical strength on a vertical double monad on a monoidal double category, and, on the other, a horizontal strength on a horizontal double monad on a horizontally monoidal double category. We prove, assuming that the vertical monad structures and monoidal structures of $\Dd$ are companion-liftable (and invertible), that a vertical strength $t$ on a vertical monad $T$ on $(\Dd, \ot, I, \alpha, \lambda, \rho)$ induces a horizontal strength $\hat t$ on a horizontal monad $\hat T$ on 
$(\Dd, \ot, I, \hat\alpha, \hat\lambda, \hat\rho)$, so that  
$\HH(\hat t)$ is a strength on the pseudomonad $\HH(\hat T)$ on the monoidal bicategory $\HH(\Dd)$. This is our \thref{strenghts thm}. 

\ssref{actions} discusses how double categorical strengths on monads induce actions on the horizontal Kleisli double category. Again, we differentiate vertical and horizontal action of a monoidal, respectively horizontally monoidal, double category on another double category. In \thref{strength-action} we show that a vertical strength $t$ on a vertical double monad $T$ on a monoidal double category $(\Dd, \ot, I, \alpha, \lambda, \rho)$, under assumptions that the involved structural vertical transformations are companion-liftable (and invertible), induces a horizontal action of 
$(\Dd, \ot, I, \hat\alpha, \hat\lambda, \hat\rho)$ on $\Kk l(\hat T)$ (passing through a horizontal strength $\hat t$ on the horizontal monad $\hat T$ on $(\Dd, \ot, I, \hat\alpha, \hat\lambda, \hat\rho)$). The corresponding claim regarding $\hat t$ 
for the underlying horizontal bicategory holds, too, so that we obtain a double categorical analogue of the bicategorical claim from \cite[Proposition 7.1]{HF2}. 

In \ssref{extensions} we study (horizontal) extensions of the canonical action of $\Dd$ on itself, which are given by an action of a (horizontally) monoidal $\Dd$ on the Kleisli double category of a horizontal double monad and a horizontal icon so that a couple of axioms hold. The above result that a vertical strength $t$ on $T$ induces a horizontal action of $\Dd$ on 
$\Kk l(\hat T)$ we extend to an extension determined by that action in \thref{strength-ext}. For a horizontal double monad $S$ we prove that there is a 1-1 correspondence between horizontal strengths $s$ on $S$ and extensions of canonical actions, 
\thref{iff}. This is a double categorical version of \cite[Theorem 1]{HF}, \cite[Theorem 7.2]{HF2}.

Lastly, we introduce bistrengths, both on a vertical and on a horizontal double monad on a corresponding monoidal double category. In 
\thref{Kl-prem} we prove, under assumptions that suitable vertical transformations are companion-liftable and invertible, that for a bistrong vertical double monad $T$ on a monoidal double category 
the horizontal Kleisli double category $\Kk l(\hat T)$ is premonoidal. We draw the natural consequences for the underlying horizontal bicategory in \coref{Kl-prem}, obtaining a relation to the analogous result from \cite[Theorem 2]{HF}, 
\cite[Proposition 5.10]{HF2}. What is remarkable in this double categorical proof, is the advantage of being able to work with identity vertical modifications 
and by the properties of companion-lifting obtain non-trivial mixed and horizontal modifications for free, which almost automatically satisfy 
the four pentagons and six triangles that one otherwise would be obliged to check in 
the context of a 
premonoidal bicategory. Instead, provided the existence of suitable companions, one only needs to verify that the 2-cells appearing in the computation are of certain type. 

\medskip

The overview of the article from Section 3 onwards is presented in the last nine paragraphs. 
The next section recollects the definitions of the ingredients of the double category of double functors, definitions and some 
properties of bicategorical adjoints and mates, as well as some basics on companions and conjoints in double categories, 
including the lifting of invertible vertical transformations to horizontal natural transformations, and the lifting of vertical identity modifications and of mixed {\em i.e.} non-globular modifications to the horizontal ones (see \prref{essence} and 
\prref{omega*}). In the Appendix we also include the following definitions related to double categories from \cite{Fem:Bif}: 
lax double quasi-functor (of two variables), a table interpreting its axioms, as well as a vertical lax transformation between two such quasi-functors.

\section{Prerequisites on double categories}

A double category is a category internal to the category of categories. A weaker notion is a pseudodouble category, which is a pseudocategory internal to the 2-category of categories. In \cite[Section~7.5]{GP:L} it is proved that any pseudodouble category is double equivalent to a double category. To simplify the writing and computation we will deal with double categories. 
For a double category $\Dd$ we denote by $\Dd_0$ the category of objects, and by $\Dd_1$ the category of morphisms. 
For horizontal 1-cells we will say shortly 1h-cells, for vertical 1-cells we will say 1v-cells, and squares we will just call 2-cells. 
In our convention the horizontal direction is weak and the vertical direction is strict. 
The underlying horizontal 2-category of $\Dd$ will be denoted by $\HH(\Dd)$. It consist of objects, 1h-cells and those 2-cells 
whose 1v-cells are identities. 
Composition of 1h-cells $A\stackrel{f}{\to}A'\stackrel{f'}{\to}A''$ we will write as $f'f$ or $[f\,\, \vert\,\, f']$ and vertical composition 
of 1v-cells $A\stackrel{u}{\to}\tilde A\stackrel{v}{\to}\tilde{\tilde A}$ as a fraction $\frac{u}{v}$. Similarly, to simplify the writing or avoid the use of too many pasting diagrams, we will denote horizontal composition of 2-cells $\alpha$ and $\beta$, 
where $\alpha$ acts first, by $[\alpha\vert\beta]$, and vertical composition of 2-cells $\alpha$ and $\gamma$, where $\alpha$ acts first, by $\frac{\alpha}{\gamma}$. In the images below we illustrate the mnemotechnical value of such notation 
$$[\alpha\vert \beta]=
\scalebox{0.8}{
\bfig
\putmorphism(0,200)(1,0)[A`B`f]{450}1a
\putmorphism(460,200)(1,0)[`C`f']{450}1a
\putmorphism(0,220)(0,-1)[\phantom{Y_2}``u]{380}1l
\putmorphism(450,220)(0,-1)[\phantom{Y_2}``u']{380}1r
\putmorphism(890,220)(0,-1)[\phantom{Y_2}``u'']{380}1r
\put(150,0){\fbox{$\alpha$}}
\put(600,0){\fbox{$\beta$}}
\putmorphism(0,-150)(1,0)[\tilde A`\tilde B`g]{450}1b
\putmorphism(460,-150)(1,0)[`\tilde C, `g']{450}1b
\efig}
\qquad\quad
\frac{\alpha}{\gamma}=
\scalebox{0.8}{
\bfig
\putmorphism(0,350)(1,0)[A`B`f]{450}1a
\putmorphism(0,370)(0,-1)[\phantom{Y_2}``u]{380}1l
\putmorphism(450,370)(0,-1)[\phantom{Y_2}``u']{380}1r
\put(150,150){\fbox{$\alpha$}}
\put(150,-210){\fbox{$\gamma$}}
\putmorphism(0,10)(0,-1)[\phantom{Y_2}``v]{380}1l
\putmorphism(450,10)(0,-1)[\phantom{Y_2}``v']{380}1r
\putmorphism(0,0)(1,0)[\tilde A`\tilde B`g]{450}1b
\putmorphism(0,-360)(1,0)[\tilde{\tilde A}`\tilde{\tilde B}.`h]{450}1b
\efig}
$$
For more on double categories we recommend \cite{GP:L, MG}.

\subsection{Some computational tools for globular 2-cells} \sslabel{strings}

For 2-cells of the form of $a$ as below we will say that they are {\em horizontally globular}, while for those of the form of $b$ 
we will say that they are {\em vertically globular}. 
$$ 
\scalebox{0.9}{
\bfig
\putmorphism(-150,170)(1,0)[A`B`f]{460}1a
\putmorphism(-150,-180)(1,0)[A`B`g]{460}1a
\putmorphism(-150,170)(0,-1)[\phantom{Y_2}``=]{350}1l
\putmorphism(310,170)(0,-1)[\phantom{Y_2}``=]{350}1r
\put(20,10){\fbox{$a$}}
\efig}
\qquad\qquad
\scalebox{0.9}{
\bfig
\putmorphism(-150,170)(1,0)[A`A`=]{460}1a
\putmorphism(-150,-180)(1,0)[\tilde A`\tilde B`=]{460}1a
\putmorphism(-150,170)(0,-1)[\phantom{Y_2}``u]{350}1l
\putmorphism(310,170)(0,-1)[\phantom{Y_2}``v]{350}1r
\put(20,0){\fbox{$b$}}
\efig}
$$ 
For computations that involve only horizontally or only vertically globular 2-cells one may use the string diagrams for 2-categories,  
which are simpler to write and easier to handle in the proofs. This way the string diagrams are to be understood as living in the 
underlying horizontal (or vertical) 2-category of the double category. In this article our string diagrams are read from top to bottom and 
from left to right. We present briefly the notation that will be used in this work. 

Let $\F$ be a lax functor and $\G$ an oplax functor between 2-categories. We depict their lax, respectively oplax structures by string diagrams in the following way:
$$
\gbeg{3}{3}
\got{1}{\F(f)} \got{3}{\F(g)} \gnl
\gwmu{3} \gnl
\gob{3}{\F(fg)}
\gend \qquad 
\gbeg{3}{3}
\got{1}{id_{\F(A)}} \gnl
\gu{1} \gnl
\gob{1}{\F(id_A)}
\gend  \qquad 
\gbeg{3}{3}
\got{3}{\G(fg)} \gnl
\gwcm{3} \gnl
\gob{1}{\G(f)} \gob{3}{\G(g)} \gnl
\gend  \qquad 
\gbeg{3}{3}
\got{1}{\G(id_A)} \gnl
\gcu{1} \gnl
\gob{1}{id_{\F(A)}}
\gend  
$$
where $f,g$ are composable 1-cells and $A$ a 0-cell in the domain 2-category. 

Recall that for a 1-cell $f \colon A\to B$ in a 2-category its \emph{left adjoint} is a 1-cell $u \colon B \to A$ together with two 
2-cells $\eta: \id_A\to u f$ and $\varepsilon: f u \to \id_B$ such that
\begin{equation*}
  \frac{[\eta\vert\Id_f]}{[\Id_f\vert\varepsilon]} = \id_f
  \qquad \qquad \text{and} \qquad \qquad
   \frac{[\Id_u\vert\eta]}{[\varepsilon\vert\Id_u]} = \id_u.
\end{equation*}

In string diagrams we will write $\eta=
\gbeg{2}{1}
\gdb \gnl
\gend$ and $\Epsilon=
\gbeg{2}{1}
\gev \gnl
\gend$, and they satisfy the laws:
$$\gbeg{3}{4}
\got{1}{} \got{3}{f} \gnl
\gdb  \gcl{1} \gnl
\gcl{1} \gev \gnl
\gob{1}{f}
\gend=\Id_f
\qquad\textnormal{and}\qquad 
\gbeg{3}{4}
\got{1}{u} \got{1}{} \gnl
\gcl{1} \gdb \gnl
\gev \gcl{1} \gnl
\gob{1}{} \gob{4}{\hspace{-0,1cm}u}
\gend=\Id_u.
$$

\medskip

It is a folklore that any equivalence between 1-cells in a 2-category can be made into an adjoint equivalence 
(for the proof see 
Theorem 3.3 of the nLab entry \\ https://ncatlab.org/nlab/show/adjoint+equivalence\#properties, 
understood as a proof in a 2-category).

\medskip

We record now two identities regarding mates that will be used in \prref{uses mates}. 
Suppose there are equivalence 1-cells $f,f'$ in a 2-category with their respective adjoint quasi-inverses $f^\bullet$ and $(f')^\bullet$.  
Let further 1-cells $g,g'$ be given so that the compositions $g'f$ and $f'g$ are possible and suppose that there is an  
invertible 2-cell $\alpha: f'g\Rightarrow g'f$ (with inverse $\alpha^{-1}: g'f\Rightarrow f'g$). 
The mates $\alpha^\bullet$ and $(\alpha^{-1})^\bullet$ of $\alpha$ and its inverse are given by 
$$\alpha^\bullet=
\gbeg{3}{5}
\got{1}{f^\bullet} \got{1}{g} \gnl
\gcl{1} \gcl{1} \gdb \gnl
\gcl{1} \glmptb \gnot{\hspace{-0,34cm}\alpha} \grmptb \gcl{2} \gnl
\gev \gcl{1} \gnl
\gob{2}{} \gob{1}{g'} \gob{1}{(f')^\bullet}
\gend
\qquad\qquad
(\alpha^{-1})^\bullet=
\gbeg{3}{5}
\gvac{2} \got{1}{g'} \got{1}{(f')^\bullet} \gnl
\gdb \gcl{1} \gcl{2} \gnl
\gcl{1} \glmptb \gnot{\hspace{-0,34cm}\alpha^{-1}} \grmptb \gnl
\gcl{1} \gcl{1} \gev \gnl
\gob{1}{f^\bullet} \gob{1}{g} 
\gend
$$
where left-most $\gbeg{2}{1}
\gdb \gnl
\gend$ stands for $\eta$ for $f'$ and the other one stands for $\Epsilon^{-1}$ for $f$, and the first appearing 
$\gbeg{2}{1}
\gev \gnl
\gend$ stands for $\Epsilon$ for $f$ and the other one for $\eta^{-1}$ for $f'$. Then given $\alpha_1=\alpha:f'g\Rightarrow g'f$ 
and $\alpha_2:fh\Rightarrow h'f'$ (with $\alpha_2^\bullet: h(f')^\bullet\Rightarrow f^\bullet h'$) one has identities
$$ 
\gbeg{4}{5}
\got{1}{g} \gvac{2} \got{1}{h} \gnl
\gcl{1} \gdb \gcl{1} \gnl
\glmptb \gnot{\hspace{-0,34cm}\alpha_1} \grmptb \glmptb \gnot{\hspace{-0,34cm}\alpha_2^\bullet} \grmptb \gnl
\gcl{1} \gcl{1} \gcl{1} \gcl{1} \gnl
\gob{1}{f} \gob{1}{g'} \gob{1}{h'} \gob{1}{f^\bullet}
\gend=
\gbeg{3}{5}
\got{1}{g} \got{1}{h} \gnl
\gcl{1} \gcl{1} \gdb \gnl
\gcl{1} \glmptb \gnot{\hspace{-0,34cm}\alpha_2} \grmptb \gcl{2} \gnl
\glmptb \gnot{\hspace{-0,34cm}\alpha_1} \grmptb \gcl{1} \gnl
\gob{1}{f} \gob{1}{g'} \gob{1}{h'} \gob{1}{f^\bullet}
\gend
\qquad\text{and}\qquad
\gbeg{4}{5}
\got{1}{f'} \got{1}{h'} \got{1}{g'} \got{1}{(f')^\bullet} \gnl 
\gcl{1} \gcl{1} \gcl{1} \gcl{1} \gnl
\glmptb \gnot{\hspace{-0,34cm}\alpha_2^{-1}} \grmptb \glmptb \gnot{\hspace{-0,34cm}(\alpha_1^{-1})^\bullet} \grmptb \gnl
\gcl{1} \gev \gcl{1} \gnl
\gob{1}{h} \gob{2}{} \gob{1}{g}
\gend
=
\gbeg{3}{5}
\got{1}{f'} \got{1}{h'} \got{1}{g'} \got{1}{(f')^\bullet} \gnl 
\glmptb \gnot{\hspace{-0,34cm}\alpha_2^{-1}} \grmptb \gcl{1} \gcl{2} \gnl
\gcl{1} \glmptb \gnot{\hspace{-0,34cm}\alpha_1^{-1}} \grmptb \gnl
\gcl{1} \gcl{1} \gev \gnl
\gob{1}{h} \gob{1}{g}
\gend
$$
and one also has 
$$
\gbeg{3}{6}
\gvac{2} \got{1}{g} \gnl
\gdb \gcl{1} \gnl
\gcl{1} \glmptb \gnot{\hspace{-0,34cm}\alpha^\bullet} \grmptb \gnl
\glmptb \gnot{\hspace{-0,34cm}\alpha^{-1}} \grmptb \gcl{1} \gnl
\gcl{1} \gev \gnl
\gob{1}{g} 
\gend
=
\Id_g
\qquad\qquad\text{and}\qquad\qquad
\gbeg{3}{6}
\gvac{2} \got{1}{g'} \gnl
\gdb \gcl{1} \gnl
\gcl{1} \glmptb \gnot{\hspace{-0,34cm}\alpha^{-1}} \grmptb \gnl
\glmptb \gnot{\hspace{-0,34cm}\alpha^\bullet} \grmptb \gcl{1} \gnl
\gcl{1} \gev \gnl
\gob{1}{g'} 
\gend
=
\Id_{g'}
$$
where in the left-hand side identity 
$\gbeg{2}{1}
\gdb \gnl
\gend$ 
denotes $\eta$ for $f$ and 
$
\gbeg{2}{1}
\gev \gnl
\gend$ denotes $\eta^{-1}$ for $f'$, and in the right-hand side identity 
$\gbeg{2}{1}
\gdb \gnl
\gend$ 
denotes $\Epsilon^{-1}$ for $f$ and 
$
\gbeg{2}{1}
\gev \gnl
\gend$ denotes $\Epsilon$ for $f'$. 
(The 2-cell $\alpha$ here will be used in \prref{uses mates} as a 2-cell component of a horizontal pseudonatural transformation 
$\alpha_{A,X,-}^3$, the 1-cells $f,f'$ as 1h-cells $\alpha_{A,X,B}^3$ and $\alpha_{A,X,B'}^3$, and the 1-cells $g,g'$ as 1h-cells 
$AB\ltimes h$ and $A(X\ltimes h)$, respectively, and similarly.)

\subsection{The double category of double functors}

We recall here the definitions of lax/pseudo double functors, horizontal and vertical (op)lax/pseudo transformations and modifications.

\begin{defn}
A {\em (horizontally) lax double functor} $F\colon\Cc\to\Dd$ between double categories is given by: 1) the data: images on objects, 1h-, 1v- and 2-cells of $\Cc$, 
globular 2-cells: 
{\em compositor} $F_{gf}\colon F(g)F(f)\Rightarrow F(gf)$ and {\em unitor} $F_A\colon 1_{F(A)}\Rightarrow F(1_A)$ in $\Dd$, 
and 2) rules (in $\Dd$): 
\begin{itemize}
\item ($\text{functoriality in vertical morphisms}$)
$$\text{{\em \axiom{lx.f.v1}}}  \quad\qquad  \frac{F(u)}{F(u')}=F(\frac{u}{u'}), \qquad\qquad 
\text{{\em \axiom{lx.f.v2}}} \quad\qquad F(1^A)=1^{F(A)};$$
\item $\text{(functoriality in squares)}$
$$\text{{\em \axiom{lx.f.s1}}}  \quad\qquad  F(\frac{\omega}{\zeta})=\frac{F(\omega)}{F(\zeta)}, \qquad\qquad 
\text{{\em \axiom{lx.f.s2}}} \quad\qquad F(Id_f)=Id_{F(f)};$$
\item $\text{(coherence with compositors and unitors)}$
$$ 
\hspace{-4cm} \text{{\em \axiom{lx.f.cmp}}}  \hspace{2cm}\qquad \frac{[F_{gf}\vert Id_{F(h)}]}{F_{h,gf}}=\frac{[Id_{F(f)}\vert F_{hg}]}{F_{hg,f}} $$
$$ \hspace{-3cm} \text{{\em \axiom{lx.f.u}}} \hspace{2cm}\qquad \frac{[F_A\vert Id_{F(f)}]}{F_{f1_A}}=\Id_{F(f)}=\frac{[Id_{F(f)}\vert F_B]}{F_{1_Bf}};$$
\item $\text{(naturality of the compositor)}$

$$ 
\text{{\em \axiom{lx.f.c-nat}}}  \qquad
\scalebox{0.82}{  
\bfig
\putmorphism(-150,500)(1,0)[F(A)`F(B)`F(f)]{600}1a
 \putmorphism(450,500)(1,0)[\phantom{F(A)}`F(C) `F(g)]{620}1a

 \putmorphism(-150,50)(1,0)[F(A')`F(B')`F(f')]{600}1a
 \putmorphism(450,50)(1,0)[\phantom{F(A)}`F(C') `F(g')]{620}1a

\putmorphism(-180,500)(0,-1)[\phantom{Y_2}``u]{450}1l
\putmorphism(450,500)(0,-1)[\phantom{Y_2}``]{450}1r
\putmorphism(300,500)(0,-1)[\phantom{Y_2}``v]{450}0r
\putmorphism(1080,500)(0,-1)[\phantom{Y_2}``w]{450}1r
\put(-10,280){\fbox{$F(\alpha)$}}
\put(640,280){\fbox{$F(\beta)$}}

\putmorphism(-150,-400)(1,0)[F(A')`F(C') `F(g'f')]{1200}1a

\putmorphism(-180,50)(0,-1)[\phantom{Y_2}``=]{450}1l
\putmorphism(1080,50)(0,-1)[\phantom{Y_3}``=]{450}1r
\put(320,-180){\fbox{$F_{g'f'}$}}

\efig}
\quad=\quad
\scalebox{0.82}{
\bfig
\putmorphism(-170,500)(1,0)[F(A)`F(B)`F(f)]{600}1a
 \putmorphism(430,500)(1,0)[\phantom{F(A)}`F(C) `F(g)]{620}1a
 \putmorphism(-150,50)(1,0)[F(A)`F(C)`F(gf)]{1200}1a

\putmorphism(-180,500)(0,-1)[\phantom{Y_2}``=]{450}1r
\putmorphism(1060,500)(0,-1)[\phantom{Y_2}``=]{450}1l
\put(350,260){\fbox{$F_{gf}$}}

\putmorphism(-150,-400)(1,0)[F(A')`F(C') `F(g'f')]{1200}1a

\putmorphism(-180,50)(0,-1)[\phantom{Y_2}``u]{450}1r
\putmorphism(1060,50)(0,-1)[\phantom{Y_3}``w]{450}1l
\put(300,-180){\fbox{$F(\beta\alpha)$}} 
\efig};$$

\item $\text{(naturality of the unitor)}$

$$\text{{\em \axiom{lx.f.u-nat}}} \qquad\quad 
\scalebox{0.86}{
\bfig
\putmorphism(-250,500)(1,0)[F(A)`F(A)` =]{550}1a
 \putmorphism(-250,50)(1,0)[F(A')`F(A')` =]{550}1a
 \putmorphism(-250,-400)(1,0)[F(A')`F(A')` F(1_{A'})]{550}1a

\putmorphism(-280,500)(0,-1)[\phantom{Y_2}``F(u)]{450}1l
 \putmorphism(-280,70)(0,-1)[\phantom{F(A)}` `=]{450}1l

\putmorphism(300,500)(0,-1)[\phantom{Y_2}``F(u)]{450}1r
\putmorphism(300,70)(0,-1)[\phantom{Y_2}``=]{450}1r
\put(-140,260){\fbox{$Id^{F(u)}$}}
\put(-100,-150){\fbox{$F_{A'}$}}
\efig}
=
\scalebox{0.86}{
\bfig
\putmorphism(-250,500)(1,0)[F(A)`F(A)` =]{550}1a
 \putmorphism(-250,50)(1,0)[F(A)`F(A)` F(1_A)]{550}1a
 \putmorphism(-250,-400)(1,0)[F(A')`F(A')` F(1_{A'})]{550}1a

\putmorphism(-280,500)(0,-1)[\phantom{Y_2}``= ]{450}1l
 \putmorphism(-280,70)(0,-1)[\phantom{F(A)}` `F(u)]{450}1l

\putmorphism(300,500)(0,-1)[\phantom{Y_2}``=]{450}1r
\putmorphism(300,70)(0,-1)[\phantom{Y_2}``F(u)]{450}1r
\put(-100,290){\fbox{$F_A$}}
\put(-140,-150){\fbox{$F(Id^u)$}}
\efig},
\qquad
$$
\end{itemize}
where $u,u'$ are composable 1v-cells, $\omega, \zeta$ vertically composable 2-cells, $\alpha, \beta$ horizontally composable 2-cells, and $f,g,h$ composable 1h-cells. 

A {\em pseudodouble functor} is a lax double functor whose compositor and unitor 2-cells are invertible. 
\end{defn}

We now define horizontal oplax and vertical lax transformations between lax double functors and their modifications. The reason why we introduce these kinds of transformations is explained at the beginning of \sslabel{quasi}.

\begin{defn} \delabel{hor nat tr}
A {\em horizontal oplax transformation} $\alpha$ between lax double functors $F,G\colon \Aa\to\Bb$ consists of the following:
\begin{enumerate}
\item for every 0-cell $A$ in $\Aa$ a 1h-cell $\alpha(A)\colon F(A)\to G(A)$ in $\Bb$,
\item for every 1v-cell $u\colon A\to A'$ in $\Aa$ a 2-cell in $\Bb$:
$$
\scalebox{0.86}{
\bfig
\putmorphism(-150,50)(1,0)[F(A)`G(A)`\alpha(A)]{560}1a
\putmorphism(-150,-320)(1,0)[F(A')`G(A')`\alpha(A')]{600}1a
\putmorphism(-180,50)(0,-1)[\phantom{Y_2}``F(u)]{370}1l
\putmorphism(410,50)(0,-1)[\phantom{Y_2}``G(u)]{370}1r
\put(30,-110){\fbox{$\alpha^u$}}
\efig}
$$
\item for every 1h-cell $f\colon A\to B$ in $\Aa$ a 2-cell in $\Bb$:
$$
\scalebox{0.86}{
\bfig
 \putmorphism(-170,500)(1,0)[F(A)`F(B)`F(f)]{540}1a
 \putmorphism(360,500)(1,0)[\phantom{F(f)}`G(B) `\alpha(B)]{560}1a
 \putmorphism(-170,120)(1,0)[F(A)`G(A)`\alpha(A)]{540}1a
 \putmorphism(360,120)(1,0)[\phantom{G(B)}`G(A) `G(f)]{560}1a
\putmorphism(-180,500)(0,-1)[\phantom{Y_2}``=]{380}1r
\putmorphism(940,500)(0,-1)[\phantom{Y_2}``=]{380}1r
\put(280,310){\fbox{$\delta_{\alpha,f}$}}
\efig}
$$
\end{enumerate}
so that the following are satisfied: 
\begin{itemize}
\item (coherence with compositors for $\delta_{\alpha,-}$): for any composable 1h-cells $f$ and $g$ in $\Aa$ the 2-cell 
$\delta_{\alpha,gf}$ satisfies: \\
{\em \axiom{h.o.t.-1}} 
$$\scalebox{0.82}{
\bfig
  \putmorphism(-750,200)(1,0)[F(A)`\phantom{F(B)}`F(f)]{600}1a
\putmorphism(-130,200)(1,0)[F(A)`F(C)`F(g)]{580}1a
\putmorphism(-730,200)(0,-1)[\phantom{Y_2}``=]{400}1r
\putmorphism(420,200)(0,-1)[\phantom{Y_2}``=]{400}1r
\putmorphism(-730,-210)(0,-1)[\phantom{Y_2}``=]{400}1r
\putmorphism(1030,-210)(0,-1)[\phantom{Y_2}``=]{400}1r
 \putmorphism(450,-210)(1,0)[F(C)`G(C) `\alpha(C)]{580}1a
  \putmorphism(-750,-210)(1,0)[F(A)`\phantom{F(B)}`F(gf)]{1200}1a
\put(-270,20){\fbox{$F_{gf}$}}
 \putmorphism(-750,-615)(1,0)[F(A)`G(A)`\alpha(A)]{620}1a
 \putmorphism(-120,-615)(1,0)[\phantom{F(B)}`G(C) `G(gf)]{1170}1a
\put(-250,-420){\fbox{$\delta_{\alpha,gf}$}}
\efig}= 
\scalebox{0.82}{
\bfig
 \putmorphism(450,150)(1,0)[F(B)`F(C) `F(g)]{680}1a
 \putmorphism(1120,150)(1,0)[\phantom{F(B)}`G(C) `\alpha(C)]{600}1a
\put(1000,-100){\fbox{$\delta_{\alpha,g}$}}

  \putmorphism(-150,-300)(1,0)[F(A)` F(B) `F(f)]{600}1a
\putmorphism(450,-300)(1,0)[\phantom{F(A)}` G(B) `\alpha(B)]{680}1a
 \putmorphism(1120,-300)(1,0)[\phantom{F(A)}`G(C) ` G(g)]{620}1a

\putmorphism(450,150)(0,-1)[\phantom{Y_2}``=]{450}1l
\putmorphism(1710,150)(0,-1)[\phantom{Y_2}``=]{450}1r

 \putmorphism(-150,-750)(1,0)[F(A)`G(A)`\alpha(A)]{600}1a
 \putmorphism(450,-750)(1,0)[\phantom{F(B)}`G(B) `G(f)]{680}1a
 \putmorphism(1120,-750)(1,0)[\phantom{F(B)}`G(C) `G(g)]{620}1a

\putmorphism(-180,-300)(0,-1)[\phantom{Y_2}``=]{450}1r
\putmorphism(1040,-300)(0,-1)[\phantom{Y_2}``=]{450}1r
\put(350,-540){\fbox{$\delta_{\alpha,f}$}}
\put(1000,-960){\fbox{$G_{gf}$}}

 \putmorphism(450,-1200)(1,0)[G(A)` G(C) `G(gf)]{1300}1a

\putmorphism(450,-750)(0,-1)[\phantom{Y_2}``=]{450}1l
\putmorphism(1750,-750)(0,-1)[\phantom{Y_2}``=]{450}1r
\efig}
$$ 
(coherence with unitors for $\delta_{\alpha,-}$): for any object $A\in\Aa$: 
$$\text{{\em \axiom{h.o.t.-2}}}  \qquad\quad
\scalebox{0.86}{
\bfig
 \putmorphism(-150,420)(1,0)[F(A)`F(A)`=]{500}1a
\putmorphism(-180,420)(0,-1)[\phantom{Y_2}``=]{370}1l
\putmorphism(320,420)(0,-1)[\phantom{Y_2}``=]{370}1r
 \putmorphism(-150,50)(1,0)[F(A)`F(A)`F(1_A)]{500}1a
 \put(-80,250){\fbox{$F_A$}} 
\putmorphism(330,50)(1,0)[\phantom{F(A)}`G(A) `\alpha(A)]{560}1a
 \putmorphism(-170,-350)(1,0)[F(A)`G(A)`\alpha(A)]{520}1a
 \putmorphism(350,-350)(1,0)[\phantom{F(A)}`G(A) `G(1_A)]{560}1a

\putmorphism(-180,50)(0,-1)[\phantom{Y_2}``=]{400}1l
\putmorphism(910,50)(0,-1)[\phantom{Y_2}``=]{400}1r
\put(240,-150){\fbox{$\delta_{\alpha,1_A}$}}
\efig}
\quad=\quad
\scalebox{0.86}{
\bfig
 \putmorphism(-150,420)(1,0)[F(A)`G(A)`\alpha(A)]{500}1a
\putmorphism(-180,420)(0,-1)[\phantom{Y_2}``=]{370}1l
\putmorphism(320,420)(0,-1)[\phantom{Y_2}``=]{370}1r
  \put(-100,230){\fbox{$\Id_{\alpha(A)}$}} 
\putmorphism(-150,50)(1,0)[F(A)` \phantom{Y_2} `\alpha(A)]{450}1a

\putmorphism(350,50)(1,0)[G(A)` G(A) `=]{470}1a
\putmorphism(330,-300)(1,0)[G(A)` G(A) `G(1_A)]{480}1b
\putmorphism(330,50)(0,-1)[\phantom{Y_2}``=]{350}1l
\putmorphism(800,50)(0,-1)[\phantom{Y_2}``=]{350}1r
\put(470,-150){\fbox{$G_A$}}
\efig}
$$

\item (coherence with vertical composition and identity for $\alpha^\bullet$): for any composable 1v-cells $u$ and $v$ in $\Aa$:
$$\text{{\em \axiom{h.o.t.-3}}} \label{h.o.t.-3} \qquad\alpha^{\frac{u}{v}}=\frac{\alpha^u}{\alpha^v}\quad\qquad\text{ and}\quad\qquad
\text{{\em \axiom{h.o.t.-4}}} \label{h.o.t.-4} \qquad\alpha^{1^A}=\Id_{\alpha(A)};$$

\item (oplax naturality of 2-cells):
for every 2-cell in $\Aa$
$\scalebox{0.86}{
\bfig
\putmorphism(-150,50)(1,0)[A` B`f]{400}1a
\putmorphism(-150,-270)(1,0)[A'`B' `g]{400}1b
\putmorphism(-170,50)(0,-1)[\phantom{Y_2}``u]{320}1l
\putmorphism(250,50)(0,-1)[\phantom{Y_2}``v]{320}1r
\put(0,-140){\fbox{$a$}}
\efig}$ 
the following identity in $\Bb$ must hold:\\
$\text{{\em \axiom{h.o.t.-5}}}$ 
$$
\scalebox{0.86}{
\bfig
\putmorphism(-150,500)(1,0)[F(A)`F(B)`F(f)]{600}1a
 \putmorphism(450,500)(1,0)[\phantom{F(A)}`G(B) `\alpha(B)]{640}1a

 \putmorphism(-150,50)(1,0)[F(A')`F(B')`F(g)]{600}1a
 \putmorphism(450,50)(1,0)[\phantom{F(A)}`G(B') `\alpha(B')]{640}1a

\putmorphism(-180,500)(0,-1)[\phantom{Y_2}``F(u)]{450}1l
\putmorphism(450,500)(0,-1)[\phantom{Y_2}``]{450}1r
\putmorphism(300,500)(0,-1)[\phantom{Y_2}``F(v)]{450}0r
\putmorphism(1100,500)(0,-1)[\phantom{Y_2}``G(v)]{450}1r
\put(0,260){\fbox{$F(a)$}}
\put(700,270){\fbox{$\alpha^v$}}

\putmorphism(-150,-400)(1,0)[F(A')`G(A') `\alpha(A')]{640}1a
 \putmorphism(450,-400)(1,0)[\phantom{A'\ot B'}` G(B') `G(g)]{680}1a

\putmorphism(-180,50)(0,-1)[\phantom{Y_2}``=]{450}1l
\putmorphism(1120,50)(0,-1)[\phantom{Y_3}``=]{450}1r
\put(320,-200){\fbox{$\delta_{\alpha,g}$}}

\efig}
\quad=\quad
\scalebox{0.86}{
\bfig
\putmorphism(-150,500)(1,0)[F(A)`F(B)`F(f)]{600}1a
 \putmorphism(450,500)(1,0)[\phantom{F(A)}`G(B) `\alpha(B)]{680}1a
 \putmorphism(-150,50)(1,0)[F(A)`G(A)`\alpha(A)]{600}1a
 \putmorphism(450,50)(1,0)[\phantom{F(A)}`G(B) `G(f)]{680}1a

\putmorphism(-180,500)(0,-1)[\phantom{Y_2}``=]{450}1r
\putmorphism(1100,500)(0,-1)[\phantom{Y_2}``=]{450}1r
\put(350,260){\fbox{$\delta_{\alpha,f}$}}
\put(650,-180){\fbox{$G(a)$}}

\putmorphism(-150,-400)(1,0)[F(A')`G(A') `\alpha(A')]{640}1a
 \putmorphism(490,-400)(1,0)[\phantom{F(A')}` G(B'). `G(g)]{640}1a

\putmorphism(-180,50)(0,-1)[\phantom{Y_2}``F(u)]{450}1l
\putmorphism(450,50)(0,-1)[\phantom{Y_2}``]{450}1l
\putmorphism(610,50)(0,-1)[\phantom{Y_2}``G(u)]{450}0l 
\putmorphism(1120,50)(0,-1)[\phantom{Y_3}``G(v)]{450}1r
\put(40,-180){\fbox{$\alpha^u$}} 
\efig}
$$
\end{itemize}
A {\em horizontal lax transformation} is defined by using the opposite direction of the 2-cells $\delta_{\alpha,f}$ in item 3 
and accommodating the axioms correspondingly. \\
A {\em horizontal pseudonatural transformation} is a horizontal oplax transformation for which the 2-cells $\delta_{\alpha,f}$ are isomorphisms. 
A {\em horizontal strict transformation} is a horizontal oplax transformation for which the 2-cells $\delta_{\alpha,f}$ are identities. 
\end{defn}

\begin{rem} \rmlabel{H trans}
Given a horizontal pseudonatural transformation $\alpha$ between lax double functors $F,G\colon \Aa\to\Bb$ we denote by 
$\HH(\alpha):\HH(F)\Rightarrow\HH(G)$ the underlying pseudonatural transformation between the underlying lax functors 
between the underlying 2-categories $\HH(\Aa)$ and $\HH(\Bb)$. Recall that a bicategorical pseudonatural transformation, in this case 
$\HH(\alpha)$, is called a {\em pseudonatural equivalence} if the 1h-cell components $\alpha(A)$ are equivalence 1-cells in $\HH(\Bb)$. 
\end{rem}

\begin{lma} \cite[Lemma 2.3]{Fem:Bif} \lelabel{vert comp hor.ps.tr.} \\
Vertical composition of two horizontal oplax transformations 
$F\stackrel{\alpha}{\Rightarrow} G \stackrel{\beta}{\Rightarrow}H$ between lax functors $F, G, H\colon\Aa\to\Bb$, 
denoted by $\frac{\alpha}{\beta}$, is given by: 
\begin{enumerate}
\item for every 0-cell $A$ in $\Aa$ a 1h-cell in $\Bb$:
$$(\frac{\alpha}{\beta})(A)=\big( F(A)\stackrel{\alpha(A)}{\longrightarrow}G(A) \stackrel{\beta(A)}{\longrightarrow} H(A) \big),$$ 
\item for every 1v-cell $u\colon A\to A'$ in $\Aa$ a 2-cell in $\Bb$:
$$(\frac{\alpha}{\beta})^u=
\scalebox{0.9}{
\bfig
\putmorphism(-150,200)(1,0)[F(A)`G(A)`\alpha(A)]{560}1a
 \putmorphism(430,200)(1,0)[\phantom{F(A)}`H(A) `\beta(A)]{620}1a

\putmorphism(-150,-170)(1,0)[F(A')`G(A')`\alpha(A')]{600}1a
\putmorphism(-180,200)(0,-1)[\phantom{Y_2}``F(u)]{370}1l
\putmorphism(380,200)(0,-1)[\phantom{Y_2}``G(u)]{370}1r
\put(0,20){\fbox{$\alpha^u$}}
 \putmorphism(460,-170)(1,0)[\phantom{F(A)}`H(A') `\beta(A')]{630}1a
\putmorphism(1060,200)(0,-1)[\phantom{Y_2}``H(u)]{370}1r
\put(670,40){\fbox{$\beta^u$}}
\efig}
$$
\item for every 1h-cell $f\colon A\to B$  in $\Aa$ a 2-cell in $\Bb$:
$$\delta_{\frac{\alpha}{\beta},f}=
\scalebox{0.9}{
\bfig

 \putmorphism(-200,300)(1,0)[F(A)`F(B)` F(f)]{650}1a
 \putmorphism(430,300)(1,0)[\phantom{F(B)}`G(B) `\alpha(B)]{700}1a

 \putmorphism(-200,-100)(1,0)[F(A)`G(A)`\alpha(A)]{650}1a
 \putmorphism(430,-100)(1,0)[\phantom{A\ot B}`G(B) `G(f)]{700}1a
 \putmorphism(1050,-100)(1,0)[\phantom{A'\ot B'}` H(B) `\beta(B)]{700}1a

\putmorphism(-230,300)(0,-1)[\phantom{Y_2}``=]{400}1r
\putmorphism(1050,300)(0,-1)[\phantom{Y_2}``=]{400}1r
\put(300,110){\fbox{$ \delta_{\alpha,f}  $}}
\put(1000,-310){\fbox{$\delta_{\beta,f}$}}

 \putmorphism(450,-500)(1,0)[G(A)` H(A) `\beta(A)]{700}1a
 \putmorphism(1080,-500)(1,0)[\phantom{A''\ot B'}`G(B). ` H(f)]{700}1a

\putmorphism(450,-100)(0,-1)[\phantom{Y_2}``=]{400}1l
\putmorphism(1750,-100)(0,-1)[\phantom{Y_2}``=]{400}1r
\efig}
$$
\end{enumerate}
\end{lma}

\begin{defn} \delabel{vlt}
A {\em vertical lax transformation} $\alpha_0$ between lax double functors $F,G\colon \Aa\to\Bb$ consists of: 
\begin{enumerate}
\item a 1v-cell $\alpha_0(A)\colon F(A)\to G(A)$ in $\Bb$ for every 0-cell $A$ in $\Aa$; 
\item 
for every 1h-cell $f\colon A\to B$ in $\Aa$ a 2-cell in $\Bb$:
$$ 
\scalebox{0.9}{
\bfig
\putmorphism(-150,180)(1,0)[F(A)`F(B)`F(f)]{560}1a
\putmorphism(-150,-190)(1,0)[G(A)`G(B)`G(f)]{600}1a
\putmorphism(-150,180)(0,-1)[\phantom{Y_2}``\alpha_0(A)]{370}1l
\putmorphism(410,180)(0,-1)[\phantom{Y_2}``\alpha_0(B)]{370}1r
\put(-30,20){\fbox{$(\alpha_0)_f$}}
\efig}
$$
\item for every 1v-cell $u\colon A\to A'$ in $\Aa$ a 2-cell in $\Bb$:
$$
\scalebox{0.9}{
\bfig
 \putmorphism(-90,500)(1,0)[F(A)`F(A) `=]{540}1a
\putmorphism(450,500)(0,-1)[\phantom{Y_2}`F(\tilde A) `F(u)]{400}1r
\putmorphism(-90,-300)(1,0)[G(\tilde A)`G(\tilde A) `=]{540}1a
\putmorphism(450,100)(0,-1)[\phantom{Y_2}``\alpha_0(\tilde A)]{400}1r
\putmorphism(-100,100)(0,-1)[\phantom{Y_2}``G(u)]{400}1l
\putmorphism(-100,500)(0,-1)[\phantom{Y_2}`G(A) `\alpha_0(A)]{400}1l
\put(80,70){\fbox{$\alpha_0^u$}}
\efig}
$$
\end{enumerate}
which need to satisfy: 
\begin{itemize}
\item (coherence with compositors for $(\alpha_0)_\bullet$): for any composable 1h-cells $f$ and $g$ in $\Aa$: \\
$$\text{{\em \axiom{v.l.t.\x 1}}}  \quad
\scalebox{0.84}{
\bfig
\putmorphism(-150,500)(1,0)[F(A)`F(B)`F(f)]{600}1a
 \putmorphism(450,500)(1,0)[\phantom{F(A)}`F(C) `F(g)]{600}1a

 \putmorphism(-150,50)(1,0)[G(A)`G(B)`G(f)]{600}1a
 \putmorphism(450,50)(1,0)[\phantom{F(A)}`G(C) `G(g)]{600}1a

\putmorphism(-180,500)(0,-1)[\phantom{Y_2}``\alpha_0(A)]{450}1l
\putmorphism(450,500)(0,-1)[\phantom{Y_2}``]{450}1r
\putmorphism(300,500)(0,-1)[\phantom{Y_2}``\alpha_0(B)]{450}0r
\putmorphism(1060,500)(0,-1)[\phantom{Y_2}``\alpha_0(C)]{450}1r
\put(0,260){\fbox{$(\alpha_0)_f$}}
\put(660,270){\fbox{$(\alpha_0)_g$}}

\putmorphism(-150,-400)(1,0)[G(A)`G(C) `G(gf)]{1220}1a

\putmorphism(-180,50)(0,-1)[\phantom{Y_2}``=]{450}1l
\putmorphism(1060,50)(0,-1)[\phantom{Y_3}``=]{450}1r
\put(320,-170){\fbox{$G_{gf}$}}

\efig}
=\quad\hspace{-0,2cm}
\scalebox{0.84}{
\bfig
\putmorphism(-160,500)(1,0)[F(A)`F(B)`F(f)]{580}1a
 \putmorphism(420,500)(1,0)[\phantom{F(A)}`F(C) `F(g)]{620}1a

 \putmorphism(-160,50)(1,0)[F(A)`F(C)`F(gf)]{1260}1a

\putmorphism(-180,500)(0,-1)[\phantom{Y_2}``=]{450}1r
\putmorphism(1060,500)(0,-1)[\phantom{Y_2}``=]{450}1l
\put(330,260){\fbox{$F_{gf}$}}
\put(300,-180){\fbox{$(\alpha_0)_{gf}$}}

\putmorphism(-180,50)(0,-1)[\phantom{Y_2}``\alpha_0(A)]{450}1r
\putmorphism(1060,50)(0,-1)[\phantom{Y_3}``\alpha_0(C)]{450}1l
\putmorphism(-150,-400)(1,0)[G(A)`G(C) `G(gf)]{1220}1a
\efig}
$$

(coherence with unitors for $(\alpha_0)_\bullet$): for any object $A$ in $\Aa$:
$$\text{{\em \axiom{v.l.t.\x 2}}}  \qquad\qquad 
\scalebox{0.86}{
\bfig
 \putmorphism(-170,420)(1,0)[F(A)`F(A)`=]{500}1a
\putmorphism(-180,420)(0,-1)[\phantom{Y_2}``=]{370}1l
\putmorphism(280,420)(0,-1)[\phantom{Y_2}``=]{370}1r
  \put(-40,250){\fbox{$F_A$}} 
\putmorphism(-170,50)(1,0)[F(A)` F(A) `F(1_A)]{450}1a

\putmorphism(-170,50)(0,-1)[\phantom{Y_2}``\alpha_0(A)]{350}1l
\putmorphism(280,50)(0,-1)[\phantom{Y_2}``\alpha_0(A)]{350}1r
\putmorphism(-170,-300)(1,0)[G(A)` G(A) `G(1_A)]{460}1b
\put(-100,-140){\fbox{$(\alpha_0)_{1_A}$}}
\efig}\quad
=
\scalebox{0.86}{
\bfig
 \putmorphism(-170,420)(1,0)[F(A)`F(A)`=]{500}1a
\putmorphism(-180,420)(0,-1)[\phantom{Y_2}``\alpha_0(A)]{370}1l
\putmorphism(280,420)(0,-1)[\phantom{Y_2}``\alpha_0(A)]{370}1r
  \put(-120,220){\fbox{$\Id^{\alpha_0(A)}$}} 
\putmorphism(-170,50)(1,0)[G(A)` G(A) `=]{450}1a

\putmorphism(-170,50)(0,-1)[\phantom{Y_2}``=]{350}1l
\putmorphism(280,50)(0,-1)[\phantom{Y_2}``=]{350}1r
\putmorphism(-170,-300)(1,0)[G(A)` G(A) `G(1_A)]{460}1b
\put(-40,-140){\fbox{$G_A$}}
\efig}
$$
\item (coherence with vertical composition for $\alpha_0^\bullet$): for any composable 1v-cells $u$ and $v$ in $\Aa$: 
$$ \text{{\em \axiom{v.l.t.\x 3}}}   \qquad\quad 
\alpha_0^{\frac{u}{u'}}=
\scalebox{0.86}{
\bfig 
 \putmorphism(-150,500)(1,0)[F(A)`F(A) `=]{500}1a
\putmorphism(-130,500)(0,-1)[\phantom{Y_2}`G(A) `\alpha_0(A)]{400}1l
\put(0,250){\fbox{$\alpha_0^u$}}
\putmorphism(-150,-300)(1,0)[G(\tilde A)`G(\tilde A) `=]{460}1a
\putmorphism(-130,110)(0,-1)[\phantom{Y_2}``G(u)]{400}1l
\putmorphism(380,500)(0,-1)[\phantom{Y_2}` `F(u)]{400}1r
\putmorphism(380,100)(0,-1)[F(\tilde A)` `\alpha_0(\tilde A)]{400}1l
\putmorphism(480,110)(1,0)[`F(\tilde A)`=]{460}1a
\putmorphism(390,-700)(1,0)[\phantom{G(A)}`G(\tilde{\tilde A})`=]{570}1a
\putmorphism(920,100)(0,-1)[\phantom{(B, \tilde A')}``F(u')]{400}1r
\putmorphism(920,-300)(0,-1)[F(\tilde{\tilde A})`` \alpha_0(\tilde{\tilde A})]{400}1r
\putmorphism(400,-300)(0,-1)[\phantom{(B, \tilde A)}`G(\tilde{\tilde A}) `G(u')]{400}1l
\put(530,-190){\fbox{$\alpha_0^{u'}$}}
\efig}
$$
(coherence with vertical identity for $\alpha^\bullet$): for any object $A$ in $\Aa$: \\
$$\text{{\em \axiom{v.l.t.\x 4}}}  \qquad\qquad \alpha_0^{1^A}=\Id^{\alpha_0(A)}  \hspace{4cm}$$
\item (lax naturality of 2-cells):
for every 2-cell in $\Aa$
$\scalebox{0.86}{
\bfig
\putmorphism(-150,50)(1,0)[A` B`f]{400}1a
\putmorphism(-150,-270)(1,0)[\tilde A` \tilde B `g]{400}1b
\putmorphism(-150,50)(0,-1)[\phantom{Y_2}``u]{320}1l
\putmorphism(250,50)(0,-1)[\phantom{Y_2}``v]{320}1r
\put(0,-140){\fbox{$a$}}
\efig}$ 
the following identity in $\Bb$ must hold:\\
$$\text{{\em \axiom{v.l.t.\x 5}}}  \quad
\scalebox{0.84}{
\bfig
 \putmorphism(-130,500)(1,0)[F(A)`F(A) `=]{560}1a
 \putmorphism(530,500)(1,0)[` `F(f)]{400}1a
\putmorphism(-130,500)(0,-1)[\phantom{Y_2}`G(A) `\alpha_0(A)]{450}1l
\put(40,50){\fbox{$\alpha_0^u$}}
\putmorphism(-130,-400)(1,0)[G(\tilde A)` `=]{470}1a
\putmorphism(-150,50)(0,-1)[\phantom{Y_2}``G(u)]{450}1l
\putmorphism(450,50)(0,-1)[\phantom{Y_2}`G(\tilde A)`\alpha_0(\tilde A)]{450}1l
\putmorphism(450,500)(0,-1)[\phantom{Y_2}`F(\tilde A) `F(u)]{450}1l
\put(600,260){\fbox{$F(a)$}}
\putmorphism(420,50)(1,0)[\phantom{(B, \tilde A)}``F(g)]{500}1a
\putmorphism(1030,50)(0,-1)[\phantom{(B, A')}`G(\tilde B)`]{450}1r
\putmorphism(1010,50)(0,-1)[\phantom{(B, A')}``\alpha_0(\tilde B)]{450}0r
\putmorphism(1030,500)(0,-1)[F(B)`F(\tilde B)`]{450}1r
\putmorphism(1010,500)(0,-1)[``F(v)]{450}0r
\putmorphism(420,-400)(1,0)[\phantom{(B, \tilde A)}``G(g)]{500}1a
\put(600,-170){\fbox{$(\alpha_0)_g$}}
\efig}
\quad\hspace{-0,1cm}=\quad\hspace{-0,14cm}
\scalebox{0.84}{
\bfig
 \putmorphism(-130,500)(1,0)[F(A)`F(B) `F(f)]{560}1a
 \putmorphism(550,500)(1,0)[` `=]{380}1a
\putmorphism(-140,500)(0,-1)[\phantom{Y_2}`G(A) `]{450}1l
\putmorphism(-120,500)(0,-1)[\phantom{Y_2}` `\alpha_0(A)]{450}0l
\put(650,50){\fbox{$\alpha_0^v$}}
\putmorphism(-130,-400)(1,0)[G(\tilde A)` `G(g)]{470}1a
\putmorphism(-140,50)(0,-1)[\phantom{Y_2}``]{450}1l
\putmorphism(-120,50)(0,-1)[\phantom{Y_2}``G(u)]{450}0l
\putmorphism(450,50)(0,-1)[\phantom{Y_2}`G(\tilde B)`G(v)]{450}1r
\putmorphism(450,500)(0,-1)[\phantom{Y_2}`G(B) `\alpha_0(B)]{450}1r
\put(40,280){\fbox{$(\alpha_0)_f$}}
\putmorphism(-150,50)(1,0)[\phantom{(B, \tilde A)}``G(f)]{500}1a
\putmorphism(1030,50)(0,-1)[\phantom{(B, A')}`G(\tilde B).`]{450}1l
\putmorphism(1050,50)(0,-1)[``\alpha_0(\tilde B)]{450}0l
\putmorphism(1030,500)(0,-1)[F(B)`F(\tilde B)`]{450}1l
\putmorphism(1050,500)(0,-1)[``F(v)]{450}0l
\putmorphism(430,-400)(1,0)[\phantom{(B, \tilde A)}``=]{470}1b
\put(0,-170){\fbox{$G(a)$}}
\efig}
$$
\end{itemize}

A {\em vertical oplax transformation} is defined by using the opposite direction of the 2-cells $\alpha_0^u$ in item 3 
and accommodating the axioms correspondingly. 

A vertical lax transformation is called a {\em vertical pseudonatural transformation} if the globular 2-cells $\alpha_0^u$ are horizontally invertible for all 1v-cells $u$ in $\Aa$. 

It is called a {\em vertical strict transformation} if all the globular 2-cells $\alpha_0^u$ are identities. A vertical strict transformation is said to be 
{\em invertible} if the 2-cells $(\alpha_0)_f$ for all 1h-cells $f$ in $\Aa$ are vertically invertible (this includes the condition that 
the 1v-cells $\alpha_0(A)$ are invertible for all $A\in\Aa$). 
\end{defn}

Similarly as in \cite[Lemma 2.6]{Fem:Bif} and \cite[Lemma 3.8]{Fem:Fil} one has:

\begin{lma} \lelabel{vert comp vert. lx tr.}
Vertical composition of two vertical strict transformations $\alpha_0\colon F\Rightarrow G\colon \Aa\to\Bb$ and 
$\beta_0\colon G\Rightarrow H\colon\Aa\to\Bb$ between lax double functors, denoted by $\frac{\alpha_0}{\beta_0}$, is given by: \\
for every 0-cell $A$ in $\Aa$ a 1v-cell on the left below, and for every 1h-cell $f\colon A\to B$ in $\Aa$ a 2-cell on the right below, both in $\Bb$:
$$
(\frac{\alpha_0}{\beta_0})(A)=
\scalebox{0.86}{
\bfig
\putmorphism(-280,500)(0,-1)[F(A)`G(A) `\alpha_0(A)]{450}1l
 \putmorphism(-280,70)(0,-1)[\phantom{F(A)}`H(A) `\beta_0(A)]{450}1l
\efig} \qquad
(\frac{\alpha_0}{\beta_0})(f)= 
\scalebox{0.9}{
\bfig
\putmorphism(-250,500)(1,0)[F(A)`F(B)` F(f)]{550}1a
 \putmorphism(-250,50)(1,0)[G(A)`G(B)` G(f)]{550}1a
 \putmorphism(-250,-400)(1,0)[H(A)`H(B).` H(f)]{550}1a

\putmorphism(-280,500)(0,-1)[\phantom{Y_2}``\alpha_0(A)]{450}1l
 \putmorphism(-280,70)(0,-1)[\phantom{F(A)}` `\beta_0(A)]{450}1l

\putmorphism(300,500)(0,-1)[\phantom{Y_2}``\alpha_0(B)]{450}1r
\putmorphism(300,70)(0,-1)[\phantom{Y_2}``\beta_0(B)]{450}1r
\put(-120,290){\fbox{$(\alpha_0)_f$}}
\put(-120,-150){\fbox{$(\beta_0)_f$}}
\efig}
$$ 
\end{lma}

We will need modifications between horizontal and vertical pseudonatural transformations. Their definition is formally the same as of 
modifications of more general transformations that we cite next from \cite[Definition 2.7]{Fem:Bif}.

\begin{defn} \delabel{modif-hv}
A modification $\Theta$ between two horizontal oplax transformations $\alpha$ and $\beta$ and two vertical lax transformations $\alpha_0$ and $\beta_0$ 
depicted below on the left, where the lax double functors $F, G, F\s', G'$ act between $\Aa\to\Bb$, is given 
by a collection of 2-cells in $\Bb$ depicted below on the right:
$$
\scalebox{0.86}{
\bfig
\putmorphism(-150,50)(1,0)[F` G`\alpha]{400}1a
\putmorphism(-150,-270)(1,0)[F'`G' `\beta]{400}1b
\putmorphism(-170,50)(0,-1)[\phantom{Y_2}``\alpha_0]{320}1l
\putmorphism(250,50)(0,-1)[\phantom{Y_2}``\beta_0]{320}1r
\put(-30,-140){\fbox{$\Theta$}}
\efig}
\qquad\qquad
\scalebox{0.86}{
\bfig
\putmorphism(-180,50)(1,0)[F(A)` G(A)`\alpha(A)]{550}1a
\putmorphism(-180,-270)(1,0)[F\s'(A)`G'(A) `\beta(A)]{550}1b
\putmorphism(-170,50)(0,-1)[\phantom{Y_2}``\alpha_0(A)]{320}1l
\putmorphism(350,50)(0,-1)[\phantom{Y_2}``\beta_0(A)]{320}1r
\put(0,-140){\fbox{$\Theta_A$}}
\efig}
$$
which satisfy the following rules: 

\medskip

\noindent {\em \axiom{m.ho-vl.-1}}  for every 1h-cell $f$, we have  
$$
\scalebox{0.86}{
\bfig
\putmorphism(-150,500)(1,0)[F(A)`F(B)`F(f)]{600}1a
 \putmorphism(450,500)(1,0)[\phantom{F(A)}`G(B) `\alpha(B)]{620}1a

 \putmorphism(-150,50)(1,0)[F\s'(A)`F\s'(B)`F\s'(f)]{600}1a
 \putmorphism(450,50)(1,0)[\phantom{F(A)}`G'(B) `\beta(B)]{620}1a

\putmorphism(-180,500)(0,-1)[\phantom{Y_2}``\alpha_0(A)]{450}1l
\putmorphism(450,500)(0,-1)[\phantom{Y_2}``]{450}1r
\putmorphism(300,500)(0,-1)[\phantom{Y_2}``\alpha_0(B)]{450}0r
\putmorphism(1080,500)(0,-1)[\phantom{Y_2}``\beta_0(B)]{450}1r
\put(0,280){\fbox{$(\alpha_0)_f$}}
\put(670,280){\fbox{$\Theta_B$}}

\putmorphism(-150,-400)(1,0)[F\s'(A)`G'(A) `\beta(A)]{600}1a
\putmorphism(510,-400)(1,0)[\phantom{Y_2}`G'(B) `G'(f)]{580}1a

\putmorphism(-180,50)(0,-1)[\phantom{Y_2}``=]{450}1l
\putmorphism(1080,50)(0,-1)[\phantom{Y_3}``=]{450}1r
\put(320,-180){\fbox{$\delta_{\beta,f}$}}

\efig}
\quad=\quad
\scalebox{0.86}{
\bfig
\putmorphism(-150,500)(1,0)[F(A)`F(B)`F(f)]{600}1a
 \putmorphism(450,500)(1,0)[\phantom{F(A)}`G(B) `\alpha(B)]{620}1a
\putmorphism(-150,50)(1,0)[F(A)`G(A) `\alpha(A)]{600}1a
\putmorphism(510,50)(1,0)[\phantom{Y_2}`G(B) `G(f)]{580}1a

\putmorphism(-180,500)(0,-1)[\phantom{Y_2}``=]{450}1r
\putmorphism(1080,500)(0,-1)[\phantom{Y_2}``=]{450}1r
\put(350,280){\fbox{$\delta_{\alpha,f}$}}

\putmorphism(-180,50)(0,-1)[\phantom{Y_2}``\alpha_0(A)]{450}1l
\putmorphism(1080,50)(0,-1)[\phantom{Y_3}``\beta_0(B)]{450}1r
\put(20,-180){\fbox{$\Theta_A$}}
\put(670,-180){\fbox{$(\beta_0)_f$}}

\putmorphism(450,50)(0,-1)[\phantom{Y_2}``]{450}1r
\putmorphism(300,50)(0,-1)[\phantom{Y_2}``\beta_0(A)]{450}0r

\putmorphism(-150,-400)(1,0)[F\s'(A)`G'(A) `\beta(A)]{600}1a
\putmorphism(510,-400)(1,0)[\phantom{Y_2}`G'(B) `G'(f)]{580}1a
\efig}
$$
and

\noindent {\em \axiom{m.ho-vl.-2}}   for every 1v-cell $u$, we have
$$
\scalebox{0.86}{
\bfig
 \putmorphism(-150,500)(1,0)[F(A)`F(A) `=]{600}1a
 \putmorphism(550,500)(1,0)[` `\alpha(A)]{400}1a
\putmorphism(-180,500)(0,-1)[\phantom{Y_2}`F\s'(A) `\alpha_0(A)]{450}1l
\put(30,50){\fbox{$\alpha_0^u$}}
\putmorphism(-150,-400)(1,0)[F\s'(\tilde A)` `=]{480}1a
\putmorphism(-180,50)(0,-1)[\phantom{Y_2}``F\s'(u)]{450}1l
\putmorphism(450,50)(0,-1)[\phantom{Y_2}`F\s'(\tilde A)` \alpha_0(\tilde A)]{450}1l
\putmorphism(450,500)(0,-1)[\phantom{Y_2}`F(\tilde A) `F(u)]{450}1l
\put(660,280){\fbox{$\alpha^u$}}
\putmorphism(450,50)(1,0)[\phantom{(B, \tilde A)}``\alpha(\tilde A)]{500}1a
\putmorphism(1070,50)(0,-1)[\phantom{(B, A')}`G'(\tilde A)`\beta_0(\tilde A)]{450}1r
\putmorphism(1070,500)(0,-1)[G(A)`G(\tilde A)`G(u)]{450}1r
\putmorphism(450,-400)(1,0)[\phantom{(B, \tilde A)}``\beta(\tilde A)]{500}1a
\put(640,-170){\fbox{$ \Theta_{\tilde A}$ } } 
\efig}\quad=\quad
\scalebox{0.86}{
\bfig
 \putmorphism(-150,500)(1,0)[F(A)`G(A) `\alpha(A)]{600}1a
 \putmorphism(450,500)(1,0)[\phantom{(B,A)}` `=]{460}1a
\putmorphism(-180,500)(0,-1)[\phantom{Y_2}`F\s'(A) `\alpha_0(A)]{450}1l
\put(650,50){\fbox{$\beta_0^u$}}
\putmorphism(-180,-400)(1,0)[F\s'(\tilde A)` `\beta(\tilde A)]{500}1a
\putmorphism(-180,50)(0,-1)[\phantom{Y_2}``F\s'(u)]{450}1l
\putmorphism(450,50)(0,-1)[\phantom{Y_2}`G'(\tilde A)`G'(u)]{450}1r
\putmorphism(450,500)(0,-1)[\phantom{Y_2}`G'(A) `\beta_0(A)]{450}1r
\put(0,260){\fbox{$\Theta_A$}}
\putmorphism(-180,50)(1,0)[\phantom{(B, \tilde A)}``\beta(A)]{500}1a
\putmorphism(1030,50)(0,-1)[\phantom{(B, A')}` G'(\tilde A). ` \beta_0(\tilde A)]{450}1r
\putmorphism(1030,500)(0,-1)[G(A)`G(\tilde A)` G(u)]{450}1r
\putmorphism(430,-400)(1,0)[\phantom{(B, \tilde A)}``=]{480}1b
\put(70,-170){\fbox{$\beta^u$}}
\efig}
$$
\end{defn}

Taking, so to say, a horizontal and a vertical restriction of modifications in \deref{modif-hv}, we obtain the definitions of:
\begin{itemize} 
\item {\em horizontal modifications} or {\em modifications between horizontal (oplax) transformations} given by families of (horizontally globular) 2-cells 
\begin{equation} \eqlabel{m-hor}
\scalebox{0.86}{
\bfig
\putmorphism(-180,50)(1,0)[F(A)` G(A)`\alpha(A)]{550}1a
\putmorphism(-180,-270)(1,0)[F\s'(A)`G'(A) `\beta(A)]{550}1b
\putmorphism(-170,50)(0,-1)[\phantom{Y_2}``=]{320}1l
\putmorphism(350,50)(0,-1)[\phantom{Y_2}``=]{320}1r
\put(0,-140){\fbox{$\Theta_A$}}
\efig}
\end{equation} 
and axioms \axiom{m.ho.-1} and \axiom{m.ho.-2} obtained from \axiomref{m.ho-vl.-1} and \axiomref{m.ho-vl.-2} by ignoring the 2-cells 
$(\alpha_0)_f, (\beta_0)_f, \alpha_0^u$ and $\beta_0^u$, and 
\item {\em vertical modifications} or {\em modifications between vertical (lax) transformations} given by families of (vertically globular) 2-cells 
\begin{equation} \eqlabel{m-vert}
\scalebox{0.86}{
\bfig
\putmorphism(-180,50)(1,0)[F(A)` G(A)`=]{550}1a
\putmorphism(-180,-270)(1,0)[F\s'(A)`G'(A) `=]{550}1b
\putmorphism(-170,50)(0,-1)[\phantom{Y_2}``\alpha_0(A)]{320}1l
\putmorphism(350,50)(0,-1)[\phantom{Y_2}``\beta_0(A)]{320}1r
\put(0,-140){\fbox{$\Theta_A$}}
\efig}
\end{equation} 
and axioms \axiom{m.vl.-1} and \axiom{m.vl.-2} obtained from \axiomref{m.ho-vl.-1} and \axiomref{m.ho-vl.-2} by ignoring the 2-cells 
$\delta_{\alpha,f}, \delta_{\beta,f}, \alpha^u$ and $\beta^u$.  
\end{itemize}

\begin{rem} \rmlabel{mod}
Given a horizontal modification $\Theta:\alpha\Rrightarrow\beta$ between horizontal pseudonatural transformations between pseudodouble functors we denote by $\HH(\Theta):\HH(\alpha)\Rrightarrow\HH(\beta)$ the underlying modification between the underlying pseudonatural transformations between the underlying pseudofunctors between the underlying 2-categories (it satisfies only the axiom \axiomref{m.ho.-1}). 
\end{rem}

We will use the following fact that is straightforwardly proved (see also \cite[Lemma 3.4]{Fem:Fil}). 

\begin{lma} \lelabel{H of omega}
Given pseudodouble functors $H,L$ and a horizontal pseudonatural transformation $\omega:F\Rightarrow G$ between lax/pseudo double functors 
so that the compositions $HF, HG, FL, GL$ make sense, and let $f$ be a 1h-cell and $u$ a 1v-cell. 
\begin{enumerate}[a)]
\item 
The 2-cells $\delta_{H(\omega),f}=\threefrac{H_{lx}}{H(\delta_{\omega,f})}{H_{oplx}}$ and $H(\omega)^u=H(\omega^u)$, 
where $H_{lx}$ and $H_{oplx}$ represent the lax and oplax functor structure of $H$ applied to the corresponding 1h-cells, respectively, 
define a horizontal pseudonatural transformation $H(\omega):HF\Rightarrow HG$ of lax/pseudo double functors. 
\item 
The 2-cells $\delta_{(\omega L),f}=\delta_{\omega,L(f)}$ and $(\omega L)^u=\omega^{L(u)}$ 
define a horizontal pseudonatural transformation $\omega L:FL\Rightarrow GL$ of lax/pseudo double functors. 
\end{enumerate}
\end{lma}

An analogous result holds for vertical pseudonatural transformations. Its proof is even simpler, as in our convention the 
vertical direction (for lax/pseudo double functors) is strict.

\begin{lma} \lelabel{H of alfa}
Given pseudodouble functors $H,L$ and a vertical pseudonatural transformation $\alpha:F\Rightarrow G$ between lax/pseudo double functors 
so that the compositions $HF, HG, FL, GL$ make sense, and let $f$ be a 1h-cell and $u$ a 1v-cell. 
\begin{enumerate}[a)]
\item 
The 2-cells $H(\alpha)_f=H(\alpha_f)$ and $H(\omega)^u=H(\omega^u)$ define a vertical pseudonatural transformation 
$H(\alpha):HF\Rightarrow HG$ of lax/pseudo double functors. 
\item 
The 2-cells $\delta_{\alpha L,f}=\delta_{\alpha, L(f)}$ and $(\alpha L)^u=\alpha^{L(u)}$ 
define a vertical pseudonatural transformation $\alpha L:FL\Rightarrow GL$ of lax/pseudo double functors. 
\end{enumerate}
\end{lma}

\subsection{Companions and conjoints}

Recall that a {\em companion} for a 1v-cell $u:A\to B$ is a 1h-cell  $\hat u:A\to B$ together with two 2-cells $\Epsilon$ and $\eta$ 
as below satisfying $[\eta\vert\Epsilon]=\Id_{\hat u}$ and $\frac{\eta}{\Epsilon}=\Id_{u}$. On the other hand, 
a {\em conjoint} for a 1v-cell $u:A\to B$ is a 1h-cell  $\check{u}:B\to A$ together with two 2-cells $\Epsilon^*$ and $\eta^*$ 
as below satisfying $[\Epsilon^*\vert\eta^*]=\Id_{\check{u}}$ and $\frac{\eta^*}{\Epsilon^*}=\Id^{u}$.
$$ 
\scalebox{0.9}{
\bfig
\putmorphism(-150,170)(1,0)[A`B`\hat{u}]{400}1a
\putmorphism(-150,-160)(1,0)[B`B`=]{400}1a
\putmorphism(-150,170)(0,-1)[\phantom{Y_2}``u]{330}1l
\putmorphism(250,170)(0,-1)[\phantom{Y_2}``=]{330}1r
\put(0,0){\fbox{$\Epsilon$}}
\efig}
\qquad
\scalebox{0.9}{
\bfig
\putmorphism(-150,170)(1,0)[A`A`=]{400}1a
\putmorphism(-150,-160)(1,0)[A`B`\hat{u}]{400}1b
\putmorphism(-150,170)(0,-1)[\phantom{Y_2}``=]{330}1l
\putmorphism(250,170)(0,-1)[\phantom{Y_2}``u]{330}1r
\put(0,10){\fbox{$\eta$}}
\efig}
\qquad\qquad
\scalebox{0.9}{
\bfig
\putmorphism(-150,170)(1,0)[A`A`=]{400}1a
\putmorphism(-150,-160)(1,0)[B`A`\check{u}]{400}1b
\putmorphism(-150,170)(0,-1)[\phantom{Y_2}``u]{330}1l
\putmorphism(250,170)(0,-1)[\phantom{Y_2}``=]{330}1r
\put(-10,-10){\fbox{$\eta^*$}}
\efig}
\qquad
\scalebox{0.9}{
\bfig
\putmorphism(-150,170)(1,0)[B`A`\check{u}]{400}1a
\putmorphism(-150,-160)(1,0)[B`B`=]{400}1b
\putmorphism(-150,170)(0,-1)[\phantom{Y_2}``=]{330}1l
\putmorphism(250,170)(0,-1)[\phantom{Y_2}``u]{330}1r
\put(-10,-10){\fbox{$\Epsilon^*$}}
\efig}
$$
Companions are unique up to a unique globular isomorphism. 
This and more properties of companions and conjoints we will recall in \ssref{liftings}, 
see also \cite[Section 1.2]{GP:Adj}, \cite[Section 3]{Shul}.
For the moment we focus on the following few facts about them.

\begin{lma} \cite[Lemma 3.20]{Shul} \lelabel{Shul1} \\
If $u$ is an invertible 1v-cell with a companion $\hat u$, then $\hat u$ is a conjoint of the inverse $u^{-1}$,  
with $\Epsilon^*_{u^{-1}}=\frac{\Epsilon_u}{\Id^{u^{-1}}}$ and $\eta^*_{u^{-1}}=\frac{\Id^{u^{-1}}}{\eta_u}$. 
\end{lma}

\begin{cor} \colabel{conj}
Let $u$ be an invertible 1v-cell. If $u^{-1}$ has a companion $\widehat{(u^{-1})}$, then $u$ has a conjoint given by $\check u=\widehat{(u^{-1})}$ with $\Epsilon^*_u=\frac{\Epsilon_{u^{-1}}}{\Id^u}$ and $\eta^*_u=\frac{\Id^u}{\eta_{u^{-1}}}$. 
\end{cor}

\begin{lma} \cite[Lemma 3.21]{Shul} \lelabel{Shul} \\
If a 1v-cell $u$ has both a companion $\hat u$ and a conjoint $\check u$, then $\check u$ is a left adjoint 1-cell to $\hat u$ 
in the underlying horizontal bicategory of the double category. If moreover $u$ is invertible, then this is an adjoint equivalence. 

The unit $\tilde\eta$ and counit $\tilde\Epsilon$ of the 
adjunction are given via $\tilde\eta=[\eta\vert\eta^*]$ and $\tilde\Epsilon=[\Epsilon^*\vert\Epsilon]$ and their inverses by 
$\tilde\eta^{-1}=\frac{[\Epsilon\vert\Epsilon^*]}{[\Id^{u^{-1}}\vert\Id^{u^{-1}}]}$ and $\tilde\Epsilon^{-1}=[\eta^*\vert\eta]$. 
\end{lma}

We record that also the following identities hold true 
\begin{equation} \eqlabel{fraction rules}
\frac{\eta_u}{\Epsilon^*_{u^{-1}}}=1, \quad \frac{\Epsilon_u}{\eta^*_{u^{-1}}}=\Id_{\hat u}, 
\quad \frac{\Epsilon^*_{u^{-1}}}{\eta}=\Id_{\hat u}, \quad \frac{\eta^*_{u^{-1}}}{\Epsilon_u}=1.
\end{equation}

We will need 
existence of companions and conjoints only for 1v-components of invertible vertical strict transformations.  
In this case by \coref{conj} it is sufficient to require the existence of companions. 
The following result will be very important for us. 


\begin{prop} \prlabel{lifting 1v to equiv}
Let $\alpha_0:F\Rightarrow G$ be a vertical pseudonatural transformation between pseudodouble functors acting between double categories 
$\Dd\to\Ee$. Suppose that the 2-cell components $(\alpha_0)_f$ of $\alpha_0$ are vertically invertible for every 1h-cell $f:A\to B$. 
\begin{enumerate}
\item The following data define a horizontal pseudonatural transformation $\alpha_1:F\Rightarrow G$:
\begin{enumerate} [a)]
\item for all 1v-cell components $\alpha_0(A)$ of $\alpha_0$ a fixed choice of companions (and conjoints) in $\Ee$, 
(we denote their companions by $\alpha_1(A)$ for every 0-cell $A$ of $\Dd$, the corresponding 2-cells by $\Epsilon^\alpha_A$ and 
$\eta^\alpha_A$, and also by $\Epsilon^*_A$ and $\eta^*_A$ the 2-cells related to conjoints of the inverse of $\alpha_0(A)$); 
\item the 2-cells 
$$\delta_{\alpha_1,f}=\quad
\scalebox{0.86}{
\bfig
 \putmorphism(-150,250)(1,0)[F(A)`\phantom{F(A)} `=]{500}1a
\put(-100,30){\fbox{$\eta^\alpha_A$}}
\putmorphism(380,250)(0,-1)[\phantom{Y_2}` `]{450}1l
\putmorphism(410,250)(0,-1)[\phantom{Y_2}` `\alpha_0(A)]{450}0l
\putmorphism(950,250)(0,-1)[\phantom{Y_2}` `]{450}1r
\putmorphism(930,250)(0,-1)[\phantom{Y_2}` `\alpha_0(B)]{450}0r
\putmorphism(350,250)(1,0)[F(A)`F(B)`F(f)]{600}1a
 \putmorphism(950,250)(1,0)[\phantom{F(A)}`G(B) `\alpha_1(B)]{600}1a
 \putmorphism(470,-200)(1,0)[`G(B)`G(f)]{500}1b
 \putmorphism(1060,-200)(1,0)[`G(B)`=]{500}1b
\putmorphism(1570,250)(0,-1)[\phantom{Y_2}``=]{450}1r
\put(530,10){\fbox{$(\alpha_0)_f$}}
\put(1280,30){\fbox{$\Epsilon^\alpha_B$}}
\putmorphism(-150,-200)(1,0)[F(A)`G(A) `\alpha_1(A)]{520}1a
\putmorphism(-150,250)(0,-1)[``=]{450}1l
\efig}
$$
in $\Ee$ for every 1h-cell $f:A\to B$ in $\Dd$;
\item the 2-cells 
$$(\alpha_1)^u=\quad
\scalebox{0.86}{
\bfig
 \putmorphism(-40,50)(1,0)[``=]{320}1b

\put(550,30){\fbox{$\delta_{\alpha_0,u}$}}
\putmorphism(-150,-400)(1,0)[F(\tilde A)`G(\tilde A) `\alpha_1(\tilde A)]{520}1a
\putmorphism(-150,50)(0,-1)[F(A')``=]{450}1l
\putmorphism(380,500)(0,-1)[\phantom{Y_2}` `F(u)]{450}1l
\putmorphism(950,500)(0,-1)[\phantom{Y_2}`G(A) `\alpha_0(A)]{450}1r
\putmorphism(380,50)(0,-1)[F(\tilde A)` `\alpha_0(\tilde A)]{450}1r
\putmorphism(350,500)(1,0)[F(A)`F(A)`=]{600}1a
 \putmorphism(950,500)(1,0)[\phantom{F(A)}`G(A) `\alpha_1(A)]{650}1a
 \putmorphism(1060,50)(1,0)[`G(A)`=]{500}1b

\putmorphism(1570,500)(0,-1)[\phantom{Y_2}``=]{450}1r
\put(1250,260){\fbox{$\Epsilon^\alpha_A$}}
\putmorphism(480,-400)(1,0)[`G(\tilde A) `=]{500}1a
\putmorphism(950,50)(0,-1)[\phantom{Y_2}``G(u)]{450}1r
\put(0,-160){\fbox{$\eta^\alpha_{\tilde A}$}}
\efig}
$$
in $\Ee$ for every 1v-cell $u:A\to \tilde A$ in $\Dd$;
\item the inverses of the 2-cells $\delta_{\alpha_1,f}$ are given by 
$$\delta_{\alpha_1,f}^{-1}=\quad
\scalebox{0.86}{
\bfig
 \putmorphism(-150,250)(1,0)[F(A)`\phantom{F(A)} `\alpha_1(A)]{500}1a
\put(-120,30){\fbox{$\Epsilon^*_A$}}
\putmorphism(400,250)(0,-1)[\phantom{Y_2}` `]{450}1l
\putmorphism(430,250)(0,-1)[\phantom{Y_2}` `\alpha_0(A)^{-1}]{450}0l
\putmorphism(950,250)(0,-1)[\phantom{Y_2}` `]{450}1r
\putmorphism(930,250)(0,-1)[\phantom{Y_2}` `\alpha_0(B)^{-1}]{450}0r
\putmorphism(350,250)(1,0)[G(A)`G(B)`G(f)]{600}1a
 \putmorphism(950,250)(1,0)[\phantom{F(A)}`G(B) `=]{600}1a
 \putmorphism(470,-200)(1,0)[`F(B)`F(f)]{500}1b
 \putmorphism(1060,-200)(1,0)[`G(B).`\alpha_1(B)]{500}1b
\putmorphism(1570,250)(0,-1)[\phantom{Y_2}``=]{450}1r
\put(530,10){\fbox{$(\alpha_0)_f^{-1}$}}
\put(1300,30){\fbox{$\eta^*_B$}} 
\putmorphism(-150,-200)(1,0)[F(A)`F(A) `=]{520}1a
\putmorphism(-150,250)(0,-1)[``=]{450}1l
\efig}
$$
\end{enumerate}
\item The 1h-cells $\alpha_1(A)$ are adjoint equivalence 1-cells in $\HH(\Ee)$. 
\item The data from 1. define a horizontal equivalence $\alpha_1:F\Rightarrow G:\Dd\to\Ee$ and a  
pseudonatural equivalence $\HH(\alpha_1):\HH(F)\Rightarrow \HH(G):\HH(\Dd)\to\HH(\Ee)$. 
\end{enumerate}
\end{prop}

\begin{proof}
The first part, up to invertibility, is proved in \cite[Proposition 4.1]{Fem:Fil}, the invertibility of $\delta_{\alpha_1,f}$ is proved in 
\cite[Theorem 4.6]{Shul}. Observe that a conjoint of $\alpha_0(A)^{-1}$ is $\alpha_1(A)$ by \leref{Shul1}. The proof of invertibility uses the identities \equref{fraction rules}. Adjoint equivalences follow from \leref{Shul}. Now the third claim is obvious. 
\qed\end{proof}

The above proposition in particular holds for invertible vertical strict transformations (recall \deref{vlt}) whose 1v-cell components 
have companions. As such transformations will be crucial in this work, we introduce some terminology. 

\begin{defn} 
Those 1v-cells that have companions we will call {\em companion-liftable} (or shortly {\em liftable}) 1v-cells. \\
Those vertical transformations all of whose 1v-cell components are liftable we will call {\em companion-liftable} (or shortly {\em liftable}) 
vertical transformations. \\
A horizontal pseudonatural equivalence $\alpha_1$ obtained in \prref{lifting 1v to equiv} from an invertible liftable vertical strict transformation $\alpha_0$ we will call a {\em companion-lift} of $\alpha_0$. 
\end{defn}

In \cite[Theorem 4.6]{Shul} it is shown that the above described assignment of bicategorical pseudonatural equivalences to invertible liftable vertical strict transformations, {\em i.e.} a companion-lifting of such vertical transformations, is functorial. 
\leref{Shul} and the above proposition, together with a reasoning that we will highlight in \prref{essence} further below, were crucial pieces on which relies the proof of Shulman's \cite[Theorem 5.1]{Shul}. Let us recall what it means for a double category to be monoidal due to Shulman, in our terminology. 

\begin{defn} \cite[\text{Definition 2.9}]{Shul} \delabel{Shul} \\ 
A {\em monoidal double category} is a double category $\Dd$ equipped with pseudodouble functors $\ot:\Dd\times\Dd\to\Dd$ and $I:*\to\Dd$ and 
invertible vertical strict transformations 
$$\alpha: \ot\comp(\Id\times\ot) \stackrel{\iso}{\to} \ot\comp(\ot\times\Id) $$
$$\lambda: \ot\comp(I\times\Id) \stackrel{\iso}{\to} \Id $$
$$\rho: \ot\comp(\Id\times I) \stackrel{\iso}{\to} \Id $$ 
satisfying the pentagonal and three triangular axioms (via identity vertical modifications). 
\end{defn} 

In any 2-category with finite products there is a notion of a pseudomonoid. Thus, a monoidal double category in other words is a pseudomonoid in the 2-category of double categories, pseudofunctors and vertical strict transformations. 

We recall now \cite[Theorem 5.1]{Shul}. ``Isofibrancy'' there supposes existence of companions for {\em all} 1v-cells in $\Dd$, 
whereas 
one indeed needs only existence of 1v-cell components of the monoidal structure $(\alpha,\lambda,\rho)$ of $\Dd$. We give an 
accordingly adapted formulation. 

\begin{thm} \cite[Theorem 5.1]{Shul} \thlabel{Shulman} \\
If the monoidal constraints $(\alpha,\lambda,\rho)$ of a monoidal pseudodouble category $\Dd$ are liftable, then its underlying horizontal bicategory $\HH(\Dd)$ is monoidal. 
\end{thm}

The proof of this theorem used a remarkable tool that we show in details in the next subsection. It will also show very important for our work, in particular in \seref{Kleisli}.

\subsection{Lifting of modifications} \sslabel{liftings}

In \prref{lifting 1v to equiv} we saw how, assuming the existence of companions, a (invertible) vertical strict transformation induces a 
(pseudonatural) horizontal transformation. For this we can say that vertical (strict) transformations lift to horizontal ones. In this subsection we will prove that a similar occurrence happens for modifications and axioms that they obey. 

\medskip

We start by recalling some more properties of companions.

\begin{lma} \cite[Lemma 3.8]{Shul} \lelabel{teta}  \\ 
Between two companions $\hat u:A\to B$ and $\hat u':A\to B$ of a 1v-cell $u:A\to B$ in $\Dd$ there is a unique globular isomorphism  
$\theta$ such that 
\begin{equation} \eqlabel{teta-property}
\bfig
 \putmorphism(-150,410)(1,0)[A`A  `=]{400}1a
\putmorphism(-160,400)(0,-1)[\phantom{Y_2}`A `=]{380}1l
\putmorphism(270,400)(0,-1)[\phantom{Y_2}`B `u]{380}1r
\putmorphism(-160,50)(0,-1)[\phantom{Y_2}`A`=]{430}1l
\putmorphism(270,50)(0,-1)[\phantom{Y_2}`B`=]{430}1r
\put(-20,240){\fbox{$\eta$}}
\putmorphism(-150,20)(1,0)[``\hat u]{380}1a
\put(-70,-160){\fbox{$\theta_{\hat u, \hat u'}$ }}
\putmorphism(-160,-350)(0,-1)[\phantom{Y_2}``u]{400}1l
\putmorphism(-150,-350)(1,0)[``\hat u']{380}1a
\putmorphism(270,-350)(0,-1)[\phantom{Y_2}`B`=]{400}1r
\putmorphism(-180,-750)(1,0)[B` `=]{440}1b
\put(-20,-560){\fbox{$\Epsilon'$ }}
\efig\quad
=\quad
\bfig
\putmorphism(-110,120)(1,0)[``=]{330}1a
\putmorphism(-160,150)(0,-1)[A`B`u]{430}1l
\putmorphism(270,150)(0,-1)[A`B.`u]{430}1r
\putmorphism(-110,-250)(1,0)[``=]{330}1a
\put(-20,-90){\fbox{$\Id^u$}}
\efig
\end{equation} 
It is given by 
\begin{equation} \eqlabel{teta}
\theta_{\hat u, \hat u'}\quad=\quad
\bfig
 \putmorphism(-150,210)(1,0)[A`A  `=]{450}1a
\putmorphism(-160,200)(0,-1)[\phantom{Y_2}`A `=]{380}1l
\putmorphism(300,200)(0,-1)[\phantom{Y_2}`B `u]{380}1l
 \putmorphism(300,210)(1,0)[\phantom{Y}`B  `\hat u']{450}1a
\putmorphism(730,200)(0,-1)[\phantom{Y_2}`B. `=]{380}1r
\putmorphism(-150,-200)(1,0)[``\hat u]{450}1a
\putmorphism(350,-200)(1,0)[``=]{330}1a
\put(-20,0){\fbox{$\eta$}}
\put(420,0){\fbox{$\Epsilon'$}}
\efig
\end{equation}
\end{lma}


\begin{lma} \lelabel{Lemma 4.8}
Let $\alpha$ be a vertical strict transformation with two different choices of companions for its 1v-cell components, giving rise to two different horizontal transformations $\hat\alpha$ and $\hat\alpha'$, as in \prref{lifting 1v to equiv}. The isomorphisms $\theta$ from 
\equref{teta} make an invertible modification $\hat\alpha\Rrightarrow\hat\alpha'$.
\end{lma}

\begin{proof}
This is proved in the proof  of \cite[Lemma 4.8]{Shul}. 
\qed\end{proof}

\medskip


\begin{lma}  \lelabel{horiz is teta}  \cite[Lemma 4.10]{Shul} \\
Let $\alpha:F\Rightarrow G$ be a vertical strict transformation between pseudodouble functors acting between double categories 
$\Dd\to\Ee$ and suppose that $u:A\to B$ has a companion $\hat u$. Then the 2-cell component 
$$\hat\alpha_{\hat u}=\quad
\scalebox{0.86}{
\bfig
 \putmorphism(-150,250)(1,0)[F(A)`\phantom{F(A)} `=]{500}1a
\put(-40,30){\fbox{$\eta$}}
\putmorphism(380,250)(0,-1)[\phantom{Y_2}` `]{450}1l
\putmorphism(410,250)(0,-1)[\phantom{Y_2}` `\alpha(A)]{450}0l
\putmorphism(950,250)(0,-1)[\phantom{Y_2}` `]{450}1r
\putmorphism(930,250)(0,-1)[\phantom{Y_2}` `\alpha(B)]{450}0r
\putmorphism(350,250)(1,0)[F(A)`F(B)`F(\hat u)]{600}1a
 \putmorphism(950,250)(1,0)[\phantom{F(A)}`G(B) `\hat\alpha(B)]{600}1a
 \putmorphism(470,-200)(1,0)[`G(B)`G(\hat u)]{500}1b
 \putmorphism(1060,-200)(1,0)[`G(B)`=]{500}1b
\putmorphism(1570,250)(0,-1)[\phantom{Y_2}``=]{450}1r
\put(560,10){\fbox{$\alpha_{\hat u}$}}
\put(1240,30){\fbox{$\Epsilon$}}
\putmorphism(-150,-200)(1,0)[F(A)`G(A) `\hat\alpha(A)]{520}1b
\putmorphism(-150,250)(0,-1)[``=]{450}1l
\efig}
$$
of the induced horizontal transformation $\hat\alpha$ is equal to $\theta_{[F(\hat u)\vert\hat\alpha(B)], [\hat\alpha(A)\vert G(\hat u)]}$, 
and in particular, it is an isomorphism.  
\end{lma}


\bigskip

Let us recall the basic algebra of companions from \cite[Section 3]{Shul}. 

\begin{prop} \prlabel{alg}
Let $\Dd$ be a double category. 
\begin{enumerate}
\item An identity 1v-cell has the identity 1h-cell as a companion. 
\item If 1v-cells $u:A\to B$ and $v:B\to C$ have companions $\hat u$ and $\hat v$, then $\frac{u}{v}$ has a companion $[\hat u\vert\hat v]$. 
\item If $u$ has three companions $\hat u, \hat u', \hat u''$, then $\theta_{\hat u,\hat u''}=
\frac{\theta_{\hat u,\hat u'}}{\theta_{\hat u',\hat u''}}$. 
\item If 1v-cells $u:A\to B$ and $v:B\to C$ have companions $\hat u, \hat u'$ and $\hat v,\hat v'$, then $[\theta_{\hat u,\hat u'} \vert 
\theta_{\hat v,\hat v'}]=\theta_{[\hat u,\hat v],[\hat u',\hat v']}$. 
\item If $u$ has a companion $\hat u$, then $\theta_{\hat u, [\id\vert\hat u]}$ and $\theta_{\hat u, [\hat u\vert\id]}$ are equal to the unit constraints $\hat u\iso [\id\vert\hat u]$ and $\hat u\iso [\hat u\vert\id]$. 
\item Let $F:\Dd\to\Ee$ be a pseudodouble functor between double categories and 
assume that $u$ has companions $\hat u, \hat u'$ in $\Dd$. Then: \vspace{-0,2cm}
\begin{enumerate} [a)]
\item $F(u)$ has companions $F(\hat u)$ and $F(\hat u')$, and 
\item $\theta_{F(\hat u),F(\hat u')}=F(\theta_{\hat u,\hat u'})$ in $\Ee$. 
\end{enumerate}
\item Assume that $\Dd$ is a monoidal double category and 
that $u:A\to\tilde A$ and $v:B\to\tilde B$ have companions $\hat u,\hat u'$ and $\hat v,\hat v'$. Then:  \vspace{-0,2cm}
\begin{enumerate} [a)]
\item $u\ot v$ has companions $\hat u\ot\hat v$ and $\hat u'\ot\hat v'$, and 
\item $\theta_{\hat u,\hat u'} \ot \theta_{\hat v,\hat v'} = \theta_{\hat u\ot\hat v,\hat u'\ot\hat v'}$. 
\end{enumerate}
%
%
\end{enumerate}
\end{prop}

We extract the essence of the mechanism used by Shulman which underlies the proof of his remarkable \cite[Theorem 5.1]{Shul} in the following proposition. Its proof is a reformulation of Shulman's one.

\begin{prop} \prlabel{essence}
Let 
$\omega$ be the identity vertical modification between two vertical composites of vertical strict transformations 
$$\omega: \,\, \threefrac{\alpha_1}{...}{\alpha_k} \,\,\, \Rrightarrow \,\,\, \threefrac{\beta_1}{...}{\beta_l} \,\,$$ 
which act between lax double functors $F\Rightarrow G:\Bb\to\Dd$ between double categories,  
so that all 1v-cell components $\alpha_1(A),...\alpha_k(A)$ and $\beta_1(A), ..., \beta_l(A)$ for 0-cells $A$ in $\Bb$, have companions in $\Dd$. 
Then: 
\begin{enumerate}
\item $\omega$ induces an invertible  
horizontal modification $\hat\omega$ between the two (vertical compositions of the) induced horizontal natural transformations; 
\item the assignment between 2-cell components of $\omega$ and $\hat\omega$ is invertible; 
\item if $\omega_1,...,\omega_m$ are 
vertical modifications with the above characteristics, then any sensible equation formed by their 2-cell components  
$\hat\omega_1(A),...,\hat\omega_m(A)$, for any 0-cell $A$ in $\Bb$, 
holds true. 
\end{enumerate}
\end{prop}

\begin{proof}
We did not specify it in the statement, but it is also possible that $\Dd$ be a monoidal double category, and that component 1v-cells of the 
$\omega(A)$'s are tensor products of other 1v-cells. This possibility will be included in the proof. 

Let $A$ be fixed. Observe that $\omega(A)$ presents the identity between its domain and codomain composite 1v-cells for every $A$. 
We will abuse notation by writing $\omega$ both for the modification and its component 2-cells $\omega(A)$. 

Let $\omega:u\Rightarrow v$ denote $\omega$ as a vertically globular 2-cell between its composite 1v-cell components $u$ and $v$. 
Then $u$ and $v$ are 1v-cell components of the composite vertical strict transformations. 
The (clearly invertible) assignment $\omega\mapsto\hat\omega$ is given by 
\begin{equation} \eqlabel{2-cell lift}
\scalebox{0.9}{
\bfig
\putmorphism(-150,170)(1,0)[A`A`=]{460}1a
\putmorphism(-150,-180)(1,0)[B`B`=]{460}1a
\putmorphism(-150,170)(0,-1)[\phantom{Y_2}``u]{350}1l
\putmorphism(310,170)(0,-1)[\phantom{Y_2}``v]{350}1r
\put(20,0){\fbox{$\omega$}}
\efig}
\quad\mapsto\quad 
\scalebox{0.9}{
\bfig
\putmorphism(-150,170)(1,0)[A`B`\hat v]{460}1a
\putmorphism(-150,-180)(1,0)[A`B`\hat u]{460}1a
\putmorphism(-150,170)(0,-1)[\phantom{Y_2}``=]{350}1l
\putmorphism(310,170)(0,-1)[\phantom{Y_2}``=]{350}1r
\put(20,0){\fbox{$\hat\omega$}}
\efig}
\quad=\quad
\scalebox{0.86}{
\bfig
 \putmorphism(-150,250)(1,0)[A`A `=]{480}1a
\put(0,30){\fbox{$\eta_u$}}
\putmorphism(340,250)(0,-1)[\phantom{Y_2}`B `]{450}1l
\putmorphism(360,250)(0,-1)[\phantom{Y_2}` `u]{450}0l
\putmorphism(810,250)(0,-1)[A`B `]{450}1r
\putmorphism(790,250)(0,-1)[` `v]{450}0r
\putmorphism(1270,250)(0,-1)[B`B`=]{450}1r
\putmorphism(350,250)(1,0)[``=]{420}1a
 \putmorphism(730,250)(1,0)[\phantom{F(A)}` `\hat v]{510}1a
 \putmorphism(-150,-200)(1,0)[A` `\hat u]{460}1a
\putmorphism(360,-200)(1,0)[``=]{440}1b
 \putmorphism(820,-200)(1,0)[``=]{440}1b
\put(530,10){\fbox{$\omega$}}
\put(980,30){\fbox{$\Epsilon_v$}}
\putmorphism(-150,250)(0,-1)[``=]{450}1l
\efig}
\end{equation} 
where $\eta_u$ is given by diagonally composing $\eta$'s for every component 1v-cell making the composite $u$, and the same for 
$\Epsilon_v$. We used here part 2. of \prref{alg}. (If some component 1v-cell in $u$ or $v$ is of the form $x\ot y$ for some 1v-cells $x,y$, then 
 by part 7a) of \prref{alg} we have $\widehat{x\ot y}=\hat x\ot\hat y$.) 

Since $u=v$ via $\omega$ by assumption, then $\hat\omega: \hat u\Rightarrow\hat v$ satisfies \equref{teta-property}, it is a $\theta$ 2-cell  and an isomorphism by \leref{teta} with inverse given by 
$$
\scalebox{0.86}{
\bfig
 \putmorphism(-150,250)(1,0)[A`A `=]{480}1a
\put(0,30){\fbox{$\eta_v$}}
\putmorphism(340,250)(0,-1)[\phantom{Y_2}`B `]{450}1l
\putmorphism(360,250)(0,-1)[\phantom{Y_2}` `v]{450}0l
\putmorphism(810,250)(0,-1)[A`B `]{450}1r
\putmorphism(790,250)(0,-1)[` `u]{450}0r
\putmorphism(1270,250)(0,-1)[B`B.`=]{450}1r
\putmorphism(350,250)(1,0)[``=]{420}1a
 \putmorphism(730,250)(1,0)[\phantom{F(A)}` `\hat u]{510}1a
 \putmorphism(-150,-200)(1,0)[A` `\hat v]{460}1a
\putmorphism(360,-200)(1,0)[``=]{440}1b
 \putmorphism(820,-200)(1,0)[``=]{440}1b
\put(470,10){\fbox{$\omega^{inv}$}}
\put(980,30){\fbox{$\Epsilon_u$}}
\putmorphism(-150,250)(0,-1)[``=]{450}1l
\efig}
$$ 
On the other hand, by \leref{Lemma 4.8} it 
induces an invertible horizontal modification between the composite horizontal natural transformations. 

If we are given any equation relating horizontally globular 2-cells $\hat\omega_1,...,\hat\omega_m$ induced from  vertically globular 2-cells 
$\omega_1,...,\omega_m$ with the above characteristics, we have as above that every $\hat\omega_i, i=1,2..,m$ is given by a $\theta$ 2-cell. 
By properties 3. and 4. of \prref{alg} we have that (combinations of horizontal and vertical) composites of 2-cells $\hat\omega_1,...,\hat\omega_m$ making the two sides of the equation are both isomorphism 2-cells of the sort of $\theta$ between their common domain and codomain. 
The uniqueness of $\theta$ implies that the equation in question holds.  
\qed\end{proof}

\bigskip

We will also need the following variation of the above claim. In it we 
treat $\omega$ 
as component 2-cells of the modification in question.

\begin{prop} \prlabel{omega*}
Let $\omega$ be as in \prref{essence} and assume moreover that at least two of 1v-cells in one of its vertical edges are non-trivial. Then: 
\begin{enumerate}
\item $\omega$ induces a 2-cell $\omega^*$ defining a modification in the sense of \deref{modif-hv};
\item there is a one-to-one correspondence between 2-cells $\omega^*$ and $\hat\omega$ from \prref{essence}, so that $\omega^*$ induces a modification if and only if so does $\hat\omega$;
\item if $\omega^*_1,...,\omega^*_m$ are 2-cells induced by vertically globular 2-cells $\omega_1,...,\omega_m$ with the characteristics as $\omega$ of this proposition, then any sensible equation formed by the 2-cells $\omega^*_1,...,\omega^*_m$ holds true. 
\end{enumerate}
\end{prop}

\begin{proof}
We illustrate the assignment from point 1) by an example: 
\begin{equation} \eqlabel{omega-inv}
\scalebox{0.9}{
\bfig
 \putmorphism(-90,500)(1,0)[A` A`=]{440}1a
\putmorphism(350,500)(0,-1)[` B'`v']{400}1r
\putmorphism(-90,-300)(1,0)[C`C `=]{440}1a
\putmorphism(350,110)(0,-1)[\phantom{Y_2}``u']{400}1r
\putmorphism(-100,110)(0,-1)[\phantom{Y_2}``v]{400}1l
\putmorphism(-100,500)(0,-1)[` B`u]{400}1l
\put(80,70){\fbox{$\omega$}}
\efig}
\quad\mapsto\quad 
\scalebox{0.9}{
\bfig
\putmorphism(-150,170)(1,0)[A`B'`\hat v']{460}1a
\putmorphism(-150,-180)(1,0)[B`C`\hat v]{460}1a
\putmorphism(-150,170)(0,-1)[\phantom{Y_2}``u]{350}1l
\putmorphism(310,170)(0,-1)[\phantom{Y_2}``u']{350}1r
\put(0,0){\fbox{$\omega^*$}}
\efig}
\quad=\quad
\scalebox{0.86}{
\bfig
 \putmorphism(-120,50)(1,0)[B`B`=]{480}1b
\put(550,30){\fbox{$\omega$}}
\putmorphism(-120,-400)(1,0)[` `\hat v]{440}1a

\putmorphism(-150,50)(0,-1)[`B`=]{450}1l
\putmorphism(380,500)(0,-1)[A``u]{450}1l
\putmorphism(850,500)(0,-1)[` `v']{450}1r
\putmorphism(380,50)(0,-1)[`C `v]{450}1r
\putmorphism(400,500)(1,0)[`A`=]{440}1a
 \putmorphism(900,500)(1,0)[` B'`\hat v']{500}1a
 \putmorphism(880,50)(1,0)[``=]{450}1b

\putmorphism(1370,500)(0,-1)[\phantom{Y_2}`B'`=]{450}1r
\put(1100,260){\fbox{$\Epsilon$}}
\putmorphism(400,-400)(1,0)[` `=]{410}1a
\putmorphism(850,50)(0,-1)[B'` C ` u']{450}1r
\put(60,-160){\fbox{$\eta$}}
\efig}
\end{equation} 
(the assumption that at least two 1v-cells in a same edge of $\omega$ are non-trivial is there only to assure that we get a non-globular 2-cell 
$\omega^*$). 
The proofs of all the three claims are similar and straightforward. We will comment on the proof of the third claim. 
The invertible assignment in part 2. is given by: 
$$\omega^*\mapsto [\eta\vert\omega^*\vert\Epsilon]\qquad\text{and}\qquad
\hat\omega\mapsto \threefrac{[\Id\vert\eta]}{\hat\omega}{[\Epsilon\vert\Id]}.$$
For the third part, write any sensible equation $E^*$ formed by $\omega^*_1,...,\omega^*_m$ in terms of $\hat\omega_1,...,\hat\omega_m$ using part 2. to obtain an equation $\tilde E$. Then $E^*$ holds true if and only if $\tilde E$ does. 
Observe that both vertical and horizontal composition of some $\omega^*_i$ and $\omega^*_j$ obtains the form 
$\fourfrac{[\Id\vert\Id\vert\eta]}{[\Id\vert\hat\omega_j]}{[\hat\omega_i\vert\Id]}{[\Epsilon\vert\Id\vert\Id]}$ in $\tilde E$ (or its symmetric diagonal version). Consequently, in the equation $\tilde E$ ``middle $\eta$'s and $\Epsilon$'s'' will cancel out and $\tilde E$ gets the form 
$\threefrac{[\Id\vert\eta]}{\hat E}{[\Epsilon\vert\Id]}$, {\em i.e.} $[\Id\vert\eta]$ and $[\Epsilon\vert\Id]$ are composed to both sides of $\hat E$, where $\hat E$ is an equation as in the third claim of \prref{essence}. We know that $\hat E$ holds true, hence we have the proof. 
\qed\end{proof}

For liftable $u,u'$ we label the invertible assignment 
\begin{equation}\eqlabel{square-lift}
\scalebox{0.9}{
\bfig
\putmorphism(-150,170)(1,0)[A`B`v]{460}1a
\putmorphism(-150,-180)(1,0)[C`D`v']{460}1a
\putmorphism(-150,170)(0,-1)[\phantom{Y_2}``u]{350}1l
\putmorphism(310,170)(0,-1)[\phantom{Y_2}``u']{350}1r
\efig}
\quad\mapsto\quad 
\scalebox{0.86}{
\bfig
 \putmorphism(-150,200)(1,0)[A`B `v]{480}1a
\putmorphism(340,200)(0,-1)[\phantom{Y_2}` `u]{400}1l
\putmorphism(780,200)(0,-1)[\phantom{Y_2}` `=]{400}1r
\putmorphism(350,200)(1,0)[`D`\hat u']{420}1a
 \putmorphism(-150,-200)(1,0)[A` C`\hat u]{460}1a
\putmorphism(360,-200)(1,0)[`D`v']{420}1a
\putmorphism(-150,200)(0,-1)[``=]{400}1l
\efig} 
\end{equation}
for future reference.

\medskip

Under assumption on existence of companions, invertibility of vertically globular 2-cells does not imply in general invertibility of squares 
obtained by the above described assignment $\omega\mapsto\omega^*$. However, we have:

\begin{lma} \lelabel{inv-spec}
If the 1v-cells $u$ and $u'$ in $\omega$ in \equref{omega-inv} and the 2-cell $\omega$ itself are invertible, an inverse of $\omega^*$ is given by 
$$
\scalebox{0.9}{
\bfig
\putmorphism(-150,170)(1,0)[B`C`\hat v]{460}1a
\putmorphism(-140,-180)(1,0)[``\hat v']{410}1b
\putmorphism(-150,170)(0,-1)[\phantom{Y_2}`A`u^{-1}]{350}1l
\putmorphism(310,170)(0,-1)[\phantom{Y_2}`B'`(u')^{-1}]{350}1r
\put(-70,-20){\fbox{$(\omega^*)^{-1}$}}
\efig}
\quad=\quad
\scalebox{0.9}{
\bfig
 \putmorphism(-90,900)(1,0)[B` B`=]{440}1a
\putmorphism(-100,900)(0,-1)[` `u^{-1}]{400}1l
\putmorphism(350,900)(0,-1)[` `u^{-1}]{400}1r
\put(0,670){\fbox{$\Id^{u^{-1}}$}}
\put(500,510){\fbox{$1$}}

 \putmorphism(370,900)(1,0)[`B `=]{440}1a
\putmorphism(790,900)(0,-1)[` `=]{800}1r

 \putmorphism(-90,500)(1,0)[` `=]{420}1a 
\putmorphism(-100,500)(0,-1)[A`B' `v']{400}1l
\putmorphism(350,500)(0,-1)[A`B `u]{400}1r
\putmorphism(-90,-270)(1,0)[` `=]{420}1a
\putmorphism(350,110)(0,-1)[\phantom{Y_2}`C`v]{380}1r
\putmorphism(-100,110)(0,-1)[\phantom{Y_2}`C`u']{380}1l
\put(40,70){\fbox{$\omega^{-1}$}}

\putmorphism(-500,110)(0,-1)[\phantom{Y_2}``=]{770}1l
\put(-350,-300){\fbox{$1$}}
\putmorphism(-100,-260)(0,-1)[\phantom{Y_2}``(u')^{-1}]{400}1l
\putmorphism(350,-260)(0,-1)[\phantom{Y_2}``(u')^{-1}]{400}1r
\put(-30,-500){\fbox{$\Id^{(u')^{-1}}$}}

\putmorphism(790,110)(0,-1)[\phantom{Y_2}``v]{380}1l
\putmorphism(1190,110)(0,-1)[C``=]{380}1r
\putmorphism(360,100)(1,0)[` `=]{420}1b
\putmorphism(800,100)(1,0)[B` `\hat v]{360}1a
\putmorphism(800,-270)(1,0)[C` C.`=]{410}1a
\put(940,-80){\fbox{$\Epsilon$}}

\putmorphism(-490,-650)(1,0)[B' ` B'`=]{380}1b
\putmorphism(-90,-650)(1,0)[`B' `=]{440}1b
\putmorphism(-470,100)(1,0)[` `=]{340}1a
\putmorphism(-890,100)(1,0)[` `\hat{v'}]{360}1b
\putmorphism(-890,500)(1,0)[` `=]{360}1a 
\putmorphism(-500,500)(0,-1)[A`B' `v']{400}1r
\putmorphism(-900,500)(0,-1)[A`A `=]{400}1l
\put(-750,260){\fbox{$\eta$}}
\efig}
$$
\end{lma}

We state for the record that although invertibility of squares implies invertibility of horizontally globular 2-cells (similar to the assignment $\omega^*\mapsto\hat\omega$ from \prref{omega*}), the converse holds under assumptions similar as in the above lemma.

\section{Premonoidal double categories and central cells}

In this section we are going to introduce a double categorical version of the notion of premonoidal bicategories 
introduced in \cite{HF1}. For reader's convenience we recall the definition of a premonoidal bicategory.

\begin{defn} \cite[Definition 6]{HF1} \delabel{premon-bicat} \\
A {\em premonoidal} bicategory is a binoidal bicategory $(\B, \ltimes,\rtimes)$ with a unit object $I\in\B$ and the following data:
\begin{enumerate}
\item pseudonatural equivalences $\lambda: I\ltimes -\Rightarrow\Id, \rho: -\rtimes I\Rightarrow\Id$ and 
$$\alpha_{-,B,C}: (-\rtimes B)\rtimes C\Rightarrow -\rtimes(B\bowtie C)$$
$$\alpha_{A,-,C}: (A\ltimes -)\rtimes C\Rightarrow A\ltimes(-\rtimes C)$$
$$\alpha_{A,B,-}: (A\bowtie B)\ltimes -\Rightarrow A\ltimes(B\ltimes -)$$
for every $A,B,C\in\Bb$, such that all the 1-cell components of these five equivalences are central 1-cells;  
\item invertible 2-cells $p_{A,B,C,D}, m_{A,B}, l_{A,B}, r_{A,B}$ such that these form modifications in each argument 
in the form of the modifications $\pi, \mu, \lambda, \rho$ from \cite{GPS}. 
\end{enumerate}
The above data is subject to the same equations between 2-cells as in a monoidal bicategory. 
\end{defn}

In our definition of a premonoidal double category we include 24 axioms which permit to describe the interrelation between the three  transformations $\alpha$. We start the section by introducing the previously necessary notions of ``binoidality'' and ``centrality'' of cells 
(we will use them in \ssref{quasi}), and we also study the relation of companions and centrality.  


\subsection{Binoidal double categories and central cells}  \sslabel{binoidal}  

We generalize the definitions of a binoidal bicategory and central 1- and 2-cells in it from \cite{HF1} to double categories. 
We will differentiate left and right central 1-cells.

\begin{defn}
We say that a double category $\Bb$ is {\em binoidal} if for all 0-cells $A,B\in\Bb$ there 
are pseudodouble functors $A\ltimes -$ and $-\rtimes B$ acting $\Bb\to\Bb$ and such that $A\ltimes B=
A\rtimes B=:A\bowtie B$. 
\end{defn}

\begin{defn}
Assume $\Bb$ is a binoidal double category. 
\begin{itemize}
\item A 1h-cell $f:A\to A'$ in $\Bb$ is said to be {\em left central} if there is a horizontal pseudonatural transformation 
$f\ltimes-:A\ltimes -\to A'\ltimes-$ such that $f\ltimes B=f\rtimes B$ for all $B\in\Bb$. 
Likewise, a 1h-cell $f:A\to A'$ in $\Bb$ is said to be {\em right central} if there is a horizontal pseudonatural transformation 
$-\rtimes f:-\ltimes A\to -\rtimes A'$ such that $B\rtimes f=B\ltimes f$ for all $B\in\Bb$. \\ 
A 1h-cell is said to be {\em central} if it is both left and right central. 
\item A 1v-cell $v:A\to\tilde A$ in $\Bb$ is said to be {\em left central} if there is a vertical pseudonatural transformation 
$v\ltimes-:A\ltimes -\to \tilde A\ltimes-$ such that $v\ltimes B=v\rtimes B$ for all $B\in\Bb$. 
Likewise, a 1v-cell $v:A\to\tilde A$ in $\Bb$ is said to be {\em right central} if there is a vertical pseudonatural transformation 
$-\rtimes v:-\rtimes A\to -\rtimes\tilde A$ such that $B\ltimes v=B\rtimes v$ for all $B\in\Bb$. \\ 
A 1v-cell is said to be {\em central} if it is both left and right central.
\end{itemize}
\end{defn}

For reader's convenience we write down the 2-cell components of the pseudonatural transformations in play. 
For a horizontal pseudonatural transformation $f\ltimes-:A\ltimes -\to A'\ltimes-$, a 1h-cell $g:B\to B'$  and 
a 1v-cell $v:A\to\tilde A$ the 2-cell components $f\ltimes-\vert_g$ and $f\ltimes-\vert_v$ of the oplax transformation structure of 
$f\ltimes-$ have the form of the left diagrams below (for the lax transformation structure of $f\ltimes-$ the first 2-cell 
component is different: it is a 2-cell obtained by reading the same upper left diagram from bottom to top). 
Likewise, the 2-cell components $-\rtimes f\vert_g$ and $-\rtimes f\vert_v$ of the oplax transformation structure of 
$-\rtimes f:-\rtimes A\to -\rtimes A'$  has the form of the right diagrams below (and the differing 2-cell component for the lax structure 
is obtained by reading the upper right diagram from bottom to top). 
\begin{equation} \eqlabel{f lr}
\scalebox{0.86}{
\bfig
 \putmorphism(-170,500)(1,0)[A\ltimes B`A\ltimes B' `A\ltimes g]{620}1a
 \putmorphism(450,500)(1,0)[\phantom{A\ot B}`A'\ltimes B' `f\ltimes B']{680}1a
 \putmorphism(-150,50)(1,0)[A\ltimes B`A'\ltimes B `f\ltimes B]{600}1a
 \putmorphism(450,50)(1,0)[\phantom{A\ot B}`A'\ltimes B' `A'\ltimes g]{680}1a

\putmorphism(-180,500)(0,-1)[\phantom{Y_2}``=]{450}1r
\putmorphism(1100,500)(0,-1)[\phantom{Y_2}``=]{450}1r
\put(350,260){\fbox{$f\ltimes-\vert_g$}}
\efig}
\qquad\qquad
\scalebox{0.86}{
\bfig
 \putmorphism(-170,500)(1,0)[B\rtimes A`B'\rtimes A `g\rtimes A]{620}1a
 \putmorphism(450,500)(1,0)[\phantom{A\ot B}`B'\rtimes A' `B'\rtimes f]{680}1a
 \putmorphism(-150,50)(1,0)[B\rtimes A`B\rtimes A' `B\rtimes f]{600}1a
 \putmorphism(450,50)(1,0)[\phantom{A\ot B}`B'\rtimes A' `g\rtimes A']{680}1a

\putmorphism(-180,500)(0,-1)[\phantom{Y_2}``=]{450}1r
\putmorphism(1100,500)(0,-1)[\phantom{Y_2}``=]{450}1r
\put(350,260){\fbox{$-\rtimes f\vert_g$}}
\efig}
\end{equation}

$$
\scalebox{0.86}{
\bfig
 \putmorphism(-150,50)(1,0)[A\ltimes B`A'\ltimes B `f\ltimes B]{640}1a
\putmorphism(-150,-400)(1,0)[A\ltimes\tilde  B`A'\ltimes\tilde  B `f\ltimes \tilde B]{640}1a
\putmorphism(-160,50)(0,-1)[\phantom{Y_2}``A\ltimes v]{450}1l
\putmorphism(450,50)(0,-1)[\phantom{Y_2}``A'\ltimes v]{450}1r
\put(-40,-170){\fbox{$f\ltimes-\vert_v$}}
\efig}
\qquad\qquad\qquad
\scalebox{0.86}{
\bfig
 \putmorphism(-150,50)(1,0)[B\rtimes A`B\rtimes A' `B\rtimes f]{640}1a
\putmorphism(-150,-400)(1,0)[\tilde B\rtimes A ` \tilde B\rtimes A' ` \tilde B\rtimes f]{640}1a
\putmorphism(-160,50)(0,-1)[\phantom{Y_2}``v\rtimes A]{450}1l
\putmorphism(450,50)(0,-1)[\phantom{Y_2}``v\rtimes A']{450}1r
\put(-40,-170){\fbox{$-\rtimes f\vert_v$}}
\efig}
$$
Similarly, for a 1v-cell $u:B\to\tilde B$ the 2-cell components of oplax transformation structures of 
$v\ltimes-:A\ltimes -\to \tilde A\ltimes-$ and $-\rtimes v:-\rtimes A\to -\rtimes\tilde A$ have the form as below (and for the lax 
structures read the respective diagrams from right to left). 
\begin{equation} \eqlabel{v lr}
\scalebox{0.86}{
\bfig
 \putmorphism(-150,500)(1,0)[A\ltimes B`A\ltimes B `=]{600}1a
\putmorphism(-180,500)(0,-1)[\phantom{Y_2}` `A\ltimes u]{450}1l 
\putmorphism(-250,500)(0,-1)[\phantom{Y_2}` A\ltimes\tilde B`]{450}0l %
\put(-80,50){\fbox{$v\ltimes-\vert_{u}$}}
\putmorphism(-170,-400)(1,0)[\tilde A\ltimes \tilde B` `=]{480}1a
\putmorphism(-180,50)(0,-1)[\phantom{Y_2}``v\ltimes \tilde B]{450}1l
\putmorphism(450,50)(0,-1)[\phantom{Y_2}`\tilde A\ltimes \tilde B`\tilde A\ltimes u]{450}1l
\putmorphism(450,500)(0,-1)[\phantom{Y_2}`\tilde A\ltimes B `v\ltimes B]{450}1l
\efig}
\qquad\qquad
\scalebox{0.86}{
\bfig
 \putmorphism(-150,500)(1,0)[B\rtimes A`B\rtimes A `=]{600}1a
\putmorphism(-180,500)(0,-1)[\phantom{Y_2}` `u\rtimes A]{450}1l 
\putmorphism(-250,500)(0,-1)[\phantom{Y_2}` \tilde B\rtimes A`]{450}0l %
\put(-80,50){\fbox{$-\rtimes v\vert_{u}$}}
\putmorphism(-170,-400)(1,0)[\tilde B\rtimes \tilde A` `=]{480}1a
\putmorphism(-180,50)(0,-1)[\phantom{Y_2}``\tilde B\rtimes v]{450}1l
\putmorphism(450,500)(0,-1)[\phantom{Y_2}`B\rtimes \tilde A `B\rtimes v]{450}1l
\putmorphism(450,50)(0,-1)[\phantom{Y_2}`\tilde B\rtimes \tilde A  ` u\rtimes \tilde A]{450}1l
\efig}
\end{equation}

$$
\hspace{0,5cm}
\scalebox{0.86}{
\bfig
 \putmorphism(-150,50)(1,0)[A\ltimes B`A\ltimes B' `A\ltimes f]{640}1a
\putmorphism(-150,-400)(1,0)[\tilde A\ltimes B ` \tilde A\ltimes B' ` \tilde A\ltimes f]{640}1a
\putmorphism(-160,50)(0,-1)[\phantom{Y_2}``v\ltimes B]{450}1l
\putmorphism(450,50)(0,-1)[\phantom{Y_2}``v\ltimes B']{450}1r
\put(-40,-170){\fbox{$v\ltimes -\vert_f$}}
\efig}
\qquad\quad
\scalebox{0.86}{
\bfig
 \putmorphism(-150,50)(1,0)[B\rtimes A`B'\rtimes A `f\rtimes A]{640}1a
\putmorphism(-150,-400)(1,0)[B\rtimes\tilde  A`B'\rtimes\tilde A `f\rtimes \tilde A]{640}1a
\putmorphism(-160,50)(0,-1)[\phantom{Y_2}``B\rtimes v]{450}1l
\putmorphism(450,50)(0,-1)[\phantom{Y_2}``B'\rtimes v]{450}1r
\put(-40,-170){\fbox{$-\rtimes v\vert_f$}}
\efig}
$$
Observe that for a left central 1h-cell $f$ the requirement $f\ltimes B=f\rtimes B$ makes that the image 
of the pseudodouble functor $-\rtimes B$ at $f$ is related to the images of the pseudodouble functor $A\ltimes-$ at $g$ and $v$ 
via the corresponding 2-cell components of $f\ltimes-$, and similarly for other central 1-cells.

\begin{defn} \delabel{central-2}
In a binoidal double category $\Bb$ a 2-cell 
$\scalebox{0.86}{
\bfig
\putmorphism(-150,50)(1,0)[A` A'`f]{400}1a
\putmorphism(-150,-270)(1,0)[\tilde A`\tilde A' `f']{400}1b
\putmorphism(-150,50)(0,-1)[\phantom{Y_2}``v]{320}1l
\putmorphism(250,50)(0,-1)[\phantom{Y_2}``v']{320}1r
\put(0,-140){\fbox{$a$}}
\efig}$ 
in $\Bb$, where $f$ and $f'$ are left central 1h-cells and $v$ and $v'$ are left central 1v-cells (with respect to horizontal and vertical transformations $f\ltimes-, f'\ltimes$ and $v\ltimes-, v'\ltimes$, respectively), is said to be {\em left central} if there is a modification $a\ltimes-$ given on 0-cells $B\in\Bb$ by 2-cells $a\ltimes B=a\rtimes B$ that satisfy: 
$$
\scalebox{0.86}{
\bfig
 \putmorphism(-170,500)(1,0)[A\ltimes B`A\ltimes B' `A\ltimes g]{620}1a
 \putmorphism(450,500)(1,0)[\phantom{A\ot B}`A'\ltimes B' `f\ltimes B']{680}1a
 \putmorphism(-150,50)(1,0)[A\ltimes B`A'\ltimes B `f\ltimes B]{600}1a
 \putmorphism(450,50)(1,0)[\phantom{A\ot B}`A'\ltimes B' `A'\ltimes g]{680}1a

\putmorphism(-180,500)(0,-1)[\phantom{Y_2}``=]{450}1r
\putmorphism(1100,500)(0,-1)[\phantom{Y_2}``=]{450}1r
\put(350,260){\fbox{$f\ltimes-\vert_g$}}

\putmorphism(-150,-400)(1,0)[\tilde A\ltimes B`\tilde A'\ltimes B `f'\ltimes B]{640}1a
 \putmorphism(450,-400)(1,0)[\phantom{A'\ot B'}` \tilde A'\ltimes B' `\tilde A'\ltimes g]{680}1a

\putmorphism(-180,50)(0,-1)[\phantom{Y_2}``v\ltimes B]{450}1l
\putmorphism(450,50)(0,-1)[\phantom{Y_2}``]{450}1r
\putmorphism(300,50)(0,-1)[\phantom{Y_2}``v'\ltimes B]{450}0r
\putmorphism(1100,50)(0,-1)[\phantom{Y_2}``]{450}1r
\putmorphism(1080,50)(0,-1)[\phantom{Y_2}``v'\ltimes B']{450}0r
\put(-20,-170){\fbox{$a\ltimes B$}}
\put(660,-170){\fbox{$v'\ltimes-\vert_g$}}

\efig}
\quad=\quad
\scalebox{0.86}{
\bfig
 \putmorphism(-170,500)(1,0)[A\ltimes B`A\ltimes B' `A\ltimes g]{620}1a
 \putmorphism(450,500)(1,0)[\phantom{A\ot B}`A'\ltimes B' `f\ltimes B']{680}1a
 \putmorphism(-150,50)(1,0)[\tilde A\ltimes B`\tilde A\ltimes B' `\tilde A\ltimes g]{600}1a
 \putmorphism(450,50)(1,0)[\phantom{A\ot B}`A'\ltimes B' `f'\ltimes B']{680}1a
\putmorphism(-180,500)(0,-1)[\phantom{Y_2}``]{450}1l
\putmorphism(-160,500)(0,-1)[\phantom{Y_2}``v\ltimes B]{450}0l
\putmorphism(450,500)(0,-1)[\phantom{Y_2}``]{450}1l
\putmorphism(610,500)(0,-1)[\phantom{Y_2}``v\ltimes B']{450}0l 
\putmorphism(1120,500)(0,-1)[\phantom{Y_3}``v'\ltimes B']{450}1r
\put(-130,270){\fbox{$v\ltimes-\vert_g$}} 
\put(650,270){\fbox{$a\ltimes B'$}}
\putmorphism(-150,-400)(1,0)[\tilde A\ltimes B`\tilde A'\ltimes B `f'\ltimes B]{640}1a
 \putmorphism(450,-400)(1,0)[\phantom{A'\ot B'}` \tilde A'\ltimes B' `\tilde A'\ltimes g]{680}1a

\putmorphism(-180,50)(0,-1)[\phantom{Y_2}``=]{450}1l
\putmorphism(1120,50)(0,-1)[\phantom{Y_3}``=]{450}1r
\put(300,-170){\fbox{$f'\ltimes-\vert_g$}}

\efig}
$$
and 
$$
\scalebox{0.86}{
\bfig
 \putmorphism(-150,500)(1,0)[A\ltimes B`A\ltimes B `=]{600}1a
 \putmorphism(450,500)(1,0)[A\ltimes B` `f\ltimes B]{480}1a
\putmorphism(-180,500)(0,-1)[\phantom{Y_2}` `A\ltimes u]{450}1l 
\putmorphism(-250,500)(0,-1)[\phantom{Y_2}` A\ltimes\tilde B`]{450}0l %
\put(-100,50){\fbox{$v\ltimes-\vert_{u}$}}
\putmorphism(-170,-400)(1,0)[\tilde A\ltimes \tilde B` `=]{480}1a
\putmorphism(-180,50)(0,-1)[\phantom{Y_2}``v\ltimes \tilde B]{450}1l
\putmorphism(450,50)(0,-1)[\phantom{Y_2}`\tilde A\ltimes \tilde B`\tilde A\ltimes u]{450}1l
\putmorphism(450,500)(0,-1)[\phantom{Y_2}`\tilde A\ltimes B `v\ltimes B]{450}1l
\put(600,260){\fbox{$a\ltimes B$}}
\putmorphism(430,50)(1,0)[\phantom{(B, \tilde A)}``f'\ltimes B]{500}1a
\putmorphism(1070,50)(0,-1)[\phantom{(B, A')}`\tilde A'\ltimes \tilde B`]{450}1r
\putmorphism(1050,50)(0,-1)[``\tilde A'\ltimes u]{450}0r
\putmorphism(1070,500)(0,-1)[A'\ltimes B`\tilde A'\ltimes B`]{450}1r
\putmorphism(1050,500)(0,-1)[``v'\ltimes B]{450}0r
\putmorphism(450,-400)(1,0)[\phantom{(B, \tilde A)}``f'\ltimes\tilde B]{500}1a
\put(540,-170){\fbox{$f'\ltimes-\vert_{u}$}}
\efig}
\quad=\quad
\scalebox{0.86}{
\bfig
 \putmorphism(-180,500)(1,0)[A\ltimes B`A'\ltimes B`f\ltimes B]{630}1a
 \putmorphism(450,500)(1,0)[\phantom{(B,A)}` `=]{520}1a
\putmorphism(-180,500)(0,-1)[\phantom{Y_2}`A\ltimes \tilde B `]{450}1l
\putmorphism(-160,500)(0,-1)[` `A\ltimes u]{450}0l
\put(620,50){\fbox{$v'\ltimes-\vert_{u}$}}
\putmorphism(-170,-400)(1,0)[\tilde A\ltimes\tilde B` `f'\ltimes \tilde B]{470}1a
\putmorphism(-180,50)(0,-1)[\phantom{Y_2}``]{450}1l
\putmorphism(-160,50)(0,-1)[\phantom{Y_2}``v\ltimes \tilde B]{450}0l
\putmorphism(450,50)(0,-1)[\phantom{Y_2}`\tilde A'\ltimes\tilde B `v'\ltimes \tilde B]{450}1r
\putmorphism(450,500)(0,-1)[\phantom{Y_2}`A'\ltimes \tilde B `A'\ltimes u]{450}1r %
\put(0,260){\fbox{$f\ltimes-\vert_{u}$}}
\putmorphism(-190,50)(1,0)[\phantom{(B, \tilde A)}``f\ltimes\tilde B]{500}1a
\putmorphism(1130,50)(0,-1)[\phantom{(B, A')}`\tilde A'\ltimes \tilde B`\tilde A'\ltimes u]{450}1r %
\putmorphism(1130,500)(0,-1)[A'\ltimes B``v'\ltimes B]{450}1r
\putmorphism(1200,500)(0,-1)[`\tilde A'\ltimes B`]{450}0r
\putmorphism(450,-400)(1,0)[\phantom{(B, \tilde A)}``=]{540}1b
\put(0,-170){\fbox{$a\ltimes \tilde B$}}
\efig}
$$
for any 1h-cell $g$ and any 1v-cell $u$ in $\Bb$. 

A {\em right central} 2-cell is defined similarly. 
A 2-cell is said to be {\em central} if it is both left and right central. 
\end{defn}

\subsection{Centrality and companions}

For concluding premonoidality of the underlying bicategory of a premonoidal double category, we will need the following results. They both hold in a binoidal double category. The first one is proved directly. 

\begin{lma}
If a 1v-cell $u$ has $\hat u$ for a companion, then $u\rtimes B$ has $\hat u\rtimes B$ for a companion. Similarly,  
$B\ltimes u$ has $B\ltimes\hat u$ for a companion.
\end{lma}

\begin{lma} \lelabel{inversely central}
Let $u$ be an invertible 1v-cell so that both $u$ and its inverse are left central with mutually inverse 2-cell components 
$u\ltimes-\vert_f$ and $u^{-1}\ltimes-\vert_f$ at 1h-cells $f$. Suppose also that $u$ has a companion and a conjoint. Then its companion $\hat u$ is left central. \\ The right-hand sided version of the claim holds, too. 
\end{lma}

\begin{proof}
By left centrality of $u$ 
we have a vertical pseudonatural transformation $u\ltimes-$. By the assumption, the 2-cell components $u\ltimes-\vert_f$ 
are invertible, so 
by \prref{lifting 1v to equiv} there is a horizontal pseudonatural transformation $\hat u\ltimes-$, hence the claim. 
\qed\end{proof}


\begin{defn} 
For an invertible 1v-cell $u$ such that both $u$ and its inverse $u^{-1}$ are left central, and so that their 2-cell components $u\ltimes-\vert_f$ and $u^{-1}\ltimes-\vert_f$ at 1h-cells $f$ are mutually inverse, we will say that it is 
{\em left inversely central}. Similarly we define {\em right inversely central} 1v-cells. Those 1v-cells which are both 
left and right inversely central we call {\em inversely central} 1v-cells. 
\end{defn}

\begin{lma} \lelabel{hat of central 2-cells}
If a 2-cell $\sigma$ is central in $\Bb$, then the 2-cell $\hat\sigma$ given as bellow is central in $\HH(\Bb)$. 
$$
\scalebox{0.9}{
\bfig
\putmorphism(-150,170)(1,0)[A`A`f]{460}1a
\putmorphism(-150,-180)(1,0)[B`B`g]{460}1a
\putmorphism(-150,170)(0,-1)[\phantom{Y_2}``u]{350}1l
\putmorphism(310,170)(0,-1)[\phantom{Y_2}``v]{350}1r
\put(20,0){\fbox{$\sigma$}}
\efig}
\quad\mapsto\quad \hat\sigma=\quad
\scalebox{0.86}{
\bfig
 \putmorphism(-150,250)(1,0)[A`A `=]{480}1a
\put(0,30){\fbox{$\eta_u$}}
\putmorphism(340,250)(0,-1)[\phantom{Y_2}`B `]{450}1l
\putmorphism(360,250)(0,-1)[\phantom{Y_2}` `u]{450}0l
\putmorphism(810,250)(0,-1)[A`B `]{450}1r
\putmorphism(790,250)(0,-1)[` `v]{450}0r
\putmorphism(1270,250)(0,-1)[B`B.`=]{450}1r
\putmorphism(350,250)(1,0)[``f]{420}1a
 \putmorphism(730,250)(1,0)[\phantom{F(A)}` `\hat v]{510}1a
 \putmorphism(-150,-200)(1,0)[A` `\hat u]{460}1a
\putmorphism(360,-200)(1,0)[``g]{440}1b
 \putmorphism(820,-200)(1,0)[``=]{440}1b
\put(530,10){\fbox{$\sigma$}}
\put(980,30){\fbox{$\Epsilon_v$}}
\putmorphism(-150,250)(0,-1)[``=]{450}1l
\efig}
$$ 
\end{lma}

\begin{proof}
To the modification condition \axiomref{m.ho-vl.-1} with $\Theta=\sigma$ and for a 1h-cell $k:D\to D'$ paste 
$\Epsilon_{\alpha_{A',B,C}\ltimes D'}$ from the right and $\eta_{\alpha_{A,B,C}\ltimes D}$ from the left. We explain the procedure on one 
side of the equation, analogous actions are done on the other side. Between the 2-cells 
$\alpha_{f,B,C}\ltimes D$ and $\alpha_{A',B,C}\ltimes -\vert_k$ insert the identity 2-cell $\frac{\eta_{\alpha_{A,B,C}\ltimes D}}
{\Epsilon_{\alpha_{A,B,C}\ltimes D}}$. 
One obtains 2-cells $\widehat{\alpha_{A',B,C}\ltimes -\vert_k}=\widehat{\alpha_{A',B,C}}\ltimes -\vert_k$ and 
$\widehat{\alpha_{f,B,C}\ltimes D}$, which is $\widehat{\alpha_{f,B,C}}\ltimes D$ conjugated by the lax and colax structure of $-\rtimes D$. 
Finally, relate the composition of components $(fB)C\ltimes-\vert_k$ and 
$\widehat{\alpha_{A',B,C}}\ltimes -\vert_k$ via \leref{vert comp hor.ps.tr.} to $\widehat{\alpha_{A',B,C}}((fB)C)\ltimes-\vert_k$. 
to get (one side of) the modification condition of $\hat\sigma\ltimes$ in $\HH(\Bb)$. Right centrality is proved in a similar way. 
\qed\end{proof}

\subsection{Premonoidal double categories}

Taking for the starting point the definition of a premonoidal bicategory from \cite{HF1}, 
one can follow two approaches to introduce a premonoidal double category. One is to add rules to the definition of a premonoidal bicategory for the 1v-cells and to accommodate the rules for 2-cells, so that forgetting the vertical direction one recovers the notion of a premonoidal bicategory. 

Another one is to follow Shulman's recipe from \cite{Shul}. Namely, to consider a premonoidal analogue of a monoidal double category $\Dd$ from his \cite[Definition 2.9]{Shul}, which is a pseudomonoid in a certain 2-category of double categories, and then assuming the existence of companions and conjoints in $\Dd$ lift isomorphism 1v-cells to equivalence 1h-cells, obtaining that the underlying horizontal bicategory 
$\HH(\Dd)$ is a proper monoidal bicategory. 

The latter is a much easier work to do, given that the first way would suppose to work with non-trivial modifications in double categories, while in the Shulman's way the (vertical) double modifications are trivial, but one still obtains non-trivial (horizontal) modifications $\pi, \mu, \lambda, \rho$ from the definition of a tricategory from \cite{GPS}. (The latter are necessary to have a monoidal bicategory, which is a one object tricategory.)
The cost one pays, though, is that this mechanism works only for those double categories whose 1v-cells have companions and conjoints. 
A synonym for such double categories is {\em framed bicategories}, \cite[Theorem A.2]{Shul:Fram}, and some examples are listed in 
\cite[Section 4.4]{Shul:Fram}.

We take the second, Shulman's approach.

\begin{defn} \delabel{unital} 
Let $\Bb$ be a binoidal double category. 
We say that the binoidal structure of $\Bb$ is {\em unital} if there exists a unit object $I$ 
and invertible vertical strict transformations $\lambda: I\ltimes -\Rightarrow\Id$ and $\rho: -\rtimes I\Rightarrow\Id$ each of whose  
1v-cell components are inversely central 1v-cells. 
\end{defn}

\begin{defn} \delabel{assoc} 
Let $\Bb$ be a binoidal double category. 
We say that the binoidal structure of $\Bb$ is {\em associative} if 
there exist invertible vertical strict transformations 
$$\alpha_{-,B,C}: (-\rtimes B)\rtimes C\Rightarrow -\rtimes(B\bowtie C)$$
$$\alpha_{A,-,C}: (A\ltimes -)\rtimes C\Rightarrow A\ltimes(-\rtimes C)$$
$$\alpha_{A,B,-}: (A\bowtie B)\ltimes -\Rightarrow A\ltimes(B\ltimes -)$$
for every $A,B,C\in\Bb$, such that the following is fulfilled  
\begin{enumerate}
\item 1v-cell components $\alpha_{A,B,C}$ of the above three vertical transformations coincide and are inversely central, 
\item and the following four pentagons, expressing equalities of vertical transformations, commute
\begin{equation*} 
\scalebox{0.84}{
\bfig 
\putmorphism(0,500)(1,0)[((-B)C)D `(-B)(CD) ` \alpha_{-B,C,D}]{2230}1a
\putmorphism(0,0)(1,0)[(-(BC))D `(-)((BC)D) ` \alpha_{-,BC,D}]{1050}1b
\putmorphism(1080,0)(1,0)[\phantom{((-B)C)D}`(-)(B(CD))` -\rtimes\alpha_{B,C,D}]{1220}1b
\putmorphism(60,500)(0,1)[`` \alpha_{-, B,C}\rtimes D]{500}1l 
\putmorphism(2210,500)(0,1)[`` \alpha_{-,B,CD}]{500}1r 
\efig}
\end{equation*}

\begin{equation*} 
\scalebox{0.84}{
\bfig 
\putmorphism(0,500)(1,0)[((A-)C)D `(A-)(CD) ` \alpha_{A-,C,D}]{2230}1a 
\putmorphism(0,0)(1,0)[(A(-C))D `A((-C)D) ` \alpha_{A,-C,D}]{1050}1b 
\putmorphism(1080,0)(1,0)[\phantom{((-B)C)D}`A(-(CD))` A\ltimes\alpha_{-,C,D}]{1220}1b 
\putmorphism(60,500)(0,1)[`` \alpha_{A, -,C}\rtimes D]{500}1l 
\putmorphism(2210,500)(0,1)[`` \alpha_{A,-,CD}]{500}1r 
\efig}
\end{equation*}

\begin{equation*} 
\scalebox{0.84}{
\bfig 
\putmorphism(0,500)(1,0)[((AB)-)D `(AB)(-D) ` \alpha_{AB,-,D}]{2230}1a 
\putmorphism(0,0)(1,0)[(A(B-))D `A((B-)D) ` \alpha_{A,B-,D}]{1050}1b 
\putmorphism(1080,0)(1,0)[\phantom{((-B)C)D}`A(B(-D))` A\ltimes\alpha_{B,-,D}]{1220}1b 
\putmorphism(60,500)(0,1)[`` \alpha_{A, B,-}\rtimes D]{500}1l 
\putmorphism(2210,500)(0,1)[`` \alpha_{A,B,-D}]{500}1r 
\efig}
\end{equation*}

\begin{equation*} 
\scalebox{0.84}{
\bfig 
\putmorphism(0,500)(1,0)[((AB)C)- `(AB)(C-) ` \alpha_{AB,C,-}]{2230}1a 
\putmorphism(0,0)(1,0)[(A(BC))- `A((BC)-) ` \alpha_{A,BC,-}]{1050}1b 
\putmorphism(1080,0)(1,0)[\phantom{((-B)C)D}`A(B(C-)).` A\ltimes\alpha_{B,C,-}]{1220}1b 
\putmorphism(60,500)(0,1)[`` \alpha_{A, B,C}\ltimes -]{500}1l 
\putmorphism(2210,500)(0,1)[`` \alpha_{A,B,C-}]{500}1r 
\efig}
\end{equation*}
\end{enumerate}
\end{defn}

\begin{defn} \delabel{premon}
Let $\Bb$ be a binoidal double category. 
We say that $\Bb$ is a {\em premonoidal} double category if its binoidal structure is unital and associative and the following six triangles, 
expressing equalities of vertical transformations, commute 
$$\scalebox{0.88}{\bfig
 \putmorphism(-180,500)(1,0)[(-I)B`(-)B`\rho_{-}\rtimes B]{680}1a
 \putmorphism(0,0)(1,1)[(-)(IB)``]{200}0a
\putmorphism(140,170)(1,1)[``-\rtimes\lambda_B]{160}1r
\putmorphism(-90,500)(0,-1)[\phantom{B\ot B}``\alpha_{-,I,B}]{480}1l
\efig}
\qquad
\scalebox{0.88}{\bfig
 \putmorphism(-180,500)(1,0)[(AI)-`(A)-`\rho_A\ltimes-]{600}1a
 \putmorphism(0,0)(1,1)[A(I-)``]{200}0a
\putmorphism(140,170)(1,1)[``A\ltimes\lambda_{-}]{160}1r
\putmorphism(-90,500)(0,-1)[\phantom{B\ot B}``\alpha_{A,I,-}]{480}1l
\efig}
$$

$$\scalebox{0.88}{\bfig
 \putmorphism(-180,500)(1,0)[(IA)-`A(-)`\lambda_A\ltimes-]{650}1a
 \putmorphism(0,0)(1,1)[I(A-)``]{200}0a
\putmorphism(140,170)(1,1)[``\lambda_{A\ltimes-}]{160}1r
\putmorphism(-90,500)(0,-1)[\phantom{B\ot B}``\alpha_{I,A,-}]{480}1l
\efig}
\qquad
\scalebox{0.88}{\bfig
 \putmorphism(-180,500)(1,0)[(I-)B`(-)B`\lambda_{-}\rtimes B]{680}1a
 \putmorphism(0,0)(1,1)[I(-B)``]{200}0a
\putmorphism(140,170)(1,1)[``\lambda_{-\rtimes B}]{160}1r
\putmorphism(-90,500)(0,-1)[\phantom{B\ot B}``\alpha_{I,-,B}]{480}1l
\efig}
$$

$$\scalebox{0.88}{\bfig
 \putmorphism(-180,500)(1,0)[(-B)I`(-)B`\rho_{-\rtimes B}]{600}1a
 \putmorphism(0,0)(1,1)[(-)(BI)``]{200}0a
\putmorphism(140,170)(1,1)[``-\rtimes\rho_B]{160}1r
\putmorphism(-90,500)(0,-1)[\phantom{B\ot B}``\alpha_{-,B,I}]{480}1l
\efig}
\qquad
\scalebox{0.88}{\bfig
 \putmorphism(-180,500)(1,0)[(A-)I`A(-)`\rho_{A\ltimes-}]{600}1a
 \putmorphism(0,0)(1,1)[A(-I).``]{200}0a
\putmorphism(140,170)(1,1)[``A\ltimes\rho_{-}]{160}1r
\putmorphism(-90,500)(0,-1)[\phantom{B\ot B}``\alpha_{A,-,I}]{480}1l
\efig}
$$
\end{defn}


\begin{rem} \rmlabel{centrality}
In the definition of a premonoidal double category, centrality of the component 1v-cells $\alpha_{A,B,C}, \lambda_A$ and 
$\rho_A$ is used in order to have well-defined invertible vertical pseudonatural transformations $-\rtimes\alpha_{B,C,D}, 
\alpha_{A,B,C}\ltimes-$ in the first and fourth pentagon in \deref{assoc}, and $-\rtimes\lambda_B, \rho_A\ltimes-, 
\lambda_A\ltimes-, -\rtimes\rho_B$ in the first three and the fifth triangle in \deref{premon}. 
Inverse centrality is included in these definitions for a reason explained in \rmref{inverse centrality and specialness}. 
\end{rem}

In the pentagons and triangles of a premonoidal double category the compositions of vertical transformations are the vertical ones, 
from \leref{vert comp vert. lx tr.}. The 2-cell $(\frac{\alpha_0}{\beta_0})_f$ for two such transformations and a 1h-cell $f$ is given by the vertical 
composition of the 2-cells $(\alpha_0)_f$ and $(\beta_0)_f$.

\begin{rem} \rmlabel{inverse centrality and specialness}
Observe that in a premonoidal double category in the above definition 1v-cell components of the five vertical strict transformations 
are invertible and inversely central. If these transformations are even liftable, then by \prref{lifting 1v to equiv} and 
\leref{inversely central} they determine bicategorical 
pseudonatural transformations in which the 1-cell components are equivalence 1-cells and central. This is what we want to obtain in the underlying horizontal 2-category $\HH(\Bb)$ (\prref{premon isofib new}).   
\end{rem}

\subsection{Other implications in a premonoidal double category}

In a premonoidal double category $\Dd$ we have: 
\begin{enumerate}
\item pseudodouble functors $A\ltimes -$ and $-\rtimes B$ such that $A\ltimes B=A\rtimes B=:A\bowtie B$; 
\item a unit object $I$
with invertible vertical strict transformations $\lambda: I\ltimes -\Rightarrow\Id$ and $\rho: -\rtimes I\Rightarrow\Id$ each of whose  
1v-cell components are central 1v-cells;
\item invertible vertical strict transformations 
$\alpha_{-,B,C}, \alpha_{A,-,C},\alpha_{A,B,-}$ for every $A,B,C\in\Dd$, such that each 1v-cell component $\alpha_{A,B,C}$ is central, 
the 24 axioms hold and the four pentagons commute;
\item the six triangles commute.
\end{enumerate} 

\bigskip

The pseudodouble functors $A\ltimes -$ and $-\rtimes B$ yield functors on the categories $\Dd_0$ of objects and $\Dd_1$ of morphisms, 
and also pseudofunctors  $A\ltimes -,-\rtimes B: \HH(\Dd)\to\HH(\Dd)$ on the underlying horizontal 2-category of $\Dd$. 

The existence of a distinguished (unit) object $I\in\Dd$ means that there is a pseudodouble functor $1\to\Dd$ from the trivial double category 
so that the image of the single object $*$ is $I$ (and the images of the remaining identity 1h- and 1v-cells and the identity 2-cell on $*$ 
are horizontal and vertical identities $1_I$ and $1^I$ on $I$ and the identity 2-cell on $I$, respectively). Then there are also functors 
$*\to\Dd_0$ and $*\to\Dd_1$ from the trivial category, and a pseudofunctor $*\to\HH(\Dd)$ from the trivial 2-category. 

Invertible vertical strict transformations are those that have vertically invertible 2-cells $(\alpha_0)_f$ for every 1h-cell $f$, 
and they satisfy the axioms \axiomref{v.l.t.\x 1}, \axiomref{v.l.t.\x 2} and \axiomref{v.l.t.\x 5} so that the latter axiom is simplified 
(the 2-cells $\alpha_0^u$ are identities), whereas the compositors in the first axiom and the unitors in the second one are invertible. Equivalently, invertible vertical strict transformations $\omega:F\to G$ between pseudodouble functors consist of two natural transformations 
$\omega_0:F_0\to G_0$ and $\omega_1:F_1\to G_1$ between the induced functors on the categories of objects and morphisms, such that the axioms 
\axiomref{v.l.t.\x 1} and \axiomref{v.l.t.\x 2} hold, and the 2-cells $(\alpha_0)_f$ are vertically invertible. 

The fact that we work with vertical transformations together with the functor $*\to\Dd_0$ produces that associative and unital binoidal structure of $\Dd$ is passed on to 
an associative and unital binoidal structure of the vertical category $\Dd_0$. For the horizontal category $\Dd_1$ we consider 
$1_I$ as the unit object with the rest of a monoidal structure given by the internal category structure of $\Dd$. 


We may conclude:

\begin{lma}
In a premonoidal double category $\Dd$ the category of objects $\Dd_0$ is premonoidal and the category of morphisms $\Dd_1$ is monoidal. 
\end{lma}

We may also prove:

\begin{prop} \prlabel{premon isofib new}
The horizontal bicategory of a premonoidal double category whose associativity and unitality constraints are liftable vertical 
transformations is a premonoidal bicategory from \deref{premon-bicat}. 
\end{prop}

\begin{proof}
We saw in \rmref{inverse centrality and specialness} that invertible liftable vertical strict transformations whose 1v-cell components are inversely central 
determine bicategorical pseudonatural equivalences whose 1-cell components are (equivalence 1-cells and) central. 
That the four pentagons and six triangles of a premonoidal double category commute means that they are given by globular vertical identity modifications between vertical compositions of invertible vertical strict transformations. Then according to \prref{essence} 
they induce horizontal {\em i.e.} bicategorical modifications between the corresponding horizontal {\em i.e.} bicategorical pseudonatural equivalences, and they satisfy in particular the equations between 2-cells in a monoidal bicategory. Thus we have the proof. 
\qed\end{proof}

For a premonoidal double category $\Bb$ whose associativity and unitality constraints are liftable vertical transformations we will 
denote by $\ul{\HH(\Bb)}$ the underlying 2-category $\HH(\Bb)$ equipped 
with bicategorical pseudonatural equivalences $\hat\alpha,\hat\lambda,\hat\rho$ induced from $\Bb$ according to \prref{lifting 1v to equiv}. 
We now know that it is premonoidal.

\section{24 axioms for an associative binoidal structure} \selabel{24 axioms}  


In a premonoidal double category a natural question arises about the coherence of any two of the three associativity constraints. This  coherence can be expressed in the form of 24 axioms that we are going to study in this section. Though, not any pair of the component 
2-cells of these structural transformations can be compared: one of the two in them appearing 1-cells should be central. Thus a coherent choice of centrality structures of 1h- and 1v-cells is necessary in order to express these axioms. The definition of a premonoidal double category does not entail such a choice, for this reason these 24 axioms do not appear in it. Center double categories of premonoidal double categories and funny functors and quasi-functors, which we will introduce in later sections, will provide a coherent choice of centrality structures for a premonoidal double category. The 24 axioms are relevant for those cases. 

\smallskip

For the interrelation between $\alpha_{-,B,C}$ and $\alpha_{A,-,C}$ we should give four axioms corresponding to compatibilities of 
their 2-cell components, 
and similarly for  the interrelations between $\alpha_{-,B,C}$ and $\alpha_{A,B,-}$, and between $\alpha_{A,-,C}$ and $\alpha_{A,B,-}$. This gives 12 axioms. 
Though, in each of these 12 axioms there appears (something that should be) a transformation of either of the following forms: 
$f\ltimes-, -\rtimes f'$ and $u\ltimes-, -\rtimes u'$, whereby we should differentiate 2-cell components $f\ltimes-\vert_{f'}$ for 
left central 1h-cells $f$, and $-\rtimes f'\vert_f$ for right central 1h-cells $f'$, and so on for the other combinations of 1-cells. 
This results in having in total 24 axioms that should hold for pairs consisting of 
{\em any left (or right) central} 1h-cell (or 1v-cell) {\em and any other} 1h-cell (or 1v-cell). Moreover, the said candidates for transformations should be transformations indeed. That the latter is fulfilled will actually be the consequence of the 24 axioms themselves.

On the other hand, one would also hope to have a coherence of structural transformations of central 1v-cell components of the three 
$\alpha$'s with respect to the binoidal structure. 
If $u$ stands for any of them, we need to know the relations between 
$u\ltimes(-,\bullet)$ and $(u\ltimes-)\ltimes\bullet$, where the places $(-,\bullet)$ can be occupied by any of the following four combinations of cells: $(g,C), (B,h), (v,C), (B, z)$ for 1h-cells $g,h$ and 1v-cells $v,z$. This and similar coherences for 
1h-cells in place of $u$ will turn out to be covered by the 24 axioms. 

\smallskip

To illustrate the above, let us have a look at the first of the 24 axioms  
\noindent \axiomref{$(f\ltimes,g,C)$}: \vspace{-0,14cm}
$$\scalebox{0.86}{
\bfig
\putmorphism(-150,500)(1,0)[``((A,g),C)]{600}1a
 \putmorphism(480,500)(1,0)[` `((f,B'),C)]{640}1a
 \putmorphism(-150,50)(1,0)[``((f,B),C)]{600}1a
 \putmorphism(470,50)(1,0)[` `((A',g),C)]{660}1a

\putmorphism(-180,500)(0,-1)[\phantom{Y_2}``=]{450}1l
\putmorphism(1100,500)(0,-1)[\phantom{Y_2}``=]{450}1r
\put(230,280){\fbox{$(f\ltimes-\vert_g,C)$}}

\putmorphism(-170,-400)(1,0)[` `(f,(B,C))]{640}1b
 \putmorphism(470,-400)(1,0)[` `(A',(g,C))]{640}1b

\putmorphism(-180,50)(0,-1)[\phantom{Y_2}``\alpha_{A,B,C}]{450}1l %
\putmorphism(450,50)(0,-1)[\phantom{Y_2}``]{450}1l
\putmorphism(660,50)(0,-1)[\phantom{Y_2}``\alpha_{A',B,C}]{450}0l 
\putmorphism(1100,50)(0,-1)[\phantom{Y_3}``\alpha_{A',B',C}]{450}1r
\put(-40,-180){\fbox{$\alpha_{f,B,C}$}} 
\put(670,-180){\fbox{$\alpha_{A',g,C}$}}
\efig}
\quad
=
\quad
\scalebox{0.86}{
\bfig
\putmorphism(-150,500)(1,0)[``((A,g),C)]{600}1a
 \putmorphism(480,500)(1,0)[` `((f,B'),C)]{640}1a

 \putmorphism(-150,50)(1,0)[``(A,(g,C))]{600}1a
 \putmorphism(450,50)(1,0)[``(f,(B',C))]{640}1a

\putmorphism(-180,500)(0,-1)[\phantom{Y_2}``\alpha_{A,B,C}]{450}1l
\putmorphism(450,500)(0,-1)[\phantom{Y_2}``]{450}1r
\putmorphism(250,500)(0,-1)[\phantom{Y_2}``\alpha_{A,B',C}]{450}0r
\putmorphism(1100,500)(0,-1)[\phantom{Y_2}``\alpha_{A',B',C}]{450}1r
\put(-90,280){\fbox{$\alpha_{A,g,C}$}}
\put(680,280){\fbox{$\alpha_{f,B',C}$}}

\putmorphism(-150,-400)(1,0)[` `(f,(B,C))]{640}1b
 \putmorphism(490,-400)(1,0)[` .`(A',(g,C))]{640}1b

\putmorphism(-180,50)(0,-1)[\phantom{Y_2}``=]{450}1l
\putmorphism(1100,50)(0,-1)[\phantom{Y_3}``=]{450}1r
\put(270,-200){\fbox{$f\ltimes-\vert_{(g,C)}$}}
\efig} 
$$
It should hold for every {\em left central 1h-cell} $f\colon A\to A'$ and {\em any 1h-cell} $g\colon B\to B'$. (We write here $(-,-)$ for 
both $-\ltimes-$ and $-\rtimes-$, which one is meant is clear from the context. 
In the rectangular diagram for the 2-cell $(f\ltimes-\vert_g,C)$ we omitted the compositor 2-cells 
$((f,B')\rtimes C)((A,g)\rtimes C) \Rightarrow (f,B')(A,g)\rtimes C$ at the top and 
$(f,B)(A',g)\rtimes C\Rightarrow ((f,B)\rtimes C)((A',g)\rtimes C)$ at the bottom, of the pseudodouble functor $-\rtimes C$). 
This axiom correlates the structural transformations $(f\ltimes-\vert_g)\rtimes C$ and $f\ltimes-\vert_{(g,C)}$ of the left central 1h-cell 
$f$, and also the 2-cell components of the vertical transformations $\alpha_{-,B,C}$ and $\alpha_{A,-,C}$ evaluated at 1h-cells. 
Moreover, observe that this axiom corresponds to the axiom \axiomref{m.ho-vl.-1}, whereas the axiom \axiomref {$(f\ltimes,v,C)$} 
(see Appendix A) 
corresponds to the axiom \axiomref{m.ho-vl.-2} for $\alpha_{f,-,C}$ to be a modification (between horizontal transformations 
$(f\ltimes-)\rtimes C$ and $f\ltimes(-\rtimes C)$ and vertical transformations $\alpha_{A,-,C}$ and $\alpha_{A',-,C}$). Though, 
we should also make sure under which conditions the former two in the parenthesis are horizontal transformations. 

We stress the fact that the 2-cell components $\alpha_{f,B,C}$ and $\alpha_{A,g,C}$ in the above axiom (and similarly in all the remaining 23 axioms) {\em
can not be compared for any 1h-cells $f$ and $g$}: $f$ should be left central in order for the 2-cells $(f\ltimes-\vert_g,C)$ and $f\ltimes-\vert_{(g,C)}$ to exist (and in general one of two determining 1-cells in any axiom must be central at an appropriate side). 

Summarizing, the $12+12$ axioms happen to correspond to the axioms \axiomref{m.ho-vl.-1} and \axiomref{m.ho-vl.-2} of certain 6+0 modifications between certain 12 horizontal and three vertical strict transformations $\alpha$. It will turn out that the same $12+12$ axioms assure that we 
indeed have those 12 horizontal and further 12 vertical pseudonatural transformations. The latter fact we show in the sequel. 

\smallskip

The complete list of the 24 axioms can be found in Appendix A. In Table \ref{table:2} we list 12 of those axioms, the 6 modifications that they 
determine, and 12 (top and bottom) horizontal transformations between which the modifications act. 
In Table \ref{table:3} further below we will list the remaining 12 axioms and 12 vertical transformations, which they interrelate.  

\begin{table}[H]
\begin{center}
\begin{tabular}{ c c c } 
Pairs of axioms & \hspace{0,2cm} that define  & \hspace{0,2cm}  candidate horiz. transf. \\ [0.5ex]
 & candidate modif. & top \hspace{2cm}  bottom \\ [0.5ex]
\hline
\axiomref{$(f\ltimes,g,C)$} and \axiomref{$(f\ltimes,v,C)$} & $\alpha_{f,-,C}$ &  $(f\ltimes-)\rtimes C$ \qquad $f\ltimes(-\rtimes C)$ \\ [1ex] 
\hdashline[0.5pt/5pt]
\axiomref{$(f\ltimes,B,h)$} and \axiomref{$(f\ltimes,B,z)$} & $\alpha_{f,B,-}$ &  $(f\ltimes B)\ltimes -$ \qquad $f\ltimes(B\ltimes -)$ \\ [1ex] 
\axiomref{$(A,g\ltimes,h)$} and \axiomref{$(A,g\ltimes,z)$}  & $\alpha_{A,g,-}$ &  $(A\ltimes g)\ltimes -$ \qquad $A\ltimes(g\ltimes -)$ \\ [1ex] 
\axiomref{$(f,\rtimes g,C)$} and \axiomref{$(u,\rtimes g,C)$}  & $\alpha_{-,g',C}$ &  $(-\rtimes g')\rtimes C$ \qquad $-\rtimes(g'\rtimes C)$ \\ [1ex] 

\axiomref{$(f,B,\rtimes h)$} and \axiomref{$(u,B,\rtimes h)$} & $\alpha_{-,B,h}$ &  $(-\rtimes B)\rtimes h$ \qquad $-\rtimes(B\ltimes h)$ \\ [1ex] 
\hdashline[0.5pt/5pt]
\axiom{$(A,g,\rtimes h)$} and \axiom{$(A,v,\rtimes h)$} & $\alpha_{A,-,h}$ &  $(A\ltimes-)\rtimes h$ \qquad $A\ltimes(-\rtimes h)$ \\ [1ex] 
\end{tabular}
\caption{Interrelations of horizontal transformations and generation of modifications}
\label{table:2}
\end{center}
\end{table}

We start by noticing that 8 out of 12 candidate horizontal transformations from the last two columns of Table \ref{table:2} are automatically horizontal transformations by \leref{H of omega}. Similarly, 8 out of 12 candidate vertical transformations in Table \ref{table:3} are automatically vertical pseudonatural transformations by \leref{H of alfa}. Namely, we have: 

\bigskip

\begin{lma} \lelabel{basic}
Let $\Bb$ be a double category with a binoidal structure. 
\begin{enumerate}
\item 
For $f,g$ left central and $g',h$ right central 1h-cells the following are horizontal pseudonatural transformations 
$$(f\ltimes -)\rtimes C, \,\,\,  \hspace{10cm} A\ltimes(-\rtimes h),$$
$$f\ltimes (-\rtimes C), \,\,\,  f\ltimes (B\ltimes-), \,\,\,  A\ltimes (g\ltimes -) \quad\text{and}\quad 
(-\rtimes g')\rtimes C, \,\,\, (-\rtimes B)\rtimes h, \,\,\, (A\ltimes-)\rtimes h.$$
\item 
For $u,v$ left central and $v',z$ right central 1v-cells the following are vertical pseudonatural transformations 
$$(u\ltimes -)\rtimes C, \,\,\,  \qquad\qquad \,\,\,  A\ltimes (v\ltimes -) \qquad\qquad\quad 
(-\rtimes v')\rtimes C, \,\,\, \qquad\qquad \,\,\, A\ltimes(-\rtimes z),$$
$$u\ltimes (-\rtimes C), \,\,\,  u\ltimes (B\ltimes-), \,\,\,  \qquad \qquad \quad\text{and}\quad \qquad 
\qquad \,\,\, (-\rtimes B)\rtimes z, \,\,\, (A\ltimes-)\rtimes z.$$
\end{enumerate}
\end{lma}

For the remaining 4 candidates of horizontal transformations in Table \ref{table:2} we have the following. Define the 2-cell components of 
\begin{equation} \eqlabel{four center}
(A\ltimes g)\ltimes -, \,\,\, (f\rtimes B)\ltimes - \quad\text{and}\quad -\rtimes(g'\rtimes C), \,\,\, -\rtimes(B\ltimes h)
\end{equation}
by expressing them out of the two axioms written in the same line in Table \ref{table:2} as the candidate transformation itself 
(this covers the 8 axioms in the middle four rows of the first column of Table \ref{table:2}). This way these component 2-cells are given in terms of the component 2-cells of the transformations 
$$A\ltimes (g\ltimes -), \,\,\, f\ltimes (B\ltimes-)\quad\text{and}\quad  (-\rtimes g')\rtimes C, \,\,\, (-\rtimes B)\rtimes h$$ 
from \leref{basic}, respectively. 
The following is verified straightforwardly: 

\begin{lma} \lelabel{middle 8 ax}
Assume that the middle 8 axioms of Table \ref{table:2} hold. 
\begin{enumerate}
\item Then \equref{four center} are horizontal pseudonatural transformations.
\item The pseudofunctors $A\ltimes-$ and $-\rtimes B$ preserve left centrality (by first four of the 8 axioms) and right centrality 
(by the other four axioms) of 1h-cells. 
\item The 8 axioms mean that $\alpha_{f,B,-}, \,\,  \alpha_{A,g,-}, \,\, \alpha_{-,g',C}, \,\, \alpha_{-,B,h}$ are modifications. 
\end{enumerate}
\end{lma}

The axioms from the first and the sixth row in Table \ref{table:2} express how the component 2-cells of 
$(f\ltimes-)\rtimes C$ and $f\ltimes(-\rtimes C)$, and of $(A\ltimes-)\rtimes h$ and $A\ltimes(-\rtimes h)$, 
respectively, are interrelated. We have:

\begin{prop} \prlabel{horiz 12 axioms}
Let $\Bb$ be a double category with a binoidal structure and three invertible vertical pseudonatural transformations 
$$\alpha_{-,B,C}: (-\rtimes B)\rtimes C\Rightarrow -\rtimes(B\bowtie C)$$
$$\alpha_{A,-,C}: (A\ltimes -)\rtimes C\Rightarrow A\ltimes(-\rtimes C)$$
$$\alpha_{A,B,-}: (A\bowtie B)\ltimes -\Rightarrow A\ltimes(B\ltimes -)$$
for every $A,B,C\in\Bb$. Assume that the 12 axioms from Table \ref{table:2} hold. Then: 
\begin{enumerate}
\item for $f,g$ left central and $g',h$ right central 1h-cells 
the following are horizontal pseudonatural transformations 
$$(f\ltimes -)\rtimes C, \,\,\,  (f\ltimes B)\ltimes-, \,\,\,  (A\ltimes g)\ltimes -,  \quad\quad\quad 
-\rtimes(g'\rtimes C), \,\,\, -\rtimes(B\rtimes h), \,\,\, A\ltimes(-\rtimes h),$$
$$f\ltimes (-\rtimes C), \,\,\,  f\ltimes (B\ltimes-), \,\,\,  A\ltimes (g\ltimes -) \quad\text{and}\quad 
(-\rtimes g')\rtimes C, \,\,\, (-\rtimes B)\rtimes h, \,\,\, (A\ltimes-)\rtimes h$$
whose component 2-cells are related via the 12 axioms;
\item the following are modifications between vertical pseudonatural transformations $\alpha_{-,B,C}, \alpha_{A,-,C}, \alpha_{A,B,-}$ 
(with suitable indexes) and horizontal pseudonatural transformations from point 1. (also indicated in Table \ref{table:2}) 
$$\alpha_{f,-,C}, \,\,\, \alpha_{f,B,-}, \,\,\,  \alpha_{A,g,-}, \,\,\, \quad \alpha_{-,g',C}, \,\,\, \alpha_{-,B,h}, \,\,\, \alpha_{A,-,h}.$$
\end{enumerate}
\end{prop}

\begin{proof}
It only remains to discuss the last statement: it is immediately seen by inspection. 
We illustrate this claim by an example: the axioms \axiomref{$(f\ltimes,g,C)$} and \axiomref{$(f\ltimes,v,C)$} formally mean that 
$$\scalebox{0.86}{
\bfig
\putmorphism(-150,50)(1,0)[` `(f\ltimes-)\rtimes C]{400}1a
\putmorphism(-150,-270)(1,0)[`` f\ltimes(-\rtimes C)]{400}1b
\putmorphism(-150,50)(0,-1)[\phantom{Y_2}``\alpha_{A,-,C}]{320}1l
\putmorphism(250,50)(0,-1)[\phantom{Y_2}``\alpha_{A',-,C}]{320}1r
\put(-100,-130){\fbox{$\alpha_{f,-,C}$}}
\efig}$$
is a modification. 
\qed\end{proof}

\bigskip

We do not get, however,  that the vertical analogue of \equref{four center} make vertical pseudonatural transformations, {\em unless 
the three vertical 
transformations $\alpha$ in \prref{horiz 12 axioms} are strict}. Only under that condition the half of the axioms in the Appendix A that 
holds for one-sided central 1v-cells 
determines interrelated vertical pseudonatural transformations. Moreover, we do not get any modification from this half of the axioms.  

As we announced, in Table \ref{table:3} are the complementing 12 axioms to those from Table \ref{table:2}. Let $u,v$ be left central and $v',z$ right central 
1v-cells. 

\begin{table}[H]
\begin{center}
\begin{tabular}{ c c c } 
Axioms & \hspace{0,2cm}   &  that define/interrelate \\ [0.5ex]
 & \hspace{0,2cm}   & \hspace{0,2cm}  vertic. transf. \\ [0.5ex]
\hline
\axiom{$(u\ltimes,g,C)$} and \axiom{$(u\ltimes,v,C)$} & & $(u\ltimes-)\rtimes C$ \qquad $u\ltimes(-\rtimes C)$ \\ [1ex] \hdashline[0.5pt/5pt]

\axiom{$(A, v\ltimes,h)$} and \axiom{$(A, v\ltimes,z)$} &  & $(A\ltimes v)\ltimes-$ \qquad $A\ltimes (v\ltimes-)$ \\ [1ex] 

\axiom{$(u\ltimes, B,h)$} and \axiom{$(u\ltimes, B,z)$} & & $(u\ltimes B)\ltimes -$ \qquad $u\ltimes(B\ltimes -)$ \\ [1ex]

\axiom{$(f,B,\rtimes z)$} and \axiom{$(u,B,\rtimes z)$} &  & $-\rtimes(B\rtimes z)$ \qquad $(-\rtimes B)\rtimes z$ \\ [1ex] 

\axiom{$(f,\rtimes v,C)$} and \axiom{$(u,\rtimes v,C)$} &  & $(-\rtimes v')\rtimes C$ \qquad $-\rtimes (v'\rtimes C)$ \\ [1ex] 
\hdashline[0.5pt/5pt]

\axiom{$(A,g,\rtimes z)$} and \axiom{$(A,v,\rtimes z)$} &  & $A\ltimes(-\rtimes z)$ \qquad $(A\ltimes-)\rtimes z$ \\ [1ex]
\end{tabular}
\caption{Definitions and interrelations of vertical transformations}
\label{table:3}
\end{center}
\end{table}

Concretely, under the above strictness assumption, by middle 8 axioms from Table \ref{table:3} 
four vertical pseudonatural transformations 
$$A\ltimes (v\ltimes -), \,\,\,\quad  u\ltimes (B\ltimes-), \,\,\,\quad  (-\rtimes B)\rtimes z, \,\,\,\quad (-\rtimes v')\rtimes C$$
from point 2. of \leref{basic} define vertical pseudonatural transformations 
\begin{equation} \eqlabel{four vertical tr}
(A\ltimes v)\ltimes-, \,\,\,\quad (u\ltimes B)\ltimes -, \,\,\,\quad -\rtimes(B\rtimes z), \qquad -\rtimes (v'\rtimes C),  
\end{equation}
while the first two and the last two axioms in this Table force the relations between vertical pseudonatural transformations
$(u\ltimes-)\rtimes C$ and $u\ltimes(-\rtimes C)$, and between $(A\ltimes-)\rtimes z$ and $A\ltimes(-\rtimes z)$, respectively. 

\smallskip

Summing up: in this subsection we proved that if the 12 axioms from Table \ref{table:2} hold for the indicated one-sided central 1h-cells, then the data in its last two columns present 12 horizontal pseudonatural transformations and the data from its middle column present 6 modifications. If additionally the 12 axioms from Table \ref{table:3} hold for the indicated one-sided central 1v-cells, then the data in its last two columns present 12 vertical pseudonatural transformations.  
Moreover, the $12+12$ axioms exhibit the interrelations between the three vertical strict transformations $\alpha$. 
(Although the globular 2-cell components of the three $\alpha$ are trivial, their presence in the diagrams in the Appendix A is clarifying.) 


\section{Funny tensor product and premonoidal double categories} \selabel{funny}


In this section we introduce two funny type of products on double categories and relate them to premonoidal double categories. 
We will call them a pure and a mixed funny product. Each of them will give rise to a representable funny multicategory. These in turn will yield monoidal categories of double categories with funny types of monoidal product. A monoid in the first one will be a strict premonoidal double category, a notion that we will make precise later on. Then we are going to construct a monoidal 
2-category of double categories with the mixed funny type of monoidal product, and a pseudomonoid in it will be a particular case of a premonoidal double category.

\subsection{Funny functors and funny product for double categories} \sslabel{funny}

We denote by $Dbl$ the category of double categories and double functors. The funny product is tied to an inner-hom in which the transformations 
are not required to obey the naturality condition. Such transformations are called unnatural transformations. 

For double categories we encounter two ways of defining horizontal unnatural transformations (the vertical unnatural transformations 
are defined then in an analogous way). On one hand, following the logic of unnatural transformations for categories \cite{PR} and 
2-categories \cite{BG}, one can define horizontal unnatural transformations as given merely by their 1h-cell components obeying no axioms. 
Thus defined unnatural transformations we will refer to as {\em purely unnatural transformations}. 
Accordingly, modifications between horizontal and analogously defined vertical purely unnatural transformations are given by 
2-cell components that are not required to fulfill any axioms. We call them {\em purely funny modifications}. 
This has for a consequence that taking for an inner-hom the double category 
$[\Aa,\Bb]_f$ of double functors $\Aa\to\Bb$, 1h- and 1v-cells horizontal and vertical purely unnatural transformations, respectively, and purely funny modifications, the funny product $\Aa\Box_f\Bb$ obtained out of $[\Aa,\Bb]_f$ will ``behave correctly'' in a sense that we will make precise in \rmref{two funny for dbl}. 

However, we are interested to build a funny type of a monoidal 2-category of double categories so that pseudomonoids in it would 
be premonoidal double categories. It will though turn out that in this way we may obtain only premonoidal double categories whose binoidal structure is given by a pair of strict double functors. (The obstacle to obtaining the most general case we will explain in \ssref{funny-premonoidal connection}.) For that purpose, we will define horizontal unnatural transformations so that they  posses square-formed 2-cells that satisfy two axioms. 
When 1v-cells in these 2-cells are identities, the 2-cells are trivial, so they do retrieve the 2-categorical notion of unnatural 
transformations. Moreover, instead of taking thus defined 
unnatural transformations in both directions, for 1v-cells in the inner-hom $[\Aa,\Bb]$ that we will consider we will take {\em vertical strict transformations}. 

We define now this second kind of unnatural horizontal transformations of lax double functors between double categories, as well as modifications between them and vertical strict transformations. Such modifications we will call {\em mixed funny modifications}. 
We will then construct a funny product $\Aa\Box\Bb$ out of $[\Aa,\Bb]$. 

\begin{defn} \delabel{hor unnat tr}
A {\em horizontal unnatural transformation} $\alpha$ between lax double functors $F,G\colon \Aa\to\Bb$ consists of the following:
\begin{enumerate}
\item for every 0-cell $A$ in $\Aa$ a 1h-cell $\alpha(A)\colon F(A)\to G(A)$ in $\Bb$,
\item for every 1v-cell $u\colon A\to A'$ in $\Aa$ a 2-cell in $\Bb$:
$$
\scalebox{0.86}{
\bfig
\putmorphism(-150,50)(1,0)[F(A)`G(A)`\alpha(A)]{560}1a
\putmorphism(-150,-320)(1,0)[F(A')`G(A')`\alpha(A')]{600}1a
\putmorphism(-180,50)(0,-1)[\phantom{Y_2}``F(u)]{370}1l
\putmorphism(410,50)(0,-1)[\phantom{Y_2}``G(u)]{370}1r
\put(30,-110){\fbox{$\alpha^u$}}
\efig}
$$
\end{enumerate}
so that the following are satisfied 
(coherence with vertical composition and identity for $\alpha^\bullet$): for any composable 1v-cells $u$ and $v$ in $\Aa$:
$$\text{{\em \axiom{h.u.t.-1}}} \label{h.u.t.-1} \qquad\alpha^{\frac{u}{v}}=\frac{\alpha^u}{\alpha^v}\quad\qquad\text{ and}\quad\qquad
\text{{\em \axiom{h.u.t.-2}}} \label{h.u.t.-2} \qquad\alpha^{1^A}=\Id_{\alpha(A)}.$$
\end{defn}

\begin{defn} \delabel{modif-mf}
A {\em mixed funny modification} $\Theta$ between two horizontal unnatural transformations $\alpha$ and $\beta$ and two vertical strict transformations $\alpha_0$ and $\beta_0$ depicted below on the left, where the lax double functors $F, G, F\s', G'$ act between $\Aa\to\Bb$, is given by a collection of 2-cells in $\Bb$ depicted below on the right:
$$
\scalebox{0.86}{
\bfig
\putmorphism(-150,50)(1,0)[F` G`\alpha]{400}1a
\putmorphism(-150,-270)(1,0)[F'`G' `\beta]{400}1b
\putmorphism(-170,50)(0,-1)[\phantom{Y_2}``\alpha_0]{320}1l
\putmorphism(250,50)(0,-1)[\phantom{Y_2}``\beta_0]{320}1r
\put(-30,-140){\fbox{$\Theta$}}
\efig}
\qquad\qquad
\scalebox{0.86}{
\bfig
\putmorphism(-180,50)(1,0)[F(A)` G(A)`\alpha(A)]{550}1a
\putmorphism(-180,-270)(1,0)[F\s'(A)`G'(A) `\beta(A)]{550}1b
\putmorphism(-170,50)(0,-1)[\phantom{Y_2}``\alpha_0(A)]{320}1l
\putmorphism(350,50)(0,-1)[\phantom{Y_2}``\beta_0(A)]{320}1r
\put(0,-140){\fbox{$\Theta_A$}}
\efig}
$$
satisfying for every 1v-cell $u$ the axiom
\noindent {\em \axiom{m.hu\x vs}} 
$$
\scalebox{0.86}{
\bfig
 \putmorphism(-150,500)(1,0)[F(A)`F(A) `=]{600}1a
 \putmorphism(550,500)(1,0)[` `\alpha(A)]{400}1a
\putmorphism(-180,500)(0,-1)[\phantom{Y_2}`F\s'(A) `\alpha_0(A)]{450}1l
\put(30,50){\fbox{$\Id$}}
\putmorphism(-150,-400)(1,0)[F\s'(\tilde A)` `=]{480}1a
\putmorphism(-180,50)(0,-1)[\phantom{Y_2}``F\s'(u)]{450}1l
\putmorphism(450,50)(0,-1)[\phantom{Y_2}`F\s'(\tilde A)` \alpha_0(\tilde A)]{450}1l
\putmorphism(450,500)(0,-1)[\phantom{Y_2}`F(\tilde A) `F(u)]{450}1l
\put(660,280){\fbox{$\alpha^u$}}
\putmorphism(450,50)(1,0)[\phantom{(B, \tilde A)}``\alpha(\tilde A)]{500}1a
\putmorphism(1070,50)(0,-1)[\phantom{(B, A')}`G'(\tilde A)`\beta_0(\tilde A)]{450}1r
\putmorphism(1070,500)(0,-1)[G(A)`G(\tilde A)`G(u)]{450}1r
\putmorphism(450,-400)(1,0)[\phantom{(B, \tilde A)}``\beta(\tilde A)]{500}1a
\put(640,-170){\fbox{$ \Theta_{\tilde A}$ } } 
\efig}\quad=\quad
\scalebox{0.86}{
\bfig
 \putmorphism(-150,500)(1,0)[F(A)`G(A) `\alpha(A)]{600}1a
 \putmorphism(450,500)(1,0)[\phantom{(B,A)}` `=]{460}1a
\putmorphism(-180,500)(0,-1)[\phantom{Y_2}`F\s'(A) `\alpha_0(A)]{450}1l
\put(650,50){\fbox{$\Id$}}
\putmorphism(-180,-400)(1,0)[F\s'(\tilde A)` `\beta(\tilde A)]{500}1a
\putmorphism(-180,50)(0,-1)[\phantom{Y_2}``F\s'(u)]{450}1l
\putmorphism(450,50)(0,-1)[\phantom{Y_2}`G'(\tilde A)`G'(u)]{450}1r
\putmorphism(450,500)(0,-1)[\phantom{Y_2}`G'(A) `\beta_0(A)]{450}1r
\put(0,260){\fbox{$\Theta_A$}}
\putmorphism(-180,50)(1,0)[\phantom{(B, \tilde A)}``\beta(A)]{500}1a
\putmorphism(1030,50)(0,-1)[\phantom{(B, A')}` G'(\tilde A). ` \beta_0(\tilde A)]{450}1r
\putmorphism(1030,500)(0,-1)[G(A)`G(\tilde A)` G(u)]{450}1r
\putmorphism(430,-400)(1,0)[\phantom{(B, \tilde A)}``=]{480}1b
\put(70,-170){\fbox{$\beta^u$}}
\efig}
$$ 
\end{defn}


The above definition of unnatural horizontal transformations is coherent with its 2-categorical version indeed. 
For identity 1v-cells $u$ the 2-cell components $\alpha^u$ of horizontal unnatural transformations are trivial by \axiomref{h.u.t.\x 2} and one recovers unnatural transformations for 2-categories. 

Let $[\Aa,\Bb]$ denote the double category of {\em strict} double functors, horizontal unnatural transformations in the sense of 
\deref{hor unnat tr}, vertical strict transformations and mixed funny modifications. By $[\Aa,\Bb]^{lx}$ (resp. $[\Aa,\Bb]^{ps}$) we denote the double category differing from the latter in that its objects are lax (resp. pseudo) double functors $\Aa\to\Bb$. The following result is straightforwardly obtained, it is a mixed funny version of \cite[Proposition 3.3]{Fem:Bif}. The reader may consult Table \ref{table:12} to keep the track of the relevant axioms (minding the fact that the horizontal transformations 
are now unnatural and that the vertical transformations are now strict).

\begin{prop} \prlabel{char lax fdf} 
A lax double functor $\F\colon\Aa\to[\Bb, \Cc]^{lx}$ of double categories consists of the following: 
\begin{enumerate}
\item two families of lax double functors 
$(-,A)\colon\Bb\to\Cc\quad\text{ and}\quad (B,-)\colon\Aa\to\Cc$ for objects $A\in\Aa, B\in\Bb$, 
such that $H(A,-)=(-, A), H(-, B)=(B,-)$ and $(-,A)\vert_B=(B,-)\vert_A=(B,A)$, and 
\item two families of 2-cells 
$$
\scalebox{0.78}{
\bfig
\putmorphism(-150,50)(1,0)[(B,A)`(B,A')`(B,K)]{600}1a
\putmorphism(-150,-400)(1,0)[(\tilde B, A)`(\tilde B,A') `(\tilde B,K)]{640}1a
\putmorphism(-180,50)(0,-1)[\phantom{Y_2}``(u,A)]{450}1l
\putmorphism(450,50)(0,-1)[\phantom{Y_2}``(u,A')]{450}1r
\put(-20,-180){\fbox{$(u, K)$}}
\efig}
\quad
\scalebox{0.78}{
\bfig
\putmorphism(-150,50)(1,0)[(B,A)`(B',A)`(k,A)]{600}1a
\putmorphism(-150,-400)(1,0)[(B, \tilde A)`(B', \tilde A) `(k,\tilde A)]{640}1a
\putmorphism(-180,50)(0,-1)[\phantom{Y_2}``(B,U)]{450}1l
\putmorphism(450,50)(0,-1)[\phantom{Y_2}``(B',U)]{450}1r
\put(0,-180){\fbox{$(k,U)$}}
\efig}
$$
and a family of (identity) 2-cells 
$$
\scalebox{0.86}{
\bfig
 \putmorphism(-150,500)(1,0)[(B,A)`(B,A) `=]{600}1a
\putmorphism(-180,500)(0,-1)[\phantom{Y_2}`(B, \tilde A) `(B,U)]{450}1l
\put(-20,50){\fbox{$(u,U)$}}
\putmorphism(-150,-400)(1,0)[(\tilde B, \tilde A)`(\tilde B, \tilde A) `=]{640}1a
\putmorphism(-180,50)(0,-1)[\phantom{Y_2}``(u,\tilde A)]{450}1l
\putmorphism(450,50)(0,-1)[\phantom{Y_2}``(\tilde B, U)]{450}1r
\putmorphism(450,500)(0,-1)[\phantom{Y_2}`(\tilde B, A) `(u,A)]{450}1r
\efig}
$$ 
in $\Cc$ determined by all 1h-cells $K\colon A\to A'$ and 1v-cells $U\colon A\to\tilde A$ in $\Aa$, and 1h-cells $k\colon B\to B'$ 
and 1v-cells $u\colon B\to\tilde B$ in $\Bb$, which satisfy the following 11 axioms from \prref{char df}: \vspace{-0,2cm}
\begin{center}
\axiomref{($u,1_A$)}, \,\, \axiomref{($1_B,U$)}, \,\, \axiomref{($1^B,K$)}, \,\, \axiomref{($k,1^A$)}, \\ 
\axiomref{($u, K'K$)}, \,\, \axiomref{($k'k, U$)}, \,\, \axiomref{($\frac{u}{u'}, K$)}, \,\, \axiomref{($k,\frac{U}{U'}$)}, \\
$(u,U)=\Id$, \,\, \axiomref{$(u,U)$-l-nat}, \,\, \axiomref{$(u,U)$-r-nat}. 
\end{center}
\end{enumerate}
\end{prop}

\begin{prop} \prlabel{char fdf} 
A double functor $\F\colon\Aa\to[\Bb, \Cc]$ is a lax double functor $\F\colon\Aa\to[\Bb, \Cc]^{lx}$ in whose data the 
double functors $(-,A)\colon\Bb\to\Cc$ and $(B,-)\colon\Aa\to\Cc$ are strict. In this case the axioms 
\axiomref{($u,1_A$)}, \, \axiomref{($1_B,U$)}, \, \axiomref{($u, K'K$)}, \, \axiomref{($k'k, U$)} 
 simplify into the expressions:
$$(u,1_A)=\Id^{(u,A)}, \,\, (1_B,U)=\Id^{(B,U)}, \,\, (u,K'K)=[(u,K)\vert(u,K')], \,\, (k'k,U)=[(k,U)\vert(k',U)].$$
\end{prop}

Analogously to \prref{char lax fdf} one characterizes a pseudodouble functor $\F\colon\Aa\to[\Bb, \Cc]^{ps}$. It satisfies formally same 
11 axioms: with the only change that the structure 2-cells of the pseudodouble functors $(-,A)\colon\Bb\to\Cc$ and $(B,-)\colon\Aa\to\Cc$ are invertible.

\begin{defn} \delabel{funny functor}
A {\em lax funny functor} $\Aa\times\Bb\to\Cc$ between double categories (or a lax double funny functor) is a 
pair of lax double functors $(-,A):\Bb\to\Cc$ and $(B,-):\Aa\to\Cc$ satisfying the conditions of \prref{char lax fdf}. 

A {\em pseudodouble funny functor} $\Aa\times\Bb\to\Cc$ is analogously defined. 

A {\em funny functor} $\Aa\times\Bb\to\Cc$ between double categories (or a double funny functor) is a 
pair of double functors $(-,A):\Bb\to\Cc$ and $(B,-):\Aa\to\Cc$ satisfying the conditions of \prref{char fdf}. 
\end{defn}

Strictly speaking the above funny functors we should call {\em mixed funny functors}, as in the vertical direction they are not determined by 
unnatural transformations. We will use the abbreviated version of the name. 

We proceed to construct a funny type of product of double categories $\Aa$ and $\Bb$. One obtains a strict funny product 
$\Aa\Box\Bb$ by considering the inner-hom $[\Bb,\Cc]$, while a lax version $\Aa\Box^{lx}\Bb$ is obtained by considering the 
inner-hom $[\Bb,\Cc]^{lx}$ (and similarly for $\Aa\Box^{ps}\Bb$). 
We will show now the construction for the first one (the strict one), as it is this one that is going to give a monoidal category structure (of funny type) to the category of double categories. 

\smallskip

To obtain a funny product $\Aa\Box\Bb$, we suppose there is a left adjoint $-\Box \Bb$ to $[\Bb,-]$. 
For this set $\Cc=\Aa\times\Bb$, consider the unit of the adjunction 
$E: \Aa\to[\Bb,\Aa\times\Bb]$ and read off the structure of the image double category $E(\Aa)(\Bb)$ in $\Aa\times\Bb$. It is a double category $\Aa\Box\Bb$ generated by certain pairs $(x,y)=:x\Box y$, where $x\in\Aa$ and $y\in\Bb$ are 0-, 1h-, 1v- or 
2-cells, that we also get as images of the above two double functors $(-,A)$ and $(B,-)$ (from the end of \deref{funny functor}). Although we do not know the above pair of functors nor $E$, knowing that 
$\Aa\Box\Bb$ is generated as the image by double functors $(-,A)$ and $(B,-)$, gives us the hint on how to define it. 

Namely, we define a funny product $\Aa\Box\Bb$ by the following generators and relations: \\ 
\ul{objects}: $A\Box B$ for objects $A\in\Aa, B\in\Bb$; \\ 
\ul{1h-cells}: $A\Box k, K\Box B$, where $k$ is a 1h-cell in $\Bb$ and $K$ a 1h-cell in $\Aa$; \\ 
\ul{1v-cells}: $A\Box u, U\Box B$ and vertical compositions of such obeying the following rules: 
\begin{equation} \eqlabel{1v-cells}
\frac{A\Box u}{A\Box u'}=A\Box \frac{u}{u'}, \quad \frac{U\Box B}{U'\Box B}=\frac{U}{U'}\Box B, \quad A\Box 1^B=1^{A\Box B}=1^A\Box B
\end{equation}
where $u,u'$ are 1v-cells of $\Bb$ and $U,U'$ 1v-cells of $\Aa$; \\
\ul{2-cells}: $A\Box\omega, \zeta\Box B$, where $\omega$ is a 2-cell in $\Bb$ and $\zeta$ is a 2-cell in $\Aa$, further 2-cells 
\begin{equation} \eqlabel{1h-1v}
\scalebox{0.86}{
\bfig
\putmorphism(-150,50)(1,0)[A\Box B`A\Box B'`A\Box k]{600}1a
\putmorphism(-150,-400)(1,0)[\tilde A\Box B`\tilde A\Box B' `\tilde A\Box k]{640}1a
\putmorphism(-180,50)(0,-1)[\phantom{Y_2}``U\Box B]{450}1l
\putmorphism(450,50)(0,-1)[\phantom{Y_2}``U\Box B']{450}1r
\put(0,-180){\fbox{$U\Box k$}}
\efig}
\qquad
\scalebox{0.86}{
\bfig
\putmorphism(-150,50)(1,0)[A\Box B`A'\Box B`K\Box B]{600}1a
\putmorphism(-150,-400)(1,0)[A\Box\tilde  B` A'\Box\tilde  B `K\Box\tilde B]{640}1a
\putmorphism(-180,50)(0,-1)[\phantom{Y_2}``A\Box u]{450}1l
\putmorphism(450,50)(0,-1)[\phantom{Y_2}``A'\Box u]{450}1r
\put(0,-180){\fbox{$K\Box u$}}
\efig}
\qquad
\scalebox{0.86}{
\bfig
 \putmorphism(-150,500)(1,0)[(B,A)`(B,A) `=]{600}1a
\putmorphism(-180,500)(0,-1)[\phantom{Y_2}`(B, \tilde A) `(B,U)]{450}1l
\put(-20,50){\fbox{$(u,U)$}}
\putmorphism(-150,-400)(1,0)[(\tilde B, \tilde A)`(\tilde B, \tilde A) `=]{640}1a
\putmorphism(-180,50)(0,-1)[\phantom{Y_2}``(u,\tilde A)]{450}1l
\putmorphism(450,50)(0,-1)[\phantom{Y_2}``(\tilde B, U)]{450}1r
\putmorphism(450,500)(0,-1)[\phantom{Y_2}`(\tilde B, A) `(u,A)]{450}1r
\efig}
\end{equation}
subject to the eleven relations: \\
\begin{equation} \eqlabel{hor-ver comp}
K\Box 1^B=Id_{K\Box B}, \,\, 1^A\Box k=Id_{A\Box k}, \,\, K\Box \frac{u}{u'}=\frac{K\Box u}{K\Box u'}, \,\, 
\frac{U}{U'}\Box k=\frac{U\Box k}{U'\Box k},
\end{equation}
\begin{equation} \eqlabel{four strictified}
1_A\Box u=\Id^{A\Box u}, \,\, U\Box 1_B=\Id^{U\Box B}, \,\, K'K\Box u=[K\Box u\vert K'\Box u], \,\, U\Box k'k=[U\Box k\vert U\Box k'];
\end{equation}
\begin{equation} \eqlabel{frac-natur}
U\Box u=\Id, \,\,\,\, \frac{A\Box\omega}{U\Box l}=\frac{U\Box k}{\tilde A\Box\omega}, \,\,\,\, \frac{K\Box u}{\zeta\Box\tilde B}=
\frac{\zeta\Box B}{L\Box u};
\end{equation}
four equations holding from the strictness of double functors $(-,A)$ and $(B,-)$: 
\begin{equation} \eqlabel{4 strictness 2-cells}
(A\Box k')(A\Box k)= A\Box (k' k), \quad (K'\Box B)(K\Box B)=(K' K)\Box B
\end{equation}
$$1_{A\Box B}= A\Box 1_B, 
\quad 1_{A\Box B}= 1_A\Box B$$
and the following ones: 
\begin{equation} \eqlabel{horiz comp 2-cells}
A\Box \omega'\omega=[A\Box\omega \vert A \Box \omega'], \quad \zeta'\zeta\Box B=[\zeta\Box B \vert\zeta' \Box B]
\end{equation}
\begin{equation} \eqlabel{vertic comp 2-cells}
A\Box\frac{\omega}{\omega'}=\frac{A\Box\omega}{A\Box\omega'}, \quad \frac{\zeta}{\zeta'}\Box B=\frac{\zeta\Box B}{\zeta'\Box B},
\end{equation}
\begin{equation} \eqlabel{Box on Id's}
A\Box\Id_k=\Id_{A\Box k}, \quad \Id_K\Box B=\Id_{K\Box B}, \quad A\Box\Id^u=\Id^{A\Box u}, \quad \Id^U\Box B=\Id^{U\Box B} .
\end{equation}
The source and target functors $s,t$ on $\Aa\Box\Bb$ are defined as in the double category $\Aa\times\Bb$, the composition functor $c$ is defined by horizontal juxtaposition of 
the corresponding 2-cells, and the unit functor $i$ is defined on generators as follows: 
$$i(A\Box B)=1_{A\Box B}, \, i(A\Box v)=1^A\Box v (=Id^{A\Box v}) \quad \text{and} \quad i(U\Box B)=U\Box 1^B (=Id^{U\Box B}).$$ 

Now it is straightforward to see that $-\Box\Bb: Dbl\to Dbl$ (resp. $\Bb\Box-: Dbl\to Dbl$) and $[\Bb, -]$ are double functors. Then by construction we have bijections
$$
Dbl(\Aa\Box\Bb,\Cc)\iso f\x Dbl(\Aa\times\Bb,\Cc)\iso Dbl(\Aa, [\Bb, \Cc]),
$$
where $Dbl(\Aa, \Bb)$ denotes a set of double functors $\Aa\to\Bb$ and $f\x Dbl(\Aa\times\Bb,\Cc)$ is the set of double funny functors.

\begin{rem} \rmlabel{functors on Box}
We defined $-\Box-:Dbl\times Dbl\to Dbl$ on objects. By the above bijection a double functor $F\Box G:\Aa\Box\Bb\to\Aa'\Box\Bb'$ is 
given via two double functors and three 2-cells satisfying eleven equations, see \prref{char fdf}. We define $(-,A):\Bb\to\Aa'\Box\Bb'$ 
and $(B,-):\Bb\to\Aa'\Box\Bb'$ by setting $(b,A)=F(A)\Box G(b)$ and $(B,a)=F(a)\Box G(B)$ for any cell $a\in\Aa$ and $b\in\Bb$. That 
these are two double functors it is proved using axioms \equref{1v-cells}, \equref{4 strictness 2-cells} -- \equref{Box on Id's} in the funny product. We introduce the following three types of 2-cells: $(u,K):=F(K)\Box G(u), \, \, (k,U):=F(k)\Box G(U), \,\, (u,U):=F(U)\Box G(u)$. 
Then the eleven axioms hold by the eleven rules \equref{hor-ver comp} -- \equref{frac-natur}. 
For so defined double functor $F\Box G$ 
we may shortly write that $(F\Box G)(a\Box b)=F(a)\Box G(b)$ for sensible cells $a\Box b\in\Aa\Box\Bb$. This clearly makes 
$-\Box-:Dbl\times Dbl\to Dbl$ a functor. 
(Observe that from the functor property of $-\Box -$ for double functors $F:\Aa\to\Aa'$ and $G:\Bb\to\Bb'$ we have in particular 
$(\Id_{\Aa'}\Box G)(F\Box \Id_{\Bb})=(F\Box \Id_{\Bb'})(\Id_{\Aa}\Box G)=F\Box G$, 
even though for 1h- and 1v-cells $f:A\to A'$ in $\Aa$ and $g:B\to B'$ in $\Bb$ the two compositions $(A'\Box g)(f\Box B)$ and 
$(f\Box B')(A\Box g)$ mapping $A\Box B\to A'\Box B'$ differ.) 
\end{rem}

It is readily seen that the above bijections are natural in all variables. 
We defer the proof that $(Dbl, \Box)$ is a monoidal category for \ssref{funny mult}, \thref{Dbl left closed mon}, then it will be 
a biclosed even symmetric monoidal category.

\begin{rem} 
It is immediate to see that there is a double functor $M:\Aa\Box\Bb\to\Aa\times\Bb$ and that $(\Id, M): (Dbl,\Box)\to(Dbl,\times)$ 
gives a symmetric oplax structure on the identity functor on $Dbl$. Moreover, it can be shown that $M$ factors as the composition 
$\Aa\Box\Bb\stackrel{P}{\to}\Aa\ot_{st}^{ps}\Bb \stackrel{Q}{\to}\Aa\times\Bb$, where $\ot_{st}^{ps}$ is the Gray tensor product 
originating from the inner-hom made of strict double functors, (horizontal and vertical) pseudonatural transformations and their modifications, and $\Aa\times\Bb$, the Cartesian monoidal product, corresponds to the inner-hom made of strict double functors, 
(horizontal and vertical) strict transformations and their modifications. The Gray tensor product $\ot_{st}^{ps}$ is a pseudonatural transformation version of the Gray product $\ot=\ot_{st}^{o\x l}$ that we studied in \cite{Fem:Gray1}. 
The above symmetric oplax structure and the factorization are analogous to the 2-categorical results of \cite[Section 2]{BG}.
\end{rem}



In the general case, where $*$ stands for any kind of funny functors (be it strict, lax, colax or pseudo), funny functors of type $*$ are 
by construction such that there is a bijection 
\begin{equation}\eqlabel{funny-basic}
f\x\Fun^*(\Aa\times\Bb,\Cc)\iso \Fun^*(\Aa, [\Bb,\Cc]^*).
\end{equation} 
Here $f\x\Fun^*(\Aa, \Bb)$ and $\Fun^*(\Aa, \Bb)$ stand for the sets of double funny functors and double functors of type $*$, respectively, 
and $[\Bb, \Cc]^*$ is the $*$-typed double functor version of $[\Bb, \Cc]$.  
Also, the induced funny product fits the bijection 
$$Dbl(\Aa\Box^*\Bb,\Cc)\iso\Fun^*(\Aa,[\Bb,\Cc]^*)$$
so that one gets 
\begin{equation} \eqlabel{no-rep}
Dbl(\Aa\Box^*\Bb,\Cc)\iso f\x\Fun^*(\Aa\times\Bb,\Cc).
\end{equation}
We state for the record that $[-,-]^*$ and $-\Box^*-$ are functors only on $Dbl$ and $Dbl_{ps}$, the categories of double categories with strict or pseudodouble functors, but they are not functors on $Dbl_{lx}$ nor $Dbl_{clx}$ 
(where morphisms are lax or colax double functors). (For the reason of this failure see \cite[Section 3.1]{Fem:Bif}.) 


\medskip

Given that $\Aa\Box^*\Bb$ is defined by generators and relations on $\Aa\times\Bb$, there is a double funny functor 
\begin{equation} \eqlabel{J}
J^*\colon\Aa\times\Bb\to\Aa\Box^*\Bb
\end{equation} 
of type $*$ given by $J^*(-, B)(a)=a\Box B, J^*(A,-)(b)=A\Box b$ for cells $a$ in $\Aa$ and $b$ in $\Bb$ and with unique 
2-cells $U\Box k, K\Box u$ satisfying the eleven axioms. 

\medskip

We end this subsection by giving another characterization of lax double funny functors analogous to that of lax double quasi-functors from \cite[Proposition 3.3]{Fem:Bif} (see \prref{char df}). 

\begin{prop} \prlabel{funny fun}
Let $\Aa,\Bb,\Cc$ be double categories. The following are equivalent:
\begin{enumerate} 
\item $H\colon \Aa\times\Bb\to\Cc$ is a lax double funny functor, 
\item there are two families of lax double functors 
$(-,A)\colon\Bb\to\Cc\quad\text{ and}\quad (B,-)\colon\Aa\to\Cc$ for objects $A\in\Aa, B\in\Bb$, 
such that $(-,A)\vert_B=(B,-)\vert_A=(B,A)$, and the following hold: 
\begin{enumerate}[(i)]
\item $(-,K)\colon (-,A)\to(-,A')$ is a horizontal unnatural transformation for each 1h-cell $K\colon A\to A'$, 
$(-,U)\colon (-,A)\to(-,\tilde A)$ is a vertical strict transformation for each 1v-cell $U\colon A\to\tilde A$ in $\Aa$,  
$(-,\zeta)$ is a mixed funny modification with respect to horizontally unnatural and vertically strict transformations for each 2-cell $\zeta$ in $\Aa$, 
and the following coincide:
$$(B,-)\vert_K=(-,K)\vert_B, \,\,\,\, (B,-)\vert_U=(-,U)\vert_B, \,\,\,\, (B,-)\vert_\zeta=(-,\zeta)\vert_B;$$
\item $(k,-)\colon (B,-)\to (B', -)$ is a horizontal unnatural transformation for each 1h-cell $k\colon B\to B'$, 
$(u,-)\colon (B,-)\to(\tilde B,-)$ is a vertical strict transformation for each 1v-cell $u\colon B\to\tilde B$ in $\Bb$, 
$(\omega,-)$ is a mixed funny modification with respect to horizontally unnatural and vertically strict transformations for each 2-cell 
$\omega$ in $\Bb$, 
and the following coincide:
$$(-,A)\vert_k=(k,-)\vert_A, \,\,\,\, (-,A)\vert_u=(u,-)\vert_A, \,\,\,\, (-,A)\vert_\omega=(\omega,-)\vert_A;$$
\item for 1h-cells $K,k$ and 1v-cells $U,u$ the following 2-cell components of the respective transformations coincide:
$(-,K)\vert_u=(u,-)\vert_K$ and $(-,U)\vert_k=(k,-)\vert_U$.
\end{enumerate}
\end{enumerate}
\end{prop}

\begin{cor} \colabel{1v-cells are central}
A lax double funny functor $H\colon \Bb\times\Bb\to\Bb$ equips all 1v-cells of $\Bb$ with a structure of central cells 
so that the square-formed 2-cell components of the vertical strict transformations $(-,U)$ and $(u,-)$ obey $(-,K)\vert_u=(u,-)\vert_K=H(K,u)$ and $(-,U)\vert_k=(k,-)\vert_U=H(U,k)$, where $(-,K), (k,-)$ are horizontal unnatural transformations. 
\end{cor}

The pseudo (double funny functor) version of the above proposition for $\Aa=\Bb=\Cc$ provides a characterization of a binoidal structure coming from a pseudodouble funny functor: 

\begin{cor} \colabel{funny yields binoidal}
There is a pseudodouble funny functor $H:\Bb\times\Bb\to\Bb$ 
if and only if 
\begin{itemize}
\item $\Bb$ is binoidal; \vspace{-0,2cm} 
\item there are 1-cells 
\vspace{-0,2cm}
$$-\rtimes k\vert_A:=A\ltimes-\vert_k, \,\,\,\, -\rtimes u\vert_A:=A\ltimes-\vert_u \qquad 
K\ltimes-\vert_B:=-\rtimes B\vert_K, \,\,\,\, U\ltimes-\vert_B:=-\rtimes B\vert_U; \vspace{-0,3cm}$$ 
\item there are 2-cells $K\ltimes-\vert_u=-\rtimes u\vert_K$ functorial both in $k$ and in $u$, and 
$U\ltimes-\vert_k=-\rtimes k\vert_U$ functorial both in $U$ and in $k$;  \vspace{-0,2cm}
\item all 1v-cells $u$ are central in $\Bb$ via vertical strict transformations $u\ltimes-$ and $-\rtimes u$; \vspace{-0,2cm}
\item there are 2-cells 
$-\rtimes \omega\vert_A:=A\ltimes-\vert_\omega$ and $\zeta\ltimes-\vert_B:=-\rtimes B\vert_\zeta$ 
natural with respect to 1v-cells (in the sense of \axiomref{m.hu\x vs}), 
\end{itemize} \vspace{-0,14cm}
for all objects $A,B$, 1h-cells $K,k$, 1v-cells $U,u$ and 2-cells $\zeta, \omega$ in $\Bb$. 
\end{cor}



\subsubsection{Purely funny product} \ssslabel{pure funny}

All the results of the above part of \ssref{funny} have their counterparts for the funny product obtained from the inner-hom $[\Aa,\Bb]_f$ whose 1h- and 1v-cells are horizontal and vertical purely unnatural transformations. In this subsection we record these results. The  obtained funny product we denote by $\Aa\Box_f\Bb$ and refer to it as a {\em purely funny product}. 

The ''pure'' version of \prref{char lax fdf} states that a lax double functor $\F\colon\Aa\to[\Bb, \Cc]_f^{lx}$ of double categories consists merely of two families of lax double functors 
$(-,A)\colon\Bb\to\Cc$ and $(B,-)\colon\Aa\to\Cc$ for objects $A\in\Aa, B\in\Bb$, 
such that $H(A,-)=(-, A), H(-, B)=(B,-)$ and $(-,A)\vert_B=(B,-)\vert_A=(B,A)$. 
This short data we may call {\em purely (lax) funny functor} $\Aa\times\Bb\to\Cc$ of double categories.

Namely, by substituting strict vertical transformations by unnatural ones, the axioms $(u,u)=\Id$, 
\axiomref{$(u,U)$-l-nat} and \axiomref{$(u,U)$-r-nat} are eliminated. Furthermore, by using purely unnatural transformations one 
eliminates square-formed 2-cells and $2+2$ axioms in both directions (\axiomref{h.u.t.-1}, \axiomref{h.u.t.-2} and their vertical counterparts) stemming from both variables in a (mixed) funny functor. In total 11 axioms from the definition of a mixed funny functor are omitted. Consequently, 11 axioms from the definition of the funny product $\Aa\Box\Bb$ are omitted to obtain the purely funny product $\Aa\Box_f\Bb$. It will have the same generators on objects and 1-cells, but its only generating 2-cells will be 
$A\Box_f\omega$ and $\zeta\Box_f B$. 

The ''pure'' versions of bijections \equref{funny-basic} and \equref{no-rep} read: 
$$pf\x\Fun^*(\Aa\times\Bb,\Cc)\iso \Fun^*(\Aa, [\Bb,\Cc]_f^*) \qquad\text{and}\qquad 
Dbl(\Aa\Box_f^*\Bb,\Cc)\iso pf\x\Fun^*(\Aa\times\Bb,\Cc).$$

Finally, \prref{funny fun} in its ``pure'' version states: 

\begin{prop} 
Let $\Aa,\Bb,\Cc$ be double categories. The following are equivalent:
\begin{enumerate} 
\item $H\colon \Aa\times\Bb\to\Cc$ is a lax double purely funny functor, 
\item there are two families of lax double functors 
$(-,A)\colon\Bb\to\Cc\quad\text{ and}\quad (B,-)\colon\Aa\to\Cc$ for objects $A\in\Aa, B\in\Bb$, 
such that $(-,A)\vert_B=(B,-)\vert_A=(B,A)$, and the following hold: 
\begin{enumerate}[(i)]
\item $(-,K)\colon (-,A)\to(-,A')$ is a horizontal purely unnatural transformation for each 1h-cell $K\colon A\to A'$, 
$(-,U)\colon (-,A)\to(-,\tilde A)$ is a vertical purely unnatural transformation for each 1v-cell $U\colon A\to\tilde A$ in $\Aa$,  
$(-,\zeta)$ is a purely funny modification for each 2-cell $\zeta$ in $\Aa$, and the following coincide:
$$(B,-)\vert_K=(-,K)\vert_B, \,\,\,\, (B,-)\vert_U=(-,U)\vert_B, \,\,\,\, (B,-)\vert_\zeta=(-,\zeta)\vert_B;$$
\item $(k,-)\colon (B,-)\to (B', -)$ is a horizontal purely unnatural transformation for each 1h-cell $k\colon B\to B'$, 
$(u,-)\colon (B,-)\to(\tilde B,-)$ is a vertical purely unnatural transformation for each 1v-cell $u\colon B\to\tilde B$ in $\Bb$, 
$(\omega,-)$ is a purely funny modification for each 2-cell $\omega$ in $\Bb$, and the following coincide:
$$(-,A)\vert_k=(k,-)\vert_A, \,\,\,\, (-,A)\vert_u=(u,-)\vert_A, \,\,\,\, (-,A)\vert_\omega=(\omega,-)\vert_A.$$
\end{enumerate}
\end{enumerate}
\end{prop}

The meaning of points (i) and (ii) in the above proposition is that the component 1-cells of horizontal and vertical 
purely unnatural transformations and the component 2-cells of purely funny modifications coincide with the corresponding 
images of the lax double functors $(-,A)\colon\Bb\to\Cc$ and $(B,-)\colon\Aa\to\Cc$. In particular, the pseudo 
(double purely funny) version of the proposition with $\Aa=\Bb=\Cc$ implies that in a premonoidal double category in which the binoidal structure comes from a pseudodouble purely funny functor, we have that there exist 1- cells 
$-\rtimes k\vert_A:=A\ltimes-\vert_k, \,\,\,\, -\rtimes u\vert_A:=A\ltimes-\vert_u, \quad 
K\ltimes-\vert_B:=-\rtimes B\vert_K, \,\,\,\, U\ltimes-\vert_B:=-\rtimes B\vert_U$ 
and 2-cells $-\rtimes \omega\vert_A:=A\ltimes-\vert_\omega$ and $\zeta\ltimes-\vert_B:=-\rtimes B\vert_\zeta$ 
for which no further laws are required. Then we may conclude:

\begin{cor} \colabel{purely funny binoidal}
Any binoidal structure in a double category $\Bb$ is given by a pseudodouble purely funny functor 
$H:\Bb\times\Bb\to\Bb$. 
\end{cor}

\subsection{Relation to the funny product for 2-categories} \sslabel{funny 2-cats}

Funny product for 2-categories was studied in \cite[Section 2]{BG}. Its construction follows an analogous process as our construction 
of the funny product for double categories, starting from an inner-hom that there was denoted by $[\A,\B]_f$ for 2-categories $\A,\B$. 
The inner-hom  $[\A,\B]_f$ has for objects 2-functors, for 1-cells transformations (consisting only of 1-cell components subject to no axioms), and for 2-cells modifications (consisting only of 2-cell components subject to no axioms).  
The obtained funny tensor product there was denoted by $\A\star\B$, it satisfies the universal property 
$$2\x\Cat(\A\star\B,\C)\iso 2\x\Cat(\A,[\B,\C]_f)$$ 
and gives a symmetric closed monoidal structure on the category $2\x\Cat$ of 2-categories. Analogously as in 
\sssref{pure funny}, 
a funny functor $\A\star\B\to\C$ between 2-categories is given merely by a pair of families of 2-functors $(-,A):\B\to\C$ and 
$(B,-):\A\to\C$ for $A\in\A,B\in\B$ satisfying $(-,A)_B=(B,-)_A=(B,A)$. An explicit description of $\A\star\B$ can be deduced from our description of 
$\A\Box\B$ for double categories by considering all vertical 1-cells to be identities. Thus it is given by the following generators and relations: \\
\ul{objects}: $A\star B$ for objects $A\in\A, B\in\B$; \\ 
\ul{1-cells}: $A\star k, K\star B$, where $k$ is a 1-cell in $\B$ and $K$ a 1-cell in $\A$; \\ 
\ul{2-cells}: $A\star\omega, \zeta\star B$, where $\omega$ is a 2-cell in $\Bb$ and $\zeta$ a 2-cell in $\A$; \\
four equations from the strictness of double functors $(-,A)$ and $(B,-)$: 
$$
(A\star k')(A\star k)= A\star (k' k), \quad 
(K'\star B)(K\star B)=(K' K)\star B
$$
$$1_{A\star B}= A\star 1_B, 
\quad 1_{A\star B}= 1_A\star B$$
and the following ones: 
$$A\star \omega'\omega=[A\star\omega \vert A \star \omega'], \quad 
\zeta'\zeta\star B=[\zeta\star B \vert\zeta' \star B]$$
$$A\star\frac{\omega}{\omega'}=\frac{A\star\omega}{A\star\omega'}, \quad \frac{\zeta}{\zeta'}\star B=\frac{\zeta\star B}{\zeta'\star B},$$
$$A\star\Id_k=\Id_{A\star k}, \quad \Id_K\star B=\Id_{K\star B}.$$
Likewise, the funny product on morphisms, {\em i.e.} on 2-functors, is defined by $(F\star G)(a\star b)=F(a)\star G(b)$ for sensible cells $a\star b\in\A\star\B$. This way we have a functor $-\star-:2\x\Cat\times 2\x\Cat\to 2\x\Cat$. 
Inspecting the above generators of $\A\star\B$ we see that $F\star G$ is indeed given by a pair of 2-functors 
$F\star\B, \A\star G$. 
This is coherent with the fact that strictly speaking $F\star G$ does not exist as a 1-cell in the funny product 2-category 
$(2\x\Cat)_2\star(2\x\Cat)_2$, where $(2\x\Cat)_2$ is the 2-category of 2-categories. 
Rather, instead there are only 1-cells $F\star\B$ and $\A\star G$. 

\begin{rem} \rmlabel{two funny for dbl}
We find that the same occurs with the purely funny product of double functors $F\Box_f G$. 
Let 
$Dbl_\bullet$ denote the 2-category of double categories, double functors and some $\bullet$-type of vertical 
transformations. In the funny product 2-category $Dbl_\bullet\star Dbl_\bullet$ 
if we set $\Aa\star\Bb:=\Aa\Box_f\Bb$ on the objects, the 1-cell defined by 
$F\star G:=F\Box_f G:\Aa\Box_f\Bb\to\Aa'\Box_f\Bb'$ for 1-cells $F,G$ presents a valid 1-cell in $Dbl_\bullet\star Dbl_\bullet$: 
from the generators of $\Aa\Box_f\Bb$ we see that $F\star G$ consists precisely of double functors $F\star\Bb$ and $\Aa\star G$. 

However, the double functor $F\Box G:\Aa\Box\Bb\to\Aa'\Box\Bb'$ from \rmref{functors on Box} 
can act on the 2-cell of the form $K\Box u$ existing in the double category $\Aa\Box\Bb$, and it is more than a pair of ``one-sided-functors''. Thus taking the mixed funny product $\Aa\Box\Bb$ 
for a candidate for objects $\Aa\star\Bb$ in $Dbl_\bullet\star Dbl_\bullet$ leads to having much richer 1-cells $F\Box G$, while   the 2-category $Dbl_\bullet\star Dbl_\bullet$ only allows for their restrictions $F\star\Bb$ and $\Aa\star G$. 
%
\end{rem}

\subsection{Funny multicategories} \sslabel{funny mult}

It was shown in \cite[Theorem 9.8]{Her} that monoidal categories are in 1-1 correspondence with representable multicategories. 
We recall the reader of the notion of a representable multicategory.

\begin{defn}
A multicategory consists of:
\begin{itemize} 
\item a collection of objects, \vspace{-0,2cm}
\item for each list of objects $a_1,...,a_n$ for $n\geq 0$ and an object $b$, a set $\M_n(a_1,...,a_n;b)$, \vspace{-0,2cm}
\item for each object $a$ an element $1_a\in\M_1(a;a)$, \vspace{-0,2cm}
\item for all lists $\crta{a_i}, i=1,...,n$, and objects $b_1,...,b_n$ and $c$ a function, called {\em substitution}:
$$\M_n(b_1,...,b_n,c)\times\Pi_{j=1}^{n}\M_{k_i}(\crta{a_i};b_i)\to\M_{\sum_{i=1}^n k_i}(\crta{a_1},...,\crta{a_n};c)$$
$$(g, f_1,...,f_n)\mapsto g\circ(f_1,...,f_n)$$
satisfying a natural associativity and two identity axioms. Here $\crta{a_i}$ is a short annotation for a list of objects, 
and the elements of $\M_n(a_1,...,a_n;b)$ are called {\em multimaps}. 
\end{itemize}
\end{defn}

\begin{defn}
A multicategory $\M$ is said to be {\em representable} if for any list of objects $\crta a=a_1,..,a_n$ there exists an object 
$m(a_1,..,a_n)$ and a multimap $j_{\crta a}:a_1,..,a_n\to m(a_1,..,a_n)$ that for every object $c$ and lists $\crta x, \crta b$ (of lengths 
$k$ and $l$, respectively) induces bijections 
$$\M_1(m(a_1,..,a_n);c)\to\M_n(a_1,..,a_n;c)$$
 and 
$$\M_{k+1+l}(\crta x,m(a_1,..,a_n),\crta b;c)\to\M_{k+n+l}(\crta x,a_1,..,a_n,\crta b;c).$$
\end{defn}

We also recall when a multicategory is closed. 

\begin{defn}
A multicategory $\M$ is said to be {\em left closed} if for all objects $b,c$ there is an object $[b,c]$ and a binary map 
$e_{b,c}:[b,c], b\to c$ so that the induced functions $\M_n(\crta a;[b,c])\to\M_{n+1}(\crta a,b;c)$ are bijections for all $n\geq 0$. 
If there exists an object $\{b,c\}$ and a binary map $\crta{e_{b,c}}:b,\{b,c\}\to c$ for which the induced functions 
$\M_n(\crta a;\{b,c\})\to\M_{n+1}(b,\crta a;c)$ are bijections for all $n\geq 0$, then $\M$ is said to be {\em right closed}. 
$\M$ is called {\em biclosed} if it is both left and right closed. 
\end{defn}

\medskip

We are going to show that (both purely and mixed) funny functors make multimaps for a multicategory. 
Moreover, these multicategories will be represented by their respective funny types of product. Such multicategories we will call {\em funny multicategories}.   

In the first six subsections of \ssref{funny mult} we will study the mixed funny case leading to the (mixed) funny monoidal product $\Box$, omitting the adjective mixed for simplicity. In \sssref{purely funny mult} we will briefly record the corresponding results for the purely funny case of $\Box_f$. The constructions that follow are a funny version of the construction that we carried out in \cite{Fem:Gray1} with quasi-functors for double categories.

\subsubsection{Ternary and $n$-ary (lax) double funny functors}

Analogously as we do in the next definition for double categories, one can define funny functors of more than two variables for categories 
and 2-categories. 

\begin{defn} \delabel{double qf 3}
A {\em lax double funny functor of $n$-variables} $H:\Aa_1\times...\times\Aa_n\to\Cc$ for $n\geq 2$ consists of binary lax double funny functors 
$$H(A_1,...,A_{i-1},\, -,\, A_{i+1},...,A_{j-1}, \, -,\, A_{j+1},..., A_n): \Aa_i\times\Aa_j\to\Cc$$
for all $i<j$ and all choices of objects $A_l\in\Aa_l, l=1,...,n$, which agree on objects as lax double functors of 1-variable (that is, they  give unambiguous lax double functors 
$$H(A_1,...,A_{i-1},\, -,\,A_{i+1},..., A_n): \Aa_i\to\Cc$$
for all $i=1,...,n$), and satisfy the axioms: \\ 
\axiom{u,v,h} 
$$\scalebox{0.86}{
\bfig
 \putmorphism(-120,500)(1,0)[` `=]{550}1a
 \putmorphism(450,500)(1,0)[` `(A,B,h)_k]{550}1a
\putmorphism(-140,520)(0,-1)[` `(A,v,C)_j]{480}1l
\put(60,50){\fbox{$\Id$}}
\putmorphism(-150,-380)(1,0)[` `=]{540}1a
\putmorphism(-140,80)(0,-1)[``(u,\tilde B,C)_i]{450}1l
\putmorphism(430,50)(0,-1)[` `(\tilde A,v,C)_j]{450}1l
\putmorphism(430,520)(0,-1)[` `(u,B,C)_i]{480}1l
\put(490,290){\fbox{$(u,B,h)_{ik}$}}
\putmorphism(430,50)(1,0)[``(\tilde A,B,h)_k]{540}1a
\putmorphism(1000,80)(0,-1)[``(\tilde A,v,C')_j]{450}1r
\putmorphism(1000,520)(0,-1)[``(u,B,C')_i]{480}1r
\putmorphism(450,-380)(1,0)[``(\tilde A,\tilde B,h)_k]{540}1b
\put(500,-190){\fbox{$(\tilde A, v,h)_{jk}$}}
\efig}
\quad=\quad
\scalebox{0.86}{
\bfig
 \putmorphism(-150,500)(1,0)[` `(A,B,h)_k]{600}1a
 \putmorphism(450,500)(1,0)[` `=]{540}1a
\putmorphism(-180,520)(0,-1)[` `(A,v,C)_j]{450}1l
\put(-120,280){\fbox{$(A,v,h)_{jk}$}}
\putmorphism(-150,-380)(1,0)[` `(\tilde A,\tilde B,h)_k]{500}1b
\putmorphism(-180,80)(0,-1)[``(u,\tilde B,C)_i]{450}1l
\putmorphism(400,80)(0,-1)[``(u,\tilde B,C')]{450}1r
\putmorphism(400,520)(0,-1)[` `(A,v,C')_j]{450}1r
\putmorphism(-150,50)(1,0)[``(A,\tilde B,h)_k]{500}1a
\putmorphism(1000,80)(0,-1)[``(\tilde A,v,C')_j]{450}1r
\putmorphism(1000,520)(0,-1)[``(u,B,C')_i]{450}1r
\putmorphism(450,-380)(1,0)[``=]{520}1b
\put(-120,-170){\fbox{$(u, \tilde B,h)_{ik}$}}
\put(600,50){\fbox{$\Id$}}
\efig}
$$
for $(u,v,h):(A,B,C)\to(\tilde A,\tilde B,C')$ in $\Aa_i\times\Aa_j\times\Aa_k$, whereby we omit writing the rest of the $n\x 3$ variables, 
and 2 similar axioms: \axiom{u,g,z} for $(u,g,z):(A,B,C)\to(\tilde A,B',\tilde C)$ and \axiom{f,v,z} for $(f,v,z):(A,B,C)\to(A',\tilde B,C')$, 
where $f,g,h$ are 1h-cells and $u,v,z$ are 1v-cells, as usual.   
\end{defn}

We illustrate the condition that the underlying funny functors of 2-variables 
give unambiguous lax double functors of 1-variables in the example of ternary funny functors $H:\Aa\times\Bb\times\Cc\to\Dd$. This means: 
\begin{equation} \eqlabel{obj}
H(A,-,-)\vert_B=H(-,B,-)\vert_A, H(A,-,-)\vert_C=H(-,-,C)\vert_A, H(-,B,-)\vert_C=H(-,-,C)\vert_B 
\end{equation} 
for $A\in\Aa,B\in\Bb,C\in\Cc$. The above definition can be restated so that an $n$-ary lax double funny functor consists of lax double funny functors of $(n\x 1)$-variables (and hence also of $k$-variables for all $2\leq k<n$) that for any choice of three variables satisfy the condition 
\equref{obj}. We will sometimes refer to the latter condition as the {\em ternary funny property (with functorialities in the third, second and first variable}, respectively.

It is clear from the definition that any kind of double funny functors $H:\Aa_1\times...\times\Aa_n\to\Cc$ for $n\geq 2$ of double categories are not double functors defined on 
$\Aa_1\times...\times\Aa_n$. Loosely speaking, they are only defined on 
cells $(a_1,...,a_n)\in\Aa_1\times...\times\Aa_n$ whereby either all $a_i$ are objects, or at most two of them are simultaneously 
1-cells (in which case one is a 1h-cell and another a 1v-cell), or a single $a_i$ is a 2-cell whereas the rest are objects. 


\medskip

Let us now introduce the sets of multimaps for our funny multicategory for double categories. 
For the class of nullary maps $\M_0^*(-;\Aa)$ we set the set of objects of $\Aa$, 
and for unary maps we set $\M_1^*(\Aa;\Bb)=\Fun^*(\Aa,\Bb)$. For the class of $n$-ary multimaps we set 
$\M_n^*(\Aa_1,...,\Aa_n;\Cc)=f_n\x\Fun^*(\Aa_1\times...\times\Aa_n,\Cc)$, the set of double funny functors of type $*$ of $n$-variables.

\subsubsection{Substitution} \ssslabel{subst}

We now define substitutions. 
Let $F^n:\Aa_1\times...\times\Aa_n\to\Bb$ denote a funny functor for $n\geq 2$, nullary for $n=0$, or unary for $n=1$. 
We separate the two extreme cases:
\begin{itemize}
\item 
for 
a unary map $G^1:\Bb\to\Cc$ and $n\geq 0$, we clearly have that $G^1\comp F^n$ is of the same kind as $F^n$; 
\item 
for $n\geq 2, k\leq n$ we set $F^n\circ_k A_k = F^n(...,A_k,...)$, which is the $k$-th $(n\x 1)$-ary multimap constituting $F^n$, 
where $A_k$ denotes the nullary multimap $(-)\to\Aa_k$ picking up the object $A_k$. 
\end{itemize}
In general, for $n\geq 1, m\geq 0$, $F^n$ as above and $F^m:\Bb_1\times...\times\Bb_m\to\Aa_k$, we set 
$$F^n\circ_k F^m(a_1,...a_{k-1}, b_1,...,b_m, a_{k+1},...,a_n)=F^n(a_1,...a_{k-1}, F^m(b_1,...,b_m),a_{k+1},...,a_n)$$ 
for any sensible combination of cells $(a_1,...,a_n, b_1,..,b_m)\in\Aa_1\times...\times\Aa_n\times\Bb_1\times...\times\Bb_m$.
When $m=1$ and $n\geq 2$, $F^n\circ_k F^1$ is a funny functor because so is $F^n$. For the rest of the cases, that is for $n,m\geq 2$, 
it remains to prove that $F^n\circ_k F^m$ gives an $(n+m-1)$-ary multimap. More generally, we prove \prref{subst-gen funny}.  
For that we will need:

\begin{lma} \lelabel{un-fun}
Let $F_1,F_2,F_3$ be unary double functors. 
\begin{enumerate}
\item For a binary funny functor $G$ the composite $G(F_1(-),F_2(-))$ is a binary funny functor. 
\item For a ternary funny functor $G$ the composite $G(F_1(-),F_2(-),F_3(-))$ is a ternary funny functor. 
\end{enumerate}
All $F_i$'s and $G$'s are considered to be of the same type $*$.
\end{lma}

\begin{proof} 
The agreement-on-objects binary property of $G(F_1(-),F_2(-))$ holds by 
that property of $G$, as $G(F_1(A),F_2(-))\vert_B=G(\crta A, -)\vert_{\crta B}=G(-,\crta B)\vert_{\crta A}
=G(F_1(-),F_2(B))\vert_A$, where $\crta A=F_1(A)$ and $\crta B=F_2(B)$. 
If $F_1$ and $F_2$ are strict double functors, we can say that the binary axioms for the composite hold because they hold for $G$. 
For non-strict $F_i$'s first observe that $G(F_1(V),F_2(v))=G(U,u)=\Id$. 
The remaining 10 axioms of the funny functor property of $G(F_1(-),F_2(-))$ 
correspond to the $2+2$ axioms of the horizontal unnatural transformations 
$G(K,-)=G(F_1(L),-)$ and $G(-,k)=G(-,F_2(l))$ and $3+3$ axioms of the vertical strict transformations $G(-,u)=G(-,F_2(v))$ and 
$G(U,-)=G(F_1(V),-)$, recall (i) and (ii) of \prref{funny fun}. 
We used notations $k,K,l,L$  and $u,U,v,V$ for the corresponding 1h- and 1v-cells, respectively.

That $G(F_1(-),F_2(-),F_3(-))$ is given by three binary funny functors follows from part 1. For \equref{obj} it is sufficient to prove functoriality at any, say first, variable. The question $G(F_1(-),F_2(B),F_3(-))_C=G(F_1(-),F_2(-),F_3(C))_B$ we can rewrite as 
 $G(F_1(-),\crta B,-)_{\crta C}=G(F_1(-),-,\crta C)_{\crta B}$, which is the same property for the ternary functor 
$G(F_1(-),-,-)$. The latter is indeed a ternary funny functor as it is given by binary funny functors in the first and any other variable by part 1. (and it was already binary funny in the last two variables). 

Ternary axioms for $G(F_1(-),F_2(-),F_3(-))$ hold because they hold for $G$. 
\qed\end{proof}

\begin{prop} \prlabel{subst-gen funny}
Given lax double unary or funny functors $F_i:\Aa_{i1}\times...\times\Aa_{im_i}\to\Bb_i$ of $m_i$-variables with $i=1,...,n, m_i\geq 1, 
n\geq 2$ and a lax double funny functor $G:\Bb_1\times...\times\Bb_n\to\Cc$ of $n$-variables. Then the composition 
$$\Pi_{j=1}^{m_1}\Aa_{1j}\times...\times\Pi_{j=1}^{m_n}\Aa_{nj}\stackrel{F_1\times...\times F_n}{\longrightarrow}\Bb_1\times...\times\Bb_n
\stackrel{G}{\to}\Cc$$
is a lax double funny functor of $m_1+..+m_n$-variables.
\end{prop}

\begin{proof}
We should check that $G(F_1\times...\times F_n)$ gives binary funny functors for any choice of two variables, and ternary funny functors 
for any choice of three variables. We can examine the binary funny property, without loss of generality, for the composite 
of the form $G(U_1(-),U_2(-),F_1(-,-), F_2(-,-))$ for unary lax double functors $U_1,U_2$. 
It is binary in the first two variables by \leref{un-fun}, and in the two variables of each 
$F_i, i=1,2$ since so is $F_i$. 
It suffices thus to consider $G(U_1(-),F_1(A,-))$ and $G(F_1(A,-),F_2(C,-))$, without loss of generality. 
But these are binary funny functors again by \leref{un-fun}. 

For the 
ternary funny property and axioms 
we can assume that neither of the $F_i$'s has more than two variables (the contrary would only be interesting to study if all the three variables of the composite in which we wish to check the ternary funny property live in a ternary $F_i$, but in that case we are done, by the property of that 
$F_i$). Also, we can assume that there are at most two unary $U_i$'s, since the ternary funny property of the composite in the variables of three unary $U_i$'s 
would hold by the part 2. of 
\leref{un-fun} in those variables. 
Then we are left with the cases: 1) $G(U_1(-),U_2(-),F(-,D))$, 2) $G(U(-), F_1(-,X), F_2(-,Y))$, 3) $G(U(-), F(-,-))$, 4) 
$G(F_1(-,X), F_2(-,Y),F_3(-,Z))$ and 5) $G(F_1(-,-), F_2(-,X))$. The cases 1), 2) and 4) hold by the part 2. of \leref{un-fun}. 
The cases 3) and 5) are actually of the same type: they both consist of a unary and a binary funny functor. Let us discuss 3). Functoriality in the first variable in the ternary funny condition holds by the binary funny property of $F(-,-)$. Functoriality in the second or third variable is the only one where there's more work to be done. We check one of them: to see if 
$G(U(A), F(-,-))\vert_B=G(U(-), F(B,-))\vert_A$ as unary functors, we evaluate them at cells $c\in\Cc$ of different order. So we need to check the equality: 
$$(G(U(A), F(-,-))\vert_B)_{\displaystyle{\vert_c}}=(G(U(-), F(B,-))\vert_A)_{\displaystyle{\vert_c}}.$$
We have: $(G(U(A), F(-,-))\vert_B)_{\displaystyle{\vert_c}}=G(U(A), F(-,-)\vert_{(B,c)})\stackrel{*}{=}G(U(-), F(B,c))\vert_A$ and \\
$(G(U(-), F(B,-))\vert_A)_{\displaystyle{\vert_c}}\stackrel{*}{=}G(U(-), F(B,-)\vert_c)\vert_A$. Observe that $F(B,-)\vert_c=F(B,c)$ and both 
equalities marked by $*$ hold true by the coincidences in (ii) of \prref{funny fun} for the binary funny functors $F$, $G(U(-),-)$ 
(with $\crta b=F(B,c)$) and $G(U(-),F(B,-))$, respectively, if $c$ is any 1-cell or a 2-cell, whereas if $c$ is an object they are true by the agreement-on-objects of the named binary funny functors. 

It remains to verify the ternary axioms. Observe that the three axioms correspond to saying that the following are mixed funny modifications: 
$(-,v,h)$ (or $(u,-,h)$) - the first axiom, $(-,g,z)$ (or $(u,g,-)$) - the second, and $(f,-,z)$ (or $(f,v,-)$) - the third, 
with usual notations. 
The first and second are true for $G(U(-), F(-,-))$ by \prref{funny fun}, as $F(v,h)$ and $F(g,z)$ 
are 2-cells. (Similarly, the second and third are true for $G(F_1(-,-), F_2(-,X))$ of case 5), since $F_1(u,g)$ and $F_1(f,v)$ are 2-cells.) 
Instead of realizing $(f,-,z)$ as a mixed funny modification in the case of $G(U(-), F(-,-))$, the third axiom we may write as the identity 
$\frac{(f,B,z)}{(f,v,\tilde C)}=\frac{(f,v,C)}{(f,\tilde B, z)}$ and check if $\frac{G(U(f),F(B,z))}{G(U(f),F(v,\tilde C))}=
\frac{G(U(f),F(v,C))}{G(U(f),F(\tilde B, z))}$ holds. By the axiom \axiomref{($k,\frac{U}{U'}$)} for $G$ 
it can be restated as $G(U(f),\frac{F(B,z)}{F(v,\tilde C)})=G(U(f),\frac{F(v,C)}{F(\tilde B, z)})$, but we know this is true, because the two fractions involving the double funny functor $F$ are equal by the identity $(u,U)=\Id$ for the double funny functor $F$.  
\qed\end{proof}

For the associativity of the substitution we find the following. Let us first discuss the cases involving nullary maps. \vspace{-0,2cm}
\begin{itemize}
\item 
When only nullary maps are composed into an $n$-ary multimap $F^n$ with $n\geq 2$ at any two variables, the associativity holds by the binary property of $F^2$ and the ternary property \equref{obj} of $F^n$ for higher $n$ that includes those two variables. \vspace{-0,2cm}
\item 
When a nullary and an $m$-ary multimap with $m\geq 1$ are composed into $F^n$, the reasoning for the associativity is similar, whereby one also applies the coincidences (i)-(iii) of \prref{funny fun}. \vspace{-0,2cm}
\item 
Another associativity that involves nullary maps is when a nullary map $A_k$ is substituted into an $m$-ary and their composite is substituted into an $n$-ary multimap, with $m,n\geq 1$. This associativity is proved in the same way as the latter one:
\begin{equation} \eqlabel{nul-ass}
F^n(...,F^m(...,A_k,...)...)=F^n(...,F^m(...,-,...)...)\vert_{(...,(...,A_k,...),...}.
\end{equation}
\end{itemize}

To analyze the other multimaps, 
take double funny functors $F^n, F^m, F^r$ of arities $n,m,r$ respectively, then we wonder if 
$(F^n\comp_i F^m)\comp_{l+i-1} F^r=F^n\comp_i (F^m\comp_l F^r)$ as double funny functors with $l\leq m$.  
Another case of associativity that can occur is that one substitutes $F_i^{m_i}$ and $F_j^{m_j}$ into $F^n$ with $i<j$ one at a time. The associativity then reads $(F^n\comp_i F_i^{m_i})\comp_{j+m_i-1} F_j^{m_j}=(F^n\comp_j F_j^{m_j})\comp_i F_i^{m_i}$. 
These two kinds of questions can be rewritten as 
$$F^n(...,F^m(...,-,...)...)\vert_{F^r(...)}=F^n(...,-,...)\vert_{F^m(...,F^r(...),...)}$$
and
$$G(...,F_i(...),...,-,...)\vert_{F_j(...)}=G(...,-,...,F_j(...),...)\vert_{F_i(...)}.$$
To check equalities of these funny functors we should check if the corresponding composites are equal for any choice of two variables 
as binary funny functors. After skipping the easy cases, similarly as in the proof of the above proposition, the non-trivial cases are when: 1) one variable 
is in $F^r$ and the other one in $F^m$ (outside of $F^r$), 2) one is in $F^r$ and the other one in $F^n$ (outside of $F^m$), 3) one is 
in $F^m$ (outside of $F^r$) and the other one in $F^n$ (outside of $F^m$), and similarly in the case of $G, F_i, F_j$. 
The equalities are proved by applying the coincidences in the items (i)-(iii) of \prref{funny fun} for 1- and 2-cells, and by applying the agreement-on-objects binary property. 





We can finally claim that double funny functors make a funny multicategory in the way we described in this subsection.

\subsubsection{Binary multimap sets as double categories} \ssslabel{mult-dc}

The bijection \equref{funny-basic} can be upgraded into an isomorphism of double categories 
\begin{equation}\eqlabel{funny-basic-d}
f\x[\Aa\times\Bb,\Cc]^*\iso [\Aa, [\Bb,\Cc]^*]^*.
\end{equation}
Let us first introduce the double category on the left-hand side. For simplicity we do it for $*=lx$. 

Objects of $f\x[\Aa\times\Bb,\Cc]^{lx}$ are lax double funny functors. Its 1h-cells are 
 horizontal unnatural transformations $\theta\colon (-,-)_1\Rightarrow (-,-)_2$ between lax double funny functors $(-,-)_1,(-,-)_2\colon 
\Aa\times\Bb\to\Cc$. They are given by: for each $A\in\Aa$ a horizontal unnatural transformation $\theta^A\colon (-,A)_1\Rightarrow(-,A)_2$ and 
for each $B\in\Bb$ a horizontal unnatural transformation $\theta^B\colon (B,-)_1\Rightarrow(B,-)_2$, both of lax double functors, such that 
$\theta^A_B=\theta^B_A$ and the following axiom holds \\
\noindent \axiom{$HUT^f$} \vspace{-0,3cm}
$$\scalebox{0.86}{
\bfig
 \putmorphism(-170,500)(1,0)[(B,A)_1`(B,A)_1 `=]{600}1a
 \putmorphism(510,500)(1,0)[\phantom{Y_2}` `\theta^{A}_B]{420}1a
\putmorphism(-200,500)(0,-1)[\phantom{Y_2}` `(B,U)_1]{450}1l
\putmorphism(-240,500)(0,-1)[\phantom{Y_2}`(B, \tilde A)_1 `]{450}0l
\put(-10,50){\fbox{$\Id$}}
\putmorphism(-190,-400)(1,0)[(\tilde B, \tilde A)_1` `=]{500}1a
\putmorphism(-200,50)(0,-1)[\phantom{Y_2}``(u,\tilde A)_1]{450}1l
\putmorphism(450,50)(0,-1)[\phantom{Y_2}`(\tilde B, \tilde A)_1`(\tilde B, U)_1]{450}1l
\putmorphism(450,500)(0,-1)[\phantom{Y_2}`(\tilde B, A)_1 `(u,A)_1]{450}1l
\put(600,280){\fbox{$\theta^A_u$}}
\putmorphism(430,50)(1,0)[\phantom{(B, \tilde A)}``\theta^{A}_{\tilde B}]{480}1a
\putmorphism(1070,50)(0,-1)[\phantom{(B, A')}`(\tilde B, \tilde A)_2`(\tilde B,U)_2]{450}1r
\putmorphism(1070,500)(0,-1)[(B, A)_2`(\tilde B, A)_2`(u,A)_2]{450}1r
\putmorphism(450,-400)(1,0)[\phantom{(B, \tilde A)}``\theta^{\tilde A}_{\tilde B}]{480}1a
\put(600,-150){\fbox{$\theta_U^{\tilde B}$}}
\efig}
\quad=\quad
\scalebox{0.86}{
\bfig
 \putmorphism(-150,500)(1,0)[(B,A)_1`(B,A)_2 `\theta^{A}_B]{600}1a
 \putmorphism(450,500)(1,0)[\phantom{(B,A)}` `=]{480}1a
\putmorphism(-180,500)(0,-1)[\phantom{Y_2}`(B, \tilde A)_1 `(B,U)_1]{450}1l
\put(0,300){\fbox{$\theta_U^B$}}
\putmorphism(-180,-400)(1,0)[(\tilde B, \tilde A)_1` `\theta^{\tilde A}_{\tilde B}]{500}1a
\putmorphism(-180,50)(0,-1)[\phantom{Y_2}``(u,\tilde A)_1]{450}1l
\putmorphism(450,50)(0,-1)[\phantom{Y_2}`(\tilde B, \tilde A)_2`(u,\tilde A)_2]{450}1r
\putmorphism(450,500)(0,-1)[\phantom{Y_2}`(B, \tilde A)_2 `(B,U)_2]{450}1r
\put(690,50){\fbox{$\Id$}}
\putmorphism(-180,50)(1,0)[\phantom{(B, \tilde A)}``\theta^{\tilde A}_B]{500}1a
\putmorphism(1070,50)(0,-1)[\phantom{(B, A')}`(\tilde B, \tilde A)_2`(\tilde B,U)_2]{450}1r
\putmorphism(1070,500)(0,-1)[(B, A)``(u,A)_2]{450}1r
\putmorphism(1110,500)(0,-1)[`(\tilde B, A)_2`]{450}0r
\putmorphism(450,-400)(1,0)[\phantom{(B, \tilde A)}``=]{480}1b
\put(0,-150){\fbox{$\theta^{\tilde A}_u$}}
\efig}
$$
for every 1v-cells $U\colon A\to \tilde A$ and $u\colon B\to\tilde B$.  

1v-cells of   $f\x[\Aa\times\Bb,\Cc]^{lx}$ are 
vertical strict transformations $\theta_0\colon (-,-)_1\Rightarrow (-,-)_2$ between lax double funny functors $(-,-)_1,(-,-)_2\colon 
\Aa\times\Bb\to\Cc$. They are given by: for each $A\in\Aa$ a vertical strict transformation $\theta_0^A\colon (-,A)_1\Rightarrow(-,A)_2$ and 
for each $B\in\Bb$ a vertical strict transformation $\theta_0^B\colon (B,-)_1\Rightarrow(B,-)_2$, both of lax double functors, such that 
$(\theta_0^A)_B=(\theta_0^B)_A$ and the following axioms hold: \\
\noindent \axiom{$VST^f_1$} \vspace{-0,3cm}
$$\scalebox{0.86}{
\bfig
 \putmorphism(-150,500)(1,0)[(B,A)_1`(B,A)_1 `=]{600}1a
 \putmorphism(450,500)(1,0)[\phantom{(B,A)_1}` `(B,K)_1]{450}1a
\putmorphism(-200,500)(0,-1)[\phantom{Y_2}`(B, A)_2 `(\theta_0^A)_B]{450}1l
\put(20,50){\fbox{$\Id$}}
\putmorphism(-170,-400)(1,0)[(\tilde B, A)_2` `=]{500}1a
\putmorphism(-200,50)(0,-1)[\phantom{Y_2}``(u,A)_2]{450}1l
\putmorphism(450,50)(0,-1)[\phantom{Y_2}`(\tilde B, A)_2 `(\theta_0^A)_{\tilde B}]{450}1l
\putmorphism(450,500)(0,-1)[\phantom{Y_2}`(\tilde B, A)_1 `(u,A)_1]{450}1l
\put(600,260){\fbox{$(u,K)_1$}}
\putmorphism(450,50)(1,0)[\phantom{(B, \tilde A)}``(\tilde B,K)_1]{500}1a
\putmorphism(1070,50)(0,-1)[\phantom{(B, A')}`(\tilde B, A')_2`(\theta^{\tilde B}_0)_{A'}]{450}1r
\putmorphism(1070,500)(0,-1)[(B, A')_1`(\tilde B, A')_1`(u,A')_1]{450}1r
\putmorphism(450,-400)(1,0)[\phantom{(B, \tilde A)}``(\tilde B,K)_2]{500}1a
\put(600,-170){\fbox{$(\theta_0^{\tilde B})_K$}}
\efig}
\quad=\quad
\scalebox{0.86}{
\bfig
 \putmorphism(-150,500)(1,0)[(B,A)_1`(B,A')_1 `(B,K)_1]{600}1a
 \putmorphism(450,500)(1,0)[\phantom{(B,A)}` `=]{480}1a
\putmorphism(-180,500)(0,-1)[\phantom{Y_2}`(B, A)_2 `(\theta^A_0)_B]{450}1l
\put(0,280){\fbox{$(\theta_0^B)_K$}}
\putmorphism(-180,-400)(1,0)[(\tilde B,A)_2` `(\tilde B,K)_2]{500}1a
\putmorphism(-180,50)(0,-1)[\phantom{Y_2}``(u,A)_2]{450}1l
\putmorphism(450,500)(0,-1)[\phantom{Y_2}`(B, A')_2 `(\theta^{A'}_0)_B]{450}1r
\putmorphism(450,50)(0,-1)[\phantom{Y_2}`(\tilde B, A')_2`(u, A')_2]{450}1r
\put(660,50){\fbox{$\Id$}}
\putmorphism(-180,50)(1,0)[\phantom{(B, \tilde A)}``(B,K)_2]{500}1a
\putmorphism(1070,50)(0,-1)[\phantom{(B, A')}`(\tilde B, A')_2`(\theta_0^{A'})_{\tilde B}]{450}1r
\putmorphism(1070,500)(0,-1)[(B, A')_1`(\tilde B, A')_1`(u,A')_1]{450}1r
\putmorphism(450,-400)(1,0)[\phantom{(B, \tilde A)}``=]{480}1b
\put(0,-170){\fbox{$(u,K)_2$}}
\efig}
$$
for every 1h-cell $K\colon A\to A'$ and 1v-cell $u\colon B\to \tilde B$, 

\pagebreak

\noindent \axiom{$VST^f_2$} \vspace{-0,6cm}

$$\scalebox{0.86}{
\bfig
 \putmorphism(-150,500)(1,0)[(B,A)_1`(B,A)_1 `=]{600}1a
 \putmorphism(450,500)(1,0)[\phantom{(B',A)_1}` `(k,A)_1]{450}1a
\putmorphism(-200,500)(0,-1)[\phantom{Y_2}`(B, A)_2 `(\theta_0^B)_A]{450}1l
\put(20,50){\fbox{$\Id$}}
\putmorphism(-170,-400)(1,0)[(B, \tilde A)_2` `=]{500}1a
\putmorphism(-200,50)(0,-1)[\phantom{Y_2}``(B,U)_2]{450}1l
\putmorphism(450,50)(0,-1)[\phantom{Y_2}`(B, \tilde A)_2 `(\theta_0^B)_{\tilde A}]{450}1l 
\putmorphism(450,500)(0,-1)[\phantom{Y_2}`(B, \tilde A)_1 `(B,U)_1]{450}1l
\put(600,260){\fbox{$(k,U)_1$}}
\putmorphism(450,50)(1,0)[\phantom{(B, \tilde A)}``(k,\tilde A)_1]{500}1a %
\putmorphism(1070,50)(0,-1)[\phantom{(B, A')}`(B', \tilde A)_2`(\theta^{\tilde A}_0)_{B'}]{450}1r
\putmorphism(1070,500)(0,-1)[(B, A')_1`(B', \tilde A)_1`(B',U)_1]{450}1r
\putmorphism(450,-400)(1,0)[\phantom{(B, \tilde A)}``(k,\tilde A)_2]{500}1a
\put(600,-170){\fbox{$(\theta_0^{\tilde A})_k$}}
\efig}
\quad=\quad
\scalebox{0.86}{
\bfig
 \putmorphism(-150,500)(1,0)[(B,A)_1`(B',A)_1 `(k,A)_1]{600}1a
 \putmorphism(450,500)(1,0)[\phantom{(B,A)}` `=]{480}1a
\putmorphism(-180,500)(0,-1)[\phantom{Y_2}`(B, A)_2 `(\theta^A_0)_B]{450}1l
\put(0,280){\fbox{$(\theta_0^A)_k$}}
\putmorphism(-180,-400)(1,0)[(\tilde B,A)_2` `(k,\tilde A)_2]{500}1a
\putmorphism(-180,50)(0,-1)[\phantom{Y_2}``(B,U)_2]{450}1l
\putmorphism(450,50)(0,-1)[\phantom{Y_2}`(B', \tilde A)_2`(B', U)_2]{450}1r
\putmorphism(450,500)(0,-1)[\phantom{Y_2}`(B', A)_2 `(\theta^{B'}_0)_A]{450}1r
\put(680,50){\fbox{$\Id$}}
\putmorphism(-180,50)(1,0)[\phantom{(B, \tilde A)}``(k,A)_2]{500}1a
\putmorphism(1070,50)(0,-1)[\phantom{(B, A')}`(B',\tilde A)_2`(\theta_0^{\tilde A})_{B'}]{450}1r
\putmorphism(1070,500)(0,-1)[(B', A)_1`(B',\tilde A)_1`(B', U)_1]{450}1r
\putmorphism(450,-400)(1,0)[\phantom{(B, \tilde A)}``=]{480}1b
\put(0,-170){\fbox{$(k,U)_2$}}
\efig}
$$
for every 1v-cell $U\colon A\to \tilde A$ and 1h-cell $k\colon B\to B'$.

Finally, 2-cells of $f\x[\Aa\times\Bb,\Cc]^{lx}$ are mixed funny modifications.  
Given horizontal unnatural transformations $\theta, \theta'$ and vertical strict transformations $\theta_0, \theta'_0$ acting between lax double funny functors 
$H_1, H_2, H_3, H_4\colon \Aa\times\Bb\to\Cc$ as in the left diagram below. Denote by $(-,A)_i\colon\Bb\to\Cc, (B,-)_i\colon\Aa\to\Cc, i=1,2,3,4$ the pairs of 
lax double functors corresponding to $H_1, H_2, H_3, H_4$, respectively. A mixed funny modification $\tau$ (on the left below) is given by a pair of mixed funny modifications 
$\tau^A, \tau^B$ acting between horizontal unnatural and vertical strict transformations among lax double functors:
\begin{equation} \eqlabel{q-modif}
\scalebox{0.86}{
\bfig
\putmorphism(-150,50)(1,0)[H_1` H_2` \theta]{400}1a
\putmorphism(-150,-270)(1,0)[H_3 ` H_4 ` \theta' ]{400}1b
\putmorphism(-170,50)(0,-1)[\phantom{Y_2}``\theta_0]{320}1l
\putmorphism(250,50)(0,-1)[\phantom{Y_2}``\theta_0']{320}1r
\put(-30,-140){\fbox{$\tau$}}
\efig}
\qquad\qquad
\scalebox{0.86}{
\bfig
\putmorphism(-180,50)(1,0)[(-,A)_1` (-,A)_2`\theta^A]{550}1a
\putmorphism(-180,-270)(1,0)[(-,A)_3`(-,A)_4 `\theta'^A]{550}1b
\putmorphism(-170,50)(0,-1)[\phantom{Y_2}``\theta_0^A]{320}1l
\putmorphism(350,50)(0,-1)[\phantom{Y_2}``\theta_0'^A]{320}1r
\put(0,-140){\fbox{$\tau^A$}}
\efig}
\qquad\qquad
\scalebox{0.86}{
\bfig
\putmorphism(-180,50)(1,0)[(B,-)_1` (B,-)_2`\theta^B]{550}1a
\putmorphism(-180,-270)(1,0)[(B,-)_3`(B,-)_4 `\theta'^B]{550}1b
\putmorphism(-170,50)(0,-1)[\phantom{Y_2}``\theta_0^B]{320}1l
\putmorphism(350,50)(0,-1)[\phantom{Y_2}``\theta_0'^B]{320}1r
\put(0,-140){\fbox{$\tau^B$}}
\efig}
\end{equation}
such that $\tau^A_B=\tau^B_A$ for every $A\in\Aa, B\in\Bb$. 

\medskip

The double isomorphism \equref{funny-basic-d} is proved analogously as the quasi-functor version of the claim from \cite[Section 4]{Fem:Bif}. 
Since double quasi-functors obey some additional axioms that do not appear for double funny functors, the proofs for double funny functors are analogous to those 
for double quasi-functors, but are to some extent simpler. 
The 1-1 correspondence on 1- and 2-cells in \equref{funny-basic-d} is proved {\em mutatis mutandi} as in Sections 4.5 and 4.8 of {\em loc.cit.}. The funny counterparts of Propositions 4.6 and 4.7 thereof read as follows: 

\begin{prop} \prlabel{equiv-modif}
Let $F,G\colon \Aa\to[\Bb,\Cc]^{lx}$ be two lax double functors with the corresponding lax double funny functors 
$(-,-)_1, (-,-)_2\colon\Aa\times\Bb\to\Cc$. 
For every $A\in\Aa$ and $B\in\Bb$ let $\alpha(A): F(A)\to G(A)$ and $\alpha(-)_B: F(-)(B)\to G(-)(B)$ be horizontal unnatural transformations 
between lax double functors. The following are equivalent: 
\begin{enumerate}
\item $\alpha_U$ of the form 
\begin{equation}\eqlabel{alfa-U}
\scalebox{0.86}{
\bfig
 \putmorphism(-200,50)(1,0)[(-,A)_1`(-,A)_2`\alpha(A)]{600}1a
\putmorphism(-200,-350)(1,0)[(-,\tilde A)_1`(-,\tilde A)_2 `\alpha(\tilde A)]{600}1a
\putmorphism(-180,50)(0,-1)[\phantom{Y_2}``(-,U)_1]{400}1l
\putmorphism(370,50)(0,-1)[\phantom{Y_2}``(-,U)_2]{400}1r
\put(0,-110){\fbox{$\alpha_U$}} 
\efig}
\end{equation}
is a mixed funny modification 
for every 1v-cell $U\colon A\to\tilde A$ in $\Aa$; 
\item the pairs $(\theta^A, \theta^B)\colon\hspace{-0,2cm}=(\alpha(A),\alpha(-)_B)$ for $A\in\Aa, B\in\Bb$ form a horizontal unnatural  transformation $\theta\colon (-,-)_1\Rightarrow(-,-)_2$ between lax double funny functors. 
\end{enumerate}
\end{prop}
 
In the above claim the 2-cell components of the mixed funny modification $\alpha_U$ (at $B$) are given by 
and correspond bijectively to the 2-cell components of the horizontal unnatural transformation $\alpha(-)_B$ (at $U$).

\begin{prop} \prlabel{horiz nat summed up}
A horizontal unnatural transformation $\alpha\colon F\Rightarrow G$ between lax double functors $F,G\colon \Aa\to[\Bb,\Cc]^{lx}$ consists of the following data:
\begin{itemize}
\item a horizontal unnatural transformation $\alpha(A)\colon F(A)\to G(A)$ between lax double functors for every $A\in\Aa$;
\item a mixed funny modification $\alpha_U$ (of the form \equref{alfa-U}) for every 1v-cell $U\colon A\to\tilde A$, 
\end{itemize} 
so that $\alpha_U$ obeys two axioms, which (after evaluation at $B\in\Bb$) yield that $\alpha(-)_B\colon F(-)(B)\to G(-)(B)$ 
is a horizontal unnatural transformation between lax double functors for every $B\in\Bb$ (by setting $\alpha(U)_B\colon\hspace{-0,2cm}=\alpha_U(B)$). 
\end{prop}

Analogously as in the quasi-functor case of \cite[Proposition 6.3]{Fem:Bif} for double funny functors one proves an isomorphism of 
double categories 
\begin{equation} \eqlabel{ot-f-d}
[\Aa\Box^*\Bb,\Cc]\iso f\x[\Aa\times\Bb,\Cc]^*. 
\end{equation}
In view of \equref{funny-basic-d} we obtain a double category isomorphism 
\begin{equation} \eqlabel{ot-corch}  
[\Aa\Box^*\Bb,\Cc]\iso[\Aa,[\Bb,\Cc]^*]^*. 
\end{equation}

\subsubsection{A characterization of ternary lax double funny functors}

In the proof of \prref{subst-gen funny} we observed that the three double funny ternary axioms are equivalent to certain mixed funny modification axioms. We 
include this finding in the following result. By mere inspection we see that there is more data equivalent to it. 

\begin{prop} \prlabel{ternary eq}
Let $F:\Aa\times\Bb\times\Cc\to\Dd$ be an assignment (possibly functorial in every variable).  
The following equivalences for $F$ to satisfy the three ternary double funny functor axioms hold true: 

$\axiomref{u,v,h}$  $\Leftrightarrow$  \axiomref{$HUT^f$} for $F(-,-,h)$ 
$\Leftrightarrow$ \hspace{0,16cm}\axiomref{$VST^f_2$} for $F(u,-,-)$ \hspace{0,1cm}
$\Leftrightarrow$ \hspace{0,16cm}\axiomref{$VST^f_2$} for $F(-,v,-)$ 

\hspace{1,2cm} $\Leftrightarrow$  $F(-,v,h)$ obeys the mixed funny modification axiom 

\hspace{1,2cm} $\Leftrightarrow$  $F(u,-,h)$ obeys the mixed funny modification axiom
\smallskip


\axiomref{u,g,z}  $\Leftrightarrow$  \axiomref{$HUT^f$} for $F(-,g,-)$ 
$\Leftrightarrow$ \hspace{0,16cm}\axiomref{$VST^f_1$} for $F(u,-,-)$ \hspace{0,1cm} 
$\Leftrightarrow$ \hspace{0,16cm}\axiomref{$VST^f_2$} for $F(-,-,z)$ 

\hspace{1,2cm} $\Leftrightarrow$  $F(-,g,z)$ obeys the mixed funny modification axiom 

\hspace{1,2cm}  $\Leftrightarrow$ $F(u,g,-)$ obeys the mixed funny modification axiom 
\smallskip

\axiomref{f,v,z}  $\Leftrightarrow$  \axiomref{$HUT^f$} for $F(f,-,-)$ 
$\Leftrightarrow$ \hspace{0,16cm}\axiomref{$VST^f_1$} for $F(-,-,z)$ \hspace{0,1cm} $\Leftrightarrow$ \hspace{0,16cm}\axiom{$VST^f_1$} 
for $F(-,v,-)$ 

\hspace{1,2cm} $\Leftrightarrow$  $F(f,-,z)$ obeys the mixed funny modification axiom 

\hspace{1,2cm} $\Leftrightarrow$ $F(f,v,-)$ obeys the mixed funny modification axiom
\smallskip

\noindent for 1h-cells $f,g,h$ and 1v-cells $u,v,z$.
\end{prop}

In the style of \prref{funny fun} we may then conclude: 

\begin{prop} \prlabel{funny ternary char}
Let $\Aa,\Bb,\Cc$ be double categories. The following are equivalent:
\begin{enumerate} 
\item $H\colon \Aa\times\Bb\times\Cc\to\Dd$ is a ternary lax double funny functor, 
\item there are three families of binary lax double funny functors 
$$(A,-,-):\Bb\times\Cc\to\Dd, (-,B,-):\Aa\times\Cc\to\Dd, (-,-,C):\Aa\times\Bb\to\Dd$$ 
such that 
$$(A,-,-)\vert_B=(-,B,-)\vert_A, \,\,(A,-,-)\vert_C=(-,-,C)\vert_A, \,\, (-,B,-)\vert_C=(-,-,C)\vert_B$$
for objects $A\in\Aa, B\in\Bb, C\in\Cc$, determining unambiguous lax double functors 
$(A,B,-):\Cc\to\Dd, (-,B,C):\Aa\to\Dd, (A,-,C):\Bb\to\Dd$ 
such that $(A,B,-)\vert_C=(-,B,C)\vert_A=(A,-,C)_B$, and the following hold for 1h-cells $f,g,h$ and 1v-cells $u,v,z$: 
\begin{enumerate}[(i)]
\item 
$$(f,-,-)\colon (A,-,-)\to(A',-,-),$$ $$(-,g,-)\colon (-,B,-)\to (-,B', -),$$ $$(-,-,h)\colon (-,-,C)\to(-,-,C')$$ 
are horizontal unnatural transformations, \\ 
$$(u,-,-)\colon (A,-,-)\to(\tilde A,-,-),$$ $$(-,v,-)\colon (-,B,-)\to(-,\tilde B,-),$$ $$(-,-,z)\colon (-,-,C)\to(-,-,\tilde C)$$
are vertical strict transformations, and 
$$(a,-,-), (-,b,-), (-,-,c)$$ 
are binary mixed funny modifications with respect to horizontally unnatural and vertically strict transformations for all 2-cells $a$ in 
$\Aa$, $b$ in $\Bb$ and $c$ in $\Cc$, 
\item the following coincide: 
$$(f,B,-)\vert_C=(f,-,C)_B=(-,B,C)\vert_f,$$ 
$$(u,B,-)\vert_C=(u,-,C)_B=(-,B,C)\vert_u, $$
$$(a,B,-)\vert_C=(a,-,C)_B=(-,B,C)\vert_a,$$

$$(A,g,-)\vert_C=(-,g,C)\vert_A=(A,-,C)\vert_g, $$
$$(A,v,-)\vert_C=(-,v,C)\vert_A=(A,-,C)\vert_v, $$
$$(A,b,-)\vert_C=(-,b,C)\vert_A=(A,-,C)\vert_b,$$

$$(A,-,h)\vert_B=(-,B,h)_A=(A,B,)\vert_h, $$
$$(A,-,z)\vert_A=(-,B,z)_B=(A,B,-)\vert_z,$$
$$(A,-,c)\vert_B=(-,B,C)_A=(A,B,-)\vert_c,$$
\item the following are mixed funny modifications 
$$(f,-,z), \,\, (u,-,h), \,\, (-,g,z), \,\, (-,v,h), \,\, (u,g,-), \,\, (f,v,-)$$
and the following 2-cell components of the respective transformations and modifications coincide:
$$(f,B,-)\vert_z=(-,B,z)\vert_f=(f,B,z), \,\,\,\, (u,B,-)\vert_h=(-,B,h)\vert_u=(u,B,h), $$ 
$$(f,-,C)\vert_v=(-,v,C)\vert_f=(f,v,C), \,\,\,\, (u,-,C)\vert_g=(-,g,C)\vert_u=(u,g,C), $$ 
$$(A,g,-)\vert_z=(A,-,z)\vert_g=(A,g,z), \,\,\,\, (A,v,-)\vert_h=(A,-,h)\vert_v=(A,v,h). $$
\end{enumerate}
\end{enumerate}
\end{prop}

In view of \prref{ternary eq} one sees that the part 2. in the above characterization contains redundant data. 
The axioms for $(a,-,-), (-,b,-), (-,-,c)$ to be binary mixed funny modifications correspond to the axioms 
\axiomref{$(u,U)$-l-nat} and \axiomref{$(u,U)$-r-nat} of the corresponding underlying component binary double funny functors 
they act upon.

\subsubsection{Transformations of ternary lax double funny functors}

We will also need horizontal unnatural transformations and vertial strict transformations of (lax) double funny functors of more than two variables. 
For more than three variables apply the definition to any triple of variables. 

\begin{defn} \delabel{hu tr of ff 3}
A {\em horizontal unnatural transformation} $\theta: H_1\Rightarrow H_2$ between lax double funny functors 
$H_1,H_2:\Aa\times\Bb\times\Cc\to\Ee$ consists of horizontal unnatural transformations 
$$\theta^{A;B}: (A,B,-)_3^1\to (A,B,-)_3^2$$
$$\theta^{B;C}: (-,B,C)_1^1\to (-,B,C)_1^2$$ 
$$\theta^{A;C}: (A,-,C)_2^1\to (A,-,C)_2^2$$ 
of lax double functors 
such that $\theta^{A;B}(C)=\theta^{B;C}(A)=\theta^{A;C}(B)$ for all $(A,B,C)\in\Aa\times\Bb\times\Cc$. 
\end{defn}

\begin{defn} \delabel{vs tr of ff 3}
A {\em vertical strict transformation} $\theta: H_1\Rightarrow H_2$ between lax double funny functors $H_1,H_2:\Aa\times\Bb\times\Cc\to\Ee$ consists 
of vertical strict transformations 
$$\theta^A: H_1(A,-,-) \Rightarrow H_2(A,-,-),$$ 
$$\theta^B: H_1(-,B,-) \Rightarrow H_2(-,B,-),$$ 
$$\theta^C: H_1(-,-,C) \Rightarrow H_2(-,-,C)$$ 
of lax double funny functors, which give unambiguous vertical strict transformations 
$$\theta^{A;B}: (A,B,-)_3^1\to (A,B,-)_3^2$$
$$\theta^{B;C}: (-,B,C)_1^1\to (-,B,C)_1^2$$ 
$$\theta^{A;C}: (A,-,C)_2^1\to (A,-,C)_2^2$$ 
of lax double functors 
so that $\theta^{A;B}(C)=\theta^{B;C}(A)=\theta^{A;C}(B)$ for all $(A,B,C)\in\Aa\times\Bb\times\Cc$, and so that 
six equalities between their structure 2-cells, on one hand, 
and the six structure 2-cells of both $H_1$ and $H_2$ from \deref{cub}, on the other hand, hold. We present these six equalities 
as vertical compositions of 2-cells together with their labels: 
$$\big((\theta^{A;C})^v, (\theta^{B;C})_f, (f,v,C)^i_{12}\big), \qquad\qquad \big((\theta^{A;B})^z, (\theta^{B;C})_f, (f,B,z)^i_{13}\big) 
\vspace{-0,2cm}$$
$$\frac{(\theta^{B;C})_f}{(f,v,C)^2}=\frac{(f,v,C)^1}{(\theta^{\tilde B;C})_f}, \qquad \qquad
\frac{(\theta^{B;C})_f}{(f,B,z)^2}=\frac{(f,B,z)^1}{(\theta^{B;\tilde C})_f}$$

$$\big((\theta^{A;B})^z, (\theta^{A;C})_g, (A,g,z)^i_{23}\big), \qquad\qquad \big((\theta^{A;C})_g, (\theta^{B;C})^u, (u,g,C)^i_{12}\big),
\vspace{-0,2cm}$$
$$\frac{(\theta^{A;C})_g}{(A,g,z)^2}=\frac{(A,g,z)^1}{(\theta^{A;\tilde C})_g}, \qquad\qquad
\frac{(\theta^{A;C})_g}{(u,g,C)^2}=\frac{(u,g,C)^1}{(\theta^{\tilde A; C})_g}$$

$$\big((\theta^{A;B})_h, (\theta^{B;C})^u, (u,B,h)^i_{13}\big), \qquad \big((\theta^{A;B})_h, (\theta^{A;C})^v, (A,v,h)^i_{23}\big)
\vspace{-0,2cm}$$
$$\frac{(\theta^{A;B})_h}{(u,B,h)^2}=\frac{(u,B,h)^1}{(\theta^{\tilde A; B})_h}, \qquad\qquad
\frac{(\theta^{A;B})_h}{(A,v,h)^2}=\frac{(A,v,h)^1}{(\theta^{A;\tilde B})_h}$$
where $(f,v,C)^i$ for $i=1,2$ presents a structure 2-cell $(f,v,C)_{12}$ from \deref{cub} for $H_1$ and $H_2$, respectively, and similarly for the remaining five 2-cells of that type. 
\end{defn}

\subsubsection{Closedness}

Knowing from \sssref{mult-dc} that $[\Bb, \Cc]^{lx}$ is a double category, we now introduce the evaluation lax double funny functors. 

\begin{prop} \prlabel{ev}
There is a lax double funny functor $ev:[\Bb, \Cc]^{lx}\times\Bb\to\Cc$ such that given any lax double funny functor 
$H:\Aa\times\Bb\to\Cc$ it is 
$$ev(H^t\times \Id_\Bb)=H$$
where $H^t:\Aa\to [\Bb, \Cc]^{lx}$ is the lax double functor corresponding to $H$ by \equref{funny-basic}.
\end{prop}

\begin{proof}
For a 0-cell $F\in [\Bb, \Cc]^{lx}$ and any cell $x$ in $\Bb$ set $ev(F,x)=F(x)$; for a 0-cell $B\in\Bb$ and 1h- or 1v-cell $\alpha$ in 
$[\Bb, \Cc]^{lx}$ set $ev(\alpha,B)=\alpha(B)$, the 1h- respectively 1v-cell component of the transformation $\alpha$ in question; for a 2-cell 
$b$ in $[\Bb, \Cc]^{lx}$ set $ev(b,B)=b(B)$, the 2-cell component of the funny modification $b$; and finally, for a 1v-cell $u$ and 1h-cell $g$ 
in $\Bb$ set $ev(\alpha,u)=\alpha^u$, respectively $ev(\alpha,g)=\alpha_g$, the 2-cell components of the horizontal unnatural, respectively 
vertical strict transformation $\alpha$.  

First condition for $ev$ to be a lax double funny functor is that $ev(F,-):\Bb\to\Cc$ and $ev(-,B):[\Bb, \Cc]^{lx}\to\Cc$ be lax double functors. 
It is clear that $ev(F,-)$ is such a functor, as so is $F$. Since both compositions of 1-cells in $[\Bb, \Cc]^{lx}$ are vertical compositions of (horizontal unnatural, resp. vertical strict) transformations, we obtain that  
$ev(-,B)$ is a strict functor (see \leref{vert comp hor.ps.tr.} and \leref{vert comp vert. lx tr.} for the vertical compositions of horizontal and vertical transformations).

The coincidences of these two lax double functors on objects, and the coincidences of 1- and 2-cells in (i)-(iii) of \prref{funny fun} 
are now trivially fulfilled. 
For a (horizontal unnatural or vertical strict) transformation $\alpha$ and a mixed funny modification 
$b$ it is clear that $ev(\alpha,-)$ is a (horizontal unnatural, resp. vertical strict) transformation, and that $ev(b,-)$ is a corresponding 
modification. In Table \ref{table:1} we argument why $ev(-,g)$ and $ev(-,v)$ are a horizontal unnatural and a vertical strict transformation, respectively. 

\begin{table}[H]
\begin{center}
\begin{tabular}{ c c } 
transformation axiom & its meaning on cells in $Lax_{hop}(\B, \C)$ on which evaluated \\ [0.5ex]
\hline
$(v.l.t.\x 1)$ for  $ev(-,v)$ & \cite[Lemma 2.3 part 2]{Fem:Fil} (vertical composition of h.u.t. $\alpha,\beta$) 
\\ [1ex]   
$\axiomref{h.u.t.\x 1}$ for  $ev(-,g)$ & \cite[Lemma 2.6 part i)]{Fem:Fil} (vertical composition of v.s.t. $\alpha_0,\beta_0$) 
\\ [1ex]   
$\axiomref{h.u.t.\x 2}$ for  $ev(-,g)$ & hold clearly \\ [1ex]    
$\axiomref{v.l.t.\x 2}$ for $ev(-,v)$ & hold clearly \\ [1ex]   
$\axiomref{v.l.t.\x 5}$ for $ev(-,v)$ & $\axiomref{m.hu\x vs}$ for $ev(b,-)$  \\ [1ex]   
\end{tabular}
\caption{Why $ev(-,g)$ is a horizontal unnatural and $ev(-,v)$ a vertical strict transformation}
\label{table:1}
\end{center}
\end{table}
For a 2-cell $\beta$ in $\Bb$ and an object $F\in [\Bb, \Cc]^{lx}$ by definition we have $F(\beta)=ev(-,\beta)$. The modification axiom 
\axiomref{m.hu\x vs} for $ev(-,\beta)$ at a 1v-cell $\alpha_0$ in $Lax_{hop}(\Bb, \Cc)$ coincides with the transformation axiom 
$(v.l.t.\x 5)$ for $\alpha_0$. Thus $ev(-,\beta)$ is a mixed funny modification. 
This finishes the proof that $ev$ is a lax double funny functor by \prref{funny fun}, and the last statement is clear. 
\qed\end{proof}

In particular, given an $m$-ary double funny functor $F^m:\Aa_1\times...\times\Aa_m\to[\Cc, \Dd]^{lx}$ with $m\geq 1$ and the unary map 
$F^1=ev_C=ev(-,C)$ for an object $C\in\Cc$, by \equref{nul-ass} we have 
\begin{equation}\eqlabel{first-last}
(ev_C\circ F^m)(...,A_k,...)=ev_C\circ F^m(...,A_k,...)
\end{equation}
as $(m\x 1)$-ary funny functors. 

\smallskip

The  following will show useful.

\begin{lma} \lelabel{binary to inner}
Let $F\in f\x[\Aa\times\Bb,[\Cc,\Dd]^*]^*$. Then: \vspace{-0,2cm} 
\begin{itemize}
\item $F(-,-)(h):\Aa\times\Bb\to\Dd$ is a horizontal binary unnatural transformation; \vspace{-0,2cm} 
\item $F(-,g)(-):\Aa\times\Cc\to\Dd$ is a horizontal binary unnatural transformation; \vspace{-0,2cm} 
\item $F(f,-)(-):\Bb\times\Cc\to\Dd$ is a horizontal binary unnatural transformation; \vspace{-0,2cm} 
\item $F(-,-)(z):\Aa\times\Bb\to\Dd$ is a vertical binary strict transformation; \vspace{-0,2cm} 
\item $F(-,v)(-):\Aa\times\Cc\to\Dd$ is a vertical binary strict transformation; \vspace{-0,2cm} 
\item $F(u,-)(-):\Bb\times\Cc\to\Dd$ is a vertical binary strict transformation. 
\end{itemize}
\end{lma}

\begin{proof}
For the first part, to prove that one has two unary horizontal unnatural transformations $F(A,-)(h)$ and $F(-,B)(h)$, 
observe that the first one is $\F_1(-)(k)=(k,-)_1$ and the second one is $\F_2(-)(k)=(k,-)_2$ in the notation of \prref{char lax fdf}, 
so that by \prref{funny fun} they are h.u.t. The agreement on objects of these h.u.t. follows by the agreement on objects of the unary components of $F$. The axiom \axiomref{$HUT^f$} holds by the property $(u,U)=\Id$ for the double funny 
functor $F$. 

The axiom \axiomref{$HUT^f$} for $F(-,g)(-)$ is the axiom \axiomref{m.ho-vl.-2} for the modification $\Theta=F(u,g)$ between 
horizontal unnatural transformations $\alpha=F(A,g)$ and $\beta=F(\tilde A,g)$ (as 1h-cells in $[\Cc,\Dd]^*$) and vertical strict transformations $\alpha_0=F(u,B)$ and $\beta=F(u,\tilde B)$ (as 1v-cells in $[\Cc,\Dd]^*$). Clearly, $F(A,g)(-):\Cc\to\Dd$ is a h.u.t 
as $F(A,g)$ is a 1h-cell in $[\Cc,\Dd]^*$. Functoriality of $F$ in the first variable yields that $F(-,g)(C):\Aa\to\Dd$ 
is a h.u.t. The agreement of these two unary h.u.t. on objects follows by \equref{first-last} with $F^m=F(-,g)$ and $A_k=C$. 

The third claim is proved analogously as the second one. 

One obtains two unary vertical strict transformations for $F(-,-)(z)$ and their agreement on objects similarly as for $F(-,-)(h)$ 
in the first part. As in \prref{ternary eq} one has that \axiomref{$HUT^f$} for $F(-,g)(-)$ (from the second part) is 
\axiomref{$VST^f_2$} for $F(-,-)(z)$. Similarly, \axiomref{$HUT^f$} for $F(f,-)(-)$ (from the third part) is 
\axiomref{$VST^f_1$} for $F(-,-)(z)$. 

The above proofs are sufficient to deduce the proofs for the remaining two claims. 
\qed\end{proof}

By iterating \equref{funny-basic-d} we obtain double category isomorphisms
\begin{equation}\eqlabel{34}
f\x[\Aa\times\Bb,[\Cc,\Dd]^*]^*\iso [\Aa, [\Bb,[\Cc,\Dd]^*]^*]^*\iso [\Aa, f\x[\Bb\times\Cc,\Dd]^*]^*. 
\end{equation}
A double funny functor from the left $F\in f\x[\Aa\times\Bb,[\Cc,\Dd]^*]^*$ consists of unary functors $F(A,-):\Bb\to [\Cc,\Dd]^*$ and 
$F(-,B):\Aa\to [\Cc,\Dd]^*$ for any $A\in\Aa, B\in\Bb$. These in turn determine unary $*$-typed double functors $F(A,-)(C):\Bb\to\Dd$, 
$F(A,-)_B:\Cc\to\Dd$, and $F(-,B)(C):\Aa\to\Dd, F(-,B)_A:\Cc\to\Dd$ respectively, whereas because of the binary funny property for $F$ the two functors $F(A,-)_B,F(-,B)_A:\Cc\to\Dd$ coincide. To prove that $F$ defines a ternary funny functor $F(-,-,-)\in[\Aa\times\Bb\times\Cc,\Dd]^{lx}$ by $F(a,b,c):=F(a,b)(c)$ for sensible cells $a\in\Aa,b\in\Bb,c\in\Cc$ (with abuse of notation), we should show that 
the above three unary functors satisfy two by two the binary funny functor property, and then that the induced three binary funny functors - that 
simultaneously are induced by $F(a,b,c)$ - 
satisfy the ternary funny property and axioms. We already have the binary funny functor $F(-,-,C):\Aa\times\Bb\to\Dd$ 
and the ternary property $F(A,-,-)_B=F(-,B,-)_A:\Cc\to\Dd$ by the properties of the original $F\in f\x[\Aa\times\Bb,[\Cc,\Dd]^*]^*$. 

Let $\crta F:\Aa\to f\x[\Bb\times\Cc,\Dd]^*$ denote the induced lax $*$-typed double functor on the right of \equref{34}. It is 
$\crta F(a)(b,c):=F(a,b)(c)$ for sensible cells $a\in\Aa,b\in\Bb,c\in\Cc$. Then 
$\crta F(A)=F(A,-,-):\Bb\times\Cc\to\Dd$ is a binary funny map for all $A\in\Aa$, so we have  $F(A,B,-)_C=F(A,-,C)_B$. Moreover, 
$F(a,B,-)_C=F(a,-,C)_B$ for any higher cell $a\in\Aa$ would mean that $\crta F(a)(B,-)_C=\crta F(a)(-,C)_B$. But this is true, since 
$\crta F$ induces a $*$-typed double functor $\crta F(-)(B,C):\Aa\to\Dd$. Thus we have the second ternary property for $F(-,-,-)$. 

The third binary property $F(A,B,-)_C=F(-,B,C)_A$ on objects now follows by \equref{first-last} with $F^1=F(-,B):\Aa\to [\Cc,\Dd]^*$ for a fixed $B\in\Bb$ and $k=1$. 
The same equality implies $F(A,-,-)_C=F(-,-,C)_A$ as unary functors, by applying $F^2=F(-,-)$ and $k=1$. 
The proof that the 11 axioms for the remaining binary funny functor $F(-,B)(-):\Aa\times\Cc\to\Dd$ hold is analogous as in the proof of \leref{un-fun}, where we apply (i) and (ii) of \prref{funny fun}. Namely, 
\axiomref{($k'k, U$)}, \axiomref{($1_B,U$)} and \axiomref{$(u,U)$-r-nat} hold by \axiomref{v.l.t.-1}, \axiomref{v.l.t.-2} and 
\axiomref{v.l.t.-5} of the vertical strict transformation $F(U,B)$ evaluated at $k'k, 1_C$ and $u$, respectively; 
\axiomref{($\frac{u}{u'}, K$)} and \axiomref{($1^B,K$)} hold by \axiomref{h.u.t.\x 1} and \axiomref{h.u.t.\x 2} of the horizontal unnatural transformation $F(k,B)$ evaluated at $\frac{u}{u'}$ and $1_C$, respectively. On the other hand, 
\axiomref{($u,1_A$)} holds since $F(1_A,B)$ as identity horizontal unnatural transformation evaluated at a 1v-cell $u$
is the identity 2-cell on $F(A,B)(u)$, and similarly for \axiomref{($k,1^A$)} use $F(1^A,B)(k)$; 
for \axiomref{($u, K'K$)} apply the composition of horizontal transformations from \leref{vert comp hor.ps.tr.} to $F(k'k,B)(v)$, 
and for \axiomref{($k,\frac{U}{U'}$)} apply the composition of vertical transformations from \leref{vert comp vert. lx tr.} to 
$F(\frac{U}{U'},B)(k)$; the axiom \axiomref{$(u,U)$-l-nat} holds by the modification axiom \axiomref{m.hu-vs} of $F(\zeta,B)$ evaluated at $u$. 
Finally, $F(U,B)(u)$ is an identity 2-cell, as $F(U,B)$ is a vertical strict transformation. 

It remains to verify the ternary axioms. 
Knowing that  $F\in f\x[\Aa\times\Bb,[\Cc,\Dd]^*]^*$, by \leref{binary to inner} we have that 
$F(-,-)(h), F(-,g)(-), F(f,-)(-)$ are binary h.u.t. The three ternary axioms follow now by \prref{ternary eq} with 
$F(f,-,-):=F(f,-)(-), F(-,g,-):=F(-,g)(-)$ and $F(-,-,h):=F(-,-)(h)$.  

\medskip

Using the identifications as above: $F(-,*,\bullet)=F(-,*)(\bullet)=\crta F(-)(*,\bullet)=\tilde F(-)(*)(\bullet)$, we 
obtain the isomorphism of double categories 
\begin{equation} \eqlabel{ternary}
f_3\x[\Aa\w\times\Bb\times\Cc,\Dd]^*\iso f\x[\Aa\w\times\Bb,[\Cc,\Dd]^*]^*.
\end{equation}

We now may prove:

\begin{thm} \thlabel{left closed multicat} 
\begin{enumerate}
\item
Given a lax double funny functor of $n$-variables $H:\Pi_{i=1}^{n-1} \Aa_i\times\Bb\to\Cc$ of type $*$, 
there is a unique lax double funny functor of $n\x 1$-variables $H^t:\Pi_{i=1}^{n-1} \Aa_i\to [\Bb,\Cc]^*$ of type $*$ 
such that $ev(H^t\times \Id_\Bb)=H$. 
\item 
The above correspondence extends to a natural isomorphism of double categories: 
$$f_{n}\x [\Pi_{i=1}^{n-1} \Aa_i\times\Bb,\Cc]^*\iso f_{n-1}\x [\Pi_{i=1}^{n-1} \Aa_i, [\Bb,\Cc]^*]^*.$$
\end{enumerate} 
\end{thm}

\begin{proof} 
We only prove the second part on the level of objects. Set $\Bb=\Aa_n$. 
Given an $n$-ary $*$-typed double funny functor $G:\Pi_{i=1}^{n-1} \Aa_i\times\Bb\to\Cc$, define $\crta G(a_1,...,a_{n-1})(a_n)
:=G(a_1,...,a_n)$ for sensible cells $a_i\in\Aa_i, i=1,..,n$. Then $\crta G$ clearly consists of binary funny functors satisfying the 
ternary funny property and axioms, as so does $G$. 

Let $\F:\Pi_{i=1}^{n-1}\Aa_i\to[\Bb,\Cc]^*$ be from the right, and define $\crta\F:\Pi_{i=1}^n\Aa_i\to\Cc$ by $\crta\F(a_1,...,a_n)=
\F(a_1,...,a_{n-1})(a_n)$ for sensible cells $a_i\in\Aa_i, i=1,..,n$. 
Let us consider the binary 
funny functor property of $\F$ at arbitrarily chosen two variables $i,j$ of $\Pi_{l=1}^{n-1}\Aa_l$, and let 
$A_i\in\Aa_i, A_j\in\Aa_j, B\in\Bb$. Then $\crta\F(-,-,B)=\F(-,-)(B):\Aa_i\times\Aa_j\to\Cc$ is a funny functor, which is a binary component  of both $\crta\F$ and $\F$. On the other hand, $\F(A_i,-)(-)\w:\Aa_j\to[\Bb,\Cc]^*$ by \equref{funny-basic-d} 
determines a double funny functor $\crta\F(A_i,-,-)\w:\Aa_j\times\Bb\to\Cc$, which is a binary component of $\crta\F$. Instead of $i$ 
we could have done the same reasoning with $j$. This way we obtain three component binary funny functors for $\crta\F$ for the variables $i,j,n$. 
%
The ternary funny functor property \equref{obj} and axioms hold as in the proof of \equref{ternary}, this finishes the proof. 
\qed\end{proof}

We have thus proved that the funny multicategory for double categories $\Dd bl^*$ with $*$-type of double functors and inner-hom $[-,-]^*$ is left closed. 

\medskip

Symmetrically to the double funny functor $ev:[\Bb, \Cc]^{lx}\times\Bb\to\Cc$ from \prref{ev} there is a double funny functor $\crta ev:
\Bb\times [\Bb, \Cc]^{lx}\to\Cc$ that induces the double category isomorphism  
\begin{equation} \eqlabel{right-cl}
f_{n}\x [\Bb\times\Pi_{i=1}^{n-1} \Aa_i,\Cc]^*\iso f_{n-1}\x [\Pi_{i=1}^{n-1} \Aa_i, [\Bb,\Cc]^*]^*
\end{equation}
implying that $\Dd bl^*$ with the same inner-hom $[-,-]^*$ is right closed.

\subsubsection{Representability} \ssslabel{rep}

In this subsection we finally prove that the funny multicategory for double categories $\Dd bl$ with {\em strict} double functors, 
{\em i.e.} with $*=st$, is representable. Observe that $\Box^*$ is not representable for $*$ other than $st$: 
representability requires that double functors from $\Aa\Box^*\Bb$ on the left and the multimaps on the right of \equref{no-rep} 
be of the same type, this forces $*=st$. 

By left and right closedness of $\Dd bl$ one can prove 
separately left and right representability of $\Box$ (the funny product 
with $*=st$):
$$
f_n\x [(\Aa\Box\Bb)\w\times\w(\Pi_{i=1}^{n-1}\Cc_i),\Dd]  
\iso f_{n+1}\x [\Aa\w\times\w \Bb\w\times\w(\Pi_{i=1}^{n-1}\Cc_i),\Dd] \quad\text{left representability}  
$$
$$
f_n\x [(\Pi_{i=1}^{n-1}\Aa_i)\times (\Bb\Box\Cc),\Dd]
\iso f_{n+1}\x [(\Pi_{i=1}^{n-1}\Aa_i)\w\times \w\Bb\w\times\w\Cc,\Dd] \quad\text{right representability} 
	$$
for $n\geq 2$. 
Namely, by multicategory-closedness \thref{left closed multicat} we can pass the $\Pi_{i=1}^{n-1}\Cc_i$'s to the codomain to get an equivalent question. Thus for left representability it is enough to prove: 
$[\Aa\Box\Bb,\Cc]\iso f\x [\Aa\times \Bb, \Cc]$ for any double category $\Cc$, but this is true by \equref{ot-f-d}. 
For the other isomorphism use \equref{right-cl} (to drop out the $\Pi_{i=1}^{n-1}\Aa_i$'s).

\bigskip

Now set 
\begin{equation} \eqlabel{acumul left}
\crta\Box^n\Aa_i=(..(\Aa_1\Box\Aa_2)\Box...\Box\Aa_{n-1})\Box\Aa_n.
\end{equation}

\begin{cor} \colabel{strict-rep} [{\em representability}]
There are double category isomorphisms:
\begin{enumerate}
\item 
$$f_n\x[(\Pi_{i=1}^{r}\Aa_i)\w\times\w(\Bb\Box\Cc)\w\times\w(\Pi_{1}^{s}\Dd_i),\Ee]  \iso 
f_{n+1}\x [(\Pi_{i=1}^{r}\Aa_i)\w\times\w(\Bb\times\Cc)\w\times\w(\Pi_{i=1}^{s}\Dd_i),\Ee] $$
natural in $\Aa_i$'s, $\Dd_i$'s and $\Ee$, with $n=r+1+s$;
\item 
$$f_n\x [(\Pi_{i=1}^{r}\Aa_i)\w\times\w(\crta\Box\Bb_i)\w\times\w(\Pi_{1}^{s}\Dd_i),\Ee]  \iso
f_p\x [(\Pi_{i=1}^{r}\Aa_i)\w\times\w(\Pi_{i=1}^k\Bb_i)\w\times\w(\Pi_{i=1}^{s}\Dd_i),\Ee] $$
natural in $\Aa_i$'s, $\Dd_i$'s and $\Ee$, with $n=r+1+s$ and $p=r+k+s$;
\item $[\crta\Box\Aa_i,\Bb]\iso f_n\x [\Pi_{i=1}^n\Aa_i,\Bb]$ \hspace{0,2cm} natural in $\Bb$. 
\end{enumerate}
\end{cor}

\begin{proof}
For the first part, by left and right representability of $\Box$ it remains only to verify the ternary property \equref{obj} and axioms of 
funny functors in both sides of 
$$f_3\x[\Aa\w\times\w(\Bb\Box\Cc)\w\times\w\Dd,\Ee]  \iso f_4\x [\Aa\w\times\w\Bb\times\Cc\w\times\Dd,\Ee] $$
for the variables living in $\Aa\times\Bb\times\Dd$, without loss of generality. The double funny functors on both sides are defined in the obvious way. The $3+3$ equalities that are to be checked are
verified by including in both sides an arbitrarily fixed object $C\in\Cc$. The functorialities in $\Aa$ and $\Dd$ at any side of the isomorphism in question follow easily by the known functorialities in $\Aa$ and $\Dd$ of the ternary property holding on the other side. Let us examine functoriality in $\Bb$ 
and let $G$ be from the left and $F$ from the right above. 
On the left-hand side it means that the equality $G(A, -\Box C, -)\vert_D=G(-, -\Box C, D)\vert_A$ holds. This equality 
holds if and only if  the equality $F(A, -, C, -)\vert_D=F(-, -, C, D)\vert_A$ is true. 
Since $F$ is a 4-ary double funny functor, we know that $F(-,-,C,-)$ is its component ternary double funny functor. 
The desired equality then holds by the agreement on objects \equref{obj}. 

The ternary axioms easily follow by the ternary axioms holding on the other side of the isomorphism in question, having in mind that 
$F(-,-,C,-)$ behaves as a ternary double funny functor.  

The second part follows from the first part by iteration, and the third one 
follows from the second part with $r=s=0$.  
\qed\end{proof}

\smallskip

The points 2. and 3. above mean that the multicategory $\Dd bl$ is representable. Then we can finally claim:

\medskip

\begin{thm} \thlabel{Dbl left closed mon}
The category $(Dbl, \Box)$ is biclosed monoidal. 
\end{thm}

We make precise how the associativity constraint for the funny monoidal product is obtained. Consider the following chain of bijections 
holding in each step by the double category isomorphism \equref{ot-f-d} (in the second and third bijection apply also left, respectively right representability):
\begin{align}
[(\Aa\Box\Bb)\Box\Cc,\Dd] & \iso f\x [(\Aa\Box\Bb)\times\Cc,\Dd] \nonumber\\
& \iso f_3\x [\Aa\times \Bb\times\Cc,\Dd] \label{eq2}\\ 
& \iso  f\x [\Aa\times (\Bb\Box\Cc),\Dd] \nonumber\\
& \iso [\Aa\Box(\Bb\Box\Cc),\Dd]. \nonumber
\end{align}
(Recall that for (left and right) representability it is necessary that the funny product is constructed via strict double functors, 
this is where we need it.) By the  Yoneda lemma  there is an isomorphism double functor 
\begin{equation} \eqlabel{alfa on funny prod of dcats}
\alpha_{\Aa,\Bb,\Cc}:(\Aa\Box\Bb)\Box\Cc\to\Aa\Box(\Bb\Box\Cc).
\end{equation}
In particular, 
setting $\Dd=\Aa\Box(\Bb\Box\Cc)$, we obtain $\alpha_{\Aa,\Bb,\Cc}$ 
as the image of the double funny functor $J_{\Aa,\Bb\Box\Cc}(\Id_\Aa\times J_{\Bb,\Cc})$ in (\ref{eq2}), with $J=J^{st}$ from \equref{J}. 
Similarly, a double functor $\Aa\Box(\Bb\Box\Cc)\to (\Aa\Box\Bb)\Box\Cc$ is obtained by setting $\Dd=(\Aa\Box\Bb)\Box\Cc$ 
as the correspondent to $J_{\Aa\Box\Bb,\Cc}(J_{\Aa,\Bb}\times\Id_\Cc)$. 

Moreover, for double functors $F,G,H$ and $\Dd=\Aa'\Box(\Bb'\Box\Cc')$ 
the two composites of double functors 
$(F\Box(G\Box H))\alpha_{\Aa,\Bb,\Cc}$ and $\alpha_{\Aa',\Bb',\Cc'}((F\Box G)\Box H)$ 
living in the left-hand side of (\ref{eq2}) 
can be shown to be equal by considering their corresponding ternary funny functors 
$\Aa\times\Bb\times\Cc\to\Aa'\Box(\Bb'\Box\Cc')$. Here the double functors $F\Box(G\Box H), (F\Box G)\Box H$ are as in \rmref{functors on Box}. Observe that both they and $\alpha_{\Aa,\Bb,\Cc}$ are given by three one-variable double functors. That $\alpha_{\Aa,\Bb,\Cc}$ obeys the pentagon it means that four pentagons commute determined by the constituent one-variable double functors, analogously as in 
\deref{assoc}. 

\subsubsection{The purely funny multicategory} \ssslabel{purely funny mult}

We finally record that all the results obtained above in this section hold in their analogous form for purely funny functors 
and the purely funny product $\Aa\Box_f\Bb$ of double categories. 

As for the notions, 
ternary purely funny functors for double categories are not required to satisfy any axioms, in particular the three axioms appearing in 
\deref{double qf 3}. Horizontal and vertical purely unnatural transformations of binary purely funny functors do not satisfy any axioms 
either (no axiom \axiomref{$HUT^f$}, nor its vertical analogues). One obtains a biclosed monoidal category $(Dbl,\Box_f)$.

\subsection{From a double funny functor to a premonoidal double category}

According to \coref{purely funny binoidal} and in light of \sssref{purely funny mult} we may say: 

\begin{prop} \prlabel{purely funny premon}
In any premonoidal double category $\Bb$ the binoidal structure is given by a pseudodouble purely funny functor 
$H:\Bb\times\Bb\to\Bb$. 
\end{prop}

On the other hand, pseudodouble (mixed) funny functors have more structure and hence induce richer binoidal structures.  
In \coref{funny yields binoidal} we saw that such a funny functor $H:\Bb\times\Bb\to\Bb$ 
induces a binoidal structure on $\Bb$ so that all 1v-cells are central, and there are 
2-cells $K\ltimes-\vert_u=-\rtimes u\vert_K$ and $U\ltimes-\vert_k=-\rtimes k\vert_U$, for all 1h-cells $K,k$ and 1v-cells $U,u$ in $\Bb$. 

To discuss associativity let us set some terminology first. 
In Tables \ref{table:2} and \ref{table:3} the 24 axioms for the associativity constraint in a premonoidal double category are listed. Each of the first columns of these Tables contains two subcolumns. The first of these four columns depends on two 1h-cells, the middle two columns depend on one 1h- and one 1v-cell, while the fourth column depends on two 1v-cells. 
Let us refer to the axioms from the last three columns as ``$0+6+6+6$ axioms'' for short, whereby the order of columns is 
respected. 

\smallskip

An associative binoidal structure on a pseudodouble funny functor $H$ supposes the existence of invertible vertical strict transformations three $\alpha$'s. Hereby 
the $\alpha$'s play the role of an associativity that acts on $H(H\times 1)\Rightarrow H(1\times H)$, which by \prref{subst-gen funny} are ternary double funny functors. So the thee $\alpha$'s are the three one-variable vertical transformations appearing in \deref{vs tr of ff 3} and they obey six identities. Because of the properties $K\ltimes-\vert_u=-\rtimes u\vert_K$ and $U\ltimes-\vert_k=-\rtimes k\vert_U$ holding for $H$ the latter six identities correspond simultaneously 
to the 
$0+6+0+0$ axioms: 
\axiomref{$(u\ltimes,g,C)$}, \axiomref{$(A, v\ltimes,h)$}, \axiomref{$(u\ltimes, B,h)$}, \axiomref{$(f,B,\rtimes z)$}, 
\axiomref{$(f,\rtimes v,C)$}, \axiomref{$(A,g,\rtimes z)$}, and to the 
$0+0+6+0$ axioms: 
\axiomref{$(f\ltimes,v,C)$}, \axiomref{$(f\ltimes,B,z)$}, \axiomref{$(A,g\ltimes,z)$}, \axiomref{$(u,\rtimes g,C)$}, 
\axiomref{$(u,B,\rtimes h)$}, \axiom{$(A,v,\rtimes h)$}. Moreover, since the structure 2-cells $(u,U)$ of the binary funny functors are trivial, the 
 $0+0+0+6$ axioms are automatically fulfilled: 
\axiomref{$(u\ltimes,v,C)$}, \axiomref{$(A, v\ltimes,z)$}, \axiomref{$(u\ltimes, B,z)$}, \axiomref{$(u,B,\rtimes z)$}, 
\axiomref{$(u,\rtimes v,C)$}, \axiomref{$(A,v,\rtimes z)$}. A similar analysis of an associativity on ternary funny functors we will carry out in more details in the proof of \thref{funny pseudomonoid}. From the above said we have:

\begin{prop} \prlabel{funny ->premon}
Let $\Bb$ be a premonidal double category. If its 
binoidal structure comes from a double funny functor, 
then it satisfies the $0+6+6+6$ axioms. 
\end{prop}


\subsection{Funny product and premonoidal double categories} \sslabel{funny-premonoidal connection} 

Binoidal structures in 1-categories were introduced in \cite{PR} for general premonoidal categories, and a strict premonoidal category was defined as a monoid in the category of categories with the funny product. A strict premonoidal category is then a 
general premonoidal category in which the monoidality constraints $\alpha,\lambda,\rho$ are trivial. We define a strict premonoidal double category in an analogous way. 

\begin{defn}
A {\em strict} premonoidal double category is a monoid in $(Dbl, \Box_f)$. 
\end{defn}

Observe that this is the setting of \sssref{purely funny mult}. 
In a monoid structure $(\Dd,M,U)$ the double functor $M:\Dd\Box_f\Dd\to\Dd$ is given by a binary double purely funny functor 
$\Dd\times\Dd\to\Dd$. 
According to \sssref{pure funny} (with $*=st$) the latter is given by a pair of families of double functors: $(-,A):\Dd\to\Dd$ and 
$(B,-):\Dd\to\Dd$ for $A,B\in\Dd$. 
Then a strict premonoidal double category $\Bb$ is 
a premonoidal double category in which the two pseudodouble functors underlying the binoidal structure of $\Bb$ are strict and the three 
$\alpha$'s, $\lambda$ and $\rho$ are identities. 

\medskip

As we announced at the beginning of this section, we are now going to construct a monoidal 
2-category of double categories with a funny type of monoidal product, so that a 
pseudomonoid in it is a particular kind of a premonoidal double category. 
Namely, observe that to recover a binoidal structure consisting of two pseudo double functors,  
one should consider the pseudo version of the funny product, according to \equref{ot-f-d}. In view of the previous paragraph, 
one would then hope that by extending the functor $-\Box^{ps}-:Dbl\times Dbl\to Dbl$ (and the natural transformations of its monoidality constraints) to a 2-functor on a suitable 2-category of double categories (and the corresponding pseudonatural transformations), one would get a monoidal 2-category so that a pseudomonoid in it would give rise to a premonoidal double category with a non-strict binoidal structure that is non-strictly associative and unital. 
However, monoidality of a 1-category, and then also of the mentioned 2-category,  
is possible to achieve only for the strict version of the funny product: recall the beginning of \sssref{rep}. In this way we will 
cover premonoidal double categories with strict binoidal structure and non-trivial monoidality constraints.

\smallskip

We denote by $Dbl_2$ the 2-category of double categories, double functors and vertical strict transformations.  
Horizontal composition of vertical pseudonatural transformations is given due to 
\cite[Lemma 3.6]{Fem:Fil} by 
$$[\alpha\vert\beta]_f= 
\scalebox{0.82}{
\bfig
\putmorphism(-250,500)(1,0)[F\s'F(A)`F\s'F(B)` F\s'F(f)]{700}1a
 \putmorphism(-250,50)(1,0)[F\s'G(A)`F\s'G(B)` F\s'G(f)]{700}1a
 \putmorphism(-250,-400)(1,0)[G'G(A)`G'G(B)` G'G(f)]{700}1a

\putmorphism(-280,500)(0,-1)[\phantom{Y_2}``F\s'(\alpha(A))]{450}1l
 \putmorphism(-280,70)(0,-1)[\phantom{F(A)}` `\beta(G(A))]{450}1l

\putmorphism(450,500)(0,-1)[\phantom{Y_2}``F\s'(\alpha(B))]{450}1r
\putmorphism(450,70)(0,-1)[\phantom{Y_2}``\beta(G(B))]{450}1r
\put(-120,290){\fbox{$F\s'(\alpha_f)$}}
\put(-80,-150){\fbox{$\beta_{G(f)}$}}
\efig}
\quad\quad
[\alpha\vert\beta]^u=
\scalebox{0.82}{
\bfig
 \putmorphism(-150,500)(1,0)[F\s'F(A)`F\s'F(A) `=]{550}1a
\putmorphism(-180,500)(0,-1)[\phantom{Y_2}`F\s'F(A') `F\s'F(u)]{450}1l
\put(-80,-160){\fbox{$F\s'(\alpha )^u$}}
\putmorphism(-150,-400)(1,0)[F\s'G(A')`F\s'G(A') `=]{550}1a
\putmorphism(-180,50)(0,-1)[\phantom{Y_2}``F\s'(\alpha (A'))]{450}1l
\putmorphism(380,500)(0,-1)[\phantom{Y_2}` `F\s'(\alpha (A))]{450}1r
\putmorphism(380,50)(0,-1)[F\s'G(A)` `F\s'G(u)]{450}1r
\putmorphism(520,60)(1,0)[`F\s'G(A)`=]{460}1a
\putmorphism(370,-850)(1,0)[\phantom{G'G(A)}``=]{440}1b
\putmorphism(960,50)(0,-1)[\phantom{(B, \tilde A')}``\beta (G(A))]{450}1r
\putmorphism(960,-400)(0,-1)[G'G(A)`G'G(A').` G'G(u)]{450}1r
\putmorphism(400,-400)(0,-1)[\phantom{(B, \tilde A)}`G'G( A') `\beta (G(A'))]{450}1l
\put(500,-630){\fbox{$\beta^{G(u)}$}}
\efig}
$$

\smallskip

We will now extend the functor $-\Box-:Dbl\times Dbl\to Dbl$ to a 2-functor. 
We can do it in two ways: one defined on the Cartesian product $Dbl_2\times Dbl_2\to Dbl_2$ and the other 
$Dbl_2\star Dbl_2\to Dbl_2$ on 
the funny product $\star$ of 2-categories from \ssref{funny 2-cats}. The former 2-functor induces the latter. 
We will define both product 2-functors for completeness, 
any of them serves for our final construction. 
We will denote both of them by $-\Box_2-$ abusing notation. 
Then we will make $(Dbl_2,-\Box_2-:Dbl_2\star Dbl_2\to Dbl_2)$ into a pseudomonoid in a funny monoidal 2-category of 
2-categories $(2\x\Cat,\star)_2$, which we will explain below. Similarly, 
$(Dbl_2,-\Box_2-:Dbl_2\times Dbl_2\to Dbl_2)$ will be a pseudomonoid in the Cartesian monoidal 2-category of 
2-categories $(2\x\Cat,\times)_2$. 
This way $(Dbl_2,-\Box_2-)$ becomes a monoidal 2-category in two ways, so that a pseudomonoid in both of them coincides and  will be the desired premonoidal double category. 


\medskip

\begin{prop} \prlabel{no funny on Cart}
The functor $-\Box-:Dbl\times Dbl\to Dbl$ extends to a 2-functor $-\Box_2-:Dbl_2\times Dbl_2\to Dbl_2$. 
\end{prop}

\begin{proof}
Let $\alpha: F\Rightarrow F':\Aa\to\Aa'$ and $\beta:G\Rightarrow G':\Bb\to\Bb'$ be vertical strict transformations of double functors. 
To define $\alpha\Box_2\beta:F\Box G\Rightarrow F'\Box G':\Aa\Box\Bb\to\Aa'\Box\Bb'$, 
we set $(\alpha\Box_2\beta)(A\Box B)=\frac{\alpha(A)\Box G(B)}{F'(A)\Box\beta(B)}$ for a 1v-cell component,  
for 1h-cells $K$ in $\Aa$ and $k$ in $\Bb$ we set $(\alpha\Box_2\beta)(K\Box B)$ and $(\alpha\Box_2\beta)(A\Box k)$ as below 
$$
\scalebox{0.86}{
\bfig
\putmorphism(-100,50)(1,0)[F(A)\Box  G(B)`F(A')\Box  G(B)`F(K)\Box  G(B)]{1050}1a
\putmorphism(-150,-400)(1,0)[F'(A)\Box  G(B)`F'(A')\Box  G(B) `F'(K)\Box  G(B)]{1040}1a
\putmorphism(-180,50)(0,-1)[\phantom{Y_2}``\alpha(A)\Box  G(B)]{450}1l
\putmorphism(900,50)(0,-1)[\phantom{Y_2}``\alpha(A')\Box  G(B)]{450}1r
\put(90,-170){\fbox{$\alpha_K\Box  G(B)$}}
%
\putmorphism(900,-370)(0,-1)[\phantom{Y_2}``F'(A')\Box  \beta(B)]{450}1r 
%
\putmorphism(-180,-370)(0,-1)[\phantom{Y_2}``F'(A)\Box  \beta(B)]{450}1l
\putmorphism(-150,-830)(1,0)[F'(A)\Box  G'(B)`F'(A')\Box  G'(B) `F'(K)\Box  G'(B)]{1140}1a
\put(90,-590){\fbox{$F'(K)\Box\beta(B)$}}
\efig}
\qquad
\scalebox{0.86}{
\bfig
\putmorphism(-100,50)(1,0)[F(A)\Box  G(B)`F(A)\Box  G(B')`1_{F(A)}\Box  G(k)]{1050}1a
\putmorphism(-150,-400)(1,0)[F'(A)\Box  G(B)`F'(A)\Box  G(B') `1_{F'(A)}\Box  G(k)]{1040}1a
\putmorphism(-180,50)(0,-1)[\phantom{Y_2}``\alpha(A)\Box  G(B)]{450}1l
\putmorphism(900,50)(0,-1)[\phantom{Y_2}``\alpha(A)\Box  G(B')]{450}1r
\put(90,-170){\fbox{$\alpha(A)\Box  G(k)$}}
\putmorphism(900,-370)(0,-1)[\phantom{Y_2}``F'(A)\Box  \beta(B')]{450}1r 
\putmorphism(-180,-370)(0,-1)[\phantom{Y_2}``F'(A)\Box  \beta(B)]{450}1l
\putmorphism(-150,-830)(1,0)[F'(A')\Box  G'(B)`F'(A)\Box  G'(B') `1_{F'(A)}\Box  G'(k)]{1140}1a
\put(90,-590){\fbox{$F'(A)\Box\beta_k$}}
\efig}
$$
respectively, for 2-cell components of $\alpha\Box_2\beta$ in $\Aa'\Box_2\Bb'$, and 
for 1v-cells $U$ in $\Aa$ and $u$ in $\Bb$ we set $(\alpha\Box_2\beta)(U\Box B)=\alpha^U\Box G(B)$ and $(\alpha\Box_2\beta)(A\Box u)=
F(A)\Box\beta^u$, which are both identity 2-cells by \equref{Box on Id's}. The axioms \axiomref{v.l.t.\x 1}, \axiomref{v.l.t.\x 2} and 
\axiomref{v.l.t.\x 5} for $\alpha\Box_2\beta$ to be a vertical strict transformation are fulfilled because of the relations 
\equref{four strictified}, \equref{Box on Id's} and \equref{frac-natur} holding in the funny product of double categories. 

The 2-functor property of $-\Box_2-$ on 2-cells of $Dbl_2\times Dbl_2$ means that the equalities 
$$\frac{\alpha\Box_2\beta}{\alpha'\Box_2\beta'}=\frac{\alpha}{\alpha'}\Box_2\frac{\beta}{\beta'}$$
$$[\alpha\vert\gamma]\Box_2[\beta\vert\delta]=[\alpha\Box_2\beta\vert\gamma\Box_2\delta]$$
$$\Id_F\Box_2\Id_G=\Id_{F\Box_2 G}$$
should hold for double functors $F,G$ and suitably composable vertical unnatural transformations $\alpha,\beta,\alpha', \beta', 
\gamma,\delta$.  
The first two equations are given in terms of vertical compositions of transformation components. The first one holds by 
\equref{frac-natur}, the second uses \equref{vertic comp 2-cells}, \equref{hor-ver comp} and \equref{frac-natur}, and 
the third equation holds by the first two identities in \equref{hor-ver comp}. 
\qed\end{proof}

We now introduce the announced 2-functor on the funny monoidal product $Dbl_2\star Dbl_2$, and then we will consider a 
``funny-pseudomonoid'', that is, 
a pseudomonoid in the funny monoidal 2-category of 2-categories $(2\x\Cat,\star)_2$. Its underlying monoidal 
category is $(2\x\Cat,\star)$ from \ssref{funny 2-cats}, and 2-cells the unnatural transformations (consisting only of 1-cell components). For unnatural transformations $\alpha,\beta$ from 2-categories $\A,\B$ we define  
$(\alpha\star B)_A=\alpha(A)\star B$ and $(A\star\beta)_B=A\star\beta(B)$ with $A\in\A, B\in\B$. Given that 
the associativity constraint $\tilde\alpha$ for the category $(2\x\Cat,\star)$ is strictly natural and satisfies the pentagon strictly, for the 
2-cell components of its extension to $(2\x\Cat,\star)_2$ we take the identities.

We proceed to construct a pseudomonoid in $(2\x\Cat,\star)_2$. We define a 2-functor 
$$-\Box_2-:Dbl_2\star Dbl_2\to Dbl_2$$ 
by giving two 2-functors $(-,\Aa), (\Bb,-):Dbl_2\to Dbl_2$ 
for any double categories $\Aa,\Bb$ such that $(-,\Aa)_\Bb=(\Bb,-)_\Aa=(\Bb,\Aa)$. We define $(-,\Aa)$ on objects $\Bb$ and 1-cells 
$G:\Bb\to\Bb'$ by $(-,\Aa)_\Bb=\Aa\Box\Bb$ and $(-,\Aa)_G=\Id_\Aa\Box G:\Aa\Box\Bb\to\Aa\Box\Bb'$,  
the funny product of double categories and double functors from \ssref{funny}, and on 2-cells $\beta$ by 
$(-,\Aa)_\beta=\Aa\Box\beta$. The transformation $\Aa\Box\beta$ is given on objects $A\Box B$ by $(\Aa\Box\beta)\vert_{A\Box B}=A\Box\beta(B)$, 
on 1h-cells $A\Box k$ and $K\Box B$ by $(\Aa\Box\beta)\vert_{A\Box k}=A\Box\beta_k$ and $(\Aa\Box\beta)\vert_{K\Box B}=K\Box\beta(B)$, 
and on 1v-cells by $(\Aa\Box\beta)\vert_{A\Box u}=A\Box\beta^u$ and $(\Aa\Box\beta)\vert_{U\Box B}=U\Box\beta(B)$, where both are identity 
2-cells. The 2-functor $(\Bb,-)$ is similarly defined. 
To prove that $\Aa\Box\beta$ is a vertical strict transformation: the axioms concerning $(\Aa\Box\beta)\vert_{A\Box k}$ 
are easily proved because of the axioms for $\beta_k$. However, the first two axioms for $(\Aa\Box\beta)\vert_{K\Box B}$ hold by the first and third axiom from \equref{four strictified}, and the third one holds by \equref{frac-natur}. 
Vertical strictness of $\Aa\Box\beta$ follows by the strictness of $\beta$ and by the first identity in \equref{frac-natur}. 
That $(-,\Aa)$ is a 2-functor is easily verified. 
Thus we obtain a 2-functor $-\Box_2-:Dbl_2\star Dbl_2\to Dbl_2$ on the funny product that we denote the same way as in 
\prref{no funny on Cart} abusing notation. 

\smallskip

Its associativity should be an invertible unnatural transformation 
$$\alpha_2:(-\Box_2-)((-\Box_2-)\star 1)\Rightarrow(-\Box_2-)(1\star(-\Box_2-))\tilde\alpha: (Dbl_2\star Dbl_2)\star Dbl_2\to 
Dbl_2$$
on 2-categories.  
Consider the underlying 2-categorical structures in the part 3. of \coref{strict-rep} and in \deref{hu tr of ff 3}. 
An unnatural transformation on a ternary funny product is given via a ternary unnatural transformation of funny functors.  
As such it is given by a triple of unnatural (2-categorical) transformations, each of which consists merely of 1-cell components in 
$Dbl_2$. The 2-functor $-\Box_2-$ is given on objects the same way as $-\Box-$, then we may set 
$\alpha_2^{\Aa;\Bb}(\Cc)=\alpha_2^{\Bb;\Cc}(\Aa)=\alpha_2^{\Aa;\Cc}(\Bb):=\alpha_{\Aa,\Bb,\Cc}$, the double functor from 
\equref{alfa on funny prod of dcats}. Since $\alpha_{\Aa,\Bb,\Cc}$ satisfies the pentagon strictly, so does $\alpha_2$. 

Let $\I$ denote the trivial double category, the terminal object in $Dbl$. Left 
unity constraint for $-\Box_2-$ is an invertible unnatural transformation $\lambda:(-\Box_2-)(\I\star\Id)\Rightarrow\Id$ given on an object 
$\Aa\in Dbl_2$ by a double functor $\lambda_\Aa:\I\Box\Aa\to\Aa$. 
One defines it by the two double functors: $(-,I):\Aa\to\Aa$ given as $\Id_\Aa$, where $I$ is the unique object of $\I$, 
and $(A,-):\I\to\Aa$ sending $I$ to $A$ and picking identity higher cells on $A$; and the 2-cells $(1^I,K):=\Id_K, (1_I,U):=\Id^U$ 
in $\Aa$ (the third one $(1^I,U)$ expressing identity on $(I,U)$ is trivial) clearly satisfy the desired eleven laws. 
It is clearly an invertible double functor and the right unity $\rho$ is defined analogously. The 2-functor $-\Box_2-$ is clearly unital on both sides with these $\lambda$ and $\rho$. The triangular axiom for $\alpha_2, \lambda, \rho$ holds strictly as in 
the category $(Dbl, \Box)$. We have obtained:

\begin{prop} \prlabel{Del_2 makes pseudomonoid}
The 2-functor $-\Box_2-:Dbl_2\star Dbl_2\to Dbl_2$ makes $Dbl_2$ a pseudomonoid in $(2\x\Cat,\star)_2$, and thus a monoidal 
2-category. 
\end{prop}

Observe that in the above proof that $\Aa\Box\beta, \alpha\Box\beta$ are vertical strict transformations, and similarly in the proof of \prref{no funny on Cart} that so is $\alpha\Box_2\beta$, we used the 2-cells $U\Box k, K\Box u$ and their axioms 
\equref{four strictified} and \equref{frac-natur} from our mixed funny product $\Aa\Box\Bb$. This testifies that it is correct to use horizontal {\em non-purely} unnatural transformations for double categories, that is, those having square-formed 2-cells and satisfying 2 axioms of \deref{hor unnat tr}: these are responsible for having the axioms \equref{four strictified} in 
$\Aa\Box\Bb$. And also, that it is correct to use {\em mixed} funny product, that is, the one 
originating from inner-hom that has unnatural transformations in the horizontal direction, and vertical strict transformations in the vertical direction: these are responsible for the axioms \equref{frac-natur} in $\Aa\Box\Bb$.

\begin{rem}
A similar construction can be carried out 
using the Cartesian monoidal 2-category $(2\x\Cat,\times)_2$, whose 2-cells are 2-transformations.  
Similarly to the above construction, an associativity for the product 2-functor $-\Box_2-:Dbl_2\times Dbl_2\to Dbl_2$ 
should be a 2-transformation 
$$\alpha_2:(-\Box_2-)((-\Box_2-)\times 1)\Rightarrow(-\Box_2-)(1\times(-\Box_2-))\tilde\alpha: (Dbl_2\times Dbl_2)\times Dbl_2\to 
Dbl_2,$$
where now $\tilde\alpha$ is the associativity constraint on the Cartesian product of 2-categories. 
The only data for $\alpha_2$ is the same double functor $\alpha_{\Aa,\Bb,\Cc}$ from 
\equref{alfa on funny prod of dcats}. The rest works the same way as above, and we have a pseudomonoid 
$(Dbl_2,-\Box_2-:Dbl_2\times Dbl_2\to Dbl_2)$ in the Cartesian monoidal 2-category $(2\x\Cat,\times)_2$. 
\end{rem}

Observe that the images of the 2-functors from \prref{no funny on Cart}  and \prref{Del_2 makes pseudomonoid} on objects 
$\Aa,\Bb$ coincide: they are both the (mixed) funny product $\Aa\Box\Bb$ of double categories. A pseudomonoid in the monoidal 
2-category $(Dbl_2, -\Box_2-:Dbl_2\times Dbl_2\to Dbl_2)$ from the above remark turns out to be the same as a pseudomonoid in the monoidal 2-category $(Dbl_2, -\Box_2-:Dbl_2\star Dbl_2\to Dbl_2)$ from \prref{Del_2 makes pseudomonoid}. 

\medskip

Recall our notion of ``$0+6+6+6$ axioms'' from \prref{funny ->premon}.

\begin{thm} \thlabel{funny pseudomonoid}
A pseudomonoid in the monoidal 2-category $(Dbl_2, -\Box_2-)$ is a premonoidal double category with a strict binoidal structure given 
by a binary double funny functor and satisfying the $0+6+6+6$ axioms. 
\end{thm}

\begin{proof}
A pseudomonoid in $(Dbl_2, -\Box_2-)$ consists of a double category $\Dd$, double functors $M:\Dd\Box_2\Dd\to\Dd$ and $U:*\to\Dd$ and 
invertible vertical strict transformations 
$$\alpha^M:M(M\Box_2 1)\Rightarrow M(1\Box_2 M)\alpha_{\Dd,\Dd,\Dd}: (\Dd\Box_2\Dd)\Box_2 \Dd\to \Dd\Box_2(\Dd\Box_2\Dd),$$
$\lambda: I\Box_2\Id_\Dd\Rightarrow\Id_\Dd, \rho: 
\Id_\Dd\Box_2 I\Rightarrow\Id_\Dd$. 
The double functor $M$ is given by a binary double funny functor $H_M:\Dd\times\Dd\to\Dd$. 
That is, by double functors: $(-,A):\Dd\to\Dd$ and 
$(B,-):\Dd\to\Dd$ for $A,B\in\Dd$ and two families of 2-cells $(u,K)$ and $(k,U)$ for 1h-cells $K,k$ and 1v-cells $U,u$ satisfying the eleven axioms of the second part of \deref{funny functor}. 
According to \coref{1v-cells are central}, $H_M$ grants 1v-cells $u$ of $\Dd$ with structures of central cells with centrality structures $(u,-)$ and $(-,u)$. 

By the third part of \coref{strict-rep} the above transformation $\alpha^M$ corresponds to an  
invertible vertical strict transformation on the ternary funny functors $\Dd\times\Dd\times\Dd\to\Dd\Box(\Dd\Box\Dd)$. This in turn 
consists by \deref{vs tr of ff 3} of three invertible vertical strict transformations 
$$\alpha_{-,B,C}^M: ((-, B), C\Rightarrow (-,(B, C))$$
$$\alpha_{A,-,C}^M: ((A, -),C)\Rightarrow (A,(-, C))$$
$$\alpha_{A,B,-}^M: ((A, B), -)\Rightarrow (A,(B,-))$$
of double functors $\Dd\to\Dd$ for every $A,B,C\in\Dd$ satisfying six identities. From what we saw above, 
1v-cell components of $\alpha^M_{A,B,C}, \lambda_A, \rho_A$ are central via $H_M$ (observe that these centrality structures 
coincide with those that the constituent binary double funny functors of $M(M\Box_2 1)$ and $M(1\Box_2 M)\alpha_{\Dd,\Dd,\Dd}$ 
deliver to $\alpha^M_{A,B,C}, \lambda_A, \rho_A$). 
Since they are also invertible, due to \leref{vert comp vert. lx tr.} they are inversely central. 
With these centrality structures, the mentioned 6 identities correspond to the 
$0+0+6+0$ axioms: 
\axiomref{$(u\ltimes,g,C)$}, \axiomref{$(A, v\ltimes,h)$}, \axiomref{$(u\ltimes, B,h)$}, \axiomref{$(f,B,\rtimes z)$}, 
\axiomref{$(f,\rtimes v,C)$}, \axiomref{$(A,g,\rtimes z)$}, as argued in the proof of \prref{funny ->premon}. 
We illustrate this by an example. The first of the 6 identities is 
$\frac{(\theta^{B;C})_f}{(f,v,C)^2}=\frac{(f,v,C)^1}{(\theta^{\tilde B;C})_f}$, where $(f,v,C)^1$ is the structure 2-cell of the binary 
funny functor $M(M\Box_2 1)$ and $(f,v,C)^2$ of $M(1\Box_2 M)\alpha_{\Dd,\Dd,\Dd}$. This corresponds to the axiom \axiomref{$(f,\rtimes v,C)$} 
(in the Appendix A) whereby the 2-cell $(f,v,C)^1$ becomes $(f\w\ltimes\w-\vert_v,C)$, and $(f,v,C)^2$ becomes 
$f\w\ltimes\w-\vert_{v\ltimes C}$ 
of the axiom \axiomref{$(f,\rtimes v,C)$}. 

On the other hand, the centrality structures obtained from $H_M$ are given via  
square-formed 2-cell components satisfying $(-,K)\vert_u=(u,-)\vert_K$ and $(-,U)\vert_k=(k,-)\vert_U$, so that the latter 
$0+0+6+0$ axioms simultaneously cover the $0+6+0+0$ axioms: 
\axiomref{$(f\ltimes,v,C)$}, \axiomref{$(f\ltimes,B,z)$}, \axiomref{$(A,g\ltimes,z)$}, \axiomref{$(u,\rtimes g,C)$}, 
\axiomref{$(u,B,\rtimes h)$}, \axiom{$(A,v,\rtimes h)$}. (Observe that in the named axioms the 1h-cells $f,g,h$ should be left/right central: in the present case $K,k$ are not central at any side, but the square-formed 2-cells relevant for the axioms exist.) Moreover, since the structure 2-cells $(u,U)$ of the binary funny functors are trivial, the $0+0+0+6$ axioms 
\axiomref{$(u\ltimes,v,C)$}, \axiomref{$(A, v\ltimes,z)$}, \axiomref{$(u\ltimes, B,z)$}, \axiomref{$(u,B,\rtimes z)$}, 
\axiomref{$(u,\rtimes v,C)$}, \axiomref{$(A,v,\rtimes z)$} are automatically fulfilled. 

That the transformation $\alpha^M$ obeys the pentagon it means that its above three components satisfy the four pentagons of \deref{assoc}, 
whereby the binoidal structure is induced by the underlying binary funny functor $H_M$ of the double functor $M:\Dd\Box\Dd\to\Dd$. Similarly, one has that $\lambda, \rho, \alpha^M$ satisfy the six triangles from \deref{premon}.  
\qed\end{proof}

\section{Double quasi-functors yield purely central binoidal double categories} \selabel{cubical1} 

In \cite{Fem:Gray1} and \cite{Fem:Bif} we generalized Gray's 2-categorical ``quasi-functors of two variables'' from \cite{Gray} 
to the double-categorical setting. In the former reference they emerged from the inner-hom studied in \cite{Gabi} whose 0-cells were double functors, while in the latter they came from a candidate for inner-hom whose 0-cells are {\em lax} double functors. 
In this article we are interested in the ``pseudo'' version of double quasi-functors. 
In \cite{GPS}, pseudodouble quasi-functors of two and more variables are called {\em cubical functors}.

In this section we investigate the relation between a binoidal structure $-\bowtie-$ and centrality of 1- and 2-cells in a double category 
$\Bb$ from \ssref{binoidal}, on one hand, and a pseudodouble quasi-functor 
structure $H$ in $\Bb$, on the other hand. 

\subsection{Pseudodouble quasi-functors and pure centrality} \sslabel{quasi}

In \cite[Section 2]{Fem:Bif} we introduced the double category $\Lax_{hop}(\Aa, \Bb)$ of lax double functors of double categories 
$\Aa\to\Bb$, horizontal oplax transformations as 1h-cells, vertical lax transformations as 1v-cells, and modifications. (Such a choice of transformations is a priori arbitrary and most general. We also explored the impact of having the transformations in the two directions of different kind. In particular, by this choice in the horizontal direction we generalized the 2-categorical constructions from 
\cite{FMS}.) In \cite[Proposition 3.3]{Fem:Bif} we characterized a lax double 
functor $\F\colon\Aa\to\Lax_{hop}(\Bb, \Cc)$ as a pair of two families of lax double functors into $\Cc$ together with 
four families of 2-cells in $\Cc$ that satisfy 20 axioms. 
The latter collection of data and axioms we call a {\em lax double quasi-functor}, 
in analogy to \cite{Gray}. For reader's convenience we include this result as \prref{char df} 
in the Appendix B. In Table 1 of \cite[Proposition 3.3]{Fem:Bif} we listed the origin of the four types of 2-cells and 
20 axioms and we labeled the axioms. Various of those axioms can be interpreted in more than one way. We include 
that Table enriched with the corresponding additional interpretations as Table \ref{table:12} in the Appendix C.

\medskip

Inspecting that Table and realizing that 
the axioms \axiomref{($k,K'K$)}, \axiomref{($k,1_A$)}, \axiomref{($k,\frac{U}{U'}$)}, \axiomref{($k,1^A$)}, \axiomref{($(k,K)$-r-nat)} 
therein mean that $(k,-)$ is a horizontal lax transformation, 
while the axioms \axiomref{($u,K'K$)}, \axiomref{($u,1_A$)}, \axiomref{($u,\frac{U}{U'}$)}, \axiomref{($u,1^A$)}, \axiomref{($(u,U)$-r-nat)} mean that $(u,-)$ is a vertical oplax transformation, and finally that the axioms \axiomref{($(k,K)$-r-nat)} and \axiomref{($(u,U)$-r-nat)} 
mean that $(\omega,-)$ is a modification (with respect to horizontally lax and vertically oplax transformations), 
one sees that the content of Corollary 3.5 of \cite{Fem:Bif} can be upgraded into an if and only if statement. 
Namely, we have: 

\begin{prop} \prlabel{quasi-fun}
Let $\Aa,\Bb,\Cc$ be double categories. The following are equivalent:
\begin{enumerate} 
\item $H\colon \Aa\times\Bb\to\Cc$ is a lax double quasi-functor, meaning that there are two families of lax double functors 
$(-,A)\colon\Bb\to\Cc\quad\text{ and}\quad (B,-)\colon\Aa\to\Cc$ for objects $A\in\Aa, B\in\Bb$, 
such that $H(A,-)=(-, A), H(-, B)=(B,-)$ and $(-,A)\vert_B=(B,-)\vert_A=(B,A)$, and 
there are four families of 2-cells 
\begin{equation} \eqlabel{four cells}
\scalebox{0.78}{
\bfig
 \putmorphism(-150,50)(1,0)[(B,A)`(B', A)`(k, A)]{600}1a
 \putmorphism(450,50)(1,0)[\phantom{A\ot B}`(B', A') `(B', K)]{680}1a
\putmorphism(-180,50)(0,-1)[\phantom{Y_2}``=]{450}1r
\putmorphism(1100,50)(0,-1)[\phantom{Y_2}``=]{450}1r
\put(330,-190){\fbox{$(k,K)$}}
 \putmorphism(-150,-400)(1,0)[(B,A)`(B,A')`(B,K)]{600}1a
 \putmorphism(450,-400)(1,0)[\phantom{A\ot B}`(B', A') `(k, A')]{680}1a
\efig}
\end{equation}

$$
\scalebox{0.78}{
\bfig
\putmorphism(-150,50)(1,0)[(B,A)`(B,A')`(B,K)]{600}1a
\putmorphism(-150,-400)(1,0)[(\tilde B, A)`(\tilde B,A') `(\tilde B,K)]{640}1a
\putmorphism(-180,50)(0,-1)[\phantom{Y_2}``(u,A)]{450}1l
\putmorphism(450,50)(0,-1)[\phantom{Y_2}``(u,A')]{450}1r
\put(-20,-180){\fbox{$(u, K)$}}
\efig}
\quad
\scalebox{0.78}{
\bfig
\putmorphism(-150,50)(1,0)[(B,A)`(B',A)`(k,A)]{600}1a
\putmorphism(-150,-400)(1,0)[(B, \tilde A)`(B', \tilde A) `(k,\tilde A)]{640}1a
\putmorphism(-180,50)(0,-1)[\phantom{Y_2}``(B,U)]{450}1l
\putmorphism(450,50)(0,-1)[\phantom{Y_2}``(B',U)]{450}1r
\put(0,-180){\fbox{$(k,U)$}}
\efig}
$$

$$
\scalebox{0.78}{
\bfig
 \putmorphism(-150,500)(1,0)[(B,A)`(B,A) `=]{600}1a
\putmorphism(-180,500)(0,-1)[\phantom{Y_2}`(B, \tilde A) `(B,U)]{450}1l
\put(-20,50){\fbox{$(u,U)$}}
\putmorphism(-150,-400)(1,0)[(\tilde B, \tilde A)`(\tilde B, \tilde A) `=]{640}1a
\putmorphism(-180,50)(0,-1)[\phantom{Y_2}``(u,\tilde A)]{450}1l
\putmorphism(450,50)(0,-1)[\phantom{Y_2}``(\tilde B, U)]{450}1r
\putmorphism(450,500)(0,-1)[\phantom{Y_2}`(\tilde B, A) `(u,A)]{450}1r
\efig}
$$ 
in $\Cc$ determined by all 1h-cells $K\colon A\to A'$ and 1v-cells $U\colon A\to\tilde A$ in $\Aa$, and 1h-cells $k\colon B\to B'$ 
and 1v-cells $u\colon B\to\tilde B$ in $\Bb$, which satisfy 20 axioms from \prref{char df}, and
\item there are two families of lax double functors 
$(-,A)\colon\Bb\to\Cc\quad\text{ and}\quad (B,-)\colon\Aa\to\Cc$ for objects $A\in\Aa, B\in\Bb$, 
such that $(-,A)\vert_B=(B,-)\vert_A=(B,A)$, and the following hold: 
\begin{enumerate}[(i)]
\item $(-,K)\colon (-,A)\to(-,A')$ is a horizontal oplax transformation for each 1h-cell $K\colon A\to A'$, 
$(-,U)\colon (-,A)\to(-,\tilde A)$ is a vertical lax transformation for each 1v-cell $U\colon A\to\tilde A$ in $\Aa$,  
$(-,\zeta)$ is a modification with respect to horizontally oplax and vertically lax transformations for each 2-cell $\zeta$ in $\Aa$, 
and the following coincide:
$$(B,-)\vert_K=(-,K)\vert_B, \,\, (B,-)\vert_U=(-,U)\vert_B, \,\, (B,-)\vert_\zeta=(-,\zeta)\vert_B;$$
\item $(k,-)\colon (B,-)\to (B', -)$ is a horizontal lax transformation for each 1h-cell $k\colon B\to B'$, 
$(u,-)\colon (B,-)\to(\tilde B,-)$ is a vertical oplax transformation for each 1v-cell $u\colon B\to\tilde B$ in $\Bb$, 
$(\omega,-)$ is a modification with respect to horizontally lax and vertically oplax transformations for each 2-cell $\omega$ in $\Bb$, 
and the following coincide:
$$(-,A)\vert_k=(k,-)\vert_A, \,\, (-,A)\vert_u=(u,-)\vert_A, \,\, (-,A)\vert_\omega=(\omega,-)\vert_A;$$
\item for 1h-cells $K,k$ the 2-cell component $(-,K)\vert_k$ of the oplax (resp. lax) structure of the horizontal transformation $(-,K)$ coincides with the 2-cell component $(k,-)\vert_K$ of the lax (resp. oplax) structure 
of the transformation $(k,-)$; 
\item for 1v-cells $U,u$ the 2-cell component $(-,U)\vert_u$ of the lax (resp. oplax) structure of the vertical transformation $(-,U)$ coincides with the 2-cell component $(u,-)\vert_U$ of the oplax (resp. lax) structure of the transformation $(u,-)$;  
\item for 1h-cells $K,k$ and 1v-cells $U,u$ the following 2-cell components of the respective transformations coincide:
$(-,K)\vert_u=(u,-)\vert_K$ and $(-,U)\vert_k=(k,-)\vert_U$.
\end{enumerate}
\end{enumerate}
\end{prop}

The modifications in item 2. (ii) above are from \deref{modif-hv}. 

\medskip

Similarly to \cite[Proposition 2.9]{Fem:Fil} and \cite[Proposition 3.3]{Fem:Bif}, one has that a {\em pseudodouble quasi-functor} $H:\Aa\times\Bb\to\Cc$ consists of two families of {\em pseudodouble functors}  
$(-,A)\colon\Bb\to\Cc$ and $(B,-)\colon\Aa\to\Cc$ with similar 20 axioms as in the cited two propositions. 

\medskip

We go back for a moment to the {\em lax} case. 
Assume now that $\Aa=\Bb=\Cc$. 
Then to have a lax double quasi-functor $H:\Bb\times\Bb\to\Bb$ it means, among other, to have 2-cells 
$(k,K)$ and $(u,U)$ in $\Bb$ from \equref{four cells}, so that $(-,K)\colon (-,A)\to(-,A')$ is a horizontal oplax transformation, $(k,-)\colon (B,-)\to (B', -)$ is a horizontal lax transformation, $(u,-)\colon (B,-)\to(\tilde B,-)$ is a vertical oplax transformation, and 
$(-,U)\colon (-,A)\to(-,\tilde A)$ is a vertical lax transformation, over lax double functors acting $\Bb\to\Bb$. Then clearly, 
if $(-,K)$ is a horizontal {\em pseudo}natural transformation for every 1h-cell $K$, necessarily so is $(k,-)$  for every 1h-cell $k$, and if 
$(-,U)$ is a vertical {\em pseudo}natural transformation for every 1v-cell $U$, necessarily so is $(u,-)$ for every 1v-cell $u$, 
as their respective component 2-cells coincide.  

\medskip
 
Let $\Pseudo_{ps}(\Bb,\Bb)$ be the double category of pseudodouble endofunctors on $\Bb$, horizontal pseudonatural transformations, 
vertical pseudonatural transformations  and modifications. (This is the fully pseudo version of the double category $\Lax_{hop}(\Bb,\Bb)$ 
from \cite[Section 2]{Fem:Bif}.) 
In view of the above said, and joining the version of \prref{char df} for pseudodouble functors $\Bb\to\Pseudo_{ps}(\Bb, \Bb)$, 
we obtain:

\begin{thm} \thlabel{lr central}
Let $\Bb$ be a double category. The following are equivalent:
\begin{enumerate}
\item there is a pseudodouble functor $\F\colon\Bb\to\Pseudo_{ps}(\Bb, \Bb)$; 
\item there is a pseudodouble quasi-functor $H:\Bb\times\Bb\to\Bb$ (with families of pseudodouble functors 
$(-,A),(B,-)\colon\Bb\to\Bb$ for $A,B\in\Bb$, four families of 2-cells satisfying 20 axioms); 
\item the following hold:
\begin{itemize}
\item $\Bb$ is binoidal with pseudodouble functors $A\ltimes -,-\rtimes B:\Bb\to\Bb$ for $A,B\in\Bb$; 
\item every 1h-cell $K$ in $\Bb$ is left central via a horizontal pseudonatural transformation $K\ltimes-$ 
and every 1h-cell $k$ is right central in $\Bb$ via a horizontal pseudonatural transformation $-\rtimes k$, 
and it is \vspace{-0,2cm} 
$$A\ltimes-\vert_k=-\rtimes k\vert_A \quad\text{and}\quad -\rtimes B\vert_K=K\ltimes-\vert_B;$$
\item every 1v-cell $U$ is left central in $\Bb$ via a vertical pseudonatural transformation $U\ltimes-$ 
and every 1v-cell $u$ is right central in $\Bb$ via a vertical pseudonatural transformation $-\rtimes u$, and 
it is \vspace{-0,2cm}
$$A\ltimes-\vert_u=-\rtimes u\vert_A \quad\text{and}\quad -\rtimes B\vert_U=U\ltimes-\vert_B;$$
\item every 2-cell $\zeta$ is left central via a modification $\zeta\ltimes-$ and every 2-cell $\omega$ is right central 
via a modification $-\rtimes\omega$, and it is \vspace{-0,2cm}
$$A\ltimes-\vert_\omega=-\rtimes\omega\vert_A \quad\text{and}\quad -\rtimes B\vert_\zeta=\zeta\ltimes-\vert_B;$$
\item it is 
$$K\ltimes-\vert_k=(-\rtimes k\vert_K)^{-1}, \qquad  U\ltimes-\vert_u=(-\rtimes u\vert_U)^{-1},$$
$$K\ltimes-\vert_u=-\rtimes u\vert_K \quad\text{and}\quad U\ltimes-\vert_k=-\rtimes k\vert_U$$ 
for all 1h-cells $K,k$ and 1v-cells $U,u$ in $\Bb$.  
\end{itemize}
In particular, with notations as above, one has: 
\begin{itemize}
\item the 2-cell component $K\ltimes-\vert_k$ of the oplax (resp. lax) structure of the horizontal transformation $K\ltimes-$ coincides with the 2-cell component at $K$ of the lax (resp. oplax) structure of the transformation $-\rtimes k$ 
(which is $(-\rtimes k\vert_K)^{-1}$ in \equref{f lr}); 
\item the 2-cell component $U\ltimes-\vert_u$ of the lax (resp. oplax) structure of the vertical transformation $U\ltimes-$ coincides with the 2-cell component at $U$ of the oplax (resp. lax) structure of the transformation $-\rtimes u$ 
(which is $(-\rtimes u\vert_U)^{-1}$ in \equref{v lr}). 
\end{itemize}
\end{enumerate}
\end{thm}

\medskip

We outline the correspondence holding in the above situation. For all $A,B\in\Bb$ it is $A\ltimes -=(-,A)$ and $-\rtimes B=(B,-)$, and moreover 
$K\ltimes-=(-,K)$ and $-\rtimes k=(k,-)$, $U\ltimes-=(-,U)$ and $-\rtimes u=(u,-)$, $\zeta\ltimes-=(-,\zeta)$ and $-\rtimes \omega=(\omega,-)$
for 1h-cells $K,k$, 1v-cells $U,u$ and 2-cells $\zeta, \omega$ in $\Bb$. 
Similarly, the four families of 2-cells in a quasi-functor correspond to the following four identities holding between centrality structural 
2-cells in the point 3. of the above theorem: 
\begin{equation} \eqlabel{4 eqs pc}
K\ltimes-\vert_k=(-\rtimes k\vert_K)^{-1}, \qquad \qquad U\ltimes-\vert_u=(-\rtimes u\vert_U)^{-1},
\end{equation}  
$$K\ltimes-\vert_u=-\rtimes u\vert_K, \qquad\qquad  U\ltimes-\vert_k=-\rtimes k\vert_U$$ 
for all 1h-cells $K,k$ and all 1v-cells $U,u$ in $\Bb$. 
Observe that in point 3. we have that every 1h-cell, every 1v-cell and every 2-cell in $\Bb$ are both left and right central, hence central, in their respective senses, with a specific structural transform. 
From point 1. we know that the assignment of those centrality structure transforms is functorial. 

\begin{defn} \delabel{pc binoidal}
A binoidal double category $\Bb$ satisfying the point 3. in the above theorem we call {\em purely central}. 
\end{defn}

\medskip


\begin{defn} \delabel{pure center}
For a purely central binoidal double category $\Bb$ a {\em pure center} double category is a pseudodouble category $\Zz_p(\Bb)$ 
that has the same objects as $\Bb$, its 1h-cells are triples 
$(f, f\ltimes-, -\rtimes f)$ consisting of all (central) 1h-cells of $\Bb$ and their two structural horizontal transformations, its 1v-cells are triples $(v, v\ltimes-, -\rtimes v)$ consisting of all (central) 1v-cells of $\Bb$ and their two structural vertical transformations, its 
2-cells are triples $(a, a\ltimes-, -\rtimes a)$ consisting of all (central) 2-cells of $\Bb$ and their two structural modifications. 
\end{defn}

It is easily seen that $\Zz_p(\Bb)$ is indeed a pseudodouble category. 
The composition of 1h-cells say $(f, f\ltimes-, -\rtimes f):A\to A'$ and $(f', f'\ltimes-, -\rtimes f'):A'\to A''$ is defined as follows: 
to the expression $\delta_{\frac{\alpha}{\beta},f}$ in \leref{vert comp hor.ps.tr.} add the 2-cell components of the 
pseudodouble functor structure of $A\ltimes -$ resp. $-\rtimes A$ (conjugate by the lax and the oplax part). Then for the associativity one applies the ``hexagon axiom'' \axiomref{lx.f.cmp}, which leads to the fact that the center double categories are indeed {\em pseudodouble} categories. Compositions of 2-cells are induced by the compositions of modifications, see Appendix A.0. 

\smallskip

We also outline that the functorial assignment of the centrality structures in a pure center double category is 
provided via the pseudodouble functors $A\ltimes -,-\rtimes B:\Bb\to\Bb$. From the perspective of point 1. in the above theorem these two 
functors are stemming from $\F(A)=A\ltimes -:\Bb\to\Bb$ and $\F(-)(B)=-\rtimes B:\Bb\to\Bb$.

\smallskip

Consequently, for $\Bb$ purely central there is a double functor $Z_p:\Bb\to\Zz_p(\Bb)$, which is identity on objects and given by 
$b\mapsto (b,b\ltimes-, -\rtimes b)$ for any 1h-, 1v- or 2-cell $b$ in $\Bb$. That it is a double functor is easily seen.   

\smallskip

Strictly speaking, a pure center double category should be written as a pair $(\Zz_p(\Bb), \F)$ (or $(\Zz_p(\Bb), Z_p)$), as its {\em functorial assignment} of centrality structures depends on the pseudodouble functor $\F$ (hence $Z_p$). (In this regard the adjective ``pure'' is 
probably not the best suiting, as observed in \cite{Phil}.) We will suppress writing the determining pseudodouble functor $\F$ (equivalently pseudodouble quasi-functor $H$), or $Z_p$, but will bear in mind its existence.

\subsection{One-sided centers of a binoidal double category}

Apart from pure centrality and pure center for a binoidal double category one can give two more equivalent interpretations to a 
pseudodouble quasi-functor and consider two other center double categories. Namely, left and right centrality structures on $\Bb$ related between 
each other as in a pure center can be seen as one-sided centrality structures on $\Bb$, by disregarding the structures on the other side. 
Accordingly, one-sided center double categories can be introduced. 

For all $A,B\in\Bb$ we keep a binoidal structure detrmined by $A\ltimes -=(-,A)$ and $-\rtimes B=(B,-)$. Then, for left centrality and left center we use: $K\ltimes-=(-,K)$, $U\ltimes-=(-,U)$ and $\zeta\ltimes-=(-,\zeta)$ for 1h-cells $K$, 1v-cells $U$ and 2-cells $\zeta$. 
However, 
whereas the four families of 2-cells of a quasi-functor were recognized in a purely central binoidal double category as the 2-cells 
\equref{4 eqs pc}, we instead consider: 
$$(k,K)=K\ltimes-\vert_k,\qquad  (u,K)=K\ltimes-\vert_u, \qquad (k,U)=U\ltimes-\vert_k, \qquad (u,U)=U\ltimes -\vert_u$$ 
and recognize that the 20 axioms mean (compare to Table \ref{table:12} of Appendix C): 
\begin{itemize}
\item every 1h-cell $K$ in $\Bb$ is left central via a horizontal pseudonatural transformation $K\ltimes-$ and it is 
$K\ltimes-\vert_B=-\rtimes B\vert_K$; 
\item every 1v-cell $U$ is left central in $\Bb$ via a vertical pseudonatural transformation $U\ltimes-$ and it is 
$U\ltimes-\vert_B=-\rtimes B\vert_U;$ 
\item every 2-cell $\zeta$ is left central via a modification $\zeta\ltimes-$  and it is $\zeta\ltimes-\vert_B=-\rtimes B\vert_\zeta;$ 
\item the axioms \axiomref{($k,K'K$)} and \axiomref{($u,K'K$)} for the invertible modification compositor 
$(K'\ltimes -)(K\ltimes -)\Rrightarrow K'K\ltimes-$ 
are saying how a horizontal pseudonatural transformation $K'K\ltimes-$, {\em i.e.} its structural 2-cells, for left centrality of the composition 1h-cell $K'K$ is given; 
\item the axioms \axiomref{($k,1_A$)} and \axiomref{($u,1_A$)} for the invertible modification unitor 
$\Id_{A\ltimes-}\Rrightarrow 1_A\ltimes-$  are saying how a horizontal pseudonatural transformation $1_A\ltimes-$, {\em i.e.} 
its structural 2-cells, for left centrality of the identity 1h-cell $1_A$ is given; 
\item the axioms \axiomref{($k,1^A$)} and \axiomref{($u,1^A$)}, stemming from the axiom \axiomref{lx.f.v2} of $\F$ when evaluated at $k$ and $u$, 
say that the structural 2-cells $1^A\ltimes-\vert_k$ and $1^A\ltimes-\vert_u$ for 1h-cells $k$ and 1v-cells $u$ for 
the vertical pseudonatural transformation $1^A\ltimes-$ are identities: $\Id_{A\ltimes k}$ and $\Id_{A\ltimes u}$, and 
\item the axioms \axiomref{($k,\frac{U}{U'}$)} and \axiomref{($u,\frac{U}{U'}$)}, stemming from the axiom \axiomref{lx.f.v1} of $\F$ 
when evaluated at $k$ and $u$, say that the structural 2-cells $\frac{U}{U'}\ltimes-\vert_k$ and $\frac{U}{U'}\ltimes-\vert_u$ for 1h-cells $k$ 
and 1v-cells $u$ for the vertical pseudonatural transformation $\frac{U}{U'}\ltimes-$ are given by 
$\frac{U}{U'}\ltimes-\vert_k=\frac{U\ltimes-\vert_k}{U'\ltimes-\vert_k}$ and 
$\frac{U}{U'}\ltimes-\vert_u=\frac{U\ltimes-\vert_u}{U'\ltimes-\vert_u}$. 
\end{itemize}

The above consideration means that there is a {\em functorial way} of assigning left centrality structural transforms to 1h- and 1v-cells of $\Bb$: the last four items express compositionality and unitality in the variable $(-)$ of the pseudodouble functor $\F(-)(\bullet)$ with image in 
$\Bb$, whereby $\F\colon\Bb\to\Pseudo_{ps}(\Bb, \Bb)$. The above consideration has its right-hand sided version.  

\begin{defn} \delabel{lr-central binoidal}
A binoidal double category $\Bb$ equipped with left centrality structures as above we call {\em left central}. \\ 
Similarly, we define {\em right central binoidal} double category. 
\end{defn}

We can now extend \thref{lr central} into: 

\begin{thm} \thlabel{summary} 
For a double category $\Bb$ the following are equivalent: 
\begin{enumerate}
\item there is a pseudodouble functor $\F\colon\Bb\to\Pseudo_{ps}(\Bb, \Bb)$; 
\item there is a pseudodouble quasi-functor $H:\Bb\times\Bb\to\Bb$; 
\item $\Bb$ is binoidal and purely central;   
\item $\Bb$ is binoidal and left central;   
\item $\Bb$ is binoidal and right central.   
\end{enumerate}
\end{thm}

\begin{rem}
The structure 2-cells and axioms in the items 4. and 5. both hold in the item 3. 
The thing is that the 8 axioms in item 4. can additionally be interpreted (see the Appendix C) 
as the axioms (h.l.t.-1) - (h.l.t.-4) of $\{-, k\}$ and \axiomref{v.o.t.\x 1} - \axiomref{v.o.t.\x 4} of $\{-, u\}$ 
(while (h.l.t.-5) of $\{-, k\}$ is contained in the modification axiom \axiomref{$(k,K)$-r-nat} for $\{\zeta, -\}$, 
and \axiomref{v.o.t.\x 5} is contained in the modification axiom \axiomref{$(u,U)$-r-nat} for $\{\zeta, -\}$). So, in the item 3. these 8 axioms  
have the meaning that $\{-, k\}$ is a horizontal pseudonatural transformation and $\{-, u\}$ is a vertical pseudonatural 
transformation defined so that $\{K,-\}_k=(\{-, k\}_K)^{-1}$, \, $\{K,-\}_u=\{-, u\}_K$ and  
$\{U,-\}_k=\{-, k\}_U$, leading to pure centrality. 
The analogous argument holds for item 5. 
\end{rem}

We may finally define left and right center double categories. 

\begin{defn} \delabel{lr center}
For a left central binoidal double category $\Bb$ a {\em left center} double category is a pseudodouble category $\Zz_l(\Bb)$ 
that has the same objects as $\Bb$, its 1h-cells are pairs  
$(f, f\ltimes-)$ consisting of all (left central) 1h-cells of $\Bb$ and their structural horizontal transformations, its 1v-cells are pairs 
$(v, v\ltimes-)$ consisting of all (left central) 1v-cells of $\Bb$ and their structural vertical transformations, its 
2-cells are pairs $(a, a\ltimes-)$ consisting of all (left central) 2-cells of $\Bb$ and their structural modifications. 

For a right central binoidal double category $\Bb$ we define similarly a {\em right center} double category $\Zz_r(\Bb)$.  
\end{defn}

If not all higher cells of $\Bb$ have their centrality structures, {\em i.e.} if one does not count with a pseudodouble functor 
$\Bb\to\Pseudo_{ps}(\Bb, \Bb)$, one can speak of a pure/left/right center double category $\Zz_\bullet(\Bb')$ for 
a double subcategory $\Bb'\subseteq\Bb$, but still one needs a pseudodouble functor $\Bb'\to\Pseudo_{ps}(\Bb', \Bb')$. 

\medskip

Given a pseudodouble functor $\F\colon\Bb\to\Pseudo_{ps}(\Bb, \Bb)$, it 
determines center double categories $\Zz_l(\Bb), \Zz_r(\Bb), \Zz_p(\Bb)$. Then there are pseudodouble functors 
\begin{equation} \eqlabel{center-pseud}
L: \Zz_l(\Bb)\to\Pseudo_{ps}(\Bb, \Bb)\quad\text{and}\quad R: \Zz_r(\Bb)\to\Pseudo_{ps}(\Bb, \Bb)
\end{equation}
(also $L_p, R_p: \Zz_p(\Bb)\to\Pseudo_{ps}(\Bb, \Bb)$) 
defined in the obvious way. (Observe that $\Zz_l(\Bb)\iso\Zz_r(\Bb)$, and that  
$L$ and $R$ are related so that $\Bb$ is purely central.) These pseudodouble functors resemble  
the known interpretation of the Drinfel`d center category of a monoidal category, or the center category of a 2-category 
from \cite{Meir} (see also Theorems 3.6 and 3.7 of \cite{FH}). Namely, the latter can be seen as $\Z_{Dr}(\B)\iso\PsNat(\Id_\B, \Id_\B)$, 
having for objects pseudonatural transformations and for morphisms modifications of the identity 2-functor. 
We obtain this in a double categorical setting by considering the endo-hom category $\HH(\Zz_p(\Bb))(I,I)$ of $\HH(\Zz_p(\Bb))$.

\subsection{Another approach to center double categories} \sslabel{not diff}

In \cite[Definition 21]{HF} a bicategory of pure maps $\C_p(\B)$ was defined for a premonoidal bicategory $\B$. 
One can consider to take a similar approach to define what we call a pure center double category for a premonoidal double category $\Bb$. 
One can even do it if $\Bb$ is only binoidal. 
Let us analyze this approach.

\begin{defn} \delabel{lr-binoidal}
Let $(\Bb, \ltimes,\rtimes)$ be a binoidal double category and $\Bb'\subseteq\Bb$ a double subcategory. Consider pseudodouble functors 
\vspace{-0,2cm} 
\begin{equation} \eqlabel{L0,R0}
L_0,R_0:\Bb'\to\Pseudo(\Bb',\Bb')
\end{equation}
given on 0-, 1h-, 1v- and 2-cells by  
$$L_0(A)=A\ltimes-, \, L_0(f)=f\ltimes -, \, L_0(u)=u\ltimes-\quad\text{and}\quad L_0(a)=a\ltimes-,$$ 
respectively, and similarly for $R_0$, whereby \vspace{-0,2cm} 
$$L_0(a)_B=-\rtimes B\vert_a, \qquad R_0(a)_B=B\ltimes -\vert_a \vspace{-0,2cm} $$
for objects $A,B\in\Bb'$ and all cells $a,b\in\Bb'$. 




When such pseudodouble functors $L_0, R_0$ exist with no further compatibilities, 
we say that $(\Bb',\ltimes,\rtimes, L_0, R_0)$ is a {\em central} binoidal double category.

When $L_0, R_0$ exist so that the four identities \equref{4 eqs pc} hold,  
we say that $(\Bb',\ltimes,\rtimes, L_0, R_0)$ is a {\em purely central} binoidal double category. 
\end{defn}

\begin{rem} \rmlabel{discuss L0,R0}
In \thref{summary} it is $\F(a)(b)=H(a,b)=a\ltimes-\vert_b=-\rtimes b\vert_a$ for all cells $a,b\in\Bb$. 
Accordingly, both $L_0$ and $R_0$ determine a quasi-functor $(H_{L_0}$ and $H_{R_0}$, respectively) 
and make $\Bb'$ purely central, left or right central, according to a perspective one takes. However, the agreement with the given binoidal structure $A\ltimes -=L_0(A)$ and $-\rtimes A=R_0(A)$ in \deref{lr-binoidal} means that a choice is fixed 
$A\ltimes -=H_{L_0}(A,-)$ and $-\rtimes A=H_{R_0}(A,-)$. Thus the latter identity equips $(\Bb',\ltimes,\rtimes, L_0, R_0)$ with a right central binoidal structure in a way opposed to the general \thref{summary}: 
whereas the relation of $R_0$ (and $H_{R_0}$) with the binoidal structure is $-\rtimes b=R_0(b)=H_{R_0}(b,-)$, in \deref{lr-central binoidal} it is $-\rtimes b=H(-,b)=\F(-)(b)$. We will keep this in mind in regard to the pseudodouble functor $R_0$. 

Regarding relations between $L_0$ and $R_0$, by \deref{lr-binoidal} one has $H_{L_0}(a,-)_B=H_{R_0}(B,-)_a$, whereas by the quasi-functor property for $H_{L_0}$ it is $H_{L_0}(a,-)_B=H_{L_0}(-,B)_a$, hence one obtains $H_{L_0}(-,B)_a=H_{R_0}(B,-)_a$. Similarly, it is 
$H_{L_0}(-,a)_B=H_{R_0}(a,-)_B$. If $L_0$ is known, the only value of $R_0$ that is not known is $H_{R_0}(f,g)$, for 1-cells $f,g$, unless 
$L_0, R_0$ make $\Bb'$ purely central. In that case one has the four identities \equref{4 eqs pc}: 
$H_{L_0}(f,g)=H_{R_0}(g,f)^{-1}, \,\, H_{L_0}(U,u)=H_{R_0}(u,U)^{-1}, \,\, H_{L_0}(K,u)=H_{R_0}(u,K), \,\, H_{L_0}(U,k)=H_{R_0}(k,U)$ 
with the usual notations of cells.  
\end{rem}

\smallskip

One can now define the pseudodouble category pure center 
$(\Zz_p(\Bb'), \ltimes,\rtimes, L_0, R_0)$, where $\Zz_p(\Bb')$ is formally as in \deref{pure center}, whereby the notion of a pure binoidal double category is as in \deref{lr-binoidal}. Observe that in view of \thref{summary} the latter notion  
of left binoidal structure is equivalent to the notion of left binoidal structure from \deref{lr-central binoidal}, whereas 
the notion of right and pure binoidal structure is ``anti-equivalent'' to that of \deref{lr-central binoidal}, in the sense discussed in the above remark. Thus the corresponding definitions of a pure center double category differ up to this subtlety.

\medskip

In \cite[Definition 7]{HF1} the authors study also the center bicategory $\Z_0(\B)$ of a premonoidal bicategory $\B$. Its 1- and 2-cells are central 1- and 2-cells of $\B$ without specifying their centrality structure. We will consider a center double category with specified centrality structures, so that it is a tuple $(\Zz(\Bb), \ltimes,\rtimes, L_0, R_0)$ where $(\Bb,\ltimes,\rtimes, 
L_0, R_0)$ is a central binoidal double category and the cells of $\Zz(\Bb)$ coincide in form with the cells of $\Zz_p(\Bb)$. We will denote by $\Z(\B)$ the bicategorical version of $\Zz(\Bb)$: it is given via pseudofunctors determining centrality structures on cells in $\B$, while $\Z_0(\B)$ does not have any functorial choice of centrality structures. 
We will further study the center double category $(\Zz(\Bb), \ltimes,\rtimes, L_0, R_0)$ in \seref{center premon}.

\section{Monoidality of purely central premonoidal double categories}   \selabel{cubical2}  

The presence of the associativity for a binoidal structure and of a functorial choice for centrality structures in a purely central premonoidal double category $\Dd$ turn out to be the just pieces of data necessary for constructing a monoidal double category structure on 
$\Dd$. We explore this fact in this section. We do it in two forms: for binoidal structures coming from pseudodouble quasi-functors, and in the form of purely central $n$-noidal structures that we introduce. Our findings are summarized in \thref{gen-1}.

\subsection{Towards monoidality of the pure center} \sslabel{towards}

In \cite{Fem:Bif} we proved a Bifunctor Theorem for lax double functors, which becomes of interest to us in this work. 
Namely, applying \cite[Proposition 5.6]{Fem:Bif} to the present setting, we have that when all the 2-cells $(u,U)$ of a pseudodouble quasi-functor $H:\Bb\times\Bb\to\Bb$ for all 1v-cells $U,u$ in $\Bb$ are identities, then $H$ induces a pseudodouble functor $P:\Bb\times\Bb\to\Bb$. 
In this case we have a candidate for a monoidal product on the double category $\Bb$. In fact,
there is a double category equivalence 
\begin{equation} \eqlabel{equiv from Fem-Gray}
\F\colon q\x\Ps_{hop}^{st}(\Bb\times\Bb,\Bb) \to \Ps_{hop}(\Bb\times\Bb,\Bb)
\end{equation}
where in the left-hand side is the double category consisting of pseudodouble quasi-functors in which the 2-cells $(u,U)$ are trivial, horizontal oplax transformations as 1h-cells, vertical lax transformations as 1v-cells, and modifications between the latter two. On the right-hand side is the double category of pseudodouble functors, 
their horizontal oplax transformations, vertical lax transformations, and modifications. (A reader interested in more details is referred to the beginning and the end of Section 5 in \cite{Fem:Bif}.) 
%

In the subsections that follow further below we will investigate conditions under which a pair $(\Bb, P)$ as above becomes a monoidal double category (recall \deref{Shul}).

\subsection{The categories of central binoidal structures} \sslabel{binoidal str}

We are going to consider three categories of central binoidal structures. First of all, 
we define the category $\Binoidal_{pc}^{st}(\Bb)$ of purely central binoidal structures. Its objects are 
purely central binoidal structures from \deref{pc binoidal}. 
Its morphisms consist of 
pairs of vertical {\em strict} transformations $\theta_l^A: A\ltimes_1-\Rightarrow A\ltimes_2-, \, \,
\theta_r^B: -\rtimes_1 B\Rightarrow-\rtimes_2 B$ for $A,B\in\Bb$ such that $\theta_l^A(B)=\theta_r^B(A)$ and which satisfy 
the analogous axioms to \axiomref{$VLT^q_1$}-\axiomref{$VLT^q_4$} from \deref{lax v tr cubical} of Appendix D, whereby 
the 2-cells $(u,U)_i, i=1,2$ are identities. 
%
%
Let $q\x\operatorname{Ps}_{vst}^{st}(\Bb\times\Bb,\Bb)$ be the (strict version of the) vertical category of the double category $q\x\Ps_{hop}^{st}(\Bb\times\Bb,\Bb)$ 
from \equref{equiv from Fem-Gray}. 
We thus obtain 
\begin{equation} \eqlabel{bin-st}
q\x\operatorname{Ps}_{hop}^{st}(\Bb\times\Bb,\Bb)\iso\Binoidal_{pc}^{st}(\Bb) \simeq \operatorname{Ps}_{hop}(\Bb\times\Bb,\Bb)
\end{equation}
where $\operatorname{Ps}_{hop}(\Bb\times\Bb,\Bb)$ is the (strict version of the) vertical category of the double category 
$\Ps_{hop}(\Bb\times\Bb,\Bb)$ from \equref{equiv from Fem-Gray}. 
Let $\G:\Binoidal_{pc}^{st}(\Bb) \to \operatorname{Ps}_{vst}(\Bb\times\Bb,\Bb)$ denote the obvious equivalence functor. 

For the other two categories of central binoidal structures we use analogous axioms to \axiomref{$VLT^q_1$}-\axiomref{$VLT^q_4$} with 
the following change. Whereas in the cited axioms the 2-cells of the form $(b,a)_1, (b,a)_2$ appear where $a\in\Aa, b\in\Bb$ are different types of 1-cells, we will differentiate ``left binoidal'' version of those axioms where instead of the latter 2-cells the 
2-cells of the form $a\ltimes_1-\vert_b, a\ltimes_2-\vert_b$ appear (here both $a,b\in\Bb$), and ``right binoidal'' version of those axioms with 2-cells of the form $-\ltimes_1 b\vert_a, -\ltimes_2 b\vert_a$. 

Now we define the category $\Binoidal_{lc}^{st}(\Bb)$ of left central binoidal structures. Its objects are 
left central binoidal structures from \deref{lr-central binoidal}. 
Its morphisms are pairs of vertical {\em strict} transformations $\theta_l^A: A\ltimes_1-\Rightarrow A\ltimes_2-, \, \,
\theta_r^B: -\rtimes_1 B\Rightarrow-\rtimes_2 B$ for $A,B\in\Bb$ such that $\theta_l^A(B)=\theta_r^B(A)$ and which satisfy the 
left binoidal version of the axioms \axiomref{$VLT^q_1$}-\axiomref{$VLT^q_4$} for left central 1h-cells $K$, left central 1v-cells $U$ and any 1h-cells $k$ and 1v-cells $U$, whereby the 2-cells $(u,U)_i, i=1,2$ are identities. 

Similarly, the category $\Binoidal_{rc}^{st}(\Bb)$ of right central binoidal structures 
is defined. Its morphisms 
$(\theta_l^A, \theta_r^B)_{A,B\in\Bb}$ with $\theta_l^A(B)=\theta_r^B(A)$ satisfy the right binoidal version of the axioms 
\axiomref{$VLT^q_1$}-\axiomref{$VLT^q_4$} for right central 1h-cells $k$, right central 1v-cells $u$ and any 1h-cells $K$ and 1v-cells 
$U$. 

\medskip

The relations from \thref{summary} can be extended functorially and added to the equivalences \equref{bin-st} so to obtain 
$$q\x\operatorname{Ps}_{hop}^{st}(\Bb\times\Bb,\Bb)\iso\Binoidal_{pc}^{st}(\Bb)\iso\Binoidal_{lc}^{st}(\Bb)\iso\Binoidal_{rc}^{st}(\Bb) \simeq \operatorname{Ps}_{hop}(\Bb\times\Bb,\Bb).$$

\subsection{Double quasi-functors with three and more variables}

We generalize 2-categorical quasi-functors in three and more variables and their quasi-natural transformations from 
\cite[Definition I.4.6]{Gray} to double categories.

\begin{defn} \delabel{cub}
A {\em pseudodouble quasi-functor} $H:\Aa\times\Bb\times\Cc\to\Ee$ (or a binary pseudodouble quasi-functor) 
consists of pseudodouble quasi-functors 
$$H(A,-,-):\Bb\times\Cc\to\Ee, \quad H(-,B,-):\Aa\times\Cc\to\Ee, \quad H(-,-,C): \Aa\times\Bb\to\Ee$$
for $(A,B,C)\in\Aa\times\Bb\times\Cc$, such that 
$$H(A,-,-)\vert_B=H(-,B,-)\vert_A, \,\, H(A,-,-)\vert_C=H(-,-,C)\vert_A, \,\, H(-,B,-)\vert_C=H(-,-,C)\vert_B,$$
and which give unambiguous pseudofunctors 
$$(A,B,-)_3:\Cc\to\Ee, \quad (A,-,C)_2:\Bb\to\Ee, \quad (-,B,C)_1: \Aa\to\Ee,$$ such that the structure 2-cells 
\equref{four cells} of the former relate in the following way:
\begin{enumerate} [i)]
\item for all 1h-cells $(f,g,h):(A,B,C)\to(A',B', C')$ in $\Aa\times\Bb\times\Cc$ it is  
$$\scalebox{0.8}{
\bfig
 \putmorphism(450,700)(1,0)[` `(A,g,C')_2]{680}1a
 \putmorphism(1140,700)(1,0)[` ` (f,B',C')_1]{650}1a

 \putmorphism(-150,250)(1,0)[``(A,B,h)_3]{600}1a
 \putmorphism(450,250)(1,0)[` `(f,B,C')_1]{680}1a
 \putmorphism(1130,250)(1,0)[` ` (A',g, C')_2]{650}1a

\putmorphism(450,700)(0,-1)[\phantom{Y_2}``=]{450}1r
\putmorphism(1750,700)(0,-1)[\phantom{Y_2}``=]{450}1r
\put(940,480){\fbox{$ (f,g,C')_{12}$}}

 \putmorphism(-150,-200)(1,0)[``(f,B,C)_1]{640}1a
 \putmorphism(460,-200)(1,0)[` `(A',B,h)_3]{680}1a

\putmorphism(-160,250)(0,-1)[\phantom{Y_2}``=]{450}1l
\putmorphism(1120,250)(0,-1)[\phantom{Y_3}``=]{450}1r
\put(310,20){\fbox{$(f,B,h)_{13}$}}

 \putmorphism(1140,-200)(1,0)[` ` (A',g,C')_2]{630}1a
\putmorphism(450,-200)(0,-1)[\phantom{Y_2}``=]{450}1r
\putmorphism(1750,-200)(0,-1)[\phantom{Y_2}``=]{450}1r

 \putmorphism(460,-650)(1,0)[` `(A',g,C)_2]{700}1a
 \putmorphism(1140,-650)(1,0)[ ` `(A',B',h)_3]{630}1a
\put(860,-430){\fbox{$ (A',g,h)_{23}$}}
\efig}
=
\scalebox{0.8}{
\bfig

 \putmorphism(-150,700)(1,0)[``(A,B,h)_3]{600}1a
 \putmorphism(450,700)(1,0)[` `(A,g,C')_2]{680}1a

 \putmorphism(-150,250)(1,0)[``(A,g,C)_2]{600}1a
 \putmorphism(450,250)(1,0)[` `(A,B',h)_3]{680}1a
 \putmorphism(1120,250)(1,0)[` `(f,B',C')_1]{650}1a

\putmorphism(-180,700)(0,-1)[\phantom{Y_2}``=]{450}1r
\putmorphism(1100,700)(0,-1)[\phantom{Y_2}``=]{450}1r
\put(280,470){\fbox{$(A,g,h)_{23}$}}
\put(910,20){\fbox{$(f,B',h)_{13}$}}

 \putmorphism(-150,-200)(1,0)[``(A,g,C)_2]{600}1a
 \putmorphism(450,-200)(1,0)[` `(f,B',C)_1]{680}1a
 \putmorphism(1100,-200)(1,0)[` ` (A',B',h)_3]{660}1a

\putmorphism(450,250)(0,-1)[\phantom{Y_2}``=]{450}1l
\putmorphism(1750,250)(0,-1)[\phantom{Y_2}``=]{450}1r
\putmorphism(-180,-200)(0,-1)[\phantom{Y_2}``=]{450}1r
\putmorphism(1100,-200)(0,-1)[\phantom{Y_2}``=]{450}1r
 \putmorphism(-150,-650)(1,0)[``(f,B,C)_1]{600}1a
 \putmorphism(580,-650)(1,0)[`` (A',g,C)_2]{540}1a
\put(280,-420){\fbox{$(f,g,C)_{12}$}}
\efig} \vspace{-0,2cm}
$$
\item for all 1v-cells $(u,v,z):(A,B,C)\to(\tilde A, \tilde B, \tilde C)$ in $\Aa\times\Bb\times\Cc$ it is  
$$
\scalebox{0.8}{
\bfig
 \putmorphism(0,480)(1,0)[` `=]{470}1a
\putmorphism(0,500)(0,-1)[` `(A,B,z)_3]{450}1l
\put(40,200){\fbox{$(v,z)_{23}$}}
\putmorphism(0,-320)(1,0)[ ` `=]{400}1a 
\putmorphism(0,90)(0,-1)[ ``(A,v,\tilde C)_2]{450}1l
\putmorphism(470,90)(0,-1)[``]{450}1l 
\putmorphism(490,90)(0,-1)[``(A,\tilde B,z)_3]{450}0l 
\putmorphism(470,500)(0,-1)[``]{450}1l 
\putmorphism(450,90)(1,0)[``=]{520}1a 
\putmorphism(970,90)(0,-1)[ ``]{450}1r
\putmorphism(950,90)(0,-1)[ ``(u,\tilde B,C)_1]{450}0r
\putmorphism(970,-300)(0,-1)[``(\tilde A,\tilde B,z)_3]{450}1r
\putmorphism(470,-730)(1,0)[ ``=]{500}1a 
\putmorphism(470,-300)(0,-1)[` `(u,\tilde B,\tilde C)_1]{450}1l 
\put(530,-430){\fbox{$(u,z)_{13}$}}

 \putmorphism(980,480)(1,0)[` `=]{470}1a %
\putmorphism(970,500)(0,-1)[` `(A,v,C)_2]{450}1l
\putmorphism(1450,500)(0,-1)[` `(u,B,C)_1]{450}1r
\putmorphism(1450,90)(0,-1)[``(\tilde A,v,C)_2]{450}1r
\putmorphism(1010,-320)(1,0)[` `=]{400}1a 
\put(1040,200){\fbox{$(u,v)_{12}$}}
\efig}
=
\scalebox{0.8}{
\bfig
 \putmorphism(-170,480)(1,0)[` `=]{480}1a
\putmorphism(-180,500)(0,-1)[\phantom{Y_2}` `(A,B,z)_3]{450}1l
\put(-120,250){\fbox{$(u,z)_{13}$}}
\putmorphism(-180,-320)(1,0)[` `=]{480}1a
\putmorphism(-180,90)(0,-1)[\phantom{Y_2}``]{450}1l
\putmorphism(-150,50)(0,-1)[\phantom{Y_2}``(u,B,\tilde C)_1]{450}0l
\putmorphism(300,90)(0,-1)[\phantom{Y_2}``]{450}1r
\putmorphism(280,90)(0,-1)[\phantom{Y_2}``(\tilde A,B,z)_3]{450}0r
\putmorphism(300,500)(0,-1)[\phantom{Y_2}` `(u,B, C)_1]{450}1r
\putmorphism(-650,90)(1,0)[ ``=]{400}1a
\putmorphism(-670,90)(0,-1)[ ``(A,v,\tilde C)_2]{450}1l
\putmorphism(-670,-300)(0,-1)[``(u,\tilde B,\tilde C)_1]{450}1l
\putmorphism(-650,-730)(1,0)[ ``=]{440}1a 
\putmorphism(-180,-300)(0,-1)[`  `(\tilde A,v,\tilde C)_2]{450}1r
\put(-600,-500){\fbox{$(u,v)_{12}$}}

\putmorphism(330,90)(1,0)[` `=]{450}1a
\putmorphism(770,90)(0,-1)[``(\tilde A,v,C)_2]{450}1r
\putmorphism(300,-300)(0,-1)[`  `]{450}1r 
\putmorphism(770,-300)(0,-1)[`  `(\tilde A,\tilde B,z)_3]{450}1r
\putmorphism(350,-730)(1,0)[` `=]{440}1a 
\put(360,-500){\fbox{$(v,z)_{23}$}}
\efig}
$$
(where we simplified the notation by writing $(v,z)_{23}$ for the 2-cell $(A,v,z)_{23}$ and so on...), 
\item for $(f,v,h):(A,B,C)\to(A',\tilde B,C')$ \vspace{-0,3cm}
$$\scalebox{0.86}{
\bfig
\putmorphism(-150,500)(1,0)[``(A,B,h)_3]{600}1a
 \putmorphism(480,500)(1,0)[` `(f,B,C')_1]{640}1a
 \putmorphism(-150,50)(1,0)[``(f,B,C)_1]{600}1a
 \putmorphism(470,50)(1,0)[` `(A',B,h)_3]{660}1a

\putmorphism(-180,500)(0,-1)[\phantom{Y_2}``=]{450}1l
\putmorphism(1100,500)(0,-1)[\phantom{Y_2}``=]{450}1r
\put(230,280){\fbox{$(f,B,h)_{13}$}}

\putmorphism(-170,-400)(1,0)[` `(f,\tilde B,C)_1]{640}1b
 \putmorphism(470,-400)(1,0)[` `(A',\tilde B,h)_3]{640}1b

\putmorphism(-180,50)(0,-1)[\phantom{Y_2}``(A,v,C)_2]{450}1l %
\putmorphism(450,50)(0,-1)[\phantom{Y_2}``]{450}1l
\putmorphism(1100,50)(0,-1)[\phantom{Y_3}``(A',v,C')_2]{450}1r
\put(-90,-180){\fbox{$(f,v,C)_{12}$}} 
\put(540,-180){\fbox{$(A',v,h)_{23}$}}
\efig}
\quad
=
\quad
\scalebox{0.86}{
\bfig
\putmorphism(-150,500)(1,0)[``(A,B,h)_3]{600}1a
 \putmorphism(480,500)(1,0)[` `(f,B,C')_1]{640}1a

 \putmorphism(-150,50)(1,0)[``(A,\tilde B,h)_3)]{600}1a
 \putmorphism(450,50)(1,0)[``(f,\tilde B,C')_1]{640}1a

\putmorphism(-180,500)(0,-1)[\phantom{Y_2}``(A,v,C)_2]{450}1l
\putmorphism(450,500)(0,-1)[\phantom{Y_2}``]{450}1r
\putmorphism(1100,500)(0,-1)[\phantom{Y_2}``(A',v,C')_2]{450}1r
\put(-90,280){\fbox{$(A,v,h)_{23}$}}
\put(540,280){\fbox{$(f,v,C')_{12}$}}

\putmorphism(-170,-400)(1,0)[` `(f,\tilde B,C)_1]{640}1b
 \putmorphism(470,-400)(1,0)[` `(A',\tilde B,h)_3]{640}1b

\putmorphism(-180,50)(0,-1)[\phantom{Y_2}``=]{450}1l
\putmorphism(1100,50)(0,-1)[\phantom{Y_3}``=]{450}1r
\put(270,-200){\fbox{$(f,\tilde B,h)_{13}$}}
\efig}
$$

and 2 similar conditions, one for $(f,g,z):(A,B,C)\to(A',B',\tilde C)$ and the other for $(u,g,h):(A,B,C)\to(\tilde A,B',C')$, 
\item for $(u,v,h):(A,B,C)\to(\tilde A,\tilde B,C')$  \vspace{-0,3cm}
$$\scalebox{0.86}{
\bfig
 \putmorphism(-120,500)(1,0)[` `=]{550}1a
 \putmorphism(450,500)(1,0)[` `(A,B,h)_3]{550}1a
\putmorphism(-140,520)(0,-1)[` `(A,v,C)_2]{480}1l
\put(-80,50){\fbox{$(u,v,C)_{12}$}}
\putmorphism(-150,-380)(1,0)[` `=]{540}1a
\putmorphism(-140,80)(0,-1)[``(u,\tilde B,C)_1]{450}1l
\putmorphism(430,50)(0,-1)[` `(\tilde A,v,C)_2]{450}1l
\putmorphism(430,520)(0,-1)[` `(u,B,C)_1]{480}1l
\put(470,290){\fbox{$(u,B,h)_{13}$}}
\putmorphism(430,50)(1,0)[``(\tilde A,B,h)_3]{540}1a
\putmorphism(1000,80)(0,-1)[``(\tilde A,v,C')_2]{450}1r
\putmorphism(1000,520)(0,-1)[``(u,B,C')_1]{480}1r
\putmorphism(450,-380)(1,0)[``(\tilde A,\tilde B,h)_3]{540}1b
\put(500,-190){\fbox{$(\tilde A, v,h)_{23}$}}
\efig}
\quad=\quad
\scalebox{0.86}{
\bfig
 \putmorphism(-150,500)(1,0)[` `(A,B,h)_3]{600}1a
 \putmorphism(450,500)(1,0)[` `=]{540}1a
\putmorphism(-180,520)(0,-1)[` `(A,v,C)_2]{450}1l
\put(-120,280){\fbox{$(A,v,h)_{23}$}}
\putmorphism(-150,-380)(1,0)[` `(\tilde A,\tilde B,h)_3]{500}1b
\putmorphism(-180,80)(0,-1)[``(u,\tilde B,C)_1]{450}1l
\putmorphism(400,80)(0,-1)[``(u,\tilde B,C')]{450}1r
\putmorphism(400,520)(0,-1)[` `(A,v,C')_2]{450}1r
\putmorphism(-150,50)(1,0)[``(A,\tilde B,h)_3]{500}1a
\putmorphism(1000,80)(0,-1)[``(\tilde A,v,C')_2]{450}1r
\putmorphism(1000,520)(0,-1)[``(u,B,C')_1]{450}1r
\putmorphism(450,-380)(1,0)[``=]{520}1b
\put(-120,-170){\fbox{$(u, \tilde B,h)_{13}$}}
\put(460,50){\fbox{$(u,v,C')_{12}$}}
\efig}
$$

and 2 similar conditions, one for $(u,g,z):(A,B,C)\to(\tilde A,B',\tilde C)$ and the other for $(f,v,z):(A,B,C)\to(A',\tilde B,C')$, 
where $f,g,h$ are 1h-cells and $u,v,z$ are 1v-cells, as usual. 
\end{enumerate}
\end{defn}

Observe that in the above definition there are $2\cdot 2\cdot 3=12$ structure 2-cells combining in $2^3=8$ ways.

\begin{defn} \delabel{vst ternary quasi}
A {\em vertical strict transformation} $\theta: H_1\Rightarrow H_2$ between pseuodouble quasi-functors $H_1,H_2:\Aa\times\Bb\times\Cc\to\Ee$ consists of vertical strict transformations 
$$\theta^A: H_1(A,-,-) \Rightarrow H_2(A,-,-),$$ 
$$\theta^B: H_1(-,B,-) \Rightarrow H_2(-,B,-),$$ 
$$\theta^C: H_1(-,-,C) \Rightarrow H_2(-,-,C)$$ 
of binary pseudodouble quasi-functors, which give unambiguous vertical strict transformations 
$$\theta^{A;B}: (A,B,-)_3^1\to (A,B,-)_3^2$$
$$\theta^{B;C}: (-,B,C)_1^1\to (-,B,C)_1^2$$ 
$$\theta^{A;C}: (A,-,C)_2^1\to (A,-,C)_2^2$$ 
of pseudodouble functors for each $(A,B,C)\in\Aa\times\Bb\times\Cc$, so that twelve equalities between their structure 2-cells, on one hand, 
and the twelve structure 2-cells of both $H_1$ and $H_2$ from \deref{cub}, on the other hand, hold. We present these twelve equalities schematically as lists consisting of those structure 2-cells which are related in the only possible way by one equation: 
$$(\theta^{A;C})_g, (\theta^{B;C})_f, (f,g,C)^i_{12}, \quad (\theta^{A;B})_h, (\theta^{B;C})_f, (f,B,h)^i_{13}, \quad 
(\theta^{A;B})_h, (\theta^{A;C})_g, (A,g,h)^i_{23}$$
$$(\theta^{A;C})^v, (\theta^{B;C})_f, (f,v,C)^i_{12}, \quad (\theta^{A;B})^z, (\theta^{B;C})_f, (f,B,z)^i_{13}, \quad 
(\theta^{A;B})^z, (\theta^{A;C})_g, (A,g,z)^i_{23}$$
$$(\theta^{A;C})_g, (\theta^{B;C})^u, (u,g,C)^i_{12}, \quad (\theta^{A;B})_h, (\theta^{B;C})^u, (u,B,h)^i_{13}, \quad 
(\theta^{A;B})_h, (\theta^{A;C})^v, (A,v,h)^i_{23}$$
$$(\theta^{A;C})^v, (\theta^{B;C})^u, (u,v,C)^i_{12}, \quad (\theta^{A;B})^z, (\theta^{B;C})^u, (u,B,z)^i_{13}, \quad 
(\theta^{A;B})^z, (\theta^{A;C})^v, (A,v,z)^i_{23}$$
where $(f,g,C)_{12}^i$ for $i=1,2$ presents a structure 2-cell $(f,g,C)_{12}$ from \deref{cub} for $H_1$ and $H_2$, respectively, and similarly for the remaining eleven 2-cells of that type. 
\end{defn}

Observe that it holds $\theta^{A;B}(C)=\theta^{B;C}(A)=\theta^{A;C}(B)$ for all $(A,B,C)\in\Aa\times\Bb\times\Cc$ in the above definition. 
Also, the first row of the axioms corresponds to the axiom \axiomref{$VLT^q_4$} applied to the three coordinates, and similarly the second, third and fourth row correspond to the axioms \axiomref{$VLT^q_2$}, \axiomref{$VLT^q_3$} and \axiomref{$VLT^q_1$}, respectively. 
Moreover, observe that the 2-cell components $(\theta^{A;B})^z, (\theta^{A;C})^v, (\theta^{B;C})^u$ are trivial, and the last three axioms 
show the way in which $(u,v,C)^1_{12}=(u,v,C)^2_{12}, \,\, (u,B,z)^1_{13}=(u,B,z)^2_{13}$ and $(A,v,z)^1_{23}=(A,v,z)^2_{23}$. 

\begin{rem}
Visually, the pasting diagrams for the 12 equations from the above definition have the same form as the diagrams of \deref{cub}. 
Indeed, in Subsection 5.6 of \cite{Fem:Gray1} we showed that there is a 1-1 correspondence between 
{\em vertical transformations of quasi-functors of two variables} and {\em quasi-functors of three variables}. We had an 
analogous situation in the double funny case, recall \prref{ternary eq}. 
\end{rem}

\bigskip

Pseudodouble quasi-functors of more than three variables and their vertical strict transformations are defined as follows. 

\begin{defn} \delabel{4-ary}
A pseudodouble quasi-functor $H:\Aa_1\times...\times\Aa_n\to\Ee$ for $n>3$ consists of pseudodouble quasi-functors of three variables 
$$H(A_1,...,A_{i-1},\, -,\, A_{i+1},...,A_{j-1}, \, -,\, A_{j+1},...,A_{k-1}, \, -,\, A_{k+1},..., A_n): 
\Aa_i\times\Aa_j\times\Aa_k\to\Ee$$
for all $i<j<k$ and all choices of objects $A_l\in\Aa_l, l=1,...,n$. 

A vertical strict transformation $\theta: H_1\Rightarrow H_2$ between pseudodouble quasi-functors 
$H_1,H_2:\Aa_1\times...\times\Aa_n\to\Ee$ for $n>3$ consists of vertical strict transformations 
$$\theta^i: H_1(A_j,-,-) \Rightarrow H_2(A_i,-,-),$$ 
$$\theta^j: H_1(-,A_j,-) \Rightarrow H_2(-,A_j,-),$$ 
$$\theta^k: H_1(-,-,A_k) \Rightarrow H_2(-,-,A_k)$$ 
of binary pseudodouble quasi-functors for all $i<j<k\leq n$ and all choices of objects $A_l\in\Aa_l, l=1,...,n$, where we omit the irrelevant variables.  
\end{defn}

\subsection{Associativity constraints} \sslabel{ass}  

We now come back to the question from the end of \ssref{towards}. We saw that a pseudodouble quasi-functor $H:\Bb\times\Bb\to\Bb$ whose 
2-cells $(u,U)$ are trivial induces a pseudodouble functor $P:\Bb\times\Bb\to\Bb$, and that there is even a 
double category equivalence 
\begin{equation} \eqlabel{F}
\F\colon q\x\Ps_{hop}^{st}(\Bb\times\Bb,\Bb) \to \Ps_{hop}(\Bb\times\Bb,\Bb)
\end{equation}
$$\hspace{1,8cm}H \mapsto P=\ot$$
(we marked the assignment on 0-cells). The pseudodouble functor $P$ is defined on a 2-cell $(\alpha,\beta)\in\Bb\times\Bb$ by 
\begin{equation} \eqlabel{P on 2-cells}
P(\alpha,\beta)\colon\hspace{-0,2cm}= 
\scalebox{0.86}{
\bfig

 \putmorphism(-520,50)(1,0)[(B,A)`(B,A')`(B,K)]{600}1a
 \putmorphism(30,50)(1,0)[\phantom{A''\ot B'}` (B', A') `(k, A')]{670}1a

\putmorphism(70,50)(0,-1)[\phantom{Y_2}`(B,\tilde A')`]{450}1r
\putmorphism(50,180)(0,-1)[``(B,U')]{450}0r
\putmorphism(650,50)(0,-1)[\phantom{Y_2}`\phantom{A''\ot B'} `(B',U')]{450}1r
\putmorphism(70,-400)(0,-1)[\phantom{Y_2}``]{450}1r
\putmorphism(70,-490)(0,-1)[\phantom{Y_2}``(u,\tilde A')]{450}0r
\putmorphism(650,-400)(0,-1)[\phantom{Y_2}``(u', \tilde A')]{450}1r

\put(210,-580){\fbox{$(\beta, \tilde A')$}}
\put(210,-210){\fbox{$(k,U')$}}

\putmorphism(-470,50)(0,-1)[\phantom{Y_2}`(B,\tilde A)`(B,U)]{450}1l
\putmorphism(-470,-400)(0,-1)[\phantom{Y_2}`(\tilde B,\tilde A)`(u,\tilde A)]{450}1l
\put(-350,-150){\fbox{$(B,\alpha)$}}
\put(-350,-650){\fbox{$(u,\tilde K)$}}

 \putmorphism(-510,-400)(1,0)[\phantom{A''\ot B'}`\phantom{A''\ot B'}`(B,\tilde K)]{600}1a
 \putmorphism(30,-400)(1,0)[\phantom{A''\ot B'}` \phantom{A''\ot B'} `]{670}1a
 \putmorphism(30,-410)(1,0)[\phantom{A''\ot B'}` (B',\tilde A')`(k, \tilde A')]{670}0a
 \putmorphism(-530,-850)(1,0)[\phantom{A''\ot B'}`(\tilde B,\tilde A')`(\tilde B,\tilde K)]{600}1b
\putmorphism(30,-850)(1,0)[\phantom{A''\ot B'}` (\tilde B',\tilde A') `(\tilde k, \tilde A')]{670}1b
\efig 
}
\end{equation} 
(mind that, as in \ssref{quasi}, we denote the two pseudodouble functors determining a pseudodouble quasi-functor $H$ 
by $(-,A)$ and $(B,-)$. ) 

For the pseudodouble functor structure $\gamma_{(K',k')(K,k)}\colon P(K',k')P(K,k)\Rightarrow P(K'K,k'k)$ and $\iota^P\colon 1_{P(A,B)}\Rightarrow P(1_{(A,B)})$ of $P$ we set \vspace{-0,4cm}
\begin{equation} \eqlabel{natur 2-cells} 
\gamma_{(K',k')(K,k)}\colon\hspace{-0,2cm}=
\scalebox{0.86}{
\bfig
\putmorphism(210,350)(1,0)[` `(k,A')]{420}1a
\putmorphism(530,350)(1,0)[\phantom{F(A)}` `(B',K')]{450}1a

\putmorphism(-190,10)(1,0)[` `(B,K)]{420}1a
\putmorphism(210,10)(1,0)[` `(B,K')]{420}1a
\putmorphism(530,10)(1,0)[\phantom{F(A)}` `(k,A'')]{450}1a
\putmorphism(900,10)(1,0)[\phantom{F(A)}` `(k',A'')]{450}1a

\putmorphism(210,350)(0,-1)[\phantom{Y_2}``=]{350}1l
\putmorphism(980,350)(0,-1)[\phantom{Y_2}``=]{350}1r

\putmorphism(-180,0)(0,-1)[\phantom{Y_2}``=]{350}1l
\putmorphism(620,0)(0,-1)[\phantom{Y_2}``=]{350}1r
\putmorphism(1340,0)(0,-1)[\phantom{Y_2}``=]{350}1r

\put(420,200){\fbox{$(k,K')$}}

 \putmorphism(-190,-320)(1,0)[``(B,K'K)]{800}1b
 \putmorphism(530,-320)(1,0)[\phantom{F(A)}` `(k'k,A'')]{800}1b

\put(-30,-180){\fbox{$(B,-)_{K'K}$}}
\put(760,-180){\fbox{$(-,A'')_{k'k}$}}
\efig}
\end{equation} 
and 
$$\iota^P_{(A,B)}\colon\hspace{-0,2cm}=
\scalebox{0.86}{
\bfig
 \putmorphism(-210,220)(1,0)[(B,A)`\phantom{F(A)} `=]{500}1a
\putmorphism(330,220)(1,0)[(B,A)`(B,A)`=]{500}1a
\putmorphism(-210,220)(0,-1)[\phantom{Y_2}``=]{370}1l
\putmorphism(320,220)(0,-1)[\phantom{Y_2}``=]{370}1l
\putmorphism(810,220)(0,-1)[\phantom{Y_2}``=]{370}1r
 \putmorphism(-210,-150)(1,0)[(B,A)`\phantom{Y_2}`(B,1_A)]{470}1a
 \put(490,50){\fbox{$\iota^B_A$}} 
\putmorphism(320,-150)(1,0)[(B,A)`(B,A) `(1_B,A)]{540}1a
\put(-60,40){\fbox{$\iota^A_B$}}
\efig}
$$
where $\iota^B_A=(B,-)_A$ and $\iota^A_B=(-,A)_B$ of $H$. 

Let us explore how associativity of $H$ is related to associativity of $P$. From now on we will write $\ot$ for $P$. 

Analogously to the proof in \cite[Section 5]{Fem:Bif}, plugging-in a third fixed variable in three different positions on the side of quasi-functors, one could first obtain an analogous double category 
$q_3\x\Ps_{hop}^{st}(\Bb\times\Bb\times\Bb,\Bb)$, and then an analogous equivalence of double categories to $\F$ in \equref{F}. 
Since here we are not interested in horizontal quasi-natural transformations nor in the modifications, we satisfy ourselves with a category equivalence. For this purpose, let $q_3\x\operatorname{Ps}_{vst}^{st}(\Bb\times\Bb\times\Bb,\Bb)$ denote the category 
of psuedodouble quasi-functors and their vertical strict transformations, and $\operatorname{Ps}_{vst}(\Bb\times\Bb\times\Bb,\Bb)$ its usual version on the Cartesian product.

\begin{thm} \thlabel{H equiv ot}
For a double category $\Bb$ there is an equivalence of categories 
\begin{equation} \eqlabel{ps-eq}
\F_3\colon q_3\x\operatorname{Ps}_{vst}^{st}(\Bb\times\Bb\times\Bb,\Bb) \to \operatorname{Ps}_{vst}(\Bb\times\Bb\times\Bb,\Bb).
\end{equation}
Similarly, there is an equivalence of categories 
$$\F_4\colon q_4\x\operatorname{Ps}_{vst}^{st}(\Bb\times\Bb\times\Bb\times\Bb,\Bb) \to 
\operatorname{Ps}_{vst}(\Bb\times\Bb\times\Bb\times\Bb,\Bb).$$
\end{thm}

\begin{proof}
The proof is analogous to that of \cite[Theorem 5.7]{Fem:Bif}. We only comment two points. The functor $\F_3(H)$ is defined for $H\in 
q_3\x\operatorname{Ps}_{vst}^{st}(\Bb\times\Bb\times\Bb,\Bb)$ by: 
$$\scalebox{0.86}{
\bfig

 \putmorphism(-480,50)(1,0)[``(f,B,C)_1]{580}1a
 \putmorphism(90,50)(1,0)[`  `(A',g,C)_2]{580}1a
 \putmorphism(700,50)(1,0)[`  `(A',B',h)_3]{620}1a

\putmorphism(90,50)(0,-1)[\phantom{Y_2}``]{450}1r
\putmorphism(700,50)(0,-1)[\phantom{Y_2}`\phantom{A''\ot B'} `]{450}1r 
\putmorphism(1310,50)(0,-1)[\phantom{Y_2}`\phantom{A''\ot B'} `(u',B',C')_1]{450}1r
\putmorphism(90,-400)(0,-1)[\phantom{Y_2}``]{450}1r
\putmorphism(700,-400)(0,-1)[\phantom{Y_2}``]{450}1r 
\putmorphism(1310,-400)(0,-1)[\phantom{Y_2}``(\tilde A',v',C')_2]{450}1r
\putmorphism(700,-800)(0,-1)[\phantom{Y_2}``]{450}1r 
\putmorphism(90,-800)(0,-1)[\phantom{Y_2}``]{450}1r
\putmorphism(1310,-800)(0,-1)[\phantom{Y_2}``(\tilde A',\tilde B',z')_3]{450}1r

\put(150,-150){\fbox{$(u',g,C)_{12}$}}
\put(150,-650){\fbox{$(\tilde A',\beta, C)_2$}}
\put(150,-1150){\fbox{$(\tilde A',\tilde g, z)_{23}$}}

\put(760,-150){\fbox{$(u',B', h)_{13}$}}
\put(760,-650){\fbox{$(\tilde A',v, h)_{23}$}}
\put(760,-1150){\fbox{$(\tilde A',\tilde B', \gamma)_3$}}

\putmorphism(-470,50)(0,-1)[\phantom{Y_2}``(u,B,C)_1]{450}1l
\putmorphism(-470,-400)(0,-1)[\phantom{Y_2}``(\tilde A,v,C)_2]{450}1l
\put(-420,-150){\fbox{$(\alpha,B,C)_1$}}
\put(-420,-650){\fbox{$(\tilde f,v,C)_{12}$}}
\put(-420,-1150){\fbox{$(\tilde f,\tilde B,z)_{13}$}}

 \putmorphism(-480,-400)(1,0)[``(\tilde f,B,C)_1]{560}1a 
 \putmorphism(80,-410)(1,0)[` `(\tilde A',g,C)_2]{620}1a
  \putmorphism(680,-400)(1,0)[` `(\tilde A',B',h)_3]{620}1a
\putmorphism(-480,-850)(1,0)[` `(\tilde f,\tilde B,C)_1]{560}1b
\putmorphism(110,-850)(1,0)[`  `(\tilde A',\tilde g, C)_2]{580}1b
\putmorphism(700,-850)(1,0)[`  `(\tilde A', \tilde B',h)_3]{580}1b

\putmorphism(-620,-1300)(1,0)[\phantom{A''\ot B'}`  `(\tilde f, \tilde B, \tilde C)_1]{700}1b
\putmorphism(-50,-1300)(1,0)[\phantom{A''\ot B'}`  `(\tilde A', \tilde g, \tilde C)_2]{720}1b
\putmorphism(550,-1300)(1,0)[\phantom{A''\ot B'}`  `(\tilde A', \tilde B', \tilde h)_3]{750}1b
\putmorphism(-470,-850)(0,-1)[\phantom{Y_2}``(\tilde A,\tilde B,z)_3]{450}1l
\efig 
}$$
with usual notations for 1-cells, and where the 2-cells are those from \deref{cub}, parts (iii) and (iv). 

Regarding the correspondence of morphisms, we use the keypoint from the proof of the binary case. Namely, from the proof of the correspondence of (horizontal and) vertical transformations on both sides of the equivalence in \cite[Section 5.4]{Fem:Bif} it is clear that the five axioms of transformations between pseudodouble functors 
(right-hand side of the equivalence) hold by: 1) the corresponding five axioms of the component transformations on the quasi-functors (left-hand side), 2) the axioms \axiomref{$VLT^q_1$}-\axiomref{$VLT^q_4$}, and 3) axioms of the quasi-functors. Then exactly the analogous happens in the current ternary case, whereby the four axioms 
\axiomref{$VLT^q_1$}-\axiomref{$VLT^q_4$} are replaced by their analogues in the three variables, which are precisely the twelve axioms of 
\deref{vst ternary quasi}. 
\qed\end{proof}

From the definition of $\F_3$ it is clear that we have:

\begin{lma}
The functor $\F_3$ preserves invertible vertical strict transformations whose 1v-cell components are inversely central 1v-cells. 
\end{lma}

As a matter of fact, due to \thref{lr central} a double category admitting a pseudodouble quasi-functor is purely central, so all 
its 1v-cells are central. 

\medskip

In \cite[Proposition 5.17]{Fem:Gray1} we generalized Gray's substitution result for quasi-functors \cite[Theorem I.4.7]{Gray} from 
2-categories to double categories and Gray-categories. We cite here the version of the result for pseudodouble quasi-functors 
for double categories.

\begin{prop} \cite[Proposition 5.17]{Fem:Gray1} \prlabel{subst-gen}
Given pseudodouble quasi-functors $F_i:\Aa_{i1}\times...\times\Aa_{im_i}\to\Bb_i$ of $m_i$-variables with $i=1,...,n, m_i\geq 2, n\geq 2$ and a quasi-functor $G:\Bb_1\times...\times\Bb_n\to\Cc$ of $n$-variables, the composition 
$$\Pi_{j=1}^{m_1}\Aa_{1j}\times...\times\Pi_{j=1}^{m_n}\Aa_{nj}\stackrel{F_1\times...\times F_n}{\longrightarrow}\Bb_1\times...\times\Bb_n
\stackrel{G}{\to}\C$$
is a pseudodouble quasi-functor of $m_1+..+m_n$-variables. 
\end{prop}

Now consider the following two diagrams 
\begin{equation} \eqlabel{teta-sigma}
\scalebox{0.84}{
\bfig 
\putmorphism(0,400)(1,0)[\Bb\times\Bb\times\Bb `\Bb\times\Bb ` 1\times H]{1000}1a
\putmorphism(50,0)(1,0)[\Bb\times\Bb `\Bb ` H]{950}1b
\putmorphism(60,400)(0,1)[``H\times 1]{380}1l 
\putmorphism(1000,400)(0,1)[`` H]{380}1r 
\put(800,250){$\Swarrow$}
\put(480,160){\fbox{$\theta$}}
\efig}
\hspace{3,5cm}
\scalebox{0.84}{
\bfig 
\putmorphism(0,400)(1,0)[\Bb\times\Bb\times\Bb `\Bb\times\Bb ` 1\times \ot]{1000}1a
\putmorphism(50,0)(1,0)[\Bb\times\Bb `\Bb. ` \ot]{950}1b
\putmorphism(60,400)(0,1)[``\ot\times 1]{380}1l 
\putmorphism(1000,400)(0,1)[`` \ot]{380}1r 
\put(800,250){$\Swarrow$}
\put(480,170){\fbox{$\Sigma$}}
\efig}
\end{equation}
According to \prref{subst-gen}, the compositions $H_1:=H(1\times H)$ and $H_2:=H(H\times 1)$ are pseudodouble quasi-functors 
$H_1,H_2: \Bb\times\Bb\times\Bb\to\Bb$. Then the compositions of the arrows in the two diagrams above present 0-cells assigned to each other via $\F_3$ in \equref{ps-eq} in the obvious way. The associativity $\theta: H_1\Rightarrow H_2$ is 
a 1-cell corresponding to $\Sigma$ in \equref{ps-eq}.

We have that $\theta: H_1\Rightarrow H_2$ is 
comprised of three unambiguous vertical strict transformations $\theta^{A;B}, \theta^{B;C}$ and $\theta^{A;C}$ 
of pseudodouble functors for each $(A,B,C)\in\Aa\times\Bb\times\Cc$, which obey the 12 axioms from \deref{vst ternary quasi}. 
Observe that the pentagon for the associativity $\theta: H_1\Rightarrow H_2$ is an equality of two vertical strict transformations of 
pseudodouble quasi-functors of four variables, one, say $\Theta^1$, is a composite of three (left and bottom arrows in the diagram below), 
and the other, say $\Theta^2$, is a composite of two (top and right arrow below): 
\begin{equation*} 
\scalebox{0.84}{
\bfig 
\putmorphism(0,500)(1,0)[H(H(H(-,-),-),-) `H(H(-,-),H(-,-)) ` \theta^{\bullet,\bullet}_{--,-,-}]{3000}1a
\putmorphism(0,0)(1,0)[H(H(-,H(-,-)),-) `H(-,H(H(-,-),-)) ` \theta^{\bullet,\bullet}_{-,--,-}]{1350}1b
\putmorphism(1380,0)(1,0)[\phantom{H(-,H(-,H(-,-)))}`H(-,H(-,H(-,-))).` -\rtimes\theta^{\bullet,\bullet}_{-,-,-}]{1620}1b
\putmorphism(60,500)(0,1)[`` \theta^{\bullet,\bullet}_{-, -,-}\rtimes -]{500}1l 
\putmorphism(3010,500)(0,1)[`` \theta^{\bullet,\bullet}_{-,-,--}]{500}1r 
\efig}
\end{equation*}
As transformations of pseudodouble quasi-functors of four variables they 
are given by four unambiguous vertical strict transformations of pseudodouble functors. Thus, there are four pentagons 
that describe the associativity $\theta: H_1\Rightarrow H_2$. 

Since $H$ gives a binoidal structure on $\Bb$ (recall \ssref{quasi}), the above three unambiguous vertical strict transformations 
$\theta^{A;B}, \theta^{B;C}$ and $\theta^{A;C}$ correspond to the vertical strict transformations $\alpha_{-,B,C}, \alpha_{A,-,C}, 
\alpha_{A,B,-}$ from \deref{assoc}, and the four pentagons for the $\theta$'s correspond to those of that same definition. 
Moreover, the 12 identities from \deref{vst ternary quasi} correspond to the 24 axioms from \seref{24 axioms} and Appendix A, 
which because of pure centrality, \thref{lr central}, come down to 12 axioms. By the comment from the end of the proof of 
\thref{H equiv ot} we know that those 12 axioms of $\theta$'s correspond to the axioms of $\Sigma$ being a vertical strict transformation.

As a consequence of \thref{H equiv ot} we have:

\begin{cor} \colabel{assoc iso}
Let $H:\Bb\times\Bb\to\Bb$ be a pseudodouble quasi-functor whose 2-cells $(u,U)$ are trivial, 
$\ot=\F(H)$ and let $\theta$ and $\Sigma$ be as in \equref{teta-sigma}. 
Then $(H,\theta:H_1\Rightarrow H_2)$ is an associative binoidal structure making $\Bb$ purely central if and only if 
$(\ot, \Sigma: \ot(\ot\times 1)\Rightarrow \ot(1\times\ot)$ is an associative product on $\Bb\times\Bb$. 
\end{cor}


\medskip

We record that in the proof of the above corollary the following is entailed:

\begin{prop} \prlabel{quasi-24 ax}
An associative binoidal structure that comes from a double quasi-functor $H$ (and makes a double category $\Bb$ premonidal) 
satisfies the 24 axioms, where the centrality structures are induced by $H$. 
\end{prop}

\subsection{Premonoidal double categories coming from quasi-functors are purely central, monoidal and satisfy the 24 axioms} 
\sslabel{quasi-premon}

For a premonoidal double category that is purely (resp. left or right) central binoidal we will say that it is purely (resp. left or right) central premonoidal. We are ready to prove:

\begin{thm} \thlabel{premon-mon}
Let $\Bb$ be a double category and bear in mind the double category equivalence \equref{F}. Then: \vspace{-0,1cm}
\begin{enumerate}
\item $\theta:H_1\Rightarrow H_2$ obeys the 24 axioms, which are indeed 12 axioms; 
\item $(\Bb,H,\theta,I)$ is a purely central premonoidal double category whereby the structure 2-cells $(u,U)$ are trivial if and only if 
$(\Bb, \ot, \Sigma, I)$ is a monoidal double category. 
\end{enumerate}
\end{thm}

\begin{proof}
First recall that a pseudodouble quasi-functor $H$ equips $\Bb$ with a binoidal structure. Because of the double category equivalence $\F$ 
we have that $(\Bb,H)$ is unital with unit $I$ and invertible vertical strict transformations $\lambda: H(I,-)\Rightarrow\Id$ 
and $\rho: H(-,I)\Rightarrow\Id$ 
if and only if $(\Bb,\ot)$ is unital with unit $I$ and invertible vertical strict transformations $\tilde\lambda: I\ot-\Rightarrow\Id$ 
and $\tilde\rho: -\ot I\Rightarrow\Id$. Moreover, the first two triangles in \deref{premon} for  
$(\Bb,H, \theta=(\alpha_{\bullet, -,-}, \alpha_{-,\bullet,-}, \alpha_{-,-,\bullet}), \lambda, \rho)$ correspond to the analogous 
triangle connecting $\Sigma, \tilde\lambda$ and $\tilde\rho$ for $(\Bb,\ot, \Sigma, \tilde\lambda, \tilde\rho)$, the second two triangles in {\em loc.cit.} for $(\Bb,H)$ correspond to the analogous triangle connecting $\Sigma,\tilde\lambda_A\ot-$ and 
$\tilde\lambda_{A\ot-}$, while the last two triangles therein correspond to the analogous triangle connecting 
$\Sigma,-\ot\tilde\rho_B$ and $\tilde\rho_{-\ot B}$.  
We finally obtain that $(\Bb,H)$ is premonoidal if and only if $(\Bb,\ot)$ is monoidal (\deref{Shul}). 
\qed\end{proof}


\bigskip

Recall from \coref{funny yields binoidal} that pseudodouble funny functors induce a binoidal structure on double categories $\Bb$ and 
equip all 1v-cells of $\Bb$ with a central structure. Though, the fact that funny functors stem from an inner-hom in which there is no naturality for horizontal transformations, makes it impossible to upgrade a possible premonoidal structure on $\Bb$ into a monoidal one 
(at least in the way studied in this section). Namely, the induction of a pseudodouble functor $\ot:\Bb\times\Bb\to\Bb$ in \equref{F} 
relies on the existence of 2-cells $(k,K)$ in the structure of a pseudodouble quasi-functor 
as it is shown in \equref{natur 2-cells}. However, precisely these 2-cells do not exist in a structure of funny functors.

\subsubsection{2-categorical case}

Similarly to the above proof, which relies on the double category equivalence \equref{F}, one proves its 2-categorical analogue, 
which relies on the 2-categorical version of \equref{F}. (The latter result is the pseudofunctor version of the result 
\cite[Theorem 5.3]{FMS} for lax functors with $\A=\B=\C$; in the introduction of \cite[Section 4]{Fem:Bif} we explained that the 
2-categories $\Dist(\B,\B,\B)$ of \cite{FMS} and $q\x\Lax_{hop}(\Bb\times\Bb,\Bb)$ of \cite{Fem:Bif} are equal.) With analogous 
notations as above one has:

\begin{thm} \thlabel{2-cat H-ot}
Let $H:\B\times\B\to\B$ be a quasi-functor of two variables (from \cite{Gray}) on a 2-category $\B$. Then $(\B, H)$ is a purely central premonoidal bicategory if and only if $(\B,\ot)$ is a monoidal bicategory. 
\end{thm}

Another way to see that the above result holds is to think of the double category equivalence \equref{F} and the double category 
equivalence functors analogous to $\F_3$ and $\F_4$ from \thref{H equiv ot}, and consider the underlying horizontal 2-categories therein. 

\medskip

In \cite[Proposition 6]{HF} it is stated that all the  
structure 2-cells $f\ltimes-\vert_g, -\rtimes g\vert_f, \alpha_{f,B,C}, \\ \alpha_{A,g,C}, \alpha_{A,B,h}, \lambda_f, \rho_f$ and 
the modification components $p_{A,B,C,D}, m_{A,B}, l_{A,B}, r_{A,B}$ of a premonoi\-dal structure of $\B$ live in $\C_p(\B)$. 
It is readily seen that an analogue of \thref{lr central} holds for bicategories. We have that the binoidal structure of a premonoidal bicategory $\B$ comes from a quasi-functor if and only if $\B$ is purely central. In this case in 
particular all 2-cells in $\B$ are central. This confirms the above claim from \cite{HF}.


\medskip


On the other hand, from the point of view of premonoidal double categories we have:

\begin{prop} \prlabel{H(B)} 
Let $\Bb$ be a premonoidal double category so that its binoidal structure comes from a pseudodouble quasi-functor $H:\Bb\times\Bb\to\Bb$ 
and assume that its associativity and unity constraints are liftable vertical transformations. 
Then: \vspace{-0,2cm}
\begin{enumerate}
\item $\ul{\HH(\Bb)}$ is a monoidal bicategory with monoidal structure $\HH(\ot)=\HH\F(H)$, \vspace{-0,2cm}
\item the structure 2-cells $f\ltimes-\vert_g, -\rtimes g\vert_f, \alpha_{f,B,C}, \alpha_{A,g,C}, \alpha_{A,B,h}, \lambda_f, \rho_f$ 
are central in $\Bb$, 
\item the above 2-cells $\omega$ induce central 2-cells $\hat\omega$ in $\ul{\HH(\Bb)}$, and moreover the induced modification components 
$p_{A,B,C,D}, m_{A,B}, l_{A,B}, r_{A,B}$ in $\ul{\HH(\Bb)}$ are central. 
\end{enumerate}
\end{prop}

\begin{proof}
From \thref{premon-mon} we know that $(\Bb, \ot, \Sigma, I)$ is a monoidal double category, where $\ot=\F(H)$ 
and $\F$ is the double 
category equivalence \equref{F}. Then by \thref{Shulman} we have that $(\ul{\HH(\Bb)},\HH(\ot))$ is a monoidal bicategory. 
The second part follows by \thref{lr central}. (In partiular, 
$(f\ltimes-\vert_g)\ltimes-$ is a modification 
according to \thref{lr central} 
coming from the structure of the ternary quasi-functor $H(H\times 1)$.) 
The first claim of 3. holds by \leref{hat of central 2-cells}, and the second one holds by \prref{essence}. 
\qed\end{proof}

\subsection{Central $n$-noidal structures} \sslabel{leading to}

In the spirit of \ssref{binoidal str}, where we established a 1-1 correspondence between pseudodouble quasi-functors and purely/left/right  central binoidal structures on a double category $\Bb$, we may extend this correspondence to more than two variables. Namely, 
observe that pseudodouble quasi-functors of three and more variables, and their transformations, from \deref{4-ary} come down to 
collections of pseudodouble quasi-functors of two variables and their transformations satisfying certain compatibilities. 
We can take a {\em mutatis mutandi} versions of \deref{cub} -- \deref{4-ary} to introduce purely central, respectively left or right central 
{\em n-noidal} structures 
and their vertical transformations, determining relations between structural centrality transforms on central cells from $\Bb$ in different variables. 

Let $3\x\noidal_{pc}^{st}(\Bb)$ and $4\x\noidal_{pc}^{st}(\Bb)$ denote the ternary and 4-ary analogue of 
the category $\Binoidal_{pc}^{st}(\Bb)$ from \equref{bin-st}. In particular, objects of $3\x\noidal_{pc}^{st}(\Bb)$ are 3-noidal structures given by three pairs of purely central binoidal structures (\deref{pc binoidal}) 
$(\{A,B\ltimes,-\}, \{A,-,\rtimes C\}), \,\, (\{-,B,\rtimes C\}, \{A\ltimes,B,-\}), \,\, (\{A\ltimes,-,C\}, \{-,\rtimes B,C\})$ 
for objects $A,B,C\in\Bb$ that agree on objects and such that the centrality structural 2-cells of these three binoidal structures obey 7 axioms in the style of \deref{cub} (the axiom in point $ii)$ corresponding to the triple $(u,v,z)$ is now trivial). 
Morphisms in $3\x\noidal_{pc}^{st}(\Bb)$ are vertical strict 
transformations of purely central 3-noidal structures, they consist of three one-variable vertical strict transformations that obey 9 axioms 
(analogous to the first 9 of the 12 axioms of \deref{vst ternary quasi} - the last 3 axioms trivially hold because the 2-cells $(u,U)$ are now trivial). We then obtain isomorphisms and equivalences of categories 
$$
q_3\x\operatorname{Ps}_{hop}^{st}(\Bb\times\Bb\times\Bb,\Bb) \iso 3\x\noidal_{pc}^{st}(\Bb)\simeq 
\operatorname{Ps}_{hop}(\Bb\times\Bb\times\Bb,\Bb)
$$
and
$$q_4\x\operatorname{Ps}_{hop}^{st}(\Bb\times\Bb\times\Bb\times\Bb,\Bb) \iso 4\x\noidal_{pc}^{st}(\Bb)
\simeq \operatorname{Ps}_{hop}(\Bb\times\Bb\times\Bb\times\Bb,\Bb),$$
extending \thref{H equiv ot}. 

\medskip


One may also consider the $n$-noidal versions of $\Binoidal_{lc}^{st}(\Bb)$ and $\Binoidal_{rc}^{st}(\Bb)$ from \ssref{binoidal str}. 
They are isomorphic to the category of purely central $n$-noidal structures. 

\medskip

For later use let us list the axioms defining left and right central 3-noidal structures. 
Observe that the 8 axioms in \deref{cub} we labeled by the triples 
$$(f,g,h), \, (u,v,z), \, (f,v,h), \, (f,g,z), \, (u,g,h), \, (u,v,h), \, (u,g,z), \, (f,v,z).$$ 
The axioms that the objects of $3\x\noidal_{lc}^{st}(\Bb)$ obey are analogous to them, whereby the axiom corresponding to $(u,v,z)$ 
is now trivial. The remaining 7 axioms we may write in an allusive way as 
$$(f\ltimes,g\ltimes,h\ltimes), \, (f\ltimes,v\ltimes,h\ltimes), \, (f\ltimes,g\ltimes,z\ltimes), \, (u\ltimes,g\ltimes,h\ltimes), 
\, (u\ltimes,v\ltimes,h\ltimes), \, (u\ltimes,g\ltimes,z\ltimes), \, (f\ltimes,v\ltimes,z\ltimes).$$ 
Similarly, the 7 axioms that the objects of $3\x\noidal_{rc}^{st}(\Bb)$ fulfill we may write as \vspace{0,16cm} \\
\hspace{0,8cm} \axiom{$\rtimes f,\rtimes g,\rtimes h$}, $(\rtimes f,\rtimes v,\rtimes h), \, (\rtimes f,\rtimes g,\rtimes z)$,  
\axiom{$\rtimes u,\rtimes g,\rtimes h$}, $(\rtimes u,\rtimes v,\rtimes h), \, (\rtimes u,\rtimes g,\rtimes z), \, 
(\rtimes f,\rtimes v,\rtimes z)$. 

\bigskip

In \thref{premon-mon} we proved a 1-1 correspondence between a premonoidal double category structure on $(\Bb,H)$ coming from a 
pseudodouble quasi-functor $H$ and a monoidal double category structure $(\Bb,\ot)$. The proof relied on the double category equivalence \equref{F} and equivalences of categories from \thref{H equiv ot}. Joining to these the above equivalences to 
$n$-noidal structures and \thref{summary}, analogously to the proofs of \coref{assoc iso} and \thref{premon-mon} we get:

\begin{thm} \thlabel{gen-1} 
The following are equivalent: 
\begin{enumerate} \vspace{-0,14cm}
\item there is a pseudodouble quasi-functor $H$ with trivial 2-cells $(u,U)$ so that $(\Bb, H)$ is a purely central premonoidal double category; \vspace{-0,14cm}
\item there is a purely central (resp. left or right central) premonoidal structure $(\Bb, \ltimes,\rtimes)$ (in the sense of 
\deref{pc binoidal} and \deref{lr-central binoidal}) with trivial 2-cells $u\ltimes-\vert_v$ (resp. $-\rtimes v\vert_u$); \vspace{-0,14cm}
\item there is a purely central premonoidal structure $(\Bb, \ltimes,\rtimes, L_0, R_0)$ with trivial 2-cells $u\ltimes-\vert_v$; 
\vspace{-0,14cm}  
\item there is a monoidal double category structure $(\Bb, \ot_r)$. 
\end{enumerate}
\end{thm}

The pair $L_0, R_0$, {\em i.e.} their quasi-functors, in point 3. are related to the first point by setting $H=H_{L_0}$ and determining 
$H_{R_0}$ according to \rmref{discuss L0,R0}. 
The index $r$ in $\ot_r$ in point 4. is allusive to {\em right}. Namely, in \equref{P on 2-cells} we defined $P=\ot$ in one of two possible ways. 
We will comment more on this in \equref{Del-r}. Similarly, analogously to \thref{2-cat H-ot} one has:

\begin{thm} \thlabel{2-cat pc-ot}
Let $\B$ be a 2-category. There is a purely central (resp. left or right central) binoidal bicategory structure $(\B, \ltimes, \rtimes)$ 
making $\B$ a premonoidal bicategory if and only if there is a monoidal bicategory structure $(\B,\ot)$. 
\end{thm}

A purely (resp. left or right) central bicategory here has an analogous meaning as a purely (resp. left or right) central double category.

\section{Centers of a premonoidal double category } \selabel{center premon}

In \seref{cubical1} we introduced pure, left and right center double categories and a center double category for a binoidal double category 
$\Bb$. In this section we introduce pure center and center for a {\em premonoidal} double category and study their relation to their bicategorical counterparts. The reason why we study separately binoidal and premonoidal case is that for centers of a premonoidal double category two additional requirements are made. One is that the associativity constraint $\alpha$ (and $\lambda, \rho$) should live in the center, and the other one is that with a functorial choice of centrality structures present in a center the appropriate axioms of the 24 axioms for $\alpha$ should hold.

\subsection{Pure center and center of a premonoidal double category}  

When $\Bb$ is premonoidal and there is a pseudodouble functor $\F\colon\Bb'\to\Pseudo_{ps}(\Bb', \Bb')$ (equivalently, a pseudodouble quasi-functor $H:\Bb'\times\Bb'\to\Bb'$) for a double subcategory $\Bb'\subseteq\Bb$, 
one has due to \prref{quasi-24 ax} that the 24 axioms for $\alpha$ hold true in $\Bb'$, where the centrality structures are induced by $H$. On the other hand, $\F$ ({\em i.e.} $H$) determines a pure/left/right center double category 
$\Zz_\bullet(\Bb')$, and we have that 12 axioms corresponding to the center $\Zz_\bullet(\Bb')$ hold. 
Namely, for the left center $\Zz_l(\Bb')$ the upper half, for the right center $\Zz_r(\Bb')$ the lower half of the 24 axioms applies, whereas 
for the pure center $\Zz_p(\Bb')$ all 24 axioms hold, but they collapse into 12 axioms. The totality of the 24 axioms says that the horizontal transformations 
$$(f\ltimes -)\rtimes C, \,\,\,  (f\ltimes B)\ltimes-, \,\,\,  (A\ltimes g)\ltimes -,  \quad\quad\quad 
-\rtimes(g'\rtimes C), \,\,\, -\rtimes(B\rtimes h), \,\,\, A\ltimes(-\rtimes h),$$
for $f,g$ left central and $g',h$ right central, are related via the $\alpha$'s to the horizontal transformations 
$$f\ltimes (-\rtimes C), \,\,\,  f\ltimes (B\ltimes-), \,\,\,  A\ltimes (g\ltimes -) \quad\text{and}\quad 
(-\rtimes g')\rtimes C, \,\,\, (-\rtimes B)\rtimes h, \,\,\, (A\ltimes-)\rtimes h$$
(recall \prref{horiz 12 axioms}), and similarly for vertical transformations (Table \ref{table:3}). In particular, the left and right centrality structures for the 1h-cells $f\rtimes B, \,\, A\ltimes g, \,\, g'\rtimes C, \,\, B\rtimes h$, 
are given via the middle 8 axioms of Table \ref{table:2} (recall \leref{middle 8 ax}). 
In terms of the quasi-functor $H$, {\em e.g.} the left centrality structure for $A\ltimes g=(g,A)$ is given by $(-,(g,A))=H(A\ltimes g,-)=
H(H(A,g),-)$, which is related via $\alpha$ to $H(A,H(g,-))=((-,g),A)$, 
and similarly for the other pairs of transformations. 
 
\medskip

Now that we have settled this, there is one more thing we should make sure regarding center double categories of a premonoidal double category 
$\Bb$. Namely, the structural 1v-cells $\alpha_{A,B,C}, \lambda_A, \rho_A$ for all $A,B,C\in\Bb$ are central, so they should live both 
in the left and the right center. In $\Zz_l(\Bb')$ we have that all 1- and 2-cells are left central and we also have pseudodouble functors 
$A\ltimes-$ and $-\rtimes B$. (Observe that $-\rtimes B\vert_A=A\ltimes-\vert_B$, while $-\rtimes B$ when evaluated at a 1- or a 2-cell 
$a\in\Bb'$ because of left centrality of $a$ is equal to $a\ltimes-\vert_B$.) Now, for a 1-cell $b\in\Bb'$ to say that it is right central we should give a corresponding transformation $-\rtimes b$. It is given by evaluation at a 0-cell $A$ and a suitable 1-cell $a$. However, we do 
not have such cells in $\Zz_l(\Bb')$. Thus, for premonoidal $(\Bb, H)$ the only center double category that we can consider is a pure one. 

\medskip

On the other hand, given pseudodouble functors $L_0,R_0\colon\Bb'\to\Pseudo_{ps}(\Bb', \Bb')$ as in \equref{L0,R0} inducing a 
central binoidal structure $(\Bb', \ltimes,\rtimes, L_0, R_0)$, similarly as in \prref{quasi-24 ax} associativity for $A\ltimes=
H_{L_0}(A,-)$ satisfies the upper half and associativity for $-\rtimes A=H_{R_0}(A,-)$ satisfies the lower half of the 24 axioms for $\alpha$. 
In the former case left centrality structures are induced from $L_0$, and in the latter right centrality structures are induced from $R_0$. 
For premonoidal $(\Bb, L_0, R_0)$ that is (purely) central binoidal we can consider a (pure) center double category $\Zz(\Bb')$ 
(resp. $\Zz_p(\Bb')$). 

\bigskip

\begin{prop} \prlabel{24 ax for L0,R0}
An associative binoidal structure $(\Bb, \ltimes,\rtimes, L_0, R_0)$ making $\Bb$ premonidal  
satisfies the 24 axioms, where the left and right centrality structures are induced by $L_0$ and $R_0$, respectively. 
\end{prop}

\begin{defn} \delabel{Z_p premon}
Let $(\Bb, \ltimes,\rtimes)$ be a premonoidal double category. 
\begin{itemize}
\item Supose that $(\Bb, \ltimes,\rtimes)$ is purely central binoidal on a double subcategory $\Bb'\subseteq\Bb$ and that the structural 
1v-cells $\alpha_{A,B,C}, \lambda_A, \rho_A$ for all $A,B,C\in\Bb$ live in $\Bb'$. In this case we say that $\Zz_p(\Bb')$ from 
\deref{pure center} is a pure center double category for $\Bb$. 
\item Supose that $(\Bb, \ltimes,\rtimes, L_0, R_0)$ is (purely) central binoidal on a double subcategory $\Bb'\subseteq\Bb$ and that the structural 1v-cells $\alpha_{A,B,C}, \lambda_A, \rho_A$ for all $A,B,C\in\Bb$ live in $\Bb'$. In this case we say that $\Zz(\Bb')$ 
(resp. $\Zz_p(\Bb')$) from \ssref{not diff} is a (pure) center double category for $\Bb$. 
\end{itemize}
\end{defn}

\bigskip

Under certain conditions there is a pseudoretraction to the pseudodouble functors \equref{center-pseud}. 
Assume that $\Bb$ is premonoidal and purely central (so that its binoidal structure comes from a pseudodouble quasi-functor $H$), that the structure vertical transformations $\alpha$'s, $\lambda$ and $\rho$ are liftable, and that the 1v-cells $\lambda_A, \rho_A$ are identities. 
Define $E:\Pseudo_{ps}(\Bb,\Bb)\to\Zz_p(\Bb)$ 
by $E(F)=F(I), \, E(\theta_h)=(\theta_h(I), \theta_h(I)\ltimes-, -\rtimes \theta_h(I))$, \, $E(\theta_v)=
(\theta_v(I), \theta_v(I)\ltimes-, -\rtimes \theta_v(I))$ and $E(\Theta)=(\Theta(I), \Theta(I)\ltimes-, -\rtimes \Theta(I))$,
for 0-, 1h-, 1v- and 2-cells $F, \theta_h, \theta_v$ and $\Theta$ from $\Pseudo_{ps}(\Bb,\Bb)$, respectively. To see that $E$ is a pseudoretraction of $L$ (and similarly $R$) observe that $\tilde\rho_f\ltimes-: (\tilde\rho_{A'}(f\ltimes I))\ltimes-\Rrightarrow (f\tilde\rho_A)\ltimes-$ is an invertible modification, where $\tilde\rho$ presents the horizontal pseudonatural equivalence obtained via 
\prref{lifting 1v to equiv} from $\rho$. This explains the level of 1h-cells, for 1v-cells the idea is similar (use directly $\rho$). 
For the level of 2-cells use the axiom \axiomref{h.o.t.-5}. 

\medskip

Whereas in the monoidal setting and even in the bicategorical setting the Drinfel'd center of a bicategory $\B$ is the {\em category} 
$\Z_{Dr}(\B)\iso\PsNat(\Id_\B, \Id_\B)$ as we commented before, in the premonoidal setting the center of a bicategory and a double category  is two-dimensional. From the point of view of a double category $\Bb$ on the endo-hom category $\HH(\Zz_p(\Bb))(I,I)$ one has an (induced) functor with an (induced) pseudoretraction under the above conditions. Still, it does not get to be an isomorphism. This illustrates the difference between premonoidal and monoidal centers. 

\subsubsection{Pure center on the underlying 2-category} 

As we mentioned, in \cite[Definition 21]{HF} a bicategory of pure maps $\C_p(\B)$ was defined for a premonoidal bicategory $\B$. 
We now want to study the relation between a pure center double category $\Zz_p(\Bb)$ of a premonoidal double category $\Bb$ and 
$\C_p(\ul{\HH(\Bb)})$, the bicategory of pure maps of the underlying premonoidal 2-category $\ul{\HH(\Bb)}$ of $\Bb$ (recall  
\prref{premon isofib new}). 

Let us spell out the definition of $\C_p(\B)$ with a small adaptation. 
According to \prref{quasi-24 ax} and \prref{24 ax for L0,R0} six bicategorical axioms analogous to 
$\alpha_{f,-,C}, \alpha_{f,B,-}, \\ \alpha_{A,g,-}, \alpha_{-,g',C}, 
\alpha_{-,B,h}, \alpha_{A,-,h}$ hold automatically in $\C_p(\B)$, whereas 
in \cite[Definition 21]{HF} it is required that $\alpha_{f,-,C}$ and $\alpha_{A,-,h}$ be modifications. 

\begin{defn}
Let $\B$ be a premonoidal bicategory and let $\B'\subseteq\B$ be a sub-bicategory such that 
\begin{enumerate}
\item there are pseudofunctors $L_0, R_0:\B'\to \operatorname{Pseudo}_{ps}(\B',\B')$, \vspace{-0,2cm}
\item structural 1-cells $\alpha_{A,B,C}, \lambda_A, \rho_A$ live in $\B'$ for all $A,B,C\in\B$,  \vspace{-0,2cm}
\item bicategorical axioms analogous to the left-most 6 axioms of Table \ref{table:2} hold, i.e. \\ $\alpha_{f,-,C}, \alpha_{f,B,-}, 
\alpha_{A,g,-}, \alpha_{-,g',C}, \alpha_{-,B,h}, \alpha_{A,-,h}$ are modifications, and \vspace{-0,2cm}
\item $f\ltimes-\vert_g=(-\rtimes g\vert_f)^{-1}$ hold for all 1-cells $f,g$. 
\end{enumerate}
\end{defn}

Here $\operatorname{Pseudo}_{ps}(\B',\B')$ is the 2-categorical analogue of $\Pseudo_{ps}(\Bb',\Bb')$. 
(It is smaller than the underlying 2-category of $\Pseudo_{ps}(\Bb',\Bb')$.) 
Also, the modifications in 3. are meant in the bicategorical sense: the associativity constraints involved are given by globular 2-cells.  

\smallskip

Observe that structural cells of the underlying horizontal bicategory $\HH(\Zz_p(\Bb))$ of the center pseudodouble category $\Zz_p(\Bb)$ of a double category $\Bb$ 
possess non-trivial 1v-cells of $\Bb$ and henceforth $\HH(\Zz_p(\Bb))$ is properly larger than the bicategory of pure maps $\C_p(\ul{\HH(\Bb)})$ of the underlying horizontal 2-category of $\Bb$. Namely, 1h-cells of $\Zz_p(\Bb)$ are triples $(f,f\ltimes-,-\rtimes f)$ where the horizontal structure transformations possess non-globular 2-cells $f\ltimes-\vert_u$ and $-\rtimes f\vert_u$. 


Let $\Bb$ be a premonoidal double category such that a premonoidal 2-category $\B$ is its underlying horizontal premonoidal 2-category, 
$\ul{\HH(\Bb)}=\B$. It is immediate to see that the points 1. and 4. for $\C_p(\B)$ follow from the definition of a pure center double category $\Zz_p(\Bb)$ from our \deref{Z_p premon}. Under the assumptions of \prref{premon isofib new} the point 2. above also holds. It remains to show that under the same conditions the bicategorical axioms for $\alpha_{f,-,C}$ and 
$\alpha_{A,-,h}$ in 3. (and for the remaining four modifications) follow from their double categorical analogues. 
(Rrecall that the axioms for $\alpha$ hold by premonoidality and existence of functorial assignment of centrality structures $L_0, R_0$, 
\prref{24 ax for L0,R0}; for $\ul{\HH(\Bb)}=\B$ we consider only the listed six axioms.) 

\medskip

In the pure center double category $\Zz_p(\Bb)$ we have that the axiom \axiomref{$(f\ltimes,g,C)$} 
(out of the 24 from Appendix A) holds, and in \prref{lifting 1v to equiv} we have a 
way how vertical strict transformations $\alpha_{-,B,C}$ and $\alpha_{A,-,C}$ determine horizontal pseudonatural transformations 
that we will denote here by $\alpha_{-,B,C}^1$ and $\alpha_{A,-,C}^2$. Then to see that the first axiom in item 3. holds we need to check the equality 
$$
\scalebox{0.86}{
\bfig
\putmorphism(-680,700)(1,0)[` `((A,g),C)]{780}1a
\putmorphism(120,700)(1,0)[` `((f,B'),C)]{800}1a
\putmorphism(-680,700)(0,-1)[` `=]{450}1l
\putmorphism(930,700)(0,-1)[` `=]{450}1r
\putmorphism(-680,250)(1,0)[` `((f,B),C)]{500}1a
\put(-130,480){\fbox{$(f\ltimes-\vert_g,C)$}}
 \putmorphism(-150,250)(1,0)[` `=]{500}1a
\put(-80,20){\fbox{$\eta_{A', B,C}$}}
\putmorphism(380,250)(0,-1)[\phantom{Y_2}` `]{450}1l
\putmorphism(410,390)(0,-1)[\phantom{Y_2}` `\alpha_{A', B,C}]{450}0l
\putmorphism(930,250)(0,-1)[\phantom{Y_2}` `]{450}1r
\putmorphism(940,390)(0,-1)[\phantom{Y_2}` `\alpha_{A', B',C}]{450}0l
\putmorphism(350,250)(1,0)[``((A',g),C)]{600}1a
 \putmorphism(950,250)(1,0)[` `\alpha^2_{A', B',C}]{550}1a
 \putmorphism(420,-200)(1,0)[``(A',(g,C))]{500}1b
 \putmorphism(1000,-200)(1,0)[``=]{500}1b
\putmorphism(1500,250)(0,-1)[``=]{450}1r
\put(510,20){\fbox{$\alpha_{A',g,C}$}}
\put(1100,30){\fbox{$\Epsilon_{A', B',C}$}}
\putmorphism(-150,-200)(1,0)[` `\alpha^2_{A', B,C}]{520}1a
\putmorphism(-150,250)(0,-1)[``=]{450}1l
\putmorphism(-150,-220)(0,-1)[` `]{450}1l
\putmorphism(-160,-280)(0,-1)[` `\alpha_{A', B,C}]{450}0r
\putmorphism(380,-220)(0,-1)[` `=]{450}1r
\put(-100,-370){\fbox{$\Epsilon_{A', B,C}$}}
\putmorphism(-150,-650)(1,0)[` `=]{520}1b

\putmorphism(-680,-200)(1,0)[` `((f,B),C)]{520}1a
\putmorphism(-680,-650)(1,0)[` `(f,(B,C))]{520}1b
\put(-550,-370){\fbox{$\alpha_{f, B,C}$}}

\putmorphism(-680,-220)(0,-1)[` `]{450}1l
\putmorphism(-700,-280)(0,-1)[` `\alpha_{A, B,C}]{450}0r

\putmorphism(-1200,-200)(1,0)[` `=]{520}1a
\putmorphism(-1200,-220)(0,-1)[` `=]{450}1l
\putmorphism(-1200,-650)(1,0)[` `\alpha^1_{A, B,C}]{520}1b
\put(-1080,-370){\fbox{$\eta_{A, B,C}$}}
\efig} 
$$

$$
=
\scalebox{0.86}{
\bfig
 \putmorphism(-150,250)(1,0)[` `=]{500}1a
\put(-80,20){\fbox{$\eta_{A, B',C}$}}
\putmorphism(380,250)(0,-1)[\phantom{Y_2}` `]{450}1l
\putmorphism(410,390)(0,-1)[\phantom{Y_2}` `\alpha_{A, B',C}]{450}0l
\putmorphism(930,250)(0,-1)[\phantom{Y_2}` `]{450}1r
\putmorphism(940,390)(0,-1)[\phantom{Y_2}` `\alpha_{A', B',C}]{450}0l
\putmorphism(350,250)(1,0)[``((f,B'),C)]{600}1a
 \putmorphism(950,250)(1,0)[` `\alpha^2_{A', B',C}]{550}1a
 \putmorphism(420,-200)(1,0)[``(f,(B',C))]{500}1b
 \putmorphism(1000,-200)(1,0)[``=]{500}1b
\putmorphism(1500,250)(0,-1)[``=]{450}1r
\put(510,20){\fbox{$\alpha_{f,B',C}$}}
\put(1100,30){\fbox{$\Epsilon_{A', B',C}$}}
\putmorphism(-150,-200)(1,0)[` `\alpha^2_{A, B',C}]{520}1a
\putmorphism(-150,250)(0,-1)[``=]{450}1l
\putmorphism(-150,-220)(0,-1)[` `]{450}1l
\putmorphism(-160,-280)(0,-1)[` `\alpha_{A, B',C}]{450}0r
\putmorphism(380,-220)(0,-1)[` `=]{450}1r
\put(-100,-370){\fbox{$\Epsilon_{A, B',C}$}}
\putmorphism(-150,-650)(1,0)[` `=]{520}1b
\putmorphism(-680,-200)(1,0)[` `((A,g),C)]{520}1a
\putmorphism(-680,-650)(1,0)[` `(A,(g,C))]{520}1b
\put(-550,-370){\fbox{$\alpha_{A, g,C}$}}
\putmorphism(-680,-220)(0,-1)[` `]{450}1l
\putmorphism(-700,-280)(0,-1)[` `\alpha_{A, B,C}]{450}0r
\putmorphism(-1200,-200)(1,0)[` `=]{520}1a
\putmorphism(-1200,-220)(0,-1)[` `=]{450}1l
\putmorphism(-1200,-650)(1,0)[` `\alpha^1_{A, B,C}]{520}1b
\put(-1080,-370){\fbox{$\eta_{A, B,C}$}}
 \putmorphism(420,-650)(1,0)[``(f,(B',C))]{500}1b
\putmorphism(-680,-650)(0,-1)[` `=]{450}1l
\putmorphism(930,-650)(0,-1)[` `=]{450}1r
\putmorphism(-680,-1100)(1,0)[` `((A,g),C)]{780}1b
\putmorphism(120,-1100)(1,0)[` .`((f,B'),C)]{800}1b
\put(-130,-880){\fbox{$f\ltimes-\vert_{(g,C)}$}}
\efig}
$$
But this is true precisely because of the axiom \axiomref{$(f\ltimes,g,C)$} and since $\frac{\eta}{\Epsilon}=\Id$ (the vertical 
composition of 2-cells). The bicategorical axiom for $\alpha_{A,-,h}$ follows similarly from the axiom \axiomref{$(A,g,\rtimes h)$}. 

\medskip

The relation between the double categorical and bicategorical pure center we can formalize as follows. 
Let $\Zz_p(\Bb)_{hm}$ denote the pure center double category in whose 2-cells $(a, a\ltimes-, -\rtimes a)$ the modifications $a\ltimes-$ and 
$-\rtimes a$ are {\em horizontal} (recall axioms \axiomref{m.ho.-1} and \axiomref{m.ho.-2}). 
For $\Bb$ such that its associativity and unity constraints $\alpha,\lambda,\rho$ are liftable vertical transformations we know from 
\prref{premon isofib new} that $\ul{\HH(\Bb)}$ is the underlying premonoidal 2-category with the induced  
pseudonatural equivalences $\hat\alpha,\hat\lambda,\hat\rho$ induced from $\Bb$ according to \prref{lifting 1v to equiv}. Similarly, 
consider $\ul{\HH(\Zz_p(\Bb)_{hm})}$ equipped with the same $\hat\alpha,\hat\lambda,\hat\rho$. There are pseudofunctors 
\begin{equation} \eqlabel{hor_p}
hor_p:\ul{\HH(\Zz_p(\Bb)_{hm})}\to C_p(\ul{\HH(\Bb)})
\end{equation} 
and
\begin{equation} \eqlabel{hor*}
hor:\ul{\HH(\Zz(\Bb)_{hm})}\to \Z(\ul{\HH(\Bb)})
\end{equation} 
that are identity on objects and send 1- and 2-cells $(b,b\ltimes-,-\rtimes b)$ to $(b,\HH(b\ltimes-),\HH(-\rtimes b))$ (for the meaning of 
$\HH$ on transformations and modifications recall \rmref{mod} and \rmref{H trans}).

\subsection{Pseudodouble functors on the center}

We have introduced two pure center double categories $\Zz_p(\Bb)$, which differ only in the way how centrality structures are given in them. 
We also introduced center double category $\Zz(\Bb)$ supposing existence of two pseudodouble functors $L_0, R_0$. 
Based on the results that we proved in \seref{24 axioms}, independently on which centrality structure one works with we may conclude:

\begin{prop} \prlabel{psfuns on centers}
In a premonoidal double category $\Bb$ 
for every $A,B\in\Bb$ there are pseudodouble functors 
$$A\ltimes-, \,\, -\rtimes B: \Zz_p(\Bb)\to\Zz_p(\Bb),$$
$$A\ltimes-, \,\, -\rtimes B: \Zz(\Bb)\to\Zz(\Bb).$$
\end{prop}

\begin{proof}
The two pseudodouble functors in both cases are clearly defined on 0-cells. 
In \prref{horiz 12 axioms} we showed that if $f$ is a left central 1h-cell, then $A\ltimes f, f\rtimes B$ are left central, 
and if $f$ is right central, then $A\ltimes f, f\rtimes B$ are right central. Then centrality of $f$ implies centrality of 
$A\ltimes f$ and $f\rtimes B$. 
Similar claims for 1v-cells we have by \leref{basic} and as explained around \equref{four vertical tr}, because of the assumption 
that the three $\alpha$'s are strict. 
Given a left central 2-cell $\sigma$ it is straightforward to check that $A\ltimes\sigma$ and $\sigma\rtimes B$ are left central, 
and if $\sigma$ is right central, that so are $A\ltimes\sigma$ and $\sigma\rtimes B$ (one uses the axiom \axiomref{v.l.t.\x 5}). 
In particular, for $\sigma$ central, both $A\ltimes\sigma$ and $\sigma\rtimes B$ are central. 
The pseudodouble functor property follows from centrality of the domain 1-cells. 
\qed\end{proof}

Up to the difference that the bicategorical center from \cite[Definition 20]{HF}, \cite[Definition 7]{HF1} does not fix centrality structures, whereas our double categorical center fixes them, the above result is an extension of \cite[Proposition 5]{HF}, 
\cite[Proposition 8]{HF1} to double categories in the following sense. 
We defined pseudonatural transformations \equref{four center} appearing in the image of the above pseudodouble functors $A\ltimes-, \,\, -\rtimes B$ on $\Zz(\Bb)$ and $\Zz_p(\Bb)$ for each $A,B\in\Bb$ via the middle 8 axioms in Table \ref{table:2}. 
Their horizontally globular 2-cell components can be expressed from the left 4 axioms of those 8 axioms, and they have the following form 
\begin{equation} \eqlabel{despejar}
(f\rtimes B)\ltimes-\vert_h=
\scalebox{0.86}{
\bfig
\putmorphism(-150,500)(1,0)[``((A,B),h)]{600}1a
 \putmorphism(480,500)(1,0)[` `((f,B),C')]{640}1a

 \putmorphism(-150,50)(1,0)[``(A,(B,h))]{600}1a
 \putmorphism(450,50)(1,0)[``(f,(B,C'))]{640}1a

\putmorphism(-180,500)(0,-1)[\phantom{Y_2}``\alpha_{A,B,C}]{450}1l
\putmorphism(450,500)(0,-1)[\phantom{Y_2}``]{450}1r
\putmorphism(250,500)(0,-1)[\phantom{Y_2}``\alpha_{A,B,C'}]{450}0r
\putmorphism(1100,500)(0,-1)[\phantom{Y_2}``\alpha_{A',B,C'}]{450}1r
\put(-90,280){\fbox{$\alpha_{A,B,h}$}}
\put(680,280){\fbox{$\alpha_{f,B,C'}$}}

\putmorphism(-150,-400)(1,0)[` `(f,(B,C))]{640}1b
 \putmorphism(490,-400)(1,0)[` `(A',(B,h))]{640}1b

\putmorphism(-180,50)(0,-1)[\phantom{Y_2}``=]{450}1l
\putmorphism(1100,50)(0,-1)[\phantom{Y_3}``=]{450}1r
\put(270,-200){\fbox{$f\ltimes-\vert_{B\ltimes h}$}}

\putmorphism(-180,-400)(0,-1)[\phantom{Y_2}``\alpha_{A,B,C}^{-1}]{450}1l
\putmorphism(450,-400)(0,-1)[\phantom{Y_3}``]{450}1l
\putmorphism(550,-400)(0,-1)[\phantom{Y_3}``\alpha_{A',B,C}^{-1}]{450}0l
\putmorphism(1100,-400)(0,-1)[\phantom{Y_3}``\alpha_{A',B,C'}^{-1}]{450}1r
\put(-90,-670){\fbox{$\alpha_{f,B,C}^{-1}$}}
\put(680,-670){\fbox{$\alpha_{A',B,h}^{-1}$}}
\putmorphism(-150,-920)(1,0)[``((f,B),C)]{600}1a
 \putmorphism(480,-920)(1,0)[` `((A',B),h)]{640}1a
\efig}
\end{equation}
and similarly for the remaining three transformations. When we pass the square-formed component 2-cells of the three vertical liftable  transformations $\alpha$ (and their inverses) appearing in these expressions to horizontal pseudonatural equivalences by 
\prref{lifting 1v to equiv} part b), we obtain exactly the same component 2-cells for the pseudonatural transformations induced by the images by $L_A$ and $R_B$ of central 1-cells $f$ in $\ul{\HH(\Bb)}$, as defined in \cite{HF1}. To see this, 
%
we will use string diagrams for 2-categories, recall \ssref{strings}. 
We will write $\beta_{A,B,C}$ for horizontal companions of the 1v-cell components $\alpha_{A,B,C}$ of the associativity 
constraint $\alpha$, 
and similarly we will denote by $\beta_{f,B,C}, \beta_{A,g,C},\beta_{A,B,h}$ the horizontally globular 2-cell components of horizontal pseudodouble transformations obtained from the 2-cell components $\alpha_{f,B,C}, \alpha_{A,g,C},\alpha_{A,B,h}$ of the corresponding 
vertical liftable transformations. We consider the 1h-cells $\beta_{A,B,C}$ as 
adjoint equivalences and denote by $\gamma_{A,B,C}$ their quasi-inverses. Suppressing the symbols $\ltimes, \rtimes$, we get 
\begin{equation} \eqlabel{hor form f} \hspace{-1,6cm}
(f\rtimes B)\ltimes -\vert_h=
\gbeg{4}{9}
\got{1}{(AB)h} \got{4}{(fB)C'} \gnl
\gcl{2} \gvac{1} \gcl{1} \gwdb{3} \gnl
\gvac{1} \gvac{1} \glmptb \gnot{\beta_{f,B,C'}} \gcmptb \grmp \gcl{5} \gnl
\glmptb \gcmp \gnot{\hspace{-0,8cm}\beta_{A,B,h}} \grmptb  \gcl{1} \gnl
\gcl{1} \gvac{1} \glmptb \gnot{f\ltimes\x\vert_{Bh}} \gcmptb \grmp \gnl 
\glmptb \gcmp \gnot{\hspace{-0,8cm}\gamma_{f,B,C}} \grmptb  \gcl{1} \gnl
\gcl{2} \gvac{1} \glmptb \gnot{\gamma_{A',B,h}} \gcmptb \grmp \gnl
\gvac{2} \gcl{1} \gwev{3} \gnl
\gvac{1} \gob{-1}{(fB)C} \gob{6}{(A'B)h}
\gend
\qquad\qquad\qquad
(A\ltimes g)\ltimes -\vert_h=
\gbeg{3}{11}
\got{1}{(AB)h} \got{4}{(Ag)C'} \gnl
\gcl{2} \gvac{1} \gcl{1} \gwdb{3} \gnl
\gvac{2} \glmptb \gnot{\beta_{A,g,C'}} \gcmptb \grmp \gcl{7} \gnl
\glmptb \gcmp \gnot{\hspace{-0,8cm}\beta_{A,B,h}} \grmptb  \gcl{1} \gnl
\gcl{3} \gvac{1} \gmu \gnl
\gvac{1} \glmp \gcmp \gnot{\hspace{-0,34cm}1(g\ltimes-\vert_h)} \gcmp \grmp \gnl
\gvac{2} \gcmu \gnl
\glmptb \gcmp \gnot{\hspace{-0,8cm}\gamma_{A,g,C}} \grmptb  \gcl{1} \gnl
\gcl{2} \glmptb \gnot{\gamma_{A,B',h}} \gcmptb \grmp \gnl
\gvac{2} \gcl{1} \gwev{3} \gnl
\gob{1}{(Ag)C} \gob{4}{(AB')h}
\gend
\end{equation}
where $\gbeg{2}{1}
\gdb \gnl
\gend$ and $\gbeg{2}{1}
\gev \gnl
\gend$ present $\eta$ and its inverse for $\alpha_{A',B,C'}$ and $\alpha_{A,B',C'}$, respectively, and 
$$\hspace{-2cm}
-\rtimes(g'\rtimes C)\vert_f=
\gbeg{3}{9}
\gvac{1} \got{-1}{f(BC)} \gvac{3} \got{5}{A'(g'C)} \gnl
\gcl{1} \gvac{1} \gdb \gvac{1} \gcl{1} \gnl
\glmptb \gnot{(\gamma_{f,B,C})^\bullet} \gcmptb \grmptb 
  \glmptb \gnot{\gamma_{A',g',C}} \gcmptb \grmptb \gnl
\gcl{1} \gvac{1} \gmu \gvac{1} \gcl{1} \gnl
\gcl{1} \glmp \gcmp \gnot{\hspace{-0,36cm}(-\rtimes g'\vert_f)1_C} \gcmp \grmp \gcl{2}  \gnl
\gcl{1} \gvac{1} \gcmu \gnl
\glmptb \gnot{(\beta_{A,g',C})^\bullet}\gcmp  \grmptb \glmptb \gnot{\beta_{f,B',C}} \gcmp \grmptb \gnl
\gcl{1} \gvac{1} \gev \gvac{1} \gcl{1} \gnl
\gvac{1} \gob{-1}{A(gC)} \gvac{3} \gob{5}{f(B'C)}
\gend
\qquad\qquad\qquad\quad
-\rtimes(B\ltimes h)\vert_f=
\gbeg{3}{9}
\gvac{1} \got{-1}{f(BC)} \gvac{3} \got{5}{A'(Bh)} \gnl
\gcl{1} \gvac{1} \gdb \gvac{1} \gcl{1} \gnl
\glmptb \gnot{(\gamma_{f,B,C})^\bullet} \gcmptb \grmptb 
  \glmptb \gnot{\gamma_{A',B,h}} \gcmptb \grmptb \gnl
\gcl{1} \gvac{1} \gmu \gvac{1} \gcl{1} \gnl
\gcl{1} \gvac{1} \glmp \gnot{-\rtimes h\vert_{fB}} \gcmp \grmp \gcl{2}  \gnl
\gcl{1} \gvac{1} \gcmu \gnl
\glmptb \gnot{(\beta_{A,B,h})^\bullet}\gcmp  \grmptb \glmptb \gnot{\beta_{f,B,C'}} \gcmp \grmptb \gnl
\gcl{1} \gvac{1} \gev \gvac{1} \gcl{1} \gnl
\gvac{1} \gob{-1}{A(Bh)} \gvac{3} \gob{5}{f(BC)}
\gend
$$
where $\gbeg{2}{1}
\gdb \gnl
\gend$ and $\gbeg{2}{1}
\gev \gnl
\gend$ in the left-hand side present $\Epsilon^{-1}$ for $\beta_{A',B,C}$ and $\Epsilon$ for $\beta_{A,B',C}$, respectively, and 
in the right-hand side $\Epsilon^{-1}$ for $\beta_{A',B,C}$ and $\Epsilon$ for $\beta_{A,B,C'}$, respectively. 
 
These four transformations are precisely pseudonatural transformations in a bicategory used in \cite[Proposition 8]{HF1} 
(two of them were not explicitly given). 

As we argued before, the first two of the above four definitions in $\ul{\HH(\Bb)}$ originate from the fact that 
$\alpha_{f,B,-}$ and $\alpha_{A,g,-}$ are modifications 
between vertical pseudonatural transformations and horizontal pseudonatural transformations. 
On the other hand, the adjoint equivalence of the companion 1h-cells makes these two definitions in $\ul{\HH(\Bb)}$ 
equivalent to the facts that 
$$\alpha_{f,X,-}:\frac{(f\ltimes X)\ltimes-}{\alpha_{A',X,-}}\Rrightarrow\frac{\alpha_{A,X,-}}{f\ltimes(X\ltimes-)}$$ 
and
$$\alpha_{X,f,-}:\frac{(X\rtimes f)\ltimes-}{\alpha_{X,A',-}}\Rrightarrow\frac{\alpha_{X,A,-}}{X\rtimes (f\ltimes-)}$$ 
are modifications of (2-categorical) pseudonatural transformations, respectively, as recorded in \cite[Proposition 5]{HF}. 

\bigskip

We may summarize the above findings in the proposition below, in which we omit the proof for 2-cells. 
Let $\LL_A,\Rr_B$ denote the restrictions of the above pseudodouble functors $A\ltimes-, \,\, -\rtimes B$ on $\Zz(\Bb)$ to 
$\ul{\HH(\Zz(\Bb)_{hm})}$, and let $L_A,R_A:\Z(\ul{\HH(\Bb)})\to\Z(\ul{\HH(\Bb)})$ denote the bicategorical counterparts of $A\ltimes-$ and 
$-\rtimes B$ analogous to those from \cite[Proposition 8]{HF2}. Let $hor:\ul{\HH(\Zz(\Bb)_{hm})}\to\Z(\ul{\HH(\Bb)})$ be as in \equref{hor*} and let 
$hor^*:\I(\ul{\HH(\Zz(\Bb)_{hm})})\to\Z(\ul{\HH(\Bb)})$ be the pseudofunctor that sends the double categorical structural data to the corresponding underlying bicategorical data in the following way. 
Here $\I(\ul{\HH(\Zz(\Bb)_{hm})})$ is the image of $\LL_A$ or $\Rr_B$, with small abuse of notation. Both pseudofunctors are identities on 
0-cells and send a 2-cell $(a,a\ltimes-,-\rtimes a)$ to the 2-cell $(a,\HH(a\ltimes-),\HH(-\rtimes a))$.  
Now, whereas $hor$ sends a 1-cell $(f,f\ltimes-,-\rtimes f)$ to the 1-cell $(f,\HH(f\ltimes-),\HH(-\rtimes f))$,  
for $hor^*((f\rtimes B,(f\rtimes B)\ltimes-,-\rtimes (f\rtimes B)))$  
instead of taking $\HH((f\rtimes B)\ltimes-)$ that is given by \equref{despejar} in which 1v-cells appear, 
one takes $\widehat{\HH}((f\rtimes B)\ltimes-)$ given by \equref{hor form f}. One proceeds similarly for $-\rtimes (f\rtimes B)$ 
and the images by $\LL_A$.

\begin{prop} \prlabel{uses mates}
Let $\Bb$ be a premonoidal double category 
whose structure vertical strict transformations are liftable. For every $A,B\in\Bb$ the following diagram of pseudofunctors between bicategories commutes 
$$\scalebox{0.86}{
\bfig
\putmorphism(-150,50)(1,0)[\ul{\HH(\Zz(\Bb)_{hm})}` \I(\ul{\HH(\Zz(\Bb)_{hm})})`\LL_A,\Rr_B]{1000}1a
\putmorphism(-150,-290)(1,0)[\Z(\ul{\HH(\Bb)})`\Z(\ul{\HH(\Bb)}). `L_A,R_A]{1000}1b
\putmorphism(-150,50)(0,-1)[\phantom{Y_2}``hor]{320}1l
\putmorphism(850,50)(0,-1)[\phantom{Y_2}``hor^*]{320}1r
\efig}$$
The analogous holds when $\Z(\ul{\HH(\Bb)})$ is replaced by $\Z_0(\ul{\HH(\Bb)})$, and also for pure (double categorical and bicategorical) centers. 
\end{prop}

\section{Monoidality of a pure center double category} \selabel{mon pure cen}

While in a premonoidal category the center subcategory is monoidal, in \cite[Section 6.5]{HF} it is explained that for a 
premonoidal bicategory $\B$ its center bicategory $\Z(\B)$ is not monoidal. 
Though, in Theorem 3 the authors prove that the corresponding bicategory of pure maps $\C_p(\B)$ is a monoidal bicategory. 

The analogous occurrences we find in the double categorical setting (the fact that our center $\Zz(\Bb)$ fixes centrality structures 
contrarily to $\Z(\B)$ does not play a role in this issue). 
In \thref{premon-mon} we proved that if the binoidal, hence premonoidal, structure of $\Bb$ comes from a quasi-functor, then $\Bb$ is purely central and monoidal. This will help us prove in \ssref{mon pc} that the pure center double category $\Zz_p^{st}(\Bb)$ is monoidal. 
Regarding a center double category $(\Zz(\Bb),\ltimes,\rtimes,L_0,R_0)$, 
the two pseudodouble functors $L_0,R_0$ give rise to two pseudodouble quasi-functors so that the binoidal structure is given by 
$A\ltimes-=A\ltimes_{L_0}-$ and $-\rtimes B=-\rtimes_{R_0} B$, whereby the 2-cells $f\ltimes_{L_0}-\vert_g$ and 
$-\rtimes_{R_0} g\vert_f$ are unrelated. In the construction of a monoidal product $\ot$ in \equref{F} we set its lax pseudodouble functor structure using $K\ltimes-\vert_k$ in \equref{natur 2-cells}, and its colax structure then uses $(K\ltimes-\vert_k)^{-1}$, 
which equals $-\rtimes k\vert_K$ because of pure centrality. 
Though, if we want to equip $\Zz(\Bb)$ with a monoidal structure, we need to give a monoidal product also on central 1h-cells 
$f$ and $g$, so we should be able to relate the 2-cell components $f\ltimes_{L_0}-\vert_g$ and $-\rtimes_{R_0} g\vert_f$, and also 
the latter with the lax and colax structure of $\ot$. Thus not having pure centrality hinders 
center double categories from being monoidal. 


In the first subsection we study necessary conditions to have a pseudodouble functor 
$\ot:\Zz_p^{st}(\Bb)\times\Zz_p^{st}(\Bb)\to\Zz_p^{st}(\Bb)$, where $\Zz_p^{st}(\Bb)$ denotes a pure center stemming from a 
pseudodouble quasi-functor whose 2-cells $(u,U)$ are trivial.

\subsection{Candidate for a monoidal product on a pure center}

Let $\Bb$ be a premonoidal double category. From \thref{gen-1} we have that in order to obtain a candidate for a monoidal product $\ot$ on 
$\Zz_p(\Bb)$, we should define a pseudodouble quasi-functor for $\Zz_p(\Bb)$ whose 2-cells $(u,U)$ are trivial. Let $\Zz_p^{st}(\Bb)$ 
denote a pure center double category stemming from a pseudodouble quasi-functor $H:\Bb\times\Bb\to\Bb$ whose 2-cells $H(U,u)$ are trivial.   
We are going to show that $H$ 
extends to a pseudodouble quasi-functor 
$\crta H:\Zz_p^{st}(\Bb)\times \Zz_p^{st}(\Bb)\to \Zz_p^{st}(\Bb)$ so that the 2-cells $\crta H(\crta U,\crta u)$ are trivial, for 
1v-cells $\crta U, \crta u\in\Zz_p^{st}(\Bb)$. 
Otherwise stated, we will use the fact that $H$ induces two pseudodouble functors as in \prref{psfuns on centers} and equipping 
$\Zz_p^{st}(\Bb)$ with thus obtained binoidal structure, we will prove that all 1- and 2-cells in $\Zz_p^{st}(\Bb)$ are left and right central in $\Zz_p^{st}(\Bb)$. The necessary four identities on 2-cells for $\crta H$ will then hold because they hold for $H$, including triviality of the 2-cells $\crta H(\crta U,\crta u)$.

\begin{rem} \rmlabel{central in center}
Before we pursue, let us see what it means for 1- and 2-cells in $\Zz_p^{st}(\Bb)$ to be left/right central in $\Zz_p^{st}(\Bb)$.
We divide this study in cases: when cells in $\Zz_l^{st}(\Bb)$ are left and right central in $\Zz_l^{st}(\Bb)$, and when cells in 
$\Zz_r^{st}(\Bb)$ are left and right central in $\Zz_r^{st}(\Bb)$. 
For example for $\Zz_r^{st}(\Bb)$: take a 1h-cell $(k, -\rtimes k)$ in $\Zz_r^{st}(\Bb)$, in order for $(k, -\rtimes k)$ to be right central in 
$\Zz_r(\Bb)$, we need a 
pseudonatural transformation $-\rtimes k: -\rtimes A\Rightarrow -\rtimes A': \Zz_r^{st}(\Bb)\to\Zz_r^{st}(\Bb)$. In particular, for a right central 1h-cell $K$ and a right central 1v-cell $U$ we want structure 2-cells $-\rtimes k\vert_K$ and $-\rtimes k\vert_U$ to live in 
$\Zz_r^{st}(\Bb)$. What we do have, by right centrality of $k$ in $\Bb$, is that $-\rtimes k: -\rtimes A\Rightarrow-\rtimes A': \Bb\to\Bb$ is a pseudonatural transformation, meaning that for $K$ and $U$ living (only) in $\Bb$ the 2-cells $-\rtimes k\vert_K$ and $-\rtimes k\vert_U$ live (only) in $\Bb$. The transformation axioms for $-\rtimes k$ in $\Zz_r(\Bb)$ will follow by the transformation axioms for $-\rtimes k$ in $\Bb$. 

For 1v-cells the situation is similar. For a 2-cell $(\sigma, -\rtimes\sigma)$ in $\Zz_r^{st}(\Bb)$, between 1h- and 1v-cells which are right central in $\Zz_r^{st}(\Bb)$, to be right central in $\Zz_r^{st}(\Bb)$ it means that $A\rtimes\sigma$ are right central 2-cells in $\Bb$ and the two modification axioms from \deref{central-2} for $\sigma$ hold for any 1h-cell $K$ and 1v-cell $U$ in $\Zz_r^{st}(\Bb)$. By right centrality of $\sigma$ in $\Bb$ we only have that $A\rtimes\sigma$ are 2-cells in $\Bb$ and $K$ and $U$ should be from $\Bb$. Though, by right center version of \prref{psfuns on centers} we already have that $A\rtimes\sigma$ are right central 2-cells, and the desired modification axioms hold in particular for $K$ and $U$ from $\Zz_r^{st}(\Bb)$. The rest of the cases for 2-cells similarly follows because of 
\prref{psfuns on centers}, so it suffices to prove only centralities of 1-cells. 
\end{rem}

\medskip

Let us consider the case that 1-cells in $\Zz_r^{st}(\Bb)$ should be right central in $\Zz_r^{st}(\Bb)$. 
For 1h-cells $(K, -\rtimes K), (k, -\rtimes k)$ and 1v-cells $(U, -\rtimes U), (u, -\rtimes u)$ in $\Zz_r^{st}(\Bb)$ we have 2-cells 
$-\rtimes k\vert_K, \,\,\, -\rtimes k\vert_U, \,\,\, -\rtimes u\vert_K, \,\,\, -\rtimes u\vert_U$ in $\Bb$. 
To simplify notation, let us denote them by $K\rtimes k, \,\,\, U\rtimes k, \,\,\, K\rtimes u, \,\,\, U\rtimes u$, respectively.   
Observe that the 2-cells $U\rtimes u$ are identities. 
By \rmref{central in center} we should prove that $K\rtimes k, \,\,\, U\rtimes k, \,\,\, K\rtimes u, \,\,\, U\rtimes u$ are right central 
2-cells in $\Bb$ (observe that here we evaluated right centrality transformations at right central 1-cells, as the Remark requires). 
To prove right centrality of $K\rtimes k$, one should prove that $-\rtimes(K\rtimes k)$ is a modification in $\Bb$. At this point we switch to the notation 
$(f,g,h):(A,B,C)\to(A',B', C')$ for 1h-cells and $(u,v,z):(A,B,C)\to(\tilde A, \tilde B,\tilde C)$ for 1v-cells, as in 
\seref{24 axioms}. To prove right centrality of $g\rtimes h$, one should prove that $-\rtimes(g\rtimes h)$ is a modification, 
between horizontal transformations $\frac{[-\rtimes gC\vert\Id]}{[\Id\vert-\rtimes B'h]}$ and 
$\frac{[-\rtimes Bh\vert\Id]}{[\Id\vert-\rtimes gC']}$ (and similarly for right centrality of $v\rtimes h$ and $g\rtimes z$). To write out the two modification conditions for $-\rtimes(g\rtimes h), \,\,\, -\rtimes(v\rtimes h), \,\,\, -\rtimes(g\rtimes z)$ one applies \leref{vert comp hor.ps.tr.} for vertical composition of horizontal transformations. 
In the case of $-\rtimes(g\rtimes h)$ the two modification conditions read 
\begin{equation} \eqlabel{modif-centr}
\end{equation}
$$\scalebox{0.86}{
\bfig
\putmorphism(-150,900)(1,0)[``f(BC)]{600}1a
 \putmorphism(380,900)(1,0)[\phantom{F(B)}` `A'(gC)]{680}1a
\put(310,660){\fbox{$-\rtimes gC\vert_f$}}

 \putmorphism(-150,450)(1,0)[``A(gC)]{600}1a
 \putmorphism(380,450)(1,0)[\phantom{F(B)}` `f(B'C)]{680}1a
 \putmorphism(1000,450)(1,0)[\phantom{F(B)}` `A'(B'h)]{550}1a 

\putmorphism(-180,900)(0,-1)[\phantom{Y_2}``=]{450}1r
\putmorphism(1060,900)(0,-1)[\phantom{Y_2}``=]{450}1r
\put(880,220){\fbox{$-\rtimes B'h\vert_f$}}
  \putmorphism(-150,0)(1,0)[` `A(gC)]{600}1a
\putmorphism(380,0)(1,0)[\phantom{F(A)}`  `A(B'h)]{650}1a
 \putmorphism(970,0)(1,0)[\phantom{F(A)}` ` f(B'C')]{550}1a

\putmorphism(450,450)(0,-1)[\phantom{Y_2}``=]{450}1l
\putmorphism(1570,450)(0,-1)[\phantom{Y_2}``=]{450}1r

\putmorphism(-180,0)(0,-1)[\phantom{Y_2}``=]{450}1l
\putmorphism(1040,0)(0,-1)[\phantom{Y_2}``=]{450}1r
\put(210,-230){\fbox{$-\rtimes (-\rtimes h\vert_g)\vert_A$}}
 \putmorphism(-220,-450)(1,0)[\phantom{F(B)}` `A(Bh)]{650}1a
 \putmorphism(380,-450)(1,0)[\phantom{F(B)}` `A(gC')]{650}1a
\efig}
\quad=
\scalebox{0.86}{
\bfig
  \putmorphism(450,900)(1,0)[` `A'(gC)]{680}1a %
\putmorphism(1000,900)(1,0)[\phantom{F(A)}`  `A'(B'h)]{680}1a %
\putmorphism(450,900)(0,-1)[\phantom{Y_2}``=]{450}1l
\putmorphism(1700,900)(0,-1)[\phantom{Y_2}``=]{450}1r
\put(810,680){\fbox{$-\rtimes (-\rtimes h\vert_g)\vert_{A'}$}}
 \putmorphism(380,450)(1,0)[\phantom{F(B)}` `A'(Bh)]{650}1a %
 \putmorphism(960,450)(1,0)[\phantom{F(B)}` `A'(gC')]{680}1a %

\putmorphism(-150,450)(1,0)[``f(BC)]{600}1a
\put(270,220){\fbox{$-\rtimes Bh\vert_f$}} %
\putmorphism(450,0)(0,-1)[\phantom{Y_2}``=]{450}1l
\putmorphism(1700,0)(0,-1)[\phantom{Y_2}``=]{450}1r

 \putmorphism(-150,0)(1,0)[``A(Bh)]{600}1a %
 \putmorphism(380,0)(1,0)[\phantom{F(B)}` `f(BC')]{680}1a %
 \putmorphism(1000,0)(1,0)[\phantom{F(B)}` `A'(gC')]{650}1a %

\putmorphism(-180,450)(0,-1)[\phantom{Y_2}``=]{450}1r
\putmorphism(1060,450)(0,-1)[\phantom{Y_2}``=]{450}1r
\put(880,-230){\fbox{$-\rtimes gC'\vert_f$}}
\putmorphism(450,900)(0,-1)[\phantom{Y_2}``=]{450}1l
\putmorphism(1700,900)(0,-1)[\phantom{Y_2}``=]{450}1r
 \putmorphism(380,-450)(1,0)[\phantom{F(B)}` `A(gC')]{650}1a
 \putmorphism(960,-450)(1,0)[\phantom{F(B)}` `f(B'C')]{680}1a
\efig}
$$ 
and 
$$\scalebox{0.86}{
\bfig
 \putmorphism(-150,450)(1,0)[``A(gC)]{600}1a
 \putmorphism(380,450)(1,0)[\phantom{F(B)}` `A(B'h)]{680}1a
\put(-140,220){\fbox{$u(gC)$}}
\put(610,220){\fbox{$u(B'h)$}}
  \putmorphism(-150,0)(1,0)[` `\tilde A(gC)]{600}1a
\putmorphism(380,0)(1,0)[\phantom{F(A)}`  `\tilde A(B'h)]{650}1a

\putmorphism(-180,450)(0,-1)[\phantom{Y_2}``u(BC)]{450}1l
\putmorphism(490,450)(0,-1)[\phantom{Y_2}``]{450}1l
\putmorphism(510,450)(0,-1)[\phantom{Y_2}``u(B'C)]{450}0l
\putmorphism(1040,450)(0,-1)[\phantom{Y_2}``u(B'C')]{450}1r

\putmorphism(-180,0)(0,-1)[\phantom{Y_2}``=]{450}1l
\putmorphism(1040,0)(0,-1)[\phantom{Y_2}``=]{450}1r
\put(210,-210){\fbox{$-\rtimes (-\rtimes h\vert_g)\vert_{\tilde A}$}}
 \putmorphism(-220,-450)(1,0)[\phantom{F(B)}` `\tilde A(Bh)]{650}1a
 \putmorphism(380,-450)(1,0)[\phantom{F(B)}` `\tilde A(gC')]{650}1a
\efig}
\quad=
\scalebox{0.86}{
\bfig
  \putmorphism(450,450)(1,0)[` `A(gC)]{600}1a %
\putmorphism(980,450)(1,0)[\phantom{F(A)}`  `A(B'h)]{680}1a %
\putmorphism(450,450)(0,-1)[\phantom{Y_2}``=]{450}1l
\putmorphism(1700,450)(0,-1)[\phantom{Y_2}``=]{450}1r
\put(780,220){\fbox{$-\rtimes (-\rtimes h\vert_g)\vert_A$}}

\putmorphism(450,0)(0,-1)[\phantom{Y_2}``u(BC)]{450}1l
\putmorphism(1000,0)(0,-1)[\phantom{Y_2}``]{450}1l
\putmorphism(980,0)(0,-1)[\phantom{Y_2}``u(BC')]{450}0r
\putmorphism(1700,0)(0,-1)[\phantom{F(A)}` ` u(B'C')]{450}1r 

 \putmorphism(360,0)(1,0)[\phantom{F(B)}` `A(Bh)]{650}1a %
 \putmorphism(920,0)(1,0)[\phantom{F(B)}` `A(gC')]{780}1a %
\put(610,-230){\fbox{$u(Bh)$}}
\put(1340,-230){\fbox{$u(gC')$}}
 \putmorphism(360,-450)(1,0)[\phantom{F(B)}` `\tilde A(Bh)]{650}1a %
 \putmorphism(920,-450)(1,0)[\phantom{F(B)}` `\tilde A(gC')]{780}1a %

\efig}
$$ 
for {\em any} 1h-cell $f$ and 1v-cell $u$ in $\Bb$, 
whereby the 2-cell $-\rtimes (-\rtimes h\vert_g)\vert_A$ is defined via the axiom \axiomref{$(A,g,\rtimes h)$} that relates it to the 2-cell 
$-\rtimes h\vert_{A\ltimes g}$. Similarly one uses the other corresponding axioms out of the 12 ones holding in $\Zz_r^{st}(\Bb)$ for the other three modification candidates. Observe that the above two conditions are formally precisely the axioms \axiomref{$\rtimes f,\rtimes g,\rtimes h$} and \axiomref{$\rtimes u,\rtimes g,\rtimes h$} of a 3-noidal right central structure, the only difference is that here $f$ and $u$ are not required to be central. To highlight this difference we write the above two axioms as $(f,\ul{\rtimes g,\rtimes h})$ and 
$(u,\ul{\rtimes g,\rtimes h})$. (As a matter of fact, for simplicity we assumed that $H$ is defined on the whole $\Bb$, so {\em any} 
$f$ and $u$ are central. Nevertheless, we proceed this meticulous style of the study for the record.) 
We conclude that the eight axioms for $-\rtimes(g\rtimes h), -\rtimes(v\rtimes h), -\rtimes(g\rtimes z)$ and $-\rtimes(v\rtimes z)$ to be modifications are the (seven) axioms from the first four lines of Table \ref{table:4} and they correspond (up to the just explained difference) to the axioms of a 3-noidal right central structure - the eighth modification axiom 
regarding $-\rtimes(v\rtimes z)\vert_u$ trivially holds, as it consists only of identity 2-cells.


\medskip

In the case where 
1-cells in $\Zz_l^{st}(\Bb)$ should be left central in $\Zz_l^{st}(\Bb)$ 
we start from left central 1-cells $K,k,U,u$, {\em i.e.} 1-cells in $\Zz_l^{st}(\Bb)$. We obtain 2-cell components 
$K\ltimes -\vert_k=K\ltimes k, \,\,\,  K\ltimes -\vert_u=K\ltimes u, \,\,\, U\ltimes -\vert_k=U\ltimes k, \,\,\, 
U\ltimes -\vert_u=U\ltimes u$ in $\Bb$, 
of the horizontal and vertical pseudonatural transformations provided by the left centrality of $K$ and $U$. In order for these 2-cells to live in $\Zz_l(\Bb)$ we obtain similar seven conditions to the above ones to have the corresponding modification conditions. These are 
the seven axioms from the last four lines of Table \ref{table:4}, and correspond to the axioms of a 3-noidal left central structure.

%

It remains to study when 1-cells in $\Zz_l^{st}(\Bb)$ are right central in $\Zz_l^{st}(\Bb)$, and analogously when 1-cells in $\Zz_r^{st}(\Bb)$ are left central in $\Zz_r(\Bb)$. Take $(k, k\ltimes-)\in\Zz_l^{st}(\Bb)$, for it to be right central in $\Zz_l^{st}(\Bb)$ we need a 
horizontal pseudonatural transformation $-\rtimes k:-\rtimes A\Rightarrow -\rtimes A'$ in $\Zz_l^{st}(\Bb)$. This means that for a 1h-cell 
$(K, K\ltimes-)\in\Zz_l^{st}(\Bb)$ we want a 2-cell $-\rtimes k\vert_K$ in $\Zz_l^{st}(\Bb)$, {\em i.e.} a modification 
$(-\rtimes k\vert_K)\ltimes-$ in $\Bb$, and a similar discussion goes for the rest of combinations of 1h- and 1v-cells. 
By the quasi-functor property for $H$ we have $-\rtimes k\vert_K=(K\ltimes-\vert_k)^{-1}$. So, we want $(K\ltimes-\vert_k)^{-1}\ltimes-$ 
to be a modification in $\Bb$. It is a modifications, if so is $(K\ltimes-\vert_k)\ltimes-$, so for right centrality in $\Zz_l^{st}(\Bb)$ 
we need the last seven axioms in Table \ref{table:4} to hold. 
Similarly, for cells in $\Zz_r^{st}(\Bb)$ to be left central in $\Zz_r^{st}(\Bb)$ we want {\em e.g.} $-\rtimes(K\ltimes-\vert_k)=
-\rtimes(-\rtimes k\vert_K)^{-1}$ to be a modification in $\Bb$ for $K,k$ right central. It is a modifications if so is 
$-\rtimes(-\rtimes k\vert_K)$. Thus for left centrality in $\Zz_r^{st}(\Bb)$ we want the first seven axioms in Table \ref{table:4} to hold.

\begin{table}[H]
\begin{center}
\begin{tabular}{ c c c c } 
For 1-cells & modification & \hspace{0,2cm} centrality for 2-cell  & \hspace{0,2cm}  the two axioms \\ [0.5ex]
\hline
right central $K,k$ & $-\rtimes(K\rtimes k)$ & right central $K\rtimes k$ &  \axiom{$(f,\ul{\rtimes g,\rtimes h})$} and 
\axiom{$(u,\ul{\rtimes g,\rtimes h})$} \\ [1ex] 
right central $U,k$ & $-\rtimes(U\rtimes k)$ & right central $U\rtimes k$ &  \axiom{$(f,\ul{\rtimes v,\rtimes h})$} and 
\axiom{$(u,\ul{\rtimes v,\rtimes h})$} \\ [1ex] 
right central $K,u$ & $-\rtimes(K\rtimes u)$ & right central $K\rtimes u$ &  \axiom{$(f,\ul{\rtimes g,\rtimes z})$} and 
\axiom{$(u,\ul{\rtimes g,\rtimes z})$} \\ [1ex] 
right central $U,u$ & $-\rtimes(U\rtimes u)$ & right central $U\rtimes u$ &  \axiom{$(f,\ul{\rtimes v,\rtimes z})$} \\ [1ex] 
\hline
left central $K,k$ & $(K\ltimes k)\ltimes-$ & left central $K\ltimes k$ & \axiom{$(\ul{f\ltimes ,g\ltimes },h)$} and 
\axiom{$(\ul{f\ltimes ,g\ltimes },z)$} \\ [1ex] 
left central $U,k$ & $(U\ltimes k)\ltimes-$ & left central $U\ltimes k$ & \axiom{$(\ul{u\ltimes ,g\ltimes },h)$} and 
\axiom{$(\ul{u\ltimes ,g\ltimes },z)$} \\ [1ex] 
left central $K,u$ & $(K\ltimes u)\ltimes-$ & left central $K\ltimes u$ & \axiom{$(\ul{f\ltimes ,v\ltimes },h)$} and 
\axiom{$(\ul{f\ltimes ,v\ltimes },z)$} \\ [1ex] 
left central $U,u$ & $(U\ltimes u)\ltimes-$ & left central $U\ltimes u$ & \axiom{$(\ul{u\ltimes ,v},h)$} \\ [1ex] 
\end{tabular}
\caption{Modification {\em i.e.} centrality conditions}
\label{table:4}
\end{center}
\end{table} \vspace{-1,24cm}
\begin{center}
(needed to have a ps.d. functor $\Zz^{st}_p(\Bb)\times\Zz^{st}_p(\Bb)\to\Zz^{st}_p(\Bb)$)
\end{center}

\smallskip

This finishes the analysis of what it takes to have a desired pseudodouble quasi-functor $\crta H:\Zz_p^{st}(\Bb)\times\Zz_p^{st}(\Bb)\to\Zz_p^{st}(\Bb)$. Hence, in order to have a pseudodouble functor $\ot:\Zz_p^{st}(\Bb)\times\Zz_p^{st}(\Bb)\to\Zz_p^{st}(\Bb)$, 
given 1h-cells $(k,k\ltimes-,-\rtimes k), 
(K,K\ltimes-,-\rtimes K)$ and 1v-cells $(u,u\ltimes-,-\rtimes u), (U,U\ltimes-,-\rtimes U)$ in $\Zz^{st}_p(\Bb)$ we want the pairs 
$(k,K), (k,U), (u,K), (u,U)$ to satisfy the above two groups of seven axioms, whereby either their left or right centralities are 
employed. 
As we said above, since $H$ exists on the whole $\Bb$,
any three 1-cells that should obey those $7\cdot 2$ axioms are both left and right central. 
Then the above axioms coincide with the seven right, respectively left 3-noidality axioms.  
However, depending on the axiom different structure 2-cells are to be employed. So for example, the axiom 
\axiomref{$(f,\ul{\rtimes g,\rtimes h})$} is about the structure 2-cell $-\rtimes(-\rtimes h\vert_g)\vert_f$, while the axiom 
\axiomref{$(\ul{f\ltimes ,g\ltimes },h)$} is about the structure 2-cell $(f\ltimes-\vert_g)\ltimes-\vert_h$. 
As a matter of fact, all these axioms hold automatically by pure centrality of $\Bb$, as we discuss next.  

\smallskip

When we proved right centrality of $K\rtimes k$ ({\em i.e.} of $g\rtimes h$), in the first  modification axiom for $-\rtimes(g\rtimes h)$ we evaluated the horizontal transformations $\frac{[-\rtimes gC\vert\Id]}{[\Id\vert-\rtimes B'h]}$ and $\frac{[-\rtimes Bh\vert\Id]}{[\Id\vert-\rtimes gC']}$ at a 1h-cell $f$. Applying pure centrality of $\Bb$, 
we have $-\rtimes Bh\vert_f=(f\ltimes-\vert_{Bh})^{-1}$ and $-\rtimes gC'\vert_f=(f\ltimes-\vert_{gC'})^{-1}$. Then by \leref{vert comp hor.ps.tr.}  
$\frac{[-\rtimes gC\vert\Id]}{[\Id\vert-\rtimes B'h]}\vert_f=(f\ltimes-\vert_{(B'\rtimes h)(g\rtimes C)})^{-1}$ and 
$\frac{[-\rtimes Bh\vert\Id]}{[\Id\vert-\rtimes gC']}\vert_f=(f\ltimes-\vert_{(g\rtimes C')(B\rtimes h)})^{-1}$. On the other hand, 
the 2-cell $A\rtimes(g\rtimes h)$ 
can be seen as $F(a)$ with $F=A\ltimes-$ and $a=g\rtimes h$, so that the first modification condition for $-\rtimes(g\rtimes h)$ actually holds by \axiomref{h.o.t.-5} of $f\ltimes-$. 
Similarly, in the second modification axiom because of pure centrality consider $[u\rtimes(gC)\vert u\rtimes(B'h)]$ as 
$(u\ltimes (B'h)(gC))^{-1}$, apply \axiomref{v.l.t.\x 5} for $u\ltimes-$ and the condition holds again by \axiomref{h.o.t.-5} but now of $-\rtimes h$. 
By similar reasoning the rest of the centrality {\em i.e.} modification proofs work.  

\medskip

Now we have, similarly as in \cite[Proposition 5.6]{Fem:Bif} and \equref{P on 2-cells}, a pseudodouble functor 
\begin{equation} \eqlabel{Del-r}
\Del_r: \Zz^{st}_p(\Bb)\times\Zz^{st}_p(\Bb)\to\Zz^{st}_p(\Bb)
\end{equation}
given as follows. For any 2-cells $(\sigma, \sigma\ltimes-, -\rtimes\sigma), (\delta, \delta\ltimes-, -\rtimes\delta)$ in $\Zz^{st}_p(\Bb)$, being 
$$
\scalebox{0.86}{
\bfig
\putmorphism(-150,175)(1,0)[A` A'`K]{450}1a
\putmorphism(-150,-175)(1,0)[\tilde A`\tilde A' `\tilde K]{440}1b
\putmorphism(-170,175)(0,-1)[\phantom{Y_2}``U]{350}1l
\putmorphism(280,175)(0,-1)[\phantom{Y_2}``U']{350}1r
\put(0,-15){\fbox{$\sigma$}}
\efig}
\quad\text{and}\quad
\scalebox{0.86}{
\bfig
\putmorphism(-150,175)(1,0)[B` B'`k]{450}1a
\putmorphism(-150,-175)(1,0)[\tilde B`\tilde B' `\tilde k]{440}1b
\putmorphism(-170,175)(0,-1)[\phantom{Y_2}``u]{350}1l
\putmorphism(280,175)(0,-1)[\phantom{Y_2}``u']{350}1r
\put(0,-15){\fbox{$\delta$}}
\efig}
$$
2-cells in $\Bb$, we define 
\begin{equation} \eqlabel{sigma ot delta}
\sigma\Del_r \delta:=
\scalebox{0.86}{
\bfig

 \putmorphism(-520,50)(1,0)[A\rtimes B`A'\rtimes B `K\rtimes B]{600}1a
 \putmorphism(60,50)(1,0)[\phantom{A'\ltimes k}` A'\ltimes B' `A'\ltimes k]{670}1a

\putmorphism(70,50)(0,-1)[\phantom{Y_2}`\tilde A'\ltimes B`]{450}1r
\putmorphism(50,180)(0,-1)[``U'\rtimes B]{450}0r
\putmorphism(650,50)(0,-1)[\phantom{Y_2}`\phantom{U'\rtimes B} `U'\rtimes B']{450}1r
\putmorphism(70,-400)(0,-1)[\phantom{Y_2}``]{450}1r
\putmorphism(70,-490)(0,-1)[\phantom{Y_2}``\tilde A'\ltimes U]{450}0r
\putmorphism(650,-400)(0,-1)[\phantom{Y_2}``\tilde A'\ltimes u']{450}1r

\put(210,-600){\fbox{$\tilde A'\ltimes \delta$}}
\put(210,-210){\fbox{$U'\rtimes k$}}

\putmorphism(-470,50)(0,-1)[\phantom{Y_2}`\tilde A\rtimes B`U\rtimes B]{450}1l
\putmorphism(-470,-400)(0,-1)[\phantom{Y_2}``\tilde A\ltimes u]{450}1l
\put(-350,-180){\fbox{$\sigma\rtimes B$}}
\put(-350,-650){\fbox{$\tilde K\rtimes u$}}

 \putmorphism(-510,-400)(1,0)[\phantom{A''\ot B'}`\phantom{A''\ot B'}`\tilde K\rtimes B]{600}1a
 \putmorphism(40,-400)(1,0)[\phantom{A''\ot B'}` \phantom{\tilde K\rtimes B} `]{670}1a
 \putmorphism(40,-410)(1,0)[\phantom{A''\ot B'}` \tilde A'\ltimes B'`\tilde A'\ltimes k]{670}0a
 \putmorphism(-530,-850)(1,0)[\tilde A\ltimes\tilde B`\tilde A'\rtimes\tilde B`\tilde K\rtimes\tilde B]{600}1b
\putmorphism(80,-850)(1,0)[\phantom{\tilde K\rtimes\tilde B}` \tilde A'\ltimes \tilde B' `\tilde A'\ltimes \tilde k]{670}1b
\efig 
}
\end{equation}
which includes the definitions 
$$K\Del_r k:=(A'\ltimes k)(K\rtimes B)\quad\text{and}\quad U\Del_r u:=\frac{U\rtimes B}{\tilde A\ltimes u},$$
for 1h-cells $K:A\to A'$ and $k:B\to B'$ and 1v-cells $U:A\to\tilde A$ and $u:B\to\tilde B$ in $\Zz^{st}_p(\Bb)$ (we suppress the structures). 
Observe that $K\Del_r k$ and $U\Del_r u$ correspond to the domain 1-cells of the right-hand sided diagrams in \equref{f lr} and
\equref{v lr}. (A pseudodouble functor $\Del_l:\Zz^{st}_p(\Bb)\times\Zz^{st}_p(\Bb)\to\Zz^{st}_p(\Bb)$ can be defined similarly 
using the domain 1-cells of the left-hand sided diagrams therein.) 

In the above considerations, instead of starting with a pseudodouble quasi-functor $H$, we could have equivalently started from a purely central binoidal structure from \deref{pc binoidal} on $\Bb$ to obtain a pseudodouble functor analogous to $\Del_l$, because of \equref{bin-st}. Then we can state:

\begin{prop} \prlabel{pc-pseudofunctor}
Let $\Bb$ be purely central premonoidal 
whose purely central binoidal structure is such that the 2-cells $-\rtimes u\vert_U$ are trivial. 
Let $\Zz_p^{st}(\Bb)$ be the corresponding pure center double category. 
The binoidal structure of $\Bb$ induces a pseudodouble functor $\Zz_p^{st}(\Bb)\times\Zz_p^{st}(\Bb)\to\Zz_p^{st}(\Bb)$. 
\end{prop}

\smallskip

In the analysis of associativity of the pseudodouble functors $\Del_l$ and $\Del_r$ additional conditions in the style of those in Table 
\ref{table:4} appear (they can be expressed in terms of a 4-noidal and a 5-noidal structure, recall \ssref{leading to}). For $\Bb$ purely central those axioms are trivially fulfilled.


\subsection{Monoidality of a pure center} \sslabel{mon pc}

For $(\Bb,\ltimes,\rtimes)$ purely central and $\G:\Binoidal_{pc}^{st}(\Bb) \to \operatorname{Ps}_{vst}(\Bb\times\Bb,\Bb)$ from \equref{bin-st}, observe that $\Del_r: \Zz^{st}_p(\Bb)\times\Zz^{st}_p(\Bb)\to\Zz^{st}_p(\Bb)$ is precisely the extension of 
$\G(\ltimes,\rtimes)=\ot_r:\Bb\times\Bb\to\Bb$ (compare \equref{Del-r} and \equref{sigma ot delta} to \equref{F} and 
\equref{P on 2-cells}). 

\begin{thm} \thlabel{gen} 
For a premonoidal and purely central double category $(\Bb, \ltimes,\rtimes)$ as in \prref{pc-pseudofunctor} its pure center 
pseudodouble category $\Zz_p^{st}(\Bb)$ is monoidal with the monoidal product $\Del_r$. 
\end{thm}

\begin{proof} 
In \thref{gen-1} we proved that $(\Bb, \ot_r)$ is a monoidal double category 
and by \prref{pc-pseudofunctor} we know that $\Del_r$ is a pseudofunctor. The monoidal structure of $\Bb$ passes then to 
$\Zz_p^{st}(\Bb)$ via the double functor $Z_p:\Bb\to\Zz_p^{st}(\Bb)$.  
\qed\end{proof}

For the underlying horizontal bicategory we have: 

\begin{thm} \thlabel{gen-H(B)}
Let $(\Bb, \ltimes,\rtimes)$ be a left/right/purely central premonoidal double category so that its associativity and unity constraints are liftable vertical transformations and so that the 2-cells $-\rtimes u\vert_U$ are trivial. Then for its underlying horizontal bicategory $\HH(\Bb)$ we have: 
\begin{enumerate}
\item $\ul{\HH(\Bb)}$ is a left/right/purely central premonoidal bicategory with structure $(\HH(\ltimes,\rtimes))$;
\item $\ul{\HH(\Bb)}$ is a monoidal bicategory with structure $\HH(\ot_r)$;
\item $\ul{\HH(\Zz_p^{st}(\Bb)_{hm})}$ is a monoidal bicategory and $C_p(\HH(\Bb))$ inherits monoidal structure from it. 
\end{enumerate}
\end{thm}

\begin{proof}
By \thref{gen-1} we have that $(\Bb, \ot_r)$ is a monoidal double category. By \thref{Shulman} then $(\HH(\Bb),\HH(\ot_r))$ is a monoidal bicategory, and by \thref{2-cat pc-ot} we obtain that $(\HH(\Bb), \HH(\ltimes,\rtimes))$ is a left/right/purely central premonoidal bicategory. 
For the last statement, by \thref{gen} we know that 
$\Zz_p^{st}(\Bb)$ is a monoidal double category, then so is $\Zz_p^{st}(\Bb)_{hm}$. By \thref{Shulman} $\ul{\HH(\Zz_p^{st}(\Bb)_{hm})}$ becomes a monoidal bicategory and $C_p(\HH(\Bb))$ inherits its monoidal structure via \equref{hor_p}. 
\qed\end{proof}

\bigskip

The following diagram illustrates the claims of \thref{gen-H(B)} and \prref{H(B)}:
$$\scalebox{0.86}{
\bfig
\putmorphism(-150,50)(1,0)[(\Bb,(\ltimes,\rtimes)_{pc})` (\Bb,\ot_r)` \thref{gen-1}]{1900}1a 
\putmorphism(-150,50)(1,0)[``\longleftrightarrow]{1900}0b 
\putmorphism(-150,-280)(1,0)[(\HH(\Bb),\HH(\ltimes,\rtimes)_{pc})` (\HH(\Bb), \HH(\ot_r)). 
   `\thref{2-cat pc-ot}, \thref{2-cat H-ot}]{1900}1b 
\putmorphism(-150,-280)(1,0)[``\longleftrightarrow]{1900}0a 
\putmorphism(-150,50)(0,-1)[\phantom{Y_2}``]{300}4l
\putmorphism(1750,50)(0,-1)[\phantom{Y_2}``\cite{Shul}]{320}1r 
\efig}$$

\medskip

We can not say that every premonoidal 2-category $\B$ comes from a premonoidal double category $\Bb$ with non-trivial 1v-cells, 
in the sense that $\ul{\HH(\Bb)}=\B$ and that the 2-categorical part of the premonoidal structure of $\Bb$ is precisely the premonoidal structure of $\B$. 
For purely central premonoidal 2-categories $\B$ coming from purely central premonoidal double categories that also satisfy the conditions of our \thref{gen-H(B)}, the fact proved in \cite[Theorem 3]{HF} that the 2-category of pure maps $\C_p(\B)$ is monoidal, can be seen as a consequence of our theorem. 

\bigskip 


We finish this first part of the paper by commenting the notion of a {\em Freyd bicategory} introduced in \cite[Definition 16]{HF1}. A  Freyd bicategory is given by the following data: 
a monoidal bicategory $(\V, \ot,I)$, a premonoidal bicategory $(\B, \ltimes, \rtimes,I)$, so that $\V$ and $\B$ have the same objects and unit $I$, and an identity-on-objects 0-strict premonoidal pseudofunctor $J:\V\to\B$, which factors strictly through the center bicategory 
$\Z_0(\B)$ of $\B$ via a binoidal pseudofunctor $J_\Z$ satisfying some axioms. (Mind the difference between the center $\Z_0(\B)$ 
from \cite{HF1} without specified centrality structures and centers $\Z(\B)$ with functorial choice of such structures.)

Let $\Bb$ be a premonoidal double category so that its binoidal structure comes from a pseudodouble quasi-functor $H:\Bb\times\Bb\to\Bb$. Let $J:\Zz_p^{st}(\Bb)_{hm}\to\Bb$ be the identity-on-objects double functor sending 1- and 2-cells $(b, b\ltimes-, -\rtimes b)\mapsto b$, 
and let $J_\Z: \ul{\HH(\Zz^{st}_p(\Bb)_{hm})}\to \Z_0(\ul{\HH(\Bb)})$ be induced from $hor: \ul{\HH(\Zz(\Bb)_{hm})}\to \Z(\ul{\HH(\Bb)})$ 
from \prref{uses mates} (by forgetting the purity condition first and the choice of centrality structures at the end). 
Then $J$ preserves strictly the structural 
1- and 2-cells (from the associativity and unity constraints). The pseudodouble category $\Zz_p^{st}(\Bb)_{hm}$ inherits 
its monoidal structure $\ot$ from $\Zz_p^{st}(\Bb)$ known from \thref{gen}.

\begin{cor} \colabel{Freyd}
Let $\Bb$ be a premonoidal double category so that its associativity and unity constraints are liftable vertical transformations and whose binoidal structure comes from a pseudodouble quasi-functor $H:\Bb\times\Bb\to\Bb$. 
Then the underlying horizontal pseudofunctor $\HH(J): \ul{\HH(\Zz^{st}_p(\Bb)_{hm})}\to\ul{\HH(\Bb)}$ determines a Freyd bicategory, with $\V=\ul{\HH(\Zz^{st}_p(\Bb)_{hm})}$, the binoidal structure on $\B=\ul{\HH(\Bb)}$ induced from $\HH(H)$, 
monoidal structure on $\V$ induced by $\HH(\ot)$, and $J_\Z$ induced by $hor$:
$$\scalebox{0.8}{
\bfig
 \putmorphism(400,700)(1,0)[\V=\ul{\HH(\Zz^{st}_p(\Bb)_{hm})} `\ul{\HH(\Bb)}=\B ` \HH(J)]{960}1a
\putmorphism(400,700)(1,-1)[\phantom{Y_2}`\Z_0(\ul{\HH(\Bb)}).`J_\Z]{460}1l
\putmorphism(1300,700)(-1,-1)[\phantom{Y_2}``forget]{450}{-1}r
\efig}
$$ 
\end{cor}


\section{Kleisli double categories and premonoidality} \selabel{Kleisli}

The last section of the paper is dedicated to the study of two aspects of Kleisli double categories in regard to strengths on their double monads. We recall them first in lower dimensions.

It is well-known in the categorical setting that strengths on a monad correspond to certain actions on the Kleisli
category of the monad, see {\em e.g.} \cite[Proposition 4.3]{MU}. As it can be seen in \cite{MS}, the latter actions are sometimes used directly to axiomatize models for effectful languages. The other feature is that for a bistrong monad on a monoidal category its Kleisli category is premonoidal. 
The Cartesian case for this can be found in \cite[Section 2]{PT}, and a recent general symmetric instance in \cite[Section 2.2]{UV}.

In one dimension higher one has pseudomonads on bicategories (see {\em e.g.} \cite{M}) and Paquet and Saville have introduced in \cite{HF2} 
strengths for pseudomonads on monoidal bicategories. 
They showed in Theorem 1 of Section 4 that strengths for a pseudomonad correspond to certain actions on its Kleisli bicategory 
and to extensions of the canonical actions of the bicategory on itself. 
Moreover, in Section 6, Theorem 2 of {\em loc.cit.} they showed that the Kleisli bicategory of a bistrong pseudomonad is premonoidal.

In this section we elaborate a double categorical treatment of the subject. We start by 
briefly recalling several notions and results from \cite{GGV}. 
These are the notion of a horizontal and a vertical monad on a double category $\Dd$, 
the fact that assuming existence of companions for certain 1v-cells, a vertical monad $T$ lifts to a horizontal monad 
$\hat T$, and the definition of a horizontal Kleisli pseudodouble category $\Kk l(S)$ for a horizontal monad $S$ on a double category $\Dd$. 
We then introduce strengths, both on a vertical monad $T$ and on a horizontal monad $S$, and prove the following results. 
In all of them $\Dd$ is a monoidal double category and suitable transformations are assumed to be liftable. 
Firstly, that a (bi)strong vertical monad $T$ induces a (bi)strong horizontal monad $\hat T$. Secondly, that 
a strength on a vertical monad $T$ induces a strength on a horizontal monad $\hat T$. Thirdly, that there is a 1-1 
correspondence between strengths on a horizontal monad $S$ and extensions of the canonical actions of $\Dd$ on itself 
(an extension is comprised of an action of $\Dd$ on the horizontal Kleisli pseudodouble category $\Kk l(S)$ 
and a double icon). Finally, that given a bistrong vertical monad $T$, the horizontal Kleisli pseudodouble category 
$\Kk l(\hat T)$ is premonoidal. 

\medskip


In this section we will heavily use the technical tools that we developed in \ssref{liftings}.

\subsection{Double monads and Kleisli double categories}

In this subsection we recall the necessary notions and results from \cite{GGV} in an abbreviated form. 
The notions in \cite{GGV} of: double functors, horizontal transformations and vertical transformations correspond to our notions of: pseudodouble functors, horizontal pseudonatural transformations and vertical strict transformations, respectively. Our formulations below are expressed accordingly in our terminology.

\begin{defn} \delabel{hor-mon}
A {\em horizontal double monad} on a double category $\Dd$ consists of: 
\begin{itemize}
\item pseudodouble functor $T:\Dd\to\Dd$;
\item horizontal pseudonatural transformations $\mu:TT\Rightarrow T$ and $\eta:\Id_{\Dd}\to T$; 
\item invertible modifications $a^T,l^T,r^T$ with respective components given by horizontally globular 2-cells: 
$$
\scalebox{0.8}{
\bfig
 \putmorphism(-170,400)(1,0)[TTT(A)`TT(A) `T\mu(A)]{620}1a
 \putmorphism(450,400)(1,0)[\phantom{A\ot B}`T(A) `\mu(A)]{620}1a
 \putmorphism(-150,50)(1,0)[TTT(A)`TT(A) `\mu(T(A))]{680}1b
 \putmorphism(530,50)(1,0)[\phantom{A\ot B}`T(A) `\mu(A)]{600}1b

\putmorphism(-180,400)(0,-1)[\phantom{Y_2}``=]{350}1r
\putmorphism(1100,400)(0,-1)[\phantom{Y_2}``=]{350}1l
\put(380,210){\fbox{$a_A^T$}}
\efig}
\quad
\scalebox{0.8}{
\bfig
 \putmorphism(-170,400)(1,0)[T(A)`TT(A) `\Epsilon(T(A))]{620}1a
 \putmorphism(450,400)(1,0)[\phantom{A\ot B}`T(A) `\mu(A)]{620}1a
 \putmorphism(-150,50)(1,0)[T(A)`T(A)`=]{1210}1b

\putmorphism(-180,400)(0,-1)[\phantom{Y_2}``=]{350}1r
\putmorphism(1100,400)(0,-1)[\phantom{Y_2}``=]{350}1l
\put(380,210){\fbox{$r_A^T$}}
\efig}
\quad
\scalebox{0.8}{
\bfig
 \putmorphism(-170,400)(1,0)[T(A)`TT(A) `T(\Epsilon(A))]{620}1a
 \putmorphism(450,400)(1,0)[\phantom{A\ot B}`T(A) `\mu(A)]{620}1a
 \putmorphism(-150,50)(1,0)[T(A)`T(A)`=]{1210}1b

\putmorphism(-180,400)(0,-1)[\phantom{Y_2}``=]{350}1r
\putmorphism(1100,400)(0,-1)[\phantom{Y_2}``=]{350}1l
\put(380,210){\fbox{$l_A^T$}}
\efig}
$$
which satisfy two axioms as in \cite[Definition 6.1]{GGV}. 
\end{itemize}
\end{defn}

\begin{defn}
A {\em vertical double monad} on a double category $\Dd$ consists of: 
\begin{itemize}
\item pseudodouble functor $T:\Dd\to\Dd$, and
\item vertical strict transformations $\mu:TT\Rightarrow T$ and $\eta:\Id_{\Dd}\to T$, 
which satisfy the usual associativity and unity laws (via identity vertical modifications). 
\end{itemize}
\end{defn}

For the next result, which is \cite[Theorem 7.4]{GGV}, we need to clarify the following. In \cite[Definition 3.6]{GGV} a vertical (strict) transformation $\alpha$ 
is said to be {\em special} if: 1) its 1v-cell components have companions, and 2) the companion transposes $\hat\alpha_f$ of the 2-cell 
components $\alpha_f$ of $\alpha$ for any 1h-cell $f$ (recall \prref{lifting 1v to equiv}) are invertible. Since we proved in our cited proposition that $\hat\alpha_f$'s are invertible if $\alpha$ is an {\em invertible} vertical strict transformation, our 
``invertible liftable vertical strict transformations'' are special in the sense of \cite{GGV}. 
Thus we get to the following formulation of \cite[Theorem 7.4]{GGV} in our terminology.

\begin{thm} \thlabel{T-lift}
Let $T:\Dd\to\Dd$ be a vertical double monad in a double category $\Dd$. Assume that its multiplication 
$\mu:TT\Rightarrow T$ and unit $\eta:\Id_{\Dd}\to T$ are invertible liftable vertical transformations. Then $(T,\mu,\eta)$ induces a horizontal double monad $(\hat T, \hat\mu,\hat\eta)$ on $\Dd$. 
\end{thm}

Although one may define a vertical Kleisli pseudodouble category of a vertical double monad $T$, there is no natural way to lift it to a horizontal Kleisli pseudodouble category of the horizontal double monad $\hat T$, which is where we want to obtain our results. For this reason we will only work with the latter version of a Kleisli pseudodouble category, that was introduced in \cite[Theorem 9.1]{GGV}. We recall it next.

\begin{thm} \cite[Theorem 9.1]{GGV} \thlabel{Kleisli}
Let $\Dd$ be a double category and $(T,\mu,\eta)$ a horizontal double monad on it. There is a pseudodouble category $\Kk l(T)$, 
called the {\em horizontal Kleisli pseudodouble category} of $\Dd$ whose objects and 1v-cells are the same as in $\Dd$, 1h-cells $A\to B$ are 
1h-cells $A\to T(B)$ in $\Dd$, and 2-cells as on the left below are the 2-cells of $\Dd$ as to the right:
$$
\bfig
\putmorphism(-110,120)(1,0)[A`B`f^K]{400}4a
\putmorphism(-120,140)(0,-1)[``u]{380}1l
\putmorphism(290,140)(0,-1)[``v]{380}1r
\putmorphism(-110,-250)(1,0)[\tilde A`\tilde B`g^K]{400}4a
\put(0,-50){\fbox{$\phi^K$}}
\efig
\qquad\qquad
\bfig
\putmorphism(-110,120)(1,0)[A`T(B)`f]{400}1a
\putmorphism(-120,140)(0,-1)[``u]{380}1l
\putmorphism(290,140)(0,-1)[``T(v)]{380}1r
\putmorphism(-110,-250)(1,0)[\tilde A`T(\tilde B).`g]{400}1a
\put(-20,-50){\fbox{$\phi$}}
\efig
$$
Vertical composition in $\Kk l(T)$ is the same as in $\Dd$, horizontal composition of Kleisli 1h-cells $f:A\to B$ and $g:B\to C$ is given by 
$$A\stackrel{f}{\to}T(B) \stackrel{T(g)}{\to}T^2(C)\stackrel{\mu(C)}{\to}T(C),$$
horizontal composition of Kleisli 2-cells $\phi$ and $\psi$ is given by 
$$
\scalebox{0.86}{
\bfig
 \putmorphism(-150,250)(1,0)[A`\phantom{F(A)} `f]{500}1a
\putmorphism(-150,250)(0,-1)[``u]{450}1l
\put(-40,30){\fbox{$\phi^K$}}
\putmorphism(380,250)(0,-1)[\phantom{Y_2}` `]{450}1l
\putmorphism(410,250)(0,-1)[\phantom{Y_2}` `T(v)]{450}0l
\putmorphism(950,250)(0,-1)[\phantom{Y_2}` `]{450}1r
\putmorphism(930,250)(0,-1)[\phantom{Y_2}` `T^2(w)]{450}0r
\putmorphism(350,250)(1,0)[T(B)`T^2(C)`T(g)]{600}1a
 \putmorphism(950,250)(1,0)[\phantom{F(A)}`T(C) `\mu(C)]{600}1a
 \putmorphism(470,-200)(1,0)[`T^2(\tilde C)`T(\tilde g)]{500}1b
 \putmorphism(1060,-200)(1,0)[`T(\tilde C).`\mu(\tilde C)]{500}1b
\putmorphism(1570,250)(0,-1)[\phantom{Y_2}``T(w)]{450}1r
\put(540,10){\fbox{$T(\psi^K)$}}
\put(1240,30){\fbox{$\mu^w$}}
\putmorphism(-150,-200)(1,0)[\tilde A`T(\tilde B) `\tilde f]{520}1b
\efig}
$$
The identity 1h-cell on $A$ is given by $\eta(A)$ in $\Dd$, and the horizontal identity 2-cell on $u:A\to\tilde A$ is 
$$
\bfig
\putmorphism(-160,120)(1,0)[A`T(A)`\eta(A)]{450}1a
\putmorphism(-160,130)(0,-1)[``u]{380}1l
\putmorphism(270,130)(0,-1)[``T(u)]{380}1r
\putmorphism(-160,-250)(1,0)[\tilde A`T(\tilde A).`\eta(\tilde A)]{450}1a
\put(-20,-30){\fbox{$\eta^u$}}
\efig
$$
\end{thm}

Observe that the underlying horizontal bicategory $\HH(\Kk l(T))$ of the Kleisli pseudodouble category coincides with the Kleisli bicategory $\Kl(\HH(T))$ of the underlying pseudomonad on the bicategory $\HH(\Dd)$ of the horizontal double monad $T$. We will also need:

\begin{defn} \cite[Definition 9.2]{GGV} \delabel{Kl embed} \\
Let $\Dd$ be a double category and $T$ a horizontal double monad on it. The
canonical embedding $K:\Dd\to\Kl(T)$ is the pseudodouble functor that is the identity on objects and vertical
1-cells, it sends a horizontal 1-cell $f: A\to B$ into a 1h-cell $K(f): A\to B$ in $\Kl(T)$ determined by $\eta_B\comp f: A\to T(B)$ in $\Dd$, and correspondingly a 2-cell $\phi$ to the horizontal composition of 2-cells 
$$
\scalebox{0.86}{
\bfig
 \putmorphism(-50,250)(1,0)[A`\phantom{F(A)} `f]{520}1a
\putmorphism(-50,-200)(1,0)[\tilde A`\tilde B `g]{460}1b
\putmorphism(-50,250)(0,-1)[``u]{450}1l
\put(90,30){\fbox{$\phi$}}
\putmorphism(410,250)(0,-1)[B` `v]{450}1l
 \putmorphism(350,250)(1,0)[\phantom{F(A)}`T(B) `\eta(B)]{600}1a
 \putmorphism(450,-200)(1,0)[`T(\tilde B).`\eta(\tilde B)]{500}1b
\putmorphism(930,250)(0,-1)[\phantom{Y_2}``T(w)]{450}1r
\put(570,30){\fbox{$\eta^v$}}
\efig}
$$
\end{defn}

\begin{prop} \cite[Proposition 9.3]{GGV}, \cite[Proposition 7.5]{CS} \prlabel{Kl compan} \\
If a 1v-cell $u$ in $\Dd$ has a companion $\hat u$ in $\Dd$, then it has a companion in $\Kl(T)$, for a 
horizontal double monad $(T,\mu,\eta)$ on $\Dd$ with structure 2-cells 
$$
e_{K(u)}=
\scalebox{0.86}{
\bfig
 \putmorphism(-50,200)(1,0)[A`\phantom{F(A)} `\hat u]{520}1a
\putmorphism(-50,-200)(1,0)[\tilde A`\tilde A `=]{460}1b
\putmorphism(-50,200)(0,-1)[``u]{400}1l
\put(90,-20){\fbox{$e_u$}}
\putmorphism(410,200)(0,-1)[\tilde A` `=]{400}1l
 \putmorphism(350,200)(1,0)[\phantom{F(A)}`T(\tilde A) `\eta(\tilde A)]{580}1a
 \putmorphism(450,-200)(1,0)[`T(\tilde A)`\eta(\tilde A)]{480}1b
\putmorphism(900,200)(0,-1)[\phantom{Y_2}``=]{400}1r
\put(570,-20){\fbox{$\Id$}}
\efig}
\qquad\qquad
\iota_{K(u)}=
\scalebox{0.86}{
\bfig
 \putmorphism(-50,200)(1,0)[A`\phantom{F(A)} `=]{520}1a
\putmorphism(-50,-200)(1,0)[A`\tilde A `\hat u]{460}1b
\putmorphism(-50,200)(0,-1)[``=]{400}1l
\put(90,-20){\fbox{$\iota_u$}}
\putmorphism(410,200)(0,-1)[A` `u]{400}1l
 \putmorphism(350,200)(1,0)[\phantom{F(A)}`T(A) `\eta(A)]{580}1a
 \putmorphism(450,-200)(1,0)[`T(\tilde A).`\eta(\tilde A)]{480}1b
\putmorphism(900,200)(0,-1)[\phantom{Y_2}``T(u)]{400}1r
\put(570,-20){\fbox{$\eta^u$}}
\efig}
$$
\end{prop}

\subsection{Strengths on double monads} \sslabel{strengths}

Pseudomonads on bicategories were introduced in \cite{M}. Based on \prref{lifting 1v to equiv} we obtain:

\begin{prop} \prlabel{str pass to bicat}
A vertical double monad $T$ on a double category $\Dd$ such that its multiplication and unit are invertible liftable vertical transformations induces a pseudomonad $\HH(\hat T)$ on the underlying bicategory $\HH(\Dd)$. 
\end{prop}

In \cite[Definition 9]{HF}, \cite[Definition 4.3]{HF2} Paquet and Saville introduced a (left) strength on a pseudomonad $T$ on a monoidal bicategory $\B$. 
We next give a definition of a strength in a double categorical setting. For this we set that a horizontal and a vertical double monad on a {\em monoidal} double category $\Dd$ are simply a horizontal and a vertical double monad on 
the underlying double category $\Dd$, respectively.

\begin{defn}
A {\em left (vertical) strength} on a vertical double monad $T$ on a monoidal double category $\Dd$ consists of:
\begin{itemize}
\item a vertical strict transformation with 1v-components \\ $t_{A,B}:A\ot T(B)\to T(A\ot B)$ for $A,B\in\Dd$; 
\item identity vertical modifications with the following identity vertically globular 2-cell components, expressing compatibility of the strength $t$ with:
 \begin{enumerate}[a)]
\item monoidal structure of $\Dd$ \vspace{-0,2cm}
$$x_A: \lambda_{T(A)}\Rightarrow\frac{t_{I,A}}{ T\lambda_A} \qquad \qquad
y_{A,B,C}:\threefrac{\alpha_{A,B,T(C)}}{At_{B,C}}{ t_{A,BC}} \Rightarrow\frac{t_{AB,C}}{T(\alpha_{A,B,C})};$$
\item monad structure of $T$ \vspace{-0,2cm} 
$$w_{A,B}: \frac{A\mu_B}{t_{A,B}}\Rightarrow\threefrac{t_{A,T(B)}}{ T(t_{A,B})}{ \mu_{A,B}} \qquad \qquad
z_{A,B}:\eta_{AB}\Rightarrow\frac{A\eta_B}{t_{A,B}}.$$
\end{enumerate}
\end{itemize}
\end{defn}

One defines a {\em right (vertical) strength} on a vertical double monad $T$ analogously, based on 
a vertical strict transformation with 1v-components $s_{A,B}:T(A)\ot B\to T(A\ot B)$. 

We can carry out the following reasoning.  

\medskip

Let $(T, \mu, \lambda)$ be a vertical double monad on a monoidal double category $\Dd$.   
Let $t$ be a left strength on $T$, and assume that the following transformations are invertible and liftable:  
strength $t$, \,\, 
$\alpha, \lambda, \rho$ of $\Dd$, and $\mu, \eta$ of $T$. 
Then by \thref{T-lift} we have a horizontal double monad $(\hat T,\hat\mu,\hat\eta)$,  
by \prref{lifting 1v to equiv} we have a horizontal natural equivalence $\hat t:-\ot \hat T(-)\to \hat T(-\ot -)$ in $\Dd$,  
by \prref{essence} we have invertible horizontal modifications 
$$\hat x_A: [\hat t_{I,A}\vert T\hat \lambda_A]\Rightarrow\hat \lambda_{T(A)} \qquad 
\hat y_{A,B,C}:[\hat t_{AB,C}\vert T(\hat \alpha_{A,B,C})]\Rightarrow[\hat \alpha_{A,B,T(C)}\vert A\hat t_{B,C}\vert \hat t_{A,BC}] 
\vspace{-0,1cm}$$
$$\hat w_{A,B}: [\hat t_{A,T(B)}\vert T(\hat t_{A,B})\vert \hat \mu_{A,B}]\Rightarrow[A\hat \mu_B\vert \hat t_{A,B}] \qquad 
\hat z_{A,B}:[A\hat \eta_B\vert \hat t_{A,B}]\Rightarrow\hat \eta_{A,B} \vspace{0,1cm}$$
in $\Dd$ satisfying any sensible equation that can be formed by them. 

\medskip

This motivates the following definition. 

\begin{defn} \delabel{hor mon D}
For a double category $\Dd$ we say that it is {\em horizontally monoidal} if there are pseudodouble functors 
$\ot:\Dd\times\Dd\to\Dd$ and $I:*\to\Dd$, horizontal equivalence transformations 
$$\alpha: \ot\comp(\Id\times\ot) \stackrel{\iso}{\to} \ot\comp(\ot\times\Id) $$
$$\lambda: \ot\comp(I\times\Id) \stackrel{\iso}{\to} \Id $$
$$\rho: \ot\comp(\Id\times I) \stackrel{\iso}{\to} \Id ,$$ 
and horizontal modifications $p,m,l,r$ whose 2-cell components satisfy axioms (TA1)-(TA3) as in \cite{GPS} (write the latter as equations of pasted horizontally globular 2-cells). 
\end{defn} 

The above definition is such that for a horizontally monoidal double category $(\Dd,\ot,\alpha,\lambda,\rho)$ 
the underlying horizontal bicategory $\HH(\Dd)$ is a monoidal bicategory. 

\medskip

A horizontal double monad on a {\em horizontally} monoidal double category $\Dd$ is a horizontal double monad on the underlying double category $\Dd$. 

\medskip

Due to \prref{lifting 1v to equiv} and \prref{essence} we have:

\begin{thm} \thlabel{two monoidalities}
A monoidal double category $(\Dd, \ot, I, \alpha,\lambda,\rho)$ in which $\alpha,\lambda,\rho$ are liftable 
vertical transformations yields a horizontally monoidal double category $(\Dd, \ot, I, \hat\alpha, \hat\lambda,\hat\rho)$. 
\end{thm}

From now on we will omit to write $I$ explicitly when referring to this result. Joining \thref{T-lift} with the latter theorem we clearly have:

\begin{prop}
Given a vertical double monad $(T, \mu,\eta)$ on a monoidal double category $(\Dd, \ot, \alpha,\lambda,\rho)$, 
and assume that $\mu,\eta, \, \alpha,\lambda,\rho$ are invertible liftable transformations. 
Then $(\hat T, \hat\mu,\hat\eta)$ is a horizontal double monad on a horizontally monoidal double category 
$(\Dd, \ot, \hat\alpha,\hat\lambda,\hat\rho)$.
\end{prop}

We finally define:

\begin{defn} \delabel{hor-strength}
A {\em left (horizontal) strength} on a horizontal double monad $T$ on a horizontally monoidal double category $\Dd$ consists of:
\begin{enumerate} [a)]
\item a horizontal pseudonatural transformation with 1h-components \\ $t_{A,B}:A\ot T(B)\to T(A\ot B)$ for $A,B\in\Dd$; 
\item invertible horizontal modifications with horizontally globular 2-cell components 
$$x_A: [t_{I,A}\vert T\lambda_A]\Rightarrow\lambda_{T(A)} \qquad 
y_{A,B,C}:[t_{AB,C}\vert T(\alpha_{A,B,C})] \Rightarrow[\alpha_{A,B,T(C)}\vert At_{B,C}\vert t_{A,BC}]$$
which fulfill three axioms as in \cite[Figure 3]{HF2} (connecting them to the horizontal modifications $p,m,l$ 
of the horizontally monoidal double category $\Dd$);
\item invertible horizontal modifications with horizontally globular 2-cell components 
$$w_{A,B}: [t_{A,T(B)}\vert T(t_{A,B})\vert \mu_{A,B}]\Rightarrow[A\mu_B\vert t_{A,B}] \qquad 
z_{A,B}:[A\eta_B\vert t_{A,B}]\Rightarrow\eta_{A,B}$$
which fulfill seven axioms as in \cite[Figure 4, p. 12]{HF} ({\em i.e.} in Figures 4 and 5 of \cite{HF2})
(three axioms for connecting $w$ and $z$ to the horizontal modifications $l^T,r^T,a^T$ of the horizontal double monad structure of $T$; two for connecting both $z$ and $w$ with $x$, and two for connecting both $z$ and $w$ with $y$). 
\end{enumerate}
\end{defn}


\medskip

For vertically globular 2-cells $\omega$ we will call the operation $\omega\mapsto\hat\omega$ as in \equref{2-cell lift} 
a {\em companion-lift of 2-cells}. 

\medskip

Coming back to our above reasoning we get to: 

\begin{thm} \thlabel{strenghts thm}
Let $(T, \mu, \eta)$ be a vertical double monad on a monoidal double category 
$(\Dd, \ot, \alpha,\lambda,\rho)$.   
Let $t$ be a left (respectively right) (vertical) strength on $T$, and assume that the following 
transformations are invertible and liftable: strength $t$, \,\, $\alpha, \lambda, \rho$ of $\Dd$, and $\mu, \eta$ of $T$. 

Then $\hat t$ is a left (respectively right) (horizontal) strength on the horizontal double monad $\hat T$ on the horizontally monoidal double category $(\Dd, \ot, \hat\alpha, \hat\lambda,\hat\rho)$. 

Moreover, $\HH(\hat t)$ is a left (respectively right) strength on the pseudomonad $\HH(\hat T)$ on the monoidal bicategory 
$\HH(\Dd)$. 
\end{thm}

\begin{proof}
For the last statement use \prref{str pass to bicat}. 
It only remains to prove that $\hat x$ and $\hat y$ obey the three axioms 
from the item $b)$ of 
\deref{hor-strength}, and $\hat w$ and $\hat y$ the seven axioms from the item c) of that Definition, since apart from modification components these axioms contain also the 2-cell components of the horizontal pseudonatural transformations $\hat t, \hat \alpha, 
\hat \lambda, \hat\mu$ and $\hat \eta$. 
Namely, the first axiom involves the 2-cell components of the modifications $\hat m, \hat x, \hat y$ and the 2-cell components 
$\hat t_{1_A,\hat\lambda_B}$ and $\hat t_{\hat\rho_A,1_B}$ of $\hat t$. The 1h-cells of their domains and codomains are all 
companions, so $\hat t_{1_A,\hat \lambda_B}$ and $\hat t_{\hat\rho_A,1_B}$ are the canonical isomorphism 2-cells of $\theta$ type (recall \leref{horiz is teta}), and the first axiom is proved to hold the same way as in the proof of \prref{essence}, part 3. The same holds for the remaining nine axioms. We only list the 2-cells appearing in the other two axioms for $\hat x$ and $\hat y$, for reader's convenience. In the second axiom there appear the modifications $\hat l, 
\hat x, \hat y$ and the 2-cell components $\hat\lambda_{t_{A,B}}$ of $\hat\lambda$ and $\hat t_{\hat\lambda_A,1_B}$ of $\hat t$. Finally, in the third axiom there appear the modifications $\hat p, \hat y$ and the 2-cell components $\hat\alpha_{1_A,1_B,t_{C,D}}$ of $\hat\alpha$ and 
$\hat t_{\alpha_{A,B,C},1_D}, \hat t_{1_A,\alpha_{B,C,D}}$ of $\hat t$. 
\qed\end{proof}

\subsection{Vertical strengths induce actions on the Kleisli double category} \sslabel{actions}

In this subsection we 
prove our double categorical version of the result known for categories and bicategories, that a strength of a double monad induces an action of a monoidal double category to the corresponding Kleisli double category of the double monad. Here we first focus on vertical strengths, 
the case of horizontal strengths we leave for \sssref{h.str.induce act}.

\subsubsection{Actions on double categories} 

As before, we may consider two kinds of actions.

\begin{defn} \delabel{ver-act} 
We say that a monoidal double category $\Dd$ acts (from the left) on a double category $\Ee$ if there is a pseudodouble functor 
$F:\Dd\times\Ee\to\Ee$, invertible vertical strict transformations with components 
$$\tilde\lambda_E: I\times E\to E\quad\text{and}\quad\tilde\alpha_{A,B,E}:(A\ot B)\times E\to A\times(B\times E)$$
with $A,B\in\Dd$ and $E\in\Ee$, and identity vertical modifications $\tilde p, \tilde l, \tilde m$, analogous to $p,l,m$ from \deref{hor mon D}.  \\
\end{defn}

This kind of action we might call a {\em vertical action}. To the contrast to it, we will differ what might be called a {\em horizontal action}. The difference is actually already anticipated by the fact that in the vertical action Shulman's kind of monoidality of the acting double category is meant, whereas in the horizontal action {\em horizontal monoidality} of the acting double category is assumed.
This ``horizontal'' action we define as follows. 

\begin{defn} \delabel{hor-act} 
By an action of a {\em horizontally monoidal} double category $\Dd$ on a double category $\Ee$ we mean the data comprised of: a 
pseudodouble functor $F:\Dd\times\Ee\to\Ee$, horizontal equivalences with components 
$$\tilde\lambda_E: I\times E\to E\quad\text{and}\quad\tilde\alpha_{A,B,E}:(A\ot B)\times E\to A\times(B\times E)$$
with $A,B\in\Dd$ and $E\in\Ee$, and horizontal modifications $\tilde p, \tilde l, \tilde m$ 
whose 2-cell components satisfy axioms analogous to (TA1)-(TA2) of \cite{GPS}. 
\end{defn}

By \thref{two monoidalities} and \prref{essence} we get:

\begin{prop} \prlabel{vert->horiz act}
Given a vertical action of a monoidal double category 
$(\Dd, \ot, \alpha,\lambda,\rho)$ by a pseudodouble functor $F$ on a double category $\Ee$, so that $\alpha,\lambda,\rho$ are liftable, 
then $F$ induces a horizontal action of the horizontally monoidal double category $(\Dd, \ot, \hat\alpha, \hat\lambda,\hat\rho)$ on $\Ee$. 
\end{prop}

\subsubsection{Vertical strengths induce actions on the Kleisli double category} \ssslabel{v.str.induce act}

Recall from \thref{Kleisli} that 1h- and 2-cells in $\Kk l(T)$  we denote by $f^K$ and $\phi^K$ and that they are given by 
 1h- and 2-cells $f$ and $\phi$ in $\Dd$, respectively. We now prove:

\begin{thm} \thlabel{strength-action}
Let $t$ be a (vertical) left strength on a vertical double monad $T$ on a monoidal double category $\Dd$. 
Assume that the vertical transformations $\alpha,\lambda,\rho$ of $\Dd$, $t$ and $\mu,\eta$ of $T$ are invertible and liftable. The following then hold. 
\begin{enumerate}
\item There is a pseudodouble functor $\rtr:\Dd\times\Kk l(\hat T)\to\Kk l(\hat T)$ defined on objects and 1v-cells by the action of the horizontally monoidal product of $\Dd$, for 1h-cells $f:A\to A'\in\Dd$ and $g^K:B\to B'\in\Kk l(\hat T)$ we define 
$$f\rtr g^K:=(A\ot B\stackrel{f\ot g}{\to}A'\ot T(B') \stackrel{\hat t_{A',B'}}{\to}T(A'\ot B'))$$
and for 2-cells $\sigma\in\Dd, \delta^K\in\Kk l(\hat T)$ we set 
$$\sigma\rtr\delta^K:=[\sigma\ot\delta \,\, \vert \,\, \hat t^{v,v'}]$$
where $v,v'$ are the right hand-side 1v-cells of $\sigma$ and $\delta$, respectively. 

\item There are invertible vertical strict (action) transformations $\tilde\lambda$ and $\tilde\alpha=\tilde\alpha^L$ with 2-cell components 
\begin{multline*} 
\scalebox{0.82}{
\bfig
\putmorphism(-110,170)(1,0)[(AB)C`T((A'B')C') `(fg)\rtr h^K]{900}1a
\putmorphism(-160,200)(0,-1)[``\tilde\alpha_{A,B,C}]{410}1l 
\putmorphism(830,200)(0,-1)[``T(\tilde\alpha_{A',B',C'})]{410}1r
\putmorphism(-110,-200)(1,0)[A(BC)`T(A'(B'C'))`f\rtr(g\rtr h^K)]{900}1a
\put(140,20){\fbox{$\tilde\alpha_{f,g,h^K}^L$}}
\efig}
\,\,:= \\
\scalebox{0.82}{
\bfig
 \putmorphism(-650,210)(1,0)[(AB)C` (A'B')T(C')`(fg)h]{660}1a 
\putmorphism(-660,200)(0,-1)[\phantom{Y_2}` `\alpha_{A,B,C}]{380}1l
\putmorphism(-50,200)(0,-1)[\phantom{Y_2}` `]{380}1r
\putmorphism(-50,200)(0,-1)[\phantom{Y_2}` `\alpha_{A',B',T(C')}]{380}0r
\putmorphism(200,210)(1,0)[\phantom{Y}` T((A'B')C')  `\hat t_{A'B',C'}]{1420}1a
\putmorphism(1600,200)(0,-1)[\phantom{Y_2}`T(A'(B'C')) `]{380}0r
\putmorphism(1530,200)(0,-1)[\phantom{Y_2}` `T(\alpha_{A',B',C'})]{380}1r
\putmorphism(-670,-200)(1,0)[A(BC) `A'(B'T(C'))`f(gh)]{640}1a 
\putmorphism(200,-200)(1,0)[``A'\hat t_{B',C'}]{380}1a
\putmorphism(780,-200)(1,0)[A'T(B'C')``\hat t_{A',B'C'}]{580}1a
\put(-500,30){\fbox{$\alpha_{f,g,h}$}}
\put(580,30){\fbox{$ y^*_{A',B',C'}$}}
\efig}
\end{multline*}
and 
$$\scalebox{0.82}{
\bfig
\putmorphism(-110,170)(1,0)[I\ot A`T(IA')`I\rtr h^K]{550}1a
\putmorphism(-110,200)(0,-1)[``\tilde\lambda_A]{400}1l 
\putmorphism(470,200)(0,-1)[``T(\tilde\lambda_{A'})]{400}1r
\putmorphism(-120,-200)(1,0)[A`T(A')`h^K]{580}1a
\put(60,-10){\fbox{$\tilde\lambda_{h^K}$}}
\efig}
\quad:=\quad
\scalebox{0.82}{
\bfig
 \putmorphism(-650,210)(1,0)[I\ot A` I\ot T(A')`I\ot h]{650}1a  
\putmorphism(-660,200)(0,-1)[\phantom{Y_2}` `\lambda_A]{380}1l
\putmorphism(-70,200)(0,-1)[\phantom{Y_2}` T(A')`]{380}1r 
\putmorphism(-80,200)(0,-1)[\phantom{Y_2}` `\lambda_{T(A')}]{380}0r 
 \putmorphism(150,210)(1,0)[\phantom{Y}` T(IA')  `\hat t_{I,A'}]{500}1a
\putmorphism(630,200)(0,-1)[\phantom{Y_2}` `T(\lambda_{A'})]{380}1r
\putmorphism(-660,-200)(1,0)[A``h]{470}1a
\putmorphism(50,-200)(1,0)[`T(A')`=]{600}1a
\put(-450,10){\fbox{$\lambda_h$}}
\put(250,0){\fbox{$x^*_{A'}$}}
\efig}
$$
for 1h-cells $f,g\in\Dd$ and $h^K\in\Kk l(\hat T)$. Here $y^*, x^*$ are the 2-cells induced from the identity vertical modification 
2-cell components of $y,x$ of the strength $t$ as described by the assignment \equref{omega-inv}. 

\item The above data define a left vertical action of $\Dd$ on $\Kk l(\hat T)$, with $\tilde p, \tilde l, \tilde m$ being 
$K(p),K(l),K(m)$, where $p,l,m$ are the identity vertical modifications from \deref{Shul}, respectively.  

\item There are horizontal equivalences $\widehat{\tilde\alpha}$ and $\widehat{\tilde\lambda}$, and 
horizontal modifications $\widehat{\tilde p}, \widehat{\tilde l}, \widehat{\tilde m}$ obeying the axioms of 
\deref{hor-act}, so that $\rtr$ induces a left horizontal action of $\Dd$ on $\Kk l(\hat T)$.
\end{enumerate}
\end{thm}

\begin{proof}
The compositor 2-cell for $\rtr$ is given by the globular 2-cell 
$$ 
\scalebox{0.82}{
\bfig
 \putmorphism(300,-450)(1,0)[A'\ot T(B')` T(A'\ot B') `\hat t_{A',B'}]{860}1a
 \putmorphism(1360,-450)(1,0)[` T(A''\ot T(B''))  `T(f'\ot g')]{730}1a

\putmorphism(340,-450)(0,-1)[\phantom{Y_2}``=]{450}1l
\putmorphism(1150,-900)(0,-1)[\phantom{Y_2}``=]{450}1l

 \putmorphism(-390,-1350)(1,0)[\phantom{A''\ot B'}` A''\ot T^2(B'') `f'f\ot T(g')g]{1460}1b 
 \putmorphism(1150,-1350)(1,0)[\phantom{A''\ot B'}` A''\ot T(B'') `A''\ot\hat\mu_{B''}]{1300}1b 
 \putmorphism(2540,-1350)(1,0)[\phantom{A''\ot B'}`  `\hat t_{A'',B''}]{1100}1b 

\put(980,-680){\fbox{$\hat t^{-1}_{f',g'}$}}
\put(280,-1150){\fbox{$lax \ot$}}

\putmorphism(-390,-900)(1,0)[A\ot B` A'\ot T(B')`f\ot g]{700}1a 
\putmorphism(-390,-900)(0,-1)[\phantom{Y_2}`A\ot B `=]{450}1r 
\putmorphism(380,-900)(1,0)[\phantom{A''\ot B'}`A''\ot T^2(B'') ` f'\ot T(g')]{810}1a
 \putmorphism(1260,-900)(1,0)[\phantom{A''\ot B'}`  `\hat t_{A'',B''}]{520}1a 

 \putmorphism(2200,-900)(1,0)[\phantom{A''\ot B'}` T^2(A''\ot B'') `T(\hat t_{A'',B''})]{770}1a 
 \putmorphism(3070,-900)(1,0)[\phantom{A''\ot B'}`T(A''\ot B'')  `\mu_{A'',B''}]{800}1a 

\putmorphism(2060,-450)(0,-1)[\phantom{Y_2}`T(A''\ot T(B'')) ` =]{450}1r 
\putmorphism(3900,-900)(0,-1)[\phantom{Y_2}`T(A''\ot B'').`=]{450}1r 
\put(2360,-1150){\fbox{$\hat w_{A'',B''}$}}
\efig}
$$
The hexagonal law for the compositor is proved as in the bicategorical case of \cite[Proposition 7.1]{HF2}. We comment that one uses the axioms 
\axiomref{h.o.t.-1} for $t$, \axiomref{m.ho-vl.-1} for $w$, and the $w-a^T$ axiom from c) of \deref{hor-strength}. 
The unitor is defined via $\hat z$ and the unitor law is proved as in the bicategorical case, too. 
The axiom \axiomref{lx.f.u-nat} reads $\frac{\hat z_{A,B}}{\hat\eta^{u,v}}=\frac{[\Id^u\ot\hat\eta^v \,\, \vert \,\, 
\hat t^{u,v}]}{\hat z_{\tilde A, \tilde B}}$. It holds true because the corresponding equality comprised of original vertical modifications and strict vertical transformations $z_{A,B}, \eta^{u,v}, z_{\tilde A, \tilde B}, \Id^u\times\eta^v, t^{u,v}$ 
(trivially) holds true. It is important to notice that the strictness of $t$ was necessary at this point. 
The rest of the pseudodouble functor axioms are proved as follows: 
\axiomref{lx.f.s1} and \axiomref{lx.f.s2} hold by \axiomref{h.o.t.-3} and \axiomref{h.o.t.-4} of $\hat t$, respectively;  
\axiomref{lx.f.c-nat} holds by: \axiomref{m.ho-vl.-1} of the modification $\hat w$, \axiomref{lx.f.c-nat} of the pseudodouble functor $\ot$ and \axiomref{h.o.t.-5} of $(\hat t)^{-1}$; and \axiomref{lx.f.u-nat}, which we commented above, holds also by \axiomref{m.ho-vl.-2} of the modification $\hat z$.

Observe that the 2-cells $\tilde\alpha_{f,g,h}^L$ and $\tilde\lambda_h$ are vertically invertible by \leref{inv-spec} and since so 
are $\alpha_{f,g,h}$ and $\lambda_h$. To check the axiom \axiomref{v.l.t.\x 1} for $\tilde\alpha$, apply \axiomref{m.ho-vl.-1} for $y^*$ 
and the axiom w-y from \deref{hor-strength} c). For \axiomref{v.l.t.\x 2}, apply \axiomref{m.ho-vl.-2} for $\hat z$, and for \axiomref{v.l.t.\x 5} for $\tilde\alpha$, apply \axiomref{v.l.t.\x 5} for $\alpha$ and \axiomref{m.ho-vl.-2} for $y^*$. Similar proof goes for $\tilde\lambda$: for 
\axiomref{v.l.t.\x 1} apply \axiomref{m.ho-vl.-1} for $x^*$ and the axiom w-x from \deref{hor-strength} c), and so on.  


The third part is clear. 
The first part of point 4. holds by \prref{lifting 1v to equiv}. 
Take the companion-lifts $\widehat{\tilde p}, \widehat{\tilde l}, \widehat{\tilde m}$ of $\tilde p, \tilde l, \tilde m$, so that by 
\prref{essence} we obtain the second claim in point 4. and hence also the whole theorem. 
\qed\end{proof}

\medskip

A right vertical strength $s$ on a vertical double monad $T$ induces in a similar way a right horizontal action of $\Dd$ on 
$\Kk l(\hat T)$ with 
\begin{multline*} 
\bfig
\putmorphism(-190,170)(1,0)[(AB)C`T((A'B')C') `(fg)\ltr h]{880}1a
\putmorphism(-240,200)(0,-1)[``\alpha_{A,B,C}]{410}1l 
\putmorphism(730,200)(0,-1)[``T(\alpha_{A',B',C'})]{410}1r
\putmorphism(-190,-200)(1,0)[A(BC)`T(A'(B'C'))`f\ltr(g\ltr h)]{880}1a
\put(90,20){\fbox{$\tilde\alpha_{f,g,h}^R$}}
\efig
\,\,:= \\
\bfig
 \putmorphism(-650,210)(1,0)[(AB)C` (T(A')B')C'`(fg)h]{660}1a 
\putmorphism(-660,200)(0,-1)[\phantom{Y_2}` `\alpha_{A,B,C}]{380}1l
\putmorphism(-50,200)(0,-1)[\phantom{Y_2}` `]{380}1r
\putmorphism(-50,200)(0,-1)[\phantom{Y_2}` `\alpha_{T(A'),B',C'}]{380}0r
\putmorphism(200,210)(1,0)[\phantom{Y}` T(A'B')C'  `\hat s_{A',B'}C']{630}1a
\putmorphism(1030,210)(1,0)[\phantom{Y}` T((A'B')C')  `\hat s_{A'B',C'}]{580}1a
\putmorphism(1600,200)(0,-1)[\phantom{Y_2}`T(A'(B'C')) `]{380}0r
\putmorphism(1530,200)(0,-1)[\phantom{Y_2}` `T(\alpha_{A',B',C'})]{380}1r
\putmorphism(-670,-200)(1,0)[A(BC) `T(A')(B'C')`f(gh)]{640}1a 
\putmorphism(200,-200)(1,0)[``\hat s_{A',B'C'}]{1180}1a
\put(-500,30){\fbox{$\alpha_{f,g,h}$}}
\put(580,30){\fbox{$ (y')^*_{A',B',C'}$}}
\efig
\end{multline*}
and 
$$\bfig
\putmorphism(-110,170)(1,0)[A\ot I`T(IA')`h\ltr I]{550}1a
\putmorphism(-110,200)(0,-1)[``\rho_A]{400}1l 
\putmorphism(470,200)(0,-1)[``T(\rho_{A'})]{400}1r
\putmorphism(-120,-200)(1,0)[A`T(A')`h]{580}1a
\put(60,-10){\fbox{$\tilde\rho_h$}}
\efig
\quad:=\quad
\bfig
 \putmorphism(-650,210)(1,0)[A\ot I ` T(A')\ot I`h\ot I]{650}1a  
\putmorphism(-660,200)(0,-1)[\phantom{Y_2}` `\rho_A]{380}1l
\putmorphism(-70,200)(0,-1)[\phantom{Y_2}` T(A')`]{380}1r 
\putmorphism(-80,200)(0,-1)[\phantom{Y_2}` `\rho_{T(A')}]{380}0r 
 \putmorphism(150,210)(1,0)[\phantom{Y}` T(IA')  `\hat s_{A',I}]{500}1a
\putmorphism(630,200)(0,-1)[\phantom{Y_2}` `T(\rho_{A'})]{380}1r
\putmorphism(-660,-200)(1,0)[A``h]{470}1a
\putmorphism(50,-200)(1,0)[`T(A')`=]{600}1a
\put(-450,10){\fbox{$\rho_h$}}
\put(230,0){\fbox{$(x')^*_{A'}$}}
\efig
$$
where $(y')^*$ and $(x')^*$ are induced by the assignment \equref{omega-inv} 
from identity vertical modification 2-cell components $y',x'$ of $s$ given by 
$$x'_A: \rho_{T(A)}\Rightarrow\frac{s_{A,I}}{ T\rho_A} \qquad \qquad
y'_{A,B,C}:\frac{\alpha_{T(A),B,C}}{s_{A,BC}}\Rightarrow\threefrac{s_{A,B}\ot C}{s_{AB,C}}{T(\alpha_{A,B,C})} .$$

\begin{cor} \colabel{strength-action-HF}
Let $t$ be a (vertical) left strength on a vertical double monad $T$ on a monoidal double category $\Dd$. 
Assume that the vertical transformations $t$, \, $\mu,\eta,$ of $T$ and $\alpha,\lambda,\rho$ of $\Dd$ are invertible and liftable. 
Then $\HH(\hat t)$ is a left strength on the pseudomonad $\HH(\hat T)$ on the monoidal bicategory $\HH(\Dd)$ and the induced pseudofunctor 
$\rtr:\HH(\Dd)\times\HH(\Kk l(\hat T))\to\HH(\Kk l(\hat T))$ defines a left action of $\HH(\Dd)$ on the 
Kleisli bicategory $\HH(\Kk l(\hat T))=\Kl(\HH(\hat T))$.
\end{cor}

In the next subsection we are going to prove 
a partial converse to \thref{strength-action}. It is a double categorical 
version of \cite[Theorem 7.2]{HF2}, where a 1-1 correspondence between left strengths on a pseudomonad $T$ on a monoidal bicategory 
$\B$, on one side, and of ``extensions of canonical actions of $\B$ on itself'', on the other, are proved. 
The latter notion comprises a left action $\B\times\B_T\to\B_T$, where $\B_T$ is the corresponding Kleisli bicategory, and an icon.

\subsection{Extensions of the canonical action} \sslabel{extensions}

In this subsection we introduce the notion of a horizontal extension of the canonical (horizontal) action of $\Dd$ on itself. 
We will first prove that vertical strengths on vertical double monads $T$ on $\Dd$ induce (horizontal strengths and then) horizontal extensions, and secondly, that horizontal strengths on the induced horizontal double monads $\hat T$ (actually, on any horizontal double monad $S$) are in one-to-one correspondence with horizontal extensions of the canonical action of (horizontally monoidal) $\Dd$ on itself. We will also explain where the difficulty lies in preventing a horizontal strength $s$ on $\hat T$ to induce a vertical strength $t$ on $T$, even though one supposes that the 1h-cells $s_{A,B}$ are companions of some 1v-cells $t_{A,B}$.

\subsubsection{Vertical strengths induce horizontal extensions}

We first prove that the results spelled out in \thref{strength-action} can be extended to having a 
horizontal extension of the canonical action. We will define the latter notion at the end of this first subsection. 

\begin{lma} \lelabel{iota-e for rtr}
In the setting of \thref{strength-action}, if 1v-cells $u,v$ have companions $\hat u, \hat v$, then a companion  $\widehat{u\w\rtr\w K(v)}$ for $u\rtr K(v)$ is given by the 1h-cell $\stackrel{\hat u\ot \hat v}{\longrightarrow}\stackrel{\tilde A\ot\eta_{\tilde B}}
{\longrightarrow}\stackrel{\hat t_{\tilde A, \tilde B}}{\longrightarrow}$ with structure 2-cells $e_{u\rtr K(v)}$ and 
$\iota_{u\rtr K(v)}$ given respectively by the following diagrams in $\Dd$
$$
\scalebox{0.86}{
\bfig
 \putmorphism(-50,200)(1,0)[AB`\phantom{F(A)} `\hat u\ot \hat v]{520}1a
 \putmorphism(370,200)(1,0)[\phantom{F(A)}`\tilde A T(\tilde B) `\tilde A\hat\eta(\tilde B)]{600}1a
 \putmorphism(1000,200)(1,0)[\phantom{F(A)}`T(\tilde A\tilde B) `\hat t_{\tilde A, \tilde B}]{580}1a
\putmorphism(-50,-200)(1,0)[\tilde A\tilde B`\tilde A\tilde B `=]{460}1b
 \putmorphism(450,-200)(1,0)[`T(\tilde A\tilde B)`\hat\eta(\tilde A\tilde B)]{1130}1b
\putmorphism(-50,200)(0,-1)[``u\ot v]{400}1l
\putmorphism(410,200)(0,-1)[\tilde A\tilde B` `=]{400}1l
\putmorphism(1570,200)(0,-1)[\phantom{Y_2}``=]{400}1r
\put(30,-20){\fbox{$e_{u\ot v}$}}
\put(870,-20){\fbox{$\hat z_{\tilde A, \tilde B}$}}
\efig}
\qquad\quad
\scalebox{0.86}{
\bfig
 \putmorphism(-50,200)(1,0)[AB`\phantom{F(A)} `=]{520}1a
 \putmorphism(350,200)(1,0)[\phantom{F(A)}`T(AB) `\hat\eta(AB)]{1130}1a
\putmorphism(-50,-200)(1,0)[AB`\tilde A\tilde B `\hat u\ot\hat v]{460}1b
 \putmorphism(450,-200)(1,0)[`T(\tilde A\tilde B)`\hat\eta(\tilde A\tilde B)]{1130}1a
\putmorphism(-50,200)(0,-1)[``=]{400}1l
\putmorphism(410,200)(0,-1)[AB` `uv]{400}1l
\putmorphism(1570,200)(0,-1)[\phantom{Y_2}``T(uv)]{400}1r
\put(30,-20){\fbox{$\iota_{u\ot v}$}}
\put(1120,-20){\fbox{$\hat\eta^{u\ot v}$}}
\putmorphism(410,-200)(0,-1)[`\tilde A\tilde B `=]{400}1l
\putmorphism(1570,-200)(0,-1)[\phantom{Y_2}``=]{400}1r
 \putmorphism(370,-600)(1,0)[\phantom{F(A)}`\tilde A T(\tilde B) `\tilde A\hat\eta(\tilde B)]{600}1a
 \putmorphism(1000,-600)(1,0)[\phantom{F(A)}`T(\tilde A\tilde B). `\hat t_{\tilde A, \tilde B}]{580}1a
\put(870,-400){\fbox{$\hat z_{\tilde A, \tilde B}^{-1}$}}
\efig}
$$
\end{lma}

\begin{proof}
The identity $\frac{\iota}{e}=\Id$ is easily proved, to prove $[\iota\vert e]=\Id$ apply the axiom \axiomref{h.o.t.\x 5} to 
$\hat\eta^{u\ot v}$ and then to $\hat\eta^{id}$, use that the structure 2-cell $\hat\eta_{\hat\eta_{\tilde A\tilde B}}$ is an identity to cancel out the $\hat z$ and its inverse. The resulting 2-cell is the identity 2-cell
$$ 
\scalebox{0.86}{
\bfig
 \putmorphism(-50,200)(1,0)[AB`\phantom{F(A)} `=]{560}1a
 \putmorphism(290,200)(1,0)[\phantom{F(A)}`T(\tilde A\tilde B) `u\w\rtr\w K(v)]{680}1a
\putmorphism(-50,-150)(1,0)[AB`\tilde A\tilde B `u\w\rtr\w K(v)]{530}1b
 \putmorphism(560,-150)(1,0)[`T(\tilde A\tilde B)`=]{430}1b
\putmorphism(-50,200)(0,-1)[``=]{350}1l
\putmorphism(960,200)(0,-1)[\phantom{Y_2}``=]{350}1r
\efig}
$$
in $\Kk l(\hat T)$ as desired.  
\qed\end{proof}

\begin{lma} \lelabel{Z cell as teta}
The 2-cell 
$$
Z_{\hat u,\hat v}:=
\scalebox{0.86}{
\bfig
 \putmorphism(-50,200)(1,0)[AB`\phantom{F(A)} `\hat u\ot \hat v]{520}1a
 \putmorphism(370,200)(1,0)[\phantom{F(A)}`\tilde A T(\tilde B) `\tilde A\hat\eta(\tilde B)]{600}1a
 \putmorphism(1000,200)(1,0)[\phantom{F(A)}`T(\tilde A\tilde B) `\hat t_{\tilde A, \tilde B}]{580}1a
\putmorphism(-50,-200)(1,0)[AB`\tilde A\tilde  B`\hat u\ot \hat v]{460}1b
 \putmorphism(450,-200)(1,0)[`T(\tilde A\tilde B)`\hat\eta(\tilde A\tilde B)]{1130}1b
\putmorphism(-50,200)(0,-1)[``=]{400}1l
\putmorphism(410,200)(0,-1)[\tilde A\tilde B` `=]{400}1l
\putmorphism(1570,200)(0,-1)[\phantom{Y_2}``=]{400}1r
\put(0,-20){\fbox{$\Id_{u\ot v}$}}
\put(870,-20){\fbox{$\hat z_{\tilde A, \tilde B}$}}
\efig}
$$
is a canonical isomorphism between two companions of the 1v-cell $u\rtr K(v)=K(u\ot v)$.
\end{lma}

\begin{proof}
One easily verifies that $Z_{\hat u,\hat v}$ obeys the condition \equref{teta-property} using $\iota_{u\rtr K(v)}$ from 
\leref{iota-e for rtr} and $e_{K(u\ot v)}$ from \prref{Kl compan}, so the claim follows by \leref{teta}.
\qed\end{proof}

\begin{defn} \delabel{vert-icon}
A {\em vertical icon} is a vertical transformation such that its 1v-cell components are all identities. 

A vertical icon that is simultaneously a strict vertical transformation we will call a {\em strict vertical icon}. 
\end{defn}

Similarly, a {\em horizontal icon} is a horizontal transformation such that the 1-cell components are all identities. 

\begin{lma} \lelabel{vert icon}
Let $\alpha: F\Rightarrow G$ be a vertical icon in a double category $\Bb$ and let $\hat\alpha$ denote the induced horizontal 
transformation (via \prref{lifting 1v to equiv}). Then: 
\begin{enumerate}
\item $\alpha$ is a strict vertical icon if and only if $\hat\alpha$ is a horizontal icon so that $\hat\alpha^u=\Id$ for 1v-cells $u$ 
if and only if $\hat\alpha$ is an icon 
in the underlying bicategory $\HH(\Bb)$;
\item for $\alpha$ as in 1. and a 1v-cell $u$ with a companion $\hat u$, the 2-cell components $\alpha_{\hat u}=\hat\alpha_{\hat u}$ coincide with the canonical isomorphism 2-cell 
$$\theta_{F(\hat u), G(\hat u)}=
\scalebox{0.86}{\bfig
 \putmorphism(-150,250)(1,0)[A`A `=]{480}1a
\put(-60,70){\fbox{$\iota$}} 
\putmorphism(340,250)(0,-1)[\phantom{Y_2}`B `]{450}1l
\putmorphism(360,250)(0,-1)[\phantom{Y_2}` `F(u)]{450}0l
\putmorphism(810,250)(0,-1)[A`B `]{450}1r
\putmorphism(790,250)(0,-1)[` `G(u)]{450}0r
\putmorphism(1270,250)(0,-1)[B`B`=]{450}1r
\putmorphism(350,250)(1,0)[``=]{420}1a
 \putmorphism(730,250)(1,0)[\phantom{F(A)}` `G(\hat u)]{510}1a
 \putmorphism(-150,-200)(1,0)[A` `F(\hat u)]{460}1b
\putmorphism(360,-200)(1,0)[``=]{440}1b
 \putmorphism(820,-200)(1,0)[``=]{440}1b
\put(500,10){\fbox{$\alpha^u$}}
\put(1080,70){\fbox{$e'$}}
\putmorphism(-150,250)(0,-1)[``=]{450}1l
\efig}$$  
from \equref{teta}, which is the companion-lift of the 2-cell component $\alpha^u=\Id$. 
\end{enumerate}
\end{lma}

\begin{proof}
For a strict vertical icon one has $F(u)=G(u)$ for 1v-cells $u$. 
Then $\Id^{F(u)}=\alpha^u=\hat\alpha^u$, and $\alpha_f=\hat\alpha_f$ for 1h-cells $f$. 
As a strict vertical transformation $\alpha$ satisfies then only the three axioms
\axiomref{v.l.t.\x 1}, \axiomref{v.l.t.\x 2} and \axiomref{v.l.t.\x 5}. Simplified by the fact that $\alpha(A)$, and hence $\hat\alpha(A)$, is identity they coincide with the axioms \axiomref{h.o.t.\x 1}, \axiomref{h.o.t.\x 2} and \axiomref{h.o.t.\x 5} for $\hat\alpha$, and with the three axioms for $\hat\alpha$ to be a pseudonatural transformation in $\HH(\Bb)$. 
The second part follows by \leref{horiz is teta} and \leref{teta}. This is an isomorphism 2-cell between two companions of $F(u)=G(u)$.
\qed\end{proof}

\begin{prop} \prlabel{icon teta}
Let $t$ be a vertical left strength on a vertical double monad $T$ on a monoidal double category $\Dd$ and 
let $K:\Dd\to\Kk l(\hat T)$ be the canonical pseudodouble functor from \deref{Kl embed}. 
In the conditions of \thref{strength-action} one has: 
\begin{enumerate}
\item $t$ induces a strict vertical icon $\theta:-\rtr K(-)\Rightarrow K(-\ot-): 
\Dd\times\Dd\to\Kk l(\hat T)$ whose invertible 2-cell components $\theta_{f,g}$ for 1h-cells $f:A\to A', g:B\to B'$ are given by 
$$\theta_{f,g}=
\scalebox{0.82}{
\bfig
 \putmorphism(-650,210)(1,0)[AB` A'B'`f\ot g]{640}1a 
\putmorphism(-660,200)(0,-1)[\phantom{Y_2}` `=]{380}1l
\putmorphism(-50,200)(0,-1)[\phantom{Y_2}` `]{380}1r
\putmorphism(-50,200)(0,-1)[\phantom{Y_2}` `=]{380}0r
\putmorphism(100,210)(1,0)[``A'\w\ot\hat\eta_{B'}]{410}1a
\putmorphism(680,210)(1,0)[A'T(B')`T(A'\ot B')`\hat t_{A',B'}]{680}1a

\putmorphism(1600,200)(0,-1)[\phantom{Y_2}``]{380}0r
\putmorphism(1370,200)(0,-1)[\phantom{Y_2}` `=]{380}1r
\putmorphism(-670,-200)(1,0)[AB `A'B'`f\ot g]{620}1a 
\putmorphism(20,-200)(1,0)[\phantom{Y}` T(A'\ot B');  `\hat\eta_{A'\ot B'}]{1360}1a

\put(-500,30){\fbox{$\Id_{f\ot g}$}}
\put(580,30){\fbox{$\hat z_{A',B'}$}}
\efig}
$$
\item for 1v-cells $u,v$ with companions $\hat u, \hat v$ the 2-cells $Z_{\hat u,\hat v}$ in $\Dd$ from \leref{Z cell as teta} present 
companion-lifts of $\theta^{u,v}=\Id$ in $\Kk l(\hat T)$. 
\end{enumerate}
\end{prop}

\begin{proof}
Recall that $K$ is identity on objects and 1v-cells, while $K(f)=\eta_{A'}\comp f$, and that the action $\rtr$ is given on objects and 
1v-cells by the horizontally monoidal product of $\Dd$. Then for $\theta$ on objects $(A,B)$ we set the identity 1v-cells on $A\ot B$, 
on 1v-cells $(u,v)$ we set $\theta^{u,v}: u\rtr K(v) \Rightarrow u\ot v$ to be the identity 2-cell in $\Kk l(\hat T)$, whereas the domain and codomain of $\theta$ at 1h-cells $(f,g)$ should be given by a horizontally globular 2-cell $\theta_{f,g}: 
f\rtr(\eta_{B'}\comp g)\Rightarrow \eta_{A'\ot B'}\comp (f\ot g)$ in $\Dd$. In \thref{strength-action} the 1h-cell 
$f\rtr(\eta_{B'}\comp g)$ is defined as $\hat t_{A',B'}\comp(f\ot(\eta_{B'}\comp g))$, so the above defined globular 2-cell makes sense as a desired 2-cell component of $\theta$, and we define $\theta$ at $(f,g)$ this way. 
By the first part of \leref{vert icon} it suffices to prove that the three axioms of a pseudonatural transformation in the bicategory 
$\HH(\Bb)$ hold for $\theta$. That they hold it was observed in the bicategorical case in \cite[Theorem 1]{HF}, \cite[Theorem 7.2]{HF2}. 
More precisely, they follow by the three strength axioms involving $z$ from \deref{hor-strength}, c) and the modification axiom for $\hat z$. 
Observe that by the second part of \leref{vert icon} it is
$$ 
\scalebox{0.86}{\bfig
 \putmorphism(-150,250)(1,0)[A`A `=]{480}4a
\put(-60,70){\fbox{$\iota_K$}} 
\putmorphism(340,250)(0,-1)[\phantom{Y_2}`B `]{450}1l
\putmorphism(360,140)(0,-1)[\phantom{Y_2}` `u\rtr K(v)]{450}0l
\putmorphism(810,250)(0,-1)[A`B `]{450}1r
\putmorphism(790,140)(0,-1)[` `K(u\ot v)]{450}0r
\putmorphism(1270,250)(0,-1)[B`B`=]{450}1r
\putmorphism(350,250)(1,0)[``=]{420}4a
 \putmorphism(730,250)(1,0)[\phantom{F(A)}` `\widehat{K(u\ot v)}]{510}4a
 \putmorphism(-150,-200)(1,0)[A` `\widehat{u\w\rtr\w K(v)}]{460}4b
\putmorphism(360,-200)(1,0)[``=]{440}4b
 \putmorphism(820,-200)(1,0)[``=]{440}4b
\put(500,10){\fbox{$\Id$}}
\put(1080,70){\fbox{$e'_K$}}
\putmorphism(-150,250)(0,-1)[``=]{450}1l
\efig}=
\scalebox{0.82}{
\bfig
 \putmorphism(-650,210)(1,0)[AB` \tilde A\tilde B`\hat u\ot \hat v]{630}1a 
\putmorphism(-660,200)(0,-1)[\phantom{Y_2}` `=]{380}1l
\putmorphism(-50,200)(0,-1)[\phantom{Y_2}` `]{380}1r
\putmorphism(-50,200)(0,-1)[\phantom{Y_2}` `=]{380}0r
\putmorphism(60,210)(1,0)[``\tilde A\w\ot\eta_{\tilde B}]{480}1a
\putmorphism(680,210)(1,0)[\tilde AT(\tilde B)`T(\tilde A\ot \tilde B)`\hat t_{\tilde A,\tilde B}]{680}1a

\putmorphism(1600,200)(0,-1)[\phantom{Y_2}``]{380}0r
\putmorphism(1370,200)(0,-1)[\phantom{Y_2}` `=]{380}1r
\putmorphism(-670,-200)(1,0)[AB `\tilde A\tilde B`\hat u\ot \hat v]{640}1a 
\putmorphism(20,-200)(1,0)[\phantom{Y}` T(\tilde A\ot \tilde B) `\eta_{\tilde A\ot \tilde B}]{1380}1a

\put(-500,30){\fbox{$\Id_{\hat u\ot \hat v}$}}
\put(580,30){\fbox{$\hat z_{\tilde A,\tilde B}$}}
\efig}
=\theta_{\hat u, \hat v}
$$
where the diagram on the left-hand side is in $\Kk l(\hat T)$. This should not be a surprise because of \leref{Z cell as teta}. 
Thus the 2-cell $Z_{\hat u,\hat v}$ in $\Dd$ presents a 2-cell in $\Kk l(\hat T)$ that is a companion-lift of $\theta^{u,v}=\Id$ in 
$\Kk l(\hat T)$. 
\qed\end{proof}

By construction of the action in \thref{strength-action} we defined the vertical modifications $\tilde p=K(p), \tilde l=K(l), \tilde m=K(m)$ 
in $\Kk l(\hat T)$ in terms of the vertical modifications $p,l,m$ from the monoidal structure of $\Dd$. We will now draw the 2-cell components of these diagrams, composing to them the identity 2-cells $\theta^{u,v}: u\rtr K(v) \Rightarrow K(u\ot v)$ in $\Kk l(\hat T)$ for the corresponding varying 1v-cells $u,v$ at suitable places. (Observe that all these 1v-cells are either identities, or 1v-cell components of $\alpha, \lambda, \rho$, which are liftable by assumption.)
To suit the space in the diagrams, instead of writing $K(u)$ for 1v-cells on the arrows of these diagrams, we will write $\crta u$. We get that 
\begin{equation} \eqlabel{l,m,p-vert-teta}
\scalebox{0.9}{
\bfig
\putmorphism(-30,500)(1,0)[` `=]{400}1a
 \putmorphism(-40,500)(0,-1)[` `\crta{\alpha_{I,A,B}}]{400}1l
 \putmorphism(-20,500)(0,-1)[` `]{400}0l
\putmorphism(-40,110)(0,-1)[\phantom{Y_2}``]{400}1l
\putmorphism(-20,110)(0,-1)[\phantom{Y_2}``\crta{\lambda_{AB}}]{400}0l
\putmorphism(350,500)(0,-1)[` `]{800}1l
\putmorphism(350,400)(0,-1)[` `\crta{\lambda_AB}]{800}0l
\put(0,240){\fbox{$K(l_{AB})$}}
\putmorphism(350,500)(1,0)[` `=]{400}1a
\putmorphism(750,500)(0,-1)[` `]{800}1l
\putmorphism(770,400)(0,-1)[` `\lambda_A\w\rtr\w\crta B]{800}0l
\putmorphism(-30,-300)(1,0)[` `=]{400}1b
\putmorphism(350,-300)(1,0)[` `=]{400}1b
\put(450,240){\fbox{$\theta^{u,v}$}}
\efig}
\qquad
\scalebox{0.9}{
\bfig
\putmorphism(-510,500)(1,0)[` `=]{420}1a
 \putmorphism(-500,500)(0,-1)[` `]{800}1r
 \putmorphism(-520,670)(0,-1)[` `\rho_A\w\rtr\w\crta B]{800}0r
\putmorphism(-510,-300)(1,0)[` `=]{420}1b
\put(-480,-90){\fbox{$(\theta^{u,v})^{-1}$}}
\putmorphism(-90,500)(1,0)[` `=]{420}1a
 \putmorphism(-100,500)(0,-1)[` `]{800}1r
 \putmorphism(-90,670)(0,-1)[` `\crta{\rho_AB}]{800}0r
\putmorphism(350,500)(0,-1)[` `\crta{\alpha_{A,I,B}}]{400}1r
\putmorphism(350,100)(0,-1)[` `]{400}1l
\putmorphism(370,40)(0,-1)[` `\crta{A\lambda_B}]{400}0l
\putmorphism(-90,-300)(1,0)[` `=]{400}1b
\put(-60,70){\fbox{K($m_{AB})$}}
\putmorphism(350,100)(1,0)[` `=]{400}1a
\putmorphism(750,100)(0,-1)[` `]{400}1l
\putmorphism(770,40)(0,-1)[` `A\w\rtr\w\crta{\lambda_B}]{400}0l
\putmorphism(350,-300)(1,0)[` `=]{420}1b
\put(450,-40){\fbox{$\theta^{u,v}$}}
\efig}
\qquad
\scalebox{0.9}{
\bfig
\putmorphism(-110,500)(1,0)[` `=]{400}1a
 \putmorphism(-120,500)(0,-1)[` `]{400}1l
 \putmorphism(-100,500)(0,-1)[` `\crta{\alpha_{AB,C,D}}]{400}0l
\putmorphism(-120,110)(0,-1)[\phantom{Y_2}``\crta{\alpha_{A,B,CD}}]{400}1l
\putmorphism(-100,110)(0,-1)[\phantom{Y_2}``]{400}0l
\putmorphism(380,500)(0,-1)[` `]{340}1l
\putmorphism(400,500)(0,-1)[` `\crta {\alpha_{A,B,C}D}]{340}0l
\putmorphism(380,230)(0,-1)[` `\crta{\alpha_{A,BC,D}}]{340}1r
\putmorphism(380,-20)(0,-1)[` `]{340}1l
\putmorphism(400,0)(0,-1)[` `\crta {A\alpha_{B,C,D}}]{340}0l
\putmorphism(-110,-300)(1,0)[` `=]{380}1b
\put(-90,70){\fbox{$K(p_{ABCD})$}}
\putmorphism(350,500)(1,0)[` `=]{400}1a
\putmorphism(750,500)(0,-1)[` `]{340}1r
\putmorphism(730,500)(0,-1)[` `\alpha\w\rtr\w\crta D]{340}0r
\putmorphism(750,-20)(0,-1)[` `]{330}1r
\putmorphism(730,-20)(0,-1)[` `A\w\rtr\w\crta\alpha]{330}0r
\put(470,340){\fbox{$\theta^{u,v}$}}
\put(470,-190){\fbox{$\theta^{u,v}$}}
\putmorphism(380,200)(1,0)[` `=]{400}1a
\putmorphism(380,-50)(1,0)[` `=]{400}1a
\putmorphism(380,-300)(1,0)[` `=]{400}1b
\efig}
\end{equation}
equal respectively $\tilde l_{AB}, \tilde m_{AB}, \tilde p_{ABCD}$.

By definition these are then respectively equal to the 2-cell components $\tilde l_{AB}, \tilde m_{AB}, \tilde p_{ABCD}$ in 
$\Kk l(\hat T)$. 
Applying companion-lifting to these three equalities of 2-cells we obtain the following equalities of the vertical compositions: 
\begin{equation} \eqlabel{p,l,m-teta}
\frac{\theta_{\hat\lambda_A, B}}{K(\hat l_{AB})}=\widehat{\tilde l_{AB}}, \qquad 
\threefrac{[\Id_{K(\hat\alpha_{A,I,B})} \vert\theta_{A,\hat\lambda_B}]}{K(\hat m_{AB})}{(\theta_{\hat\rho_A, B})^{-1}}=\widehat{\tilde m_{AB}}, \qquad 
\frac{[\theta_{K(\hat\alpha_{A,B,C}),D} \vert K(\hat\alpha_{A,BC,D})\vert\theta_{A,K(\hat\alpha_{B,C,D})}]}{K(\hat p_{ABCD})}=
\widehat{\tilde p_{ABCD}}
\end{equation}
whereby note that $\widehat{K(u)}=K(\hat u)$ for 1v-cells $u$ by \prref{Kl compan}, and that the companion-lift of $\theta^{u,v}$ in 
 $\Kk l(\hat T)$ is $\theta_{\hat u, \hat v}$ by the second part of \prref{icon teta}. 

\medskip

In the setting of  \thref{strength-action} we have two more identities that hold true. One connects $\hat\alpha, \widehat{\tilde\alpha}, \theta$, and another $\hat\lambda, \widehat{\tilde\lambda}, \theta$. We discuss them next. 

\smallskip

We comment first in more detail why the equality of 2-cells 
\begin{equation} \eqlabel{lambda-teta}
\frac{[\theta_{1_I,f}\vert\Id_{K(\hat\lambda_{A'})}]}{K(\hat\lambda_f)}=\widehat{\tilde\lambda_{K(f)}}: \quad
\widehat{\tilde\lambda_{A'}}(I\rtr K(f))\Rightarrow K(f)\widehat{\tilde\lambda_A}
\end{equation} 
holds, for any 1h-cell $f$ in $\Dd$, whereby $K(\hat\lambda_A)=\widehat{\tilde\lambda_A}$ (by \prref{Kl compan} and 
by construction $K(\lambda_A)=\tilde\lambda_A$.) 
The 2-cell $\theta_{1_I,f}$ is comprised basically of the 2-cell $\hat z_{I,A}$. Next, $K(\hat\lambda_f)$ is a companion-lift in 
$\Kk l(\hat T)$ of the 2-cell $K(\lambda_f)$, which is the horizontal composition $[\lambda_f\vert\hat\eta^{\lambda_A}]$. 
Finally, $\widehat{\tilde\lambda_{K(f)}}$ is a companion-lift in $\Kk l(\hat T)$ of the 2-cell $\tilde\lambda_{K(f)}$ in 
$\Kk l(\hat T)$, which by \thref{strength-action} is given by the horizontal composition $[\lambda_{K(f)}\vert x^*_A]$, whereby 
$\lambda_{K(f)}$ consists of the composition $[\lambda_f\vert\lambda_{\hat\eta_A}]$. (All the 2-cells in these horizontal compositions are square-formed.) After canceling out the equal parts stemming from companion-lifting (the structure 
2-cells $\iota$ and $e$) and the 2-cell $\lambda_f$, for the desired equality to hold it remains to prove 
\begin{equation} \eqlabel{z-x}
\scalebox{0.86}{
\bfig
 \putmorphism(370,400)(1,0)[\phantom{F(A)}`I T(A) `I\hat\eta(A)]{600}1a
 \putmorphism(1000,400)(1,0)[\phantom{F(A)}`T(IA) `\hat t_{I, A}]{580}1a
 \putmorphism(450,00)(1,0)[`T(IA)`\hat\eta(IA)]{1130}1b
\putmorphism(410,400)(0,-1)[IA`IA `=]{400}1l
\putmorphism(1570,400)(0,-1)[\phantom{Y_2}``=]{400}1r
\put(870,180){\fbox{$\hat z_{I, A}$}}
\putmorphism(410,0)(0,-1)[`A `\lambda_A]{400}1l
 \putmorphism(450,-400)(1,0)[`T(IA)`\hat\eta(IA)]{1130}1b
\putmorphism(1570,0)(0,-1)[\phantom{Y_2}``T(\lambda_A)]{400}1r
\put(870,-280){\fbox{$\hat\eta^{\lambda_A}$}}
\efig}
=
\scalebox{0.86}{
\bfig
\putmorphism(340,250)(0,-1)[IA`A `\lambda_A]{450}1l
\putmorphism(810,250)(0,-1)[IT(A)`T(A) `]{450}1r
\putmorphism(790,250)(0,-1)[` `\lambda_{T(A)}]{450}0r
\putmorphism(1270,250)(0,-1)[T(IA)`T(A).`T(\lambda_A)]{450}1r
\putmorphism(370,250)(1,0)[``I\hat\eta_A]{330}1a
 \putmorphism(900,250)(1,0)[` `\hat t_{I,A}]{240}1a
\putmorphism(360,-200)(1,0)[``\hat\eta_A]{370}1b
 \putmorphism(900,-200)(1,0)[``=]{270}1b
\put(450,10){\fbox{$\lambda_{\hat\eta_A}$}}
\put(1080,70){\fbox{$x^*$}}
\efig}
\end{equation}

\begin{rem}
Note that both $K(\hat\lambda_f)$ and $\widehat{\tilde\lambda_{K(f)}}$ contain $\lambda_f$ as a component in the composition constituting them, and that it is being cancelled out. Thus the question of liftability does not affect the 1h-cell $f$. 
On the other hand, $\theta_{1_I,f}=[\Id_{1_I\ot f}\vert\hat z_{I,A'}]$, so that the only non-trivial 2-cell depending on $f$ in 
\equref{lambda-teta} is $\lambda_f$ that we already discussed.
\end{rem}

To obtain \equref{z-x} start from the trivially holding identity between vertical modifications $z$ and $x$ that gives rise to the $z\x x$ axiom holding in \deref{hor-strength} for $\hat T$. Namely, we have
\begin{equation} \eqlabel{vert z-x}
\scalebox{0.9}{
\bfig
\putmorphism(-30,500)(1,0)[\phantom{IA}`IA `=]{400}1a
 \putmorphism(-40,500)(0,-1)[IA`A `\lambda_A]{400}1l
\putmorphism(-40,110)(0,-1)[\phantom{Y_2}`T(A)`\eta_A]{500}1l
\putmorphism(380,500)(0,-1)[` T(IA)`]{560}1l
\putmorphism(400,500)(0,-1)[` `\eta_A]{600}0l
\putmorphism(750,500)(0,-1)[`IT(A) `I\eta_A]{300}1r
\putmorphism(750,230)(0,-1)[`T(IA) `t_{I,A}]{280}1r
\putmorphism(380,-50)(0,-1)[` `T(\lambda_A)]{340}1r
\putmorphism(-30,-400)(1,0)[\phantom{T(A)}` T(A)`=]{380}1a
\put(20,100){\fbox{$\eta^{\lambda_A}$}}
\putmorphism(400,500)(1,0)[`IA `=]{380}1a
\put(470,320){\fbox{$z_{I,A}$}}
\putmorphism(500,-50)(1,0)[` `=]{140}1a
\efig}
=
\scalebox{0.9}{
\bfig
\putmorphism(-30,500)(1,0)[\phantom{IA}`IA `=]{400}1a
 \putmorphism(-40,500)(0,-1)[IA`A `\lambda_A]{400}1l
\putmorphism(-40,110)(0,-1)[\phantom{Y_2}`T(A)`\eta_A]{500}1l
\putmorphism(380,500)(0,-1)[`IT(A) `I\eta_A]{300}1l
\putmorphism(750,230)(0,-1)[`T(IA) `t_{I,A}]{300}1r
\putmorphism(380,230)(0,-1)[` `\lambda_{T(A)}]{640}1l
\putmorphism(-30,-400)(1,0)[\phantom{T(A)}` T(A)`=]{380}1a
\put(40,100){\fbox{$\lambda^{\eta_A}$}}
\putmorphism(750,-40)(0,-1)[` `T(\lambda_A)]{380}1r
\put(470,-190){\fbox{$x_A$}}
\putmorphism(500,200)(1,0)[`IT(A) `=]{300}1a
\putmorphism(440,-400)(1,0)[` T(A).`=]{300}1a
\efig}
\end{equation}
Then compose on the left the companion structure 2-cell $\iota_{\eta_A}$ and on the right $e_{I\eta_A}$ and $e_{t_{I,A}}$, 
and also insert the identity 2-cell $\frac{\iota_{\eta_A}}{e_{\eta_A}}=\Id$ between $z$ and $\eta^{\lambda_A}$ on the left-hand 
side to obtain \equref{z-x}. 

\medskip

By similar arguments, using trivial identity holding between vertical modifications $z$ and $y$ that gives rise to the $z\x y$ axiom in 
\deref{hor-strength}, one proves the equality of 2-cells
\begin{equation} \eqlabel{alfa-teta}
\frac{[\theta_{fg,h}\vert\Id_{K(\hat\alpha_{A',B',C'})}]}{K(\hat\alpha_{f,g,h})}=
\threefrac{\widehat{\tilde\alpha_{f,g,K(h)}}}{[\Id_{\hat\alpha_{A,B,C}}\vert f\rtr\theta_{g,h}]}{[\Id_{\hat\alpha_{A,B,C}}\vert 
\theta_{f,gh}]}: \quad \hat\alpha_{A',B',C'}\big( (fg)K(h)\big) \Rightarrow K(f(gh))\hat\alpha_{A,B,C},
\end{equation} 
whereby $K(\hat\alpha_{A,B,C})=\widehat{\tilde\alpha_{A,B,C}}$. 

\medskip

Given that all the 2-cells appearing in equations \equref{p,l,m-teta}, \equref{lambda-teta}, \equref{alfa-teta} 
are horizontally globular, the diagrammatic representation of the equations corresponding to them in $\HH(\Kk l(\hat T))$ 
is precisely the one appearing in \cite[Definition 11]{HF} where ``an extension of the canonical action of a bicategory on itself'' is defined (in \cite[Definition 20]{HF1} it is called ``0-strict morphism of actions''). This brings us to the following. 

\begin{defn}
Let $T$ be a horizontal double monad on a horizontally monoidal double category $(\Dd, \ot, \alpha, \lambda, \rho; p,l,m)$. 
A {\em horizontal extension of the canonical action of $\Dd$ on itself} is a pair $(\rtr, \theta)$, where \vspace{-0,2cm}
\begin{enumerate}
\item $\rtr:\Dd\times\Kk l(T)\to\Kk l(T)$ 
is a horizontal action with structure pseudonatural transformations $\tilde\alpha, \tilde\lambda$ and structure horizontal modifications 
$\tilde p, \tilde m, \tilde l$, and \vspace{-0,2cm}
\item $\theta:-\rtr K(-)\Rightarrow K(-\ot-)$ is a horizontal icon, 
\end{enumerate} 
so that $\theta_{1_A,1_B}=\rtr^0_{A,B}$ (the unitor of the action on $A\rtr B$) and the identities 
\begin{equation} \eqlabel{hor-p,l,m-teta}
\frac{\theta_{\lambda_A, B}}{K(l_{AB})}=\tilde l_{AB}, \qquad 
\threefrac{[\Id_{K(\alpha_{A,I,B})} \vert\theta_{A,\lambda_B}]}{K(m_{AB})}{(\theta_{\rho_A, B})^{-1}}=\tilde m_{AB}, \qquad 
\frac{[\theta_{K(\alpha_{A,B,C}),D} \vert K(\alpha_{A,BC,D})\vert\theta_{A,K(\alpha_{B,C,D})}]}{K(p_{ABCD})}=
\tilde p_{ABCD}
\end{equation}
\begin{equation} \eqlabel{hor-lambda-teta}
\frac{[\theta_{1_I,f}\vert\Id_{K(\lambda_{A'})}]}{K(\lambda_f)}=\tilde\lambda_{K(f)}: \quad
\tilde\lambda_{A'}(I\rtr K(f))\Rightarrow K(f)\tilde\lambda_A
\end{equation} 
\begin{equation} \eqlabel{hor-alfa-teta}
\frac{[\theta_{fg,h}\vert\Id_{K(\alpha_{A',B',C'})}]}{K(\alpha_{f,g,h})}=
\threefrac{\tilde\alpha_{f,g,K(h)}}{[\Id_{\alpha_{A,B,C}}\vert f\rtr\theta_{g,h}]}{[\Id_{\alpha_{A,B,C}}\vert 
\theta_{f,gh}]}: \quad \alpha_{A',B',C'}\big( (fg)K(h)\big) \Rightarrow K(f(gh))\alpha_{A,B,C}
\end{equation} 
hold, whereby $K(\lambda_A)=\tilde\lambda_A$ and $K(\alpha_{A,B,C})=\tilde\alpha_{A,B,C}$. 
\end{defn}

Recall from \leref{vert icon} that $\theta$ is a strict vertical icon if and only if $\hat\theta$ 
is a horizontal icon with trivial (square-formed) 2-cell components $\hat\theta^u$. 
Our above findings we can then formulate into the following claim. 

\begin{thm} \thlabel{strength-ext}
Let $T$ be a vertical double monad on a monoidal double category $\Dd$ and 
assume that the vertical transformations $\alpha,\lambda,\rho$ of $\Dd$, and $\mu,\eta$ of $T$ are liftable. 
A vertical liftable (left) strength $t$ on $T$ that is invertible as a vertical transformation 
induces a horizontal extension of the canonical action of $\Dd$ on itself, 
where the action and the icon are given as in \thref{strength-action} and \prref{icon teta}. 
\end{thm}

\begin{rem}
One can also define a vertical extension for a vertical double monad on a monoidal double category as a pair consisting of 
a vertical action and a strict vertical icon, so that the three identities \equref{l,m,p-vert-teta} hold and $\theta_{1_A,1_B}=
\rtr^0_{A,B}$. 
Nevertheless, 
for a general action not necessarily induced by a strength, it is not clear if one can lift a vertical extension into a horizontal one. 
Namely, if for liftable $\alpha, \lambda, \rho$ of $\Dd$ and $\mu,\eta$ of $T$ we had a vertical action (and hence extension 
$(\rtr,\theta)$), we would get an induced horizontal icon $\hat\theta$, 
but we would not be able to recover \equref{lambda-teta} and \equref{alfa-teta} as lifts of vertically globular 2-cells. This is so 
because we can not 
factor $\tilde\alpha_{f,g,h}, \tilde\lambda_f$ and $\theta_{f,g}$ to be able to compare them with 
$\alpha_{f,g,h}, \lambda_f$ and $\theta_{A,B}$ and eventually cancel out $\alpha_{f,g,h}$ and $\lambda_f$ in those expressions to get to a companion-lift of some trivially holding identity among vertical structures, as we did in \thref{strength-action}. 
\end{rem}

In a bicategorical case, in Table 1 on p.16 of \cite{FH} the correspondence is listed between the underlying 
axioms of left strengths, on one hand, and of extensions of the canonical actions, on the other hand. We present that Table here: 
observe that the axioms in the left column are both bicategorical and double categorical, as they are expressed in terms of globular 2-cells. 
\begin{table}[H]
\begin{center}
\begin{tabular}{ c c } 
(horizontal) strength & (horizontal) extension \\ [0.5ex]
\hline
three axioms for $x,y$ & three axioms \equref{hor-p,l,m-teta} \\ [0.3ex]   
in \deref{hor-strength}, b) & \\ [1ex] 
first three axioms for $w,z$ & pseudofunctor axioms for $\rtr$ \\ [0.3ex] 
in \deref{hor-strength}, c) &   \\ [1ex] 
middle two axioms relating $w,z$ & (h.) ps.nat.tr. axioms for $\tilde\lambda$ of $\rtr$ \\ [0.3ex]  
with $x$ in \deref{hor-strength}, c)  & \\ [1ex]
last two axioms relating $w,z$ & (h.) ps.nat. tr. axioms for $\tilde\alpha$ of $\rtr$ \\ [0.3ex]   
with $y$ in \deref{hor-strength}, c) & \\ [1ex]
\end{tabular}
\caption{Correspondence between axioms of (h.) strengths and (h.) extensions}
\label{table:6}
\end{center}
\end{table} \vspace{-0,4cm}
\noindent As far as the right-hand side column applied to the double categorical situation concerns, the axioms in the first row are globular, 
whereas pseudodouble functor axioms and 
axioms of horizontal pseudonatural transformations additionally involve vertical 1-cells. 
This added complexity in the axioms with respect to the bicategorical situation 
corresponds to the same sort of change in the axioms for the underlying given data on the side of a strength: 
horizontal transformations (\axiomref{h.o.t.-3} - \axiomref{h.o.t.-5} for $\hat t$) and horizontal modifications (\axiomref{m.ho-vl.-1} and \axiomref{m.ho-vl.-2} for $\hat w, \hat z$ and for $x^*, y^*$ {\em i.e.} $\hat x, \hat y$). 
Consequently, the correspondence listed in the above Table applies equally to the double categorical case. 

\smallskip

Let us summarize the structure of our above proof that a vertical strength yields a (horizontal) extension, in order to 
conclude what the gain is of working with a {\em vertical} structure in a double categorical setting and by assuming the existence of certain companions. 
Namely, starting from a {\em vertical} strength we obtained: 1) a {\em horizontal} strength (due to \thref{strenghts thm}) - this is a double categorical version of a bicategorical strength, 2) the vertical action - as in \thref{strength-action}, 3) the horizontal action (in 
\thref{strength-action}, obtained by lifting the vertical action), 4) an icon - as in \prref{icon teta} and 5) the five extension axioms -  we obtained them for free 
by companion-lifting of trivially holding identities between vertical structures (so they follow essentially by \prref{essence}).

\subsubsection{Horizontal strengths induce horizontal actions on the Kleisli double category} \ssslabel{h.str.induce act}

Let $s$ be a horizontal strength on a horizontal double monad $S$ on a horizontally monoidal double category $(\Dd, \ot, \alpha, 
\lambda, \rho)$. By \leref{vert icon} and the first part of \prref{icon teta} we indeed have: 

\begin{prop} \prlabel{horiz str-icon}
A horizontal left strength $s$ on a horizontal double monad $S$ on a horizontally monoidal double category $\Dd$ 
induces a horizontal icon $\theta:-\rtr K(-)\Rightarrow K(-\ot-): 
\Dd\times\Dd\to\Kk l(S)$ whose invertible 2-cell components $\theta_{f,g}$ for 1h-cells $f:A\to A', g:B\to B'$ are given by 
$$\theta_{f,g}=
\scalebox{0.82}{
\bfig
 \putmorphism(-650,210)(1,0)[AB` A'B'`f\ot g]{640}1a 
\putmorphism(-660,200)(0,-1)[\phantom{Y_2}` `=]{380}1l
\putmorphism(-50,200)(0,-1)[\phantom{Y_2}` `]{380}1r
\putmorphism(-50,200)(0,-1)[\phantom{Y_2}` `=]{380}0r
\putmorphism(100,210)(1,0)[``A'\w\ot\eta_{B'}]{410}1a
\putmorphism(680,210)(1,0)[A'S(B')`S(A'\ot B')`s_{A',B'}]{680}1a

\putmorphism(1600,200)(0,-1)[\phantom{Y_2}``]{380}0r
\putmorphism(1370,200)(0,-1)[\phantom{Y_2}` `=]{380}1r
\putmorphism(-670,-200)(1,0)[AB `A'B'`f\ot g]{620}1a 
\putmorphism(20,-200)(1,0)[\phantom{Y}` S(A'\ot B')  `\eta_{A'\ot B'}]{1360}1a

\put(-500,30){\fbox{$\Id_{f\ot g}$}}
\put(580,30){\fbox{$z_{A',B'}$}}
\efig}
$$
and the square-formed 2-cell components $\theta^{u,v}=\Id^{u\ot v}$ for 1v-cells $u,v$. 
\end{prop}

By substituting in the proof of \thref{strength-action} the companion-lifts $\hat T, \hat t, \hat\alpha, \hat\lambda, 
y^*, x^*$ by their horizontal analogues $S,s,\alpha, \lambda, y, x$, that proof is easily adapted to the proof of the following fact:

\begin{thm} \thlabel{hor-str-act}
A horizontal strength $s$ on a horizontal double monad $S$ on a horizontally monoidal double category $(\Dd, \ot, \alpha, \lambda, \rho)$ induces an action $\rtr:\Dd\times\Kk l(S)\to\Kk l(S)$.  
\end{thm}
 
\begin{proof}
The action functor $\rtr$ in the horizontal context is defined similarly as in \thref{strength-action}. 
The action constraints $\tilde\alpha, \tilde\lambda$ 
are defined on objects and 1v-cells as $\alpha$ and $\lambda$, and on 1h-cells as  
$$\tilde\alpha_{f,g,h}=\frac{[\Id_{(fg)h}\,\vert\,y_{A',B',C'}]}{[\alpha_{f,g,h}\,\vert\,\Id_{[A's_{B',C'}\vert s_{A',B'C'}]}]}
\quad\text{ and }\quad
\tilde\lambda_h=\frac{[\Id_{Ih}\vert  x_{A'}]}{\lambda_h},$$
whereby now $\alpha_{f,g,h}$ and $\lambda_h$ are globular 2-cells. To prove that $\tilde\alpha$ is a horizontal 
transformation, observe that the axioms \axiomref{h.o.t.-1} and \axiomref{h.o.t.-2} hold as in the bicategorical case of \cite{FH}, 
and that \axiomref{h.o.t.-3} and \axiomref{h.o.t.-4} hold because they hold for $\alpha$. One only needs to work out the axiom 
\axiomref{h.o.t.-5}. It is proved to hold applying the modification axiom \axiomref{m.ho-vl.-2} for $y$ and \axiomref{h.o.t.-5} for $\alpha$. 
The proof that $\tilde\lambda$ is a horizontal transformation is similar. 

The horizontal modifications $\tilde l, \tilde m, \tilde p$ are defined as in \equref{hor-p,l,m-teta} with $\theta$ from 
\prref{horiz str-icon}, whereby now $l,m,p$ are horizontal modifications of the horizontally monoidal structure on $\Dd$. The modification axiom \axiomref{m.ho-vl.-1} holds as in the bicategorical case. For the other modification axiom \axiomref{m.ho-vl.-2}, for $\tilde l$ use \axiomref{h.o.t.-5} for $\theta$ and \axiomref{m.ho-vl.-2} of the horizontal modification $l$. The proofs for $\tilde m$ and $\tilde p$ work by similar arguments. 

The remaining work required to prove \thref{hor-str-act} 
is the verification of the action axioms for the horizontal modifications $\tilde l, 
\tilde m, \tilde p$. This verification coincides with the one carried out in the bicategorical case in \cite[Proposition 2]{HF}, 
\cite[Proposition 7.1]{HF2}, given that those axioms are comprised of globular 2-cells, which are the 2-cell components of horizontal modifications. 
\qed\end{proof}

Note that for vertical $T, t$ in the setting of \thref{strength-action} and such that $S=\hat T, s=\hat t$ 
the horizontal action structures 
defined in the above theorem are companion-lifts of the type \equref{square-lift} of 
the corresponding action structures from \thref{strength-action}. 
In other words, the horizontal action that $\hat t$ induces as in \thref{hor-str-act} 
is the companion-lift $(\rtr; \widehat{\tilde\alpha}, \widehat{\tilde\lambda})$ of the vertical action 
$(\rtr; \tilde\alpha, \tilde\lambda)$ obtained in \thref{strength-action}.

\subsubsection{Particularities of actions on the Kleisli double category} 

Given that cells in a Kleisli double category are represented via cells in the underlying double category, it is necessary to 
settle down some technical details regarding actions $\rtr$ of $\Dd$ on $\Kk l(T)$ for a horizontal double monad $T$. 
Let us first fix some notation. 
The lax pseudodouble functor structure of the action $\rtr$ we will denote by $\rtr^2, \rtr^0$. 
Recall that for 1h-cells in $\Kk l(T)$ we write 
$f^K:A\to A'$ 
and that they are given by 1h-cells $f:A\to T(A')$ in $\Dd$. For the 1h-cell unit on an object $A$ in $\Kk l(T)$ we will write 
$1^K_A: A\to A$ and we know that it is given by $\eta_A:A\to T(A)$ in $\Dd$. In contrast, we will denote by 
$1_{T(A)}^\bullet:T(A)\to A$ a 1h-cell in $\Kk l(T)$ 
given by $1_{T(A)}$ in $\Dd$. 
Then for 1h-cells $f:A\to A'$ in $\Dd$ and $g^K:B\to B'$ in $\Kk l(T)$ we have that $(f\rtr g^K)^K: A\rtr B\to A'\rtr B'$ in $\Kk l(T)$ 
is given by a 1h-cell $A\rtr B\to T(A'\rtr B')$ in $\Dd$ that can be understood as the composition 
\begin{equation} \eqlabel{rtr explain}
f\rtr g^K=\big(A\rtr B \stackrel{f\rtr g}{\longrightarrow} A'\rtr T(B') \stackrel{1_{A'}\rtr 1^\bullet_{T(B')}}{\xrightarrow{\hspace{1cm}}} T(A'\rtr B')\big).
\end{equation}
Observe that in the component $f\rtr g$ the 1h-cells $f$ and $g$ in $\Dd$ act componentwise, whereas $1_{A'}\rtr 1^\bullet_{T(B')}$ is 
a 1h-cell in $\Dd$ determining a 1h-cell $(1_{A'}\rtr 1^\bullet_{T(B')})^K: A'\rtr T(B')\to A'\rtr B'$ in $\Kk l(T)$  
acting componentwise. 
(By the notation $1_{A'}\rtr 1^\bullet_{T(B')}$ we keep an analogy with the notation in \cite[Proposition 4.3]{UV} and \cite[Theorem 1]{FH}.)

\begin{rem}
Furthermore, observe that the component $f\rtr g$ in \equref{rtr explain} 
behaves {\em in a strict way}, in the sense that $f'f\rtr g'g=(f'\rtr g')(f\rtr g): 
A\rtr B \stackrel{f\rtr g}{\longrightarrow} A'\rtr T(B') \stackrel{f'\rtr T(g')}{\longrightarrow} A''\rtr T^2(B'') 
\stackrel{A''\mu_{B''}}{\longrightarrow} A''\rtr T(B'')$, and we have that $1_A\rtr 1_{T(B)}:A\rtr T(B)\to A\rtr T(B)$ is precisely 
$1_{A\rtr T(B)}$ (set $g^K=1^\bullet_{T(B)}$ in \equref{rtr explain}). 
\end{rem}

Similarly, for 2-cells 
$\scalebox{0.86}{
\bfig
\putmorphism(-150,150)(1,0)[A` A'`f]{450}1a
\putmorphism(-150,-150)(1,0)[\tilde A`\tilde A' `\tilde f]{440}1b
\putmorphism(-150,150)(0,-1)[\phantom{Y_2}``u]{300}1l
\putmorphism(280,150)(0,-1)[\phantom{Y_2}``u']{300}1r
\put(0,-25){\fbox{$\sigma$}}
\efig}$
in $\Dd$ and 
$\scalebox{0.86}{
\bfig
\putmorphism(-150,150)(1,0)[B` B'`g^K]{450}4a
\putmorphism(-150,-150)(1,0)[\tilde B`\tilde B' `\tilde g^K]{440}4b
\putmorphism(-150,150)(0,-1)[\phantom{Y_2}``v]{300}1l
\putmorphism(280,150)(0,-1)[\phantom{Y_2}``v']{300}1r
\put(0,-25){\fbox{$\delta^K$}}
\efig}$ 
in $\Kk l(T)$ we have that the 2-cell $(\sigma\rtr\delta^K)^K$ is given as the composition 
\begin{equation} \eqlabel{s-u,v}
\sigma\rtr\delta^K=
\scalebox{0.86}{
\bfig
\putmorphism(40,250)(0,-1)[` `u\rtr u']{400}1l
\putmorphism(640,250)(0,-1)[` `]{400}1r
\putmorphism(620,250)(0,-1)[` `v\rtr T(v')]{400}0r
\putmorphism(1700,250)(0,-1)[T(A'\rtr B')``T(v\rtr v')]{400}1r
\putmorphism(70,250)(1,0)[A\rtr B`A'\rtr T(B')`f\rtr g]{630}1a
 \putmorphism(900,250)(1,0)[` `1_{A'}\rtr 1^\bullet_{T(B')}]{600}1a
\putmorphism(60,-150)(1,0)[\tilde A\rtr\tilde B`\tilde A'\rtr T(\tilde B') `f'\rtr g']{670}1b
 \putmorphism(940,-150)(1,0)[`T(\tilde A'\rtr\tilde B')`1_{\tilde A'}\rtr 1^\bullet_{T(\tilde B')}]{750}1b
\put(200,50){\fbox{$\sigma\rtr\delta$}}
\put(1080,40){\fbox{$\Id^v\rtr\Id^{T(v')\bullet}$}}
\efig}
\end{equation}  
in $\Dd$, where $\Id^{T(v)\bullet}: T(v)\Rightarrow v$ in $\Kk l(T)$ is a 2-cell given by the identity 2-cell $\Id^{T(v)}$ in $\Dd$.  

\medskip

For the compositor and unitor we have 
\begin{equation} \eqlabel{compositor}
\scalebox{0.8}{
\bfig
 \putmorphism(60,200)(1,0)[AB` T(A'B') `f\w\rtr\w g^K]{650}1a 
 \putmorphism(710,200)(1,0)[\phantom{A''\ot B'}` T^2(A''B'') `T(f'\w\rtr\w g'^K)]{980}1a 
 \putmorphism(1720,200)(1,0)[\phantom{A''\ot B'}`T(A''B'')  `\mu_{A'',B''}]{800}1a 
\putmorphism(60,200)(0,-1)[\phantom{Y_2}``=]{400}1l
\putmorphism(2540,200)(0,-1)[\phantom{Y_2}``=]{400}1r 
 \putmorphism(60,-200)(1,0)[AB` T(A''B'') `f'f\w\rtr\w g'^Kg^K]{2460}1b 
\put(1160,-20){\fbox{$\rtr^2_{f'f,g'^Kg^K}$}}
\putmorphism(60,-200)(0,-1)[\phantom{Y_2}``=]{400}1l
\putmorphism(2540,-200)(0,-1)[\phantom{Y_2}``=]{400}1r 
 \putmorphism(60,-620)(1,0)[AB` A''\w\rtr\w T(B'') `f'f\w\rtr\w g'g]{1160}1b 
 \putmorphism(1430,-600)(1,0)[` T(A''B'') `1_{A''}\w\rtr\w 1^\bullet_{B''}]{1090}1b 
\put(1240,-470){\fbox{$\equref{rtr explain}$}}
\efig}
\end{equation}
and
\begin{equation} \eqlabel{unitor}
\scalebox{0.8}{
\bfig
 \putmorphism(370,200)(1,0)[\phantom{F(A)}`A T(B) ` 1_A\rtr\eta(B)]{600}1a
 \putmorphism(1000,200)(1,0)[\phantom{F(A)}`T(AB) `1_A\rtr 1_{T(B)}^\bullet]{620}1a
 \putmorphism(470,-100)(1,0)[`T(AB)`1_A\rtr 1_B^K]{1110}1b
 \putmorphism(470,-500)(1,0)[`T(AB).`\eta(AB)]{1150}1b
\putmorphism(410,200)(0,-1)[AB`AB `=]{300}1l
\putmorphism(410,-100)(0,-1)[`AB `=]{400}1l
\putmorphism(1610,200)(0,-1)[\phantom{Y_2}``=]{300}1r
\putmorphism(1610,-100)(0,-1)[\phantom{Y_2}``=]{400}1r
\put(870,0){\fbox{$\equref{rtr explain}$}}
\put(870,-370){\fbox{$\rtr^0_{A,B}$}}
\efig}
\end{equation}
The reader may have noticed that the action constructed in \thref{strength-action} fits the above presented framework.

\subsubsection{From an extension to a horizontal strength}

In this final part of the current subsection we are going to prove that a (horizontal) extension induces a horizontal strength.  
We will also explain why we can not construct a vertical strength even though if we start from a monoidal double category with liftable 
monoidality constraints $\alpha, \lambda, \rho$, and a vertical double monad with liftable monad structures $\mu, \eta$. 

For this reason we start from a horizontally monoidal double category $(\Dd, \ot, \alpha, \lambda, \rho; \\ p,l,m)$ and a horizontal double monad 
$(T, \mu, \eta)$. We assume that there is a horizontal extension of the canonical action of $\Dd$ on itself 
given by: \vspace{-0,2cm}
\begin{itemize}
\item a horizontal action, that is a pseudodouble functor $\rtr:\Dd\times\Kk l(T)\to\Kk l(T)$ with 
horizontal equivalences $\tilde\lambda,\tilde\alpha$, and horizontal modifications $\tilde p, \tilde l, \tilde m$ 
satisfying the appropriate axioms, and \vspace{-0,2cm} 
\item a horizontal icon $\theta:-\rtr K(-)\Rightarrow K(-\ot-)$ \vspace{-0,2cm}
\end{itemize}
so that $\theta_{1_A,1_B}=\rtr^0_{A,B}$, the identities \equref{hor-p,l,m-teta}, \equref{hor-lambda-teta} and \equref{hor-alfa-teta} 
hold, and it is $K(\alpha_{A,B,C})=\tilde\alpha_{A,B,C}$ and $K(\lambda_A)=\tilde\lambda_A$.

The construction of a horizontal transformation $s$ out of the above pair $(\rtr, \theta)$ is technically involved. 
In the categorical case in 
\cite{UV} it is only stated how morphisms $s_{A,B}$ are defined, and in the bicategorical case in \cite[Theorem 1]{FH} apart from 
the definition of 1-cells $s_{A,B}$, bijective interrelation between the axioms of an extension and a strength is recorded (see Table 
\ref{table:6}). Given the complexity of the construction, we 
explain for the record how the structure 2-cells for a horizontal strength $s$ in 
a double categorical case are given. 

\smallskip

We define a horizontal transformation $s$ with 1h-cell components $s_{A,B}:A\ot T(B)\to T(A\ot B)$ 
given by $s_{A,B}=1_A\rtr 1_{ T(B)}^\bullet$. For 1h-cells $f:A\to A', g:B\to B'$ in $\Dd$, set $g^K:=g$ for a cell in 
$\Kk l( T)$ determined by $\eta_{B'}g$ in $\Dd$. Then compose $\rtr^2_{f1_A,g^K1_{T(B)}^\bullet}$ in \equref{compositor} with
$(\rtr^2_{1_{A'}f,1_{T(B')}^\bullet T(g^K)})^{-1}$ - let us call this a {\em braiding of the action} - to get a 2-cell 
\begin{multline*}
\mu_{A'B'}\big( T(1_{A'}\w\rtr\w 1^\bullet_{T(B')})  T(A'\rtr\eta_{T(B')})  T(f\rtr g)\big)(1_A\w\rtr\w 1^\bullet_{T(B)})\Rightarrow \\
\mu_{A'B'} T(1_{A'}\w\rtr\w 1_{T(B')}^\bullet)\big((1_{A'}\w\rtr\w 1^\bullet_{T^2(B')})(A'\rtr\eta_{T(B')})(f\rtr  T(g))\big)  
\end{multline*}
whereby we applied \equref{rtr explain} at two places. We apply it once again in the codomain to get 
$(1_{A'}\rtr 1^\bullet_{T^2(B')})(A'\rtr\eta_{T(B')})=1_{A'}\rtr 1_{T(B')}^K$, and compose the domain from above with (the inverse of) 
\equref{unitor} in the form of a composite 2-cell 
$ T(\eta_{A'B'})\Rightarrow  T(1_{A'}\rtr 1^\bullet_{T(B')})  T(A'\rtr\eta_{T(B')})$. We 
find ourselves in a diagram chasing, where we apply then a braiding of the action in the codomain to use an isomorphism  
$ T(1_{A'}\rtr 1_{T(B')}^\bullet)(1_{A'}\rtr 1_{T(B')}^K)\iso T(1_{A'}\rtr 1_{B'}^K)(1_{A'}\rtr 1_{T(B')}^\bullet)$ 
(this is obtained by composing $\rtr^2_{1_A1_A,1_{T(B)}^\bullet 1_{T(B)}^K}$ with $(\rtr^2_{1_A1_A,1_B^K 1_{T(B)}^\bullet})^{-1}$). We then compose 
with the 2-cell $ T(\rtr^0_{A',B'})$ to ``substitute'' the border 1h-cell $ T(1_{A'}\rtr 1_{B'}^K)$ by $ T(\eta_{A'B'})$. 
Now we may compose with the double monad structure 2-cell $l^{ T}$ in the codomain and with its inverse in the domain to ``substitute'' 
the border 1h-cell $\mu_{A'B'} T(\eta_{A'B'})$ by the identity on $A'\rtr B'$, so that we end up with an invertible composite 2-cell that we take to be  
$$s_{f,g}:
 T(f\rtr g)(1_A\rtr 1^\bullet_{T(B)})\Rightarrow 
(1_{A'}\rtr 1_{T(B')}^\bullet)(f\rtr T(g)).$$

For the 2-cell component $s^{u,v}$ for 1v-cells $u,v$ we take the square 2-cell $\Id^u\rtr\Id^{T(v)\bullet}$ as in the right-hand side of \equref{s-u,v}. The axioms \axiomref{h.o.t.-3} and \axiomref{h.o.t.-4} concerning $s^{u,v}$ then clearly hold. 
For the axiom \axiomref{h.o.t.-5}: since $s_{f,g}$ is basically constructed as a braiding of the action, this axiom follows 
by the naturality of the action functor \axiomref{lx.f.c-nat} - one applies it twice in an appropriate way, being that the braiding is comprised of the composition of the lax and the colax structure of the action functor. The axioms 
\axiomref{h.o.t.-1} and \axiomref{h.o.t.-2} are proved as in the bicategorical case. In conclusion, $s$ makes a horizontal pseudonatural transformation. 

For the 2-cell 
$z_{A,B}: s_{A,B}\comp(A\rtr\eta_B)\Rightarrow\eta_{AB}$ we set the globular isomorphism $\theta_{1_A,1_B}=\rtr^0_{A,B}$, 
see also \equref{unitor}.

To construct $x_A$, substitute $K(f)$ in \equref{lambda-teta} by $1^\bullet_{T(A)}$ and observe that 
${\tilde\lambda}_A(1_I\rtr 1^\bullet_{T(A)})^K
\iso  T(\lambda_A)(1_I\rtr 1^\bullet_{T(A)})$ as 1h-cells in $\Dd$ via the isomorphism 2-cell $l^{ T}$, and similarly, 
$1^\bullet_{T(A)}\tilde\lambda_{T(A)}\iso\lambda_{T(A)}$. Then, 
given that 
a 2-cell component of the horizontal modification $x$ is of the form $ T(\lambda_{A})\comp s_{I,A}\Rightarrow\lambda_{T(A)}$, we set $x_A$ to be $\tilde\lambda_{1^\bullet_{T(A)}}$ up to the conjugation by $l^{ T}$ ({\em i.e.} its inverse), similarly as in the construction of $s_{f,g}$ above. 

To define  the horizontal modification $y$ with 2-cell components mapping \\
$y_{A,B,C}:[s_{AB,C}\vert T(\alpha_{A,B,C})] \Rightarrow[\alpha_{A,B,T(C)}\vert As_{B,C}\vert s_{A,BC}]$ we define 
$$y_{A,B,C}:=\scalebox{0.82}{
\bfig
\putmorphism(-110,170)(1,0)[(AB)T(C)` T((AB)C) `(1_A1_B)\rtr 1^K_{T(C)}]{1200}1a
\putmorphism(1290,170)(1,0)[` T(A(BC)) `T(\alpha_{A,B,C})]{1410}1a
\putmorphism(-160,180)(0,-1)[``=]{380}1l 
\putmorphism(2690,180)(0,-1)[``]{380}1r
\putmorphism(-110,-200)(1,0)[(AB)T(C) ` A(BT(C))`\alpha_{A,B,T(C)}]{1000}1a
\putmorphism(1090,-200)(1,0)[`T(A(BC))`1_A\rtr(1_B\rtr 1^K_{T(C)})]{1700}1a
\put(690,0){\fbox{$\tilde\alpha_{1_A,1_B, 1^K_{T(C)}}$}} 
\putmorphism(890,-200)(0,-1)[``]{380}1r
\putmorphism(2690,-200)(0,-1)[``]{380}1r
\putmorphism(790,-600)(1,0)[A(BT(C))`A T(BC)`A(1_B\w\rtr\w 1^K_{T(C)})]{1000}1b
\putmorphism(1950,-600)(1,0)[`T(A(BC)).`1_A\w\rtr\w 1^K_{T(BC)}]{750}1b
\put(1690,-400){\fbox{$\iso$}} 
\efig}
$$

Finally, we define the 2-cell components for a modification $w$ by  
$$w_A:=\scalebox{0.8}{
\bfig
 \putmorphism(60,200)(1,0)[A T^2(B)` T(AT(B)) `1_A\w\rtr\w 1^\bullet_{T^2(B)}]{900}1a 
 \putmorphism(980,200)(1,0)[\phantom{A''\ot B'}` T^2(A B) `T(1_A\w\rtr\w 1^\bullet_{T(B)})]{950}1a 
 \putmorphism(1900,200)(1,0)[\phantom{A''\ot B'}`T(A B)  `\mu_{A,B}]{600}1a 
\putmorphism(0,200)(0,-1)[\phantom{Y_2}``=]{400}1l
\putmorphism(2540,200)(0,-1)[\phantom{Y_2}`T(A  B)`=]{400}1r 
 \putmorphism(-50,-200)(1,0)[A T^2(B)`  `1_A\w\rtr\w (1^\bullet_{T(B)}1^\bullet_{T^2(B)})]{2470}1b 
\putmorphism(0,-200)(0,-1)[\phantom{Y_2}``=]{400}1l
\putmorphism(2540,-200)(0,-1)[\phantom{Y_2}`T(A  B).`=]{400}1r 
 \putmorphism(-50,-600)(1,0)[A T^2(B)` A T(B) `A \mu_{B}]{1300}1b 
 \putmorphism(1220,-600)(1,0)[\phantom{A''\ot B'}`  `1_A\w\rtr\w 1^\bullet_{T(B)}]{1180}1b 
\put(1160,-20){\fbox{$\rtr^2_{1_A1_A,1_{T(B)}^\bullet 1_{T^2(B)}^\bullet}$}}
\put(1340,-470){\fbox{\equref{rtr explain}}}
\efig}$$

\medskip

The above horizontal action, that is, the pseudodouble functor $\rtr:\Dd\times\Kk l( T)\to\Kk l( T)$ with 
horizontal equivalences $\tilde\lambda,\tilde\alpha$, and horizontal modifications $\tilde p, \tilde l, \tilde m$, 
determines an action $\HH(\rtr):\HH(\Dd) \times\HH(\Kk l( T))\to\HH(\Kk l( T))$ of the underlying bicategory $\HH(\Dd)$ on 
$\HH(\Kk l( T))=\Kl(\HH( T))$. 
Also, the above horizontal icon $\theta:-\rtr K(-)\Rightarrow K(-\ot-)$ induces a bicategorical icon, that abusing 
notation we will denote the same way, $\theta:-\rtr K(-)\Rightarrow K(-\ot-)$. 
Finally, the identities \equref{hor-p,l,m-teta}, \equref{hor-lambda-teta} and \equref{hor-alfa-teta} already have a bicategorical  form so that we have an extension of the canonical action of $\HH(\Dd)$ on itself. By \cite[Theorem 1]{HF}, 
\cite[Theorem 7.2]{HF2} this extension determines a 
left strength $s$ on $\HH(T)$, whereby we spelled out the data $s_{f,g}, z,x,y,w$ for $s$ explicitly (observe that they are all globular 2-cells). Since a bicategorical strength satisfies the same axioms as a horizontal double categorical strength 
(in the form of 2-cell components of horizontal modifications, see \deref{hor-strength}, b), c)) and we proved above that $s$ is a horizontal pseudonatural transformation, we have that $s$ is a left horizontal strength on $ T$. 

\smallskip

In this subsection we have thus proved that given a horizontal monad $T$ on a double category $\Dd$, a horizontal extension 
$\rtr:\Dd\times\Kk l( T)\to\Kk l( T)$ (with a horizontal icon $\theta:-\rtr K(-)\Rightarrow K(-\ot-)$) of the canonical action of $\Dd$ on itself induces a left horizontal strength on $ T$. In view of \thref{hor-str-act} we may claim: 

\begin{thm} \thlabel{iff}
Let $(\Dd, \ot, \alpha, \lambda, \rho)$ be a horizontally monoidal double category and 
let $S$ be a horizontal double monad on $\Dd$. 
There is a one-to-one correspondence betwen horizontal strengths on $S$ and extensions of the canonical action of $\Dd$ on itself (given by horizontal actions $\rtr:\Dd\times\Kk l(S)\to\Kk l(S)$ and horizontal icons $\theta:-\rtr K(-)\Rightarrow K(-\ot-)$). 
\end{thm}

The above result obviously applies to a monoidal double category with liftable $\alpha, \lambda, \rho$ and a vertical double monad $T$ 
with invertible liftable transformations $\mu, \eta$ considering $S=\hat T$. However, although in \thref{strength-action} we constructed a horizontal extension out of a vertical strength on 
such $T$, the converse 
is not possible to construct, as we can not obtain a liftable {\em strict} invertible vertical transformation $t$ such that 
$\hat t=s$,  
where $s$ is a horizontal pseudonatural transformation obtained in this subsection. Assuming that the 1h-cells 
$s_{A,B}$ are companions of some 1v-cells $t_{A,B}$, similarly as in part 1. of \prref{lifting 1v to equiv} 
one can construct an invertible vertical transformation $t$ so that $\hat t=s$ (see \cite[Proposition 4.2]{Fem:Fil}). Though, 
there is no way to obtain the strictness for $t$, necessary for $t$ to be a vertical strength.

Likewise, horizontal pseudonatural transformations $\tilde\alpha, \tilde\lambda$ in a horizontal action do not yield the 
corresponding invertible vertical {\em strict} transformations, so we can not obtain a vertical action out of a horizontal one. 
An easy but insignificant part that works is that by assuming that the horizontal icon $\theta$ in a horizontal extension has trivial 
square-formed 2-cell components $\theta^{u,v}$, by \leref{vert icon} it determines a vertical icon $-\rtr K(-)\Rightarrow K(-\ot-)$.

\subsection{Bistrong monads yield premonoidal Kleisli double categories}

As several times so far, we can define three versions of the notion of a bistrength: a vertical one, a so-to-say ``square'' or ``mixed'' one, and 
a horizontal one, where the vertical one induces the other two, assuming that both left and right strength are given by invertible liftable  vertical transformations. 

Given a monoidal double category $(\Dd,\ot,\alpha,\lambda,\rho)$ and a vertical double monad $T$ on $\Dd$. Assume that $t$ 
and $s$ are a left and a right vertical strength on $T$ and that they are invertible liftable vertical transformations.  
We say that $(T,t,s)$ is {\em vertically bistrong} if there is an identity vertical modification $q$, and similarly that $(T,\hat t,\hat s)$ is (mixed) {\em bistrong} if there is a modification in the sense of \deref{modif-hv} with components 
$$
\scalebox{0.86}{
\bfig
 \putmorphism(-470,500)(1,0)[(AT(B))C`T(AB)C`\hat t_{A,B}\ot C]{800}1a
 \putmorphism(400,500)(1,0)[\phantom{F(f)}`T((AB)C) `\hat s_{AB,C}]{700}1a
 \putmorphism(-470,120)(1,0)[A(T(B)C) ` AT(BC) `A\ot \hat s_{B,C}]{800}1a
 \putmorphism(380,120)(1,0)[\phantom{G(B)}`T(A(BC)) `\hat t_{A,BC}]{740}1a
\putmorphism(-480,500)(0,-1)[\phantom{Y_2}``\alpha]{380}1l
\putmorphism(1140,500)(0,-1)[\phantom{Y_2}``T(\alpha)]{380}1r
\put(250,310){\fbox{$q_{A,B,C}^*$}}
\efig}
$$
(we denote it by $q^*$ alluding to \prref{omega*}), which satisfies the two axioms analogous to those of \cite[Figure 6]{HF2} 
(they stand for compatibility of $q^*$ with $z^*$ and $w^*$). 

It is immediate that a (vertically or mixed) bistrong vertical monad $T$ induces a bistrong horizontal monad $\hat T$ with the defining horizontal modification $\hat q$ satisfying the same axioms as in \cite[Figure 6]{HF2}.

\bigskip

Let a bistrong vertical double monad $(T,t,s)$ on a monoidal double category $(\Dd,\ot, \\ \alpha,\lambda,\rho)$ be given, 
and assume that $\alpha,\lambda,\rho, \,\, t,s,\,\, \mu,\eta$ are invertible liftable  
vertical transformations. By \thref{strength-action} we have that 
$t$ and similarly $s$ induce a left and a right horizontal actions $\rtr:\Dd\times\Kk l(\hat T)\to\Kk l(\hat T)$ and 
$\ltr:\Kk l(\hat T)\times\Dd\to\Kk l(\hat T)$, respectively. The accompanying horizontal equivalences are induced by invertible vertical strict transformations that we denote as
$$\tilde\alpha_{A,B,E}^L:(A\ot B)\rtr E\to A\rtr(B\rtr E), \qquad \tilde\alpha_{W,X,Y}^R:(W\ltr X)\ltr Y\to W\ltr(X\ot Y)$$
$$\tilde\lambda_E: I\rtr E\to E\qquad\text{and}\qquad \tilde\rho_W: W\ltr I\to W$$
where $A,B,X,Y\in\Dd$ and $E,W\in\Kk l(\hat T)$, with their respective horizontal modifications $\tilde p^L, \tilde m^L, \tilde l$ and 
$\tilde p^R, \tilde m^R, \tilde r$. 

\medskip

We will say that these two (underlying vertical) actions are compatible if 
\begin{itemize} 
\item $\tilde\alpha_{A,B,C}^L=\tilde\alpha_{A,B,C}^R$ for $A,B,C\in\Dd$, and 
\item $\tilde p^L=\tilde p^R$ and $\tilde m^L=\tilde m^R$.
\end{itemize}

We now may prove:

\begin{thm} \thlabel{Kl-prem}
Let $(T,t,s)$ be a bistrong vertical double monad on a monoidal double category $\Dd$ so that 
the vertical transformations $\alpha, \lambda, \rho, \,\, t,s, \,\, \mu,\eta$ are invertible and liftable. Then the Kleisli pseudodouble category 
$\Kk l(\hat T)$ is premonoidal. 
\end{thm}

\begin{proof}
We define a binoidal structure on $\Kk l(\hat T)$ by setting 
$$A\ltimes g^K:=(A\ot B\stackrel{A\ot g}{\to}A\ot T(B') \stackrel{\hat t_{A,B'}}{\to}T(A\ot B'))$$
$$f^K\rtimes B:=(A\ot B\stackrel{f\ot B}{\to}T(A')\ot B \stackrel{\hat s_{A',B}}{\to}T(A'\ot B))$$
for 1h-cells $f^K:A\to A'$ and $g^K:B\to B'$ in $\Kk l(\hat T)$. That thus defined assignments $A\ltimes-, -\rtimes B: 
\Kk l(\hat T)\to\Kk l(\hat T)$ determine pseudodouble functors is proved analogously as for $\rtr$ of \thref{strength-action}. 
(We also have invertible vertical strict transformations $\tilde\alpha^L, \tilde\alpha^R$ and $\tilde\lambda, \tilde\rho$, as we commented above.) 
We now define three $\overline\alpha$'s as in \deref{assoc} which should live in $\Kk l(\hat T)$. We set: 
$$\overline\alpha_{-,B,C}:=\tilde\alpha^R_{-,B,C}: (-\rtimes B)\rtimes C\Rightarrow -\rtimes(B\bowtie C)$$
$$\overline\alpha_{A,B,-}:=\tilde\alpha^L_{A,B,-}: (A\bowtie B)\ltimes -\Rightarrow A\ltimes(B\ltimes -)$$
$$\hspace{-1cm}\overline\alpha_{A,-,C}:
(A\ltimes -)\rtimes C\Rightarrow A\ltimes(-\rtimes C)$$
where the latter we define by $\overline\alpha_{A,g,C}:=[\alpha_{A,g,C}\,\vert\, q^*_{A,B',C}]$ and 
$\overline\alpha_{A,v,C}:=\alpha_{A,v,C}$ at a 1h-cell $g$ and a 1v-cell $v$. 
We already know that $\tilde\alpha^R_{-,B,C}$ and $\tilde\alpha^L_{A,B,-}$ are invertible vertical strict transformations, and similarly so is 
$[\alpha_{A,-,C}\,\vert\, q^*_{A,-,C}]$. For 
$\overline\lambda: I\ltimes -\Rightarrow\Id$ and $\overline\rho: -\rtimes I\to \Id$ 
we set to be $\tilde\lambda: I\rtr -\Rightarrow\Id$ and $\tilde\rho: -\ltr I\Rightarrow \Id$. 

Observe that 1v-cell components of the three $\overline\alpha$'s are all equal to $\alpha_{A,B,C}$ (of $\Dd$), thus the first condition 
for the left and the right actions of $\Dd$ on $\Kk l(\hat T)$ to be compatible is fulfilled. We also have that the 1v-cell
components of $\overline\lambda$ and $\overline\rho$ are equal to the 1v-cells $\lambda_A$ and $\rho_A$. We have that they and 
$\alpha_{A,B,C}$ are invertible, we are now going to make them inversely central. We define
\begin{multline*} 
\bfig
\putmorphism(-110,170)(1,0)[``(AB)C\ltimes k^K]{570}1a
\putmorphism(-160,200)(0,-1)[``\alpha_{A,B,C}\ltimes D]{430}1l 
\putmorphism(480,200)(0,-1)[``\alpha_{A,B,C}\ltimes D']{430}1r
\putmorphism(-110,-200)(1,0)[``A(BC)\ltimes k^K]{570}1a
\put(-120,20){\fbox{$\alpha_{A,B,C}\ltimes-\vert_{k^K}$}}
\efig
\,\,:= \\
\bfig
 \putmorphism(-650,210)(1,0)[((AB)C)D` `(AB)C\ot k]{750}1a 
\putmorphism(-660,200)(0,-1)[\phantom{Y_2}` `\alpha_{A,B,C}D]{380}1l
\putmorphism(100,200)(0,-1)[\phantom{Y_2}` `]{380}1r
\putmorphism(80,200)(0,-1)[\phantom{Y_2}` `\alpha_{A,B,C}T(D')]{380}0r
\putmorphism(100,210)(1,0)[\phantom{Y}` T(((AB)C)D')  `\hat t_{(AB)C,D'}]{1150}1a
\putmorphism(1230,200)(0,-1)[\phantom{Y_2}`T((A(BC))D') `]{380}0r
\putmorphism(1130,200)(0,-1)[\phantom{Y_2}` `T(\alpha_{A,B,C}D')]{380}1r
\putmorphism(-670,-200)(1,0)[(A(BC))D ``A(BC)\ot k]{730}1a 
\putmorphism(150,-200)(1,0)[``\hat t_{A(BC),D'}]{780}1a
\put(-600,0){\fbox{$\alpha_{1_A,1_B,1_C}\ot\Id_k$}}
\put(650,0){\fbox{$\hat t^{\alpha_{A,B,C},1_{D'}}$}}
\efig
\end{multline*}
and 
\begin{multline*} 
\bfig
\putmorphism(-110,170)(1,0)[``f^K\rtimes(BC)D]{570}1a
\putmorphism(-160,200)(0,-1)[``A\rtimes\alpha_{B,C,D}]{430}1l 
\putmorphism(490,200)(0,-1)[``A'\rtimes\alpha_{B,C,D}]{430}1r
\putmorphism(-110,-200)(1,0)[``f^K\rtimes B(CD)]{570}1a
\put(-120,20){\fbox{$-\rtimes\alpha_{B,C,D}\vert_{f^K}$}}
\efig
\,\,:=\\
\bfig
 \putmorphism(-650,210)(1,0)[((AB)C)D` `f\ot(BC)D]{750}1a 
\putmorphism(-660,200)(0,-1)[\phantom{Y_2}` `A\alpha_{B,C,D}]{380}1l
\putmorphism(100,200)(0,-1)[\phantom{Y_2}` `]{380}1r 
\putmorphism(80,200)(0,-1)[\phantom{Y_2}` `T(A')\alpha_{B,C,D}]{380}0r 
 \putmorphism(100,210)(1,0)[\phantom{Y}` T(((AB)C)D')  `\hat s_{A',(BC)D}]{1150}1a
\putmorphism(1230,200)(0,-1)[\phantom{Y_2}`T((A(BC))D') `]{380}0r
\putmorphism(1130,200)(0,-1)[\phantom{Y_2}` `T(A'\alpha_{B,C,D})]{380}1r
\putmorphism(-650,-200)(1,0)[(A(BC))D``f\ot B(CD)]{730}1a
\putmorphism(150,-200)(1,0)[``\hat s_{A',B(CD)}]{780}1a
\put(-600,10){\fbox{$\Id_f\ot\alpha_{1_B,1_C,1_D}$}}
\put(650,0){\fbox{$\hat s^{1_{A'},\alpha_{B,C,D}}$}}
\efig
\end{multline*}
(the right-most 2-cells of the form $\hat t^u$ and $\hat s^u$ are obtained from \prref{lifting 1v to equiv}, part 1. c)). 
It is directly proved that these make $\alpha_{A,B,C}\ltimes-$ and $-\rtimes\alpha_{B,C,D}$ invertible vertical strict transformations, 
so that the 1-cell $\alpha_{A,B,C}$ is central (one uses \axiomref{h.o.t.-5} for $\hat t$ and $\hat s$ and modification axioms for $\hat w$ and 
$\hat z$). 
They are obviously inversely central. 
The construction and proof of inverse centrality for $\tilde\lambda_A$ and $\tilde\rho_A$ are similar.


It remains to check the four pentagons and six triangles for a premonoidal structure. We first observe that 2-cell components of the three 
$\overline\alpha$'s comprising the pentagons are given by:
$$\overline\alpha_{A,B,h}=\tilde\alpha_{A,B,h}^L= 
\scalebox{0.82}{
\bfig
 \putmorphism(-650,210)(1,0)[(AB)C` (AB)T(C')`(AB)h]{660}1a 
\putmorphism(-660,200)(0,-1)[\phantom{Y_2}` `\alpha_{A,B,C}]{380}1l
\putmorphism(-50,200)(0,-1)[\phantom{Y_2}` `]{380}1r
\putmorphism(-50,200)(0,-1)[\phantom{Y_2}` `\alpha_{A,B,T(C')}]{380}0r
\putmorphism(200,210)(1,0)[\phantom{Y}` T((AB)C')  `\hat t_{AB,C'}]{1420}1a
\putmorphism(1600,200)(0,-1)[\phantom{Y_2}`T(A(BC')) `]{380}0r
\putmorphism(1530,200)(0,-1)[\phantom{Y_2}` `T(\alpha_{A,B,C'})]{380}1r
\putmorphism(-670,-200)(1,0)[A(BC) `A(BT(C'))`A(Bh)]{640}1a 
\putmorphism(200,-200)(1,0)[``A\hat t_{B,C'}]{380}1a
\putmorphism(780,-200)(1,0)[AT(BC')``\hat t_{A,BC'}]{580}1a
\put(-500,30){\fbox{$\alpha_{A,B,h}$}}
\put(580,30){\fbox{$ y^*_{A,B,C'}$}}
\efig}
$$

\vspace{0,2cm}

$$\overline\alpha_{A,g,C}=
\scalebox{0.82}{
\bfig
 \putmorphism(-650,210)(1,0)[(Ag)C` (AT(B'))C`(Ag)C]{660}1a 
\putmorphism(-660,200)(0,-1)[\phantom{Y_2}` `\alpha_{A,B,C}]{380}1l
\putmorphism(-50,200)(0,-1)[\phantom{Y_2}` `]{380}1r
\putmorphism(-50,200)(0,-1)[\phantom{Y_2}` `\alpha_{A,T(B'),C}]{380}0r
\putmorphism(-670,-200)(1,0)[A(BC) `A(T(B')C)`A(gC)]{640}1a 

\putmorphism(200,210)(1,0)[\phantom{Y}` T(AB')C`\hat t_{A,B'}\ot C]{630}1a
\putmorphism(1030,210)(1,0)[\phantom{Y}` T((AB')C) `\hat s_{AB',C}]{580}1a
\putmorphism(200,-200)(1,0)[``A\ot \hat s_{B',C}]{380}1a
\putmorphism(780,-200)(1,0)[AT(B'C)``\hat t_{A,B'C}]{580}1a
\putmorphism(1600,200)(0,-1)[\phantom{Y_2}`T(A(B'C)) `]{380}0r
\putmorphism(1530,200)(0,-1)[\phantom{Y_2}` `T(\alpha_{A,B',C})]{380}1r

\put(-500,30){\fbox{$\alpha_{A,g,C}$}}
\put(580,30){\fbox{$q^*_{A,B',C}$}}
\efig}
$$

\vspace{0,2cm}

$$\overline\alpha_{f,B,C}=\tilde\alpha_{f,B,C}^R=
\scalebox{0.82}{\bfig
 \putmorphism(-650,210)(1,0)[(AB)C` (T(A')B)C`(fB)C]{660}1a 
\putmorphism(-660,200)(0,-1)[\phantom{Y_2}` `\alpha_{A,B,C}]{380}1l
\putmorphism(-50,200)(0,-1)[\phantom{Y_2}` `]{380}1r
\putmorphism(-50,200)(0,-1)[\phantom{Y_2}` `\alpha_{T(A'),B,C}]{380}0r
\putmorphism(200,210)(1,0)[\phantom{Y}` T(A'B)C  `\hat s_{A',B}C]{630}1a
\putmorphism(1030,210)(1,0)[\phantom{Y}` T((A'B)C)  `\hat s_{A'B,C}]{580}1a
\putmorphism(1600,200)(0,-1)[\phantom{Y_2}`T(A'(BC)) `]{380}0r
\putmorphism(1530,200)(0,-1)[\phantom{Y_2}` `T(\alpha_{A',B,C})]{380}1r
\putmorphism(-670,-200)(1,0)[A'(BC) `T(A')(BC)`f(BC)]{640}1a 
\putmorphism(200,-200)(1,0)[``\hat s_{A',BC}]{1180}1a
\put(-500,30){\fbox{$\alpha_{f,B,C}$}}
\put(580,30){\fbox{$ (y')^*_{A',B,C}$}}
\efig}
$$
for 1h-cells $f,g,h$ as usual. 
Before going on, at this point we note the following. Recall from the end of the proof of \thref{strength-action} that both $\tilde p^L$ 
and $\tilde p^R$ we took to be the identity vertical modification $p$ from $\Dd$, and similarly, for $\tilde m^L$ and $\tilde m^R$ 
we took to be $m$ from $\Dd$. Thus also the second condition for the left and the right actions of $\Dd$ on $\Kk l(\hat T)$ to be compatible
is fulfilled. 

Now, each of the four pentagons becomes an equation of vertical compositions of horizontal compositions of 2-cell components, where  these horizontal compositions have some of the following four forms:     
$$[\alpha_{f,B,C}\ot D\,\,\vert\,\, (y')^*_{A',B,C}\ot D\,\,\vert\,\,\hat s^{\alpha_{A',B,C},1_D}]$$
$$[\alpha_{f,BC,D}\,\,\vert\,\, (y')^*_{A',BC,D}]$$
$$[\Id_f\ot\alpha_{1_B,1_C,1_D}\,\,\vert\,\, \hat s^{1_{A'},\alpha_{B,C,D}}]$$
$$[\alpha_{fB,C,D}\,\,\vert\,\, \alpha_{\hat s_{A',B},C,D}\,\,\vert\,\, (y')^*_{A'B,C,D}]$$
with the following variations: \vspace{-0,2cm}
\begin{itemize}
\item instead of $-\ot D$ there appears $A\ot-$, \vspace{-0,2cm}
\item instead of $(y')^*$ there appears $y^*$ or $q^*$, \vspace{-0,2cm}
\item instead of $\hat s$ there appears $\hat t$, \vspace{-0,2cm}
\item instead of $\alpha_{f,BC,D}$ there can be any of: $\alpha_{f,B,CD}, \alpha_{A,g,CD}, 
\alpha_{AB,h,D}, \alpha_{A,BC,k}, \alpha_{AB,C,k}$, and \vspace{-0,2cm}
\item instead of $\alpha_{fB,C,D}$ there can be any of: 
$\alpha_{A,gC,D}, \alpha_{Ag,C,D}, \alpha_{A,Bh,D}, \alpha_{A,B,hD}, \alpha_{A,B,Ck}$ \vspace{-0,2cm}
\end{itemize}
(with 1h-cells $f:A\to A', g:B\to B', h:C\to C', k:D\to D'$ in $\Dd$). 
In the above described vertical compositions at some places there also appear 2-cells (compositors) for the lax and the colax pseudodouble 
functor structure of $\ot$. 

Although it would be very tedious to check that such pentagon-equations hold, the argument that we have at hand is simple. At the first place, in each of the four pentagon-equations, to the most left, on both sides of the equations, in all the rows of the pasted diagrams there is a 2-cell component of $\alpha$. When we first compose vertically these 2-cell components of $\alpha$ (and then we compose this vertical composition horizontally with the rest of the pasting diagrams, {\em i.e.} 2-cells, in both sides of the four pentagon-equations, which is allowed by the strict interchange law holding in double categories), we find the single pentagon-equation for $\alpha$ of $\Dd$ that we know that holds. Then we cancel out the first columns in both sides of every of the four pentagon-equations. 

In the remaining part, there appear (vertical and horizontal) compositions of 2-cells of the form 
$(y')^*_{A',B,C}\ot D, \,\, \hat s^{\alpha_{A',B,C},1_D}, \,\,  (y')^*_{A',BC,D}, \,\, \alpha_{\hat s_{A',B},C,D}$ 
(and their announced variations). 
We make several observations. \vspace{-0,2cm}
\begin{enumerate}
\item The 2-cell components of the modifications $y^*, (y')^*, q^*$ originate from 2-cell components of the identity vertical modifications $y,y',q$ and are in 1-1 correspondence with the 2-cell components of the corresponding horizontal modifications $\hat y, \hat y', \hat q$ (see \prref{omega*}). \vspace{-0,2cm} 
\item The 2-cells of the form $\hat s^u, \hat t^u$ are in 1-1 correspondence with the vertically globular 2-cells 
$s^u, t^u$ (this is similar to \prref{omega*}), and these in turn are in 1-1 correspondence with the horizontally globular 2-cells $\hat s_{\hat u}, \hat t_{\hat u}$ (as in \prref{essence}).\vspace{-0,2cm}  
\item The 2-cells of the form $\alpha_{\hat u}$ are in 1-1 correspondence with the horizontally globular 2-cells 
$\hat\alpha_{\hat u}$. \vspace{-0,2cm}
\end{enumerate}
Now, by \leref{horiz is teta} the 2-cells $\hat s_{\hat u}, \hat t_{\hat u}$ and $\hat\alpha_{\hat u}$ (as in 2. and 3.) are of $\theta$ type, 
and by \leref{Lemma 4.8} so are the 2-cell components of $\hat y, \hat y', \hat q$ (from 1.). 

At last, the 2-cells $(y')^*\ot D, y^*\ot D, A\ot q^*$ ({\em i.e.} the corresponding 2-cell components of them) actually appear vertically composed 
to the compositors of the lax, colax and both lax and colax pseudodouble functor structures of $-\ot D$ {\em i.e.} $A\ot-$, respectively. 
As a matter of fact, all of the appearances of the compositors that we mentioned earlier above can be attached to some of these 2-cells (or their variations), and for these vertical composite 2-cells we have:

\begin{lma}
Let $\omega_1^*, \omega_2^*, \omega_3^*$ be three 2-cells as below induced as in \prref{omega*} with liftable 1v-cells $l,r$ and 
so that all the indicated 1-cells are non-trivial.  
$$
\scalebox{0.86}{
\bfig
 \putmorphism(-170,500)(1,0)[``\hat u]{500}1a
 \putmorphism(250,500)(1,0)[\phantom{F(f)}``\hat v]{550}1a
 \putmorphism(-170,120)(1,0)[``\hat w]{1000}1a
\putmorphism(-180,500)(0,-1)[\phantom{Y_2}``l]{380}1l
\putmorphism(820,500)(0,-1)[\phantom{Y_2}``r]{380}1r
\put(260,310){\fbox{$\omega_1^*$}}
\efig}
\qquad\scalebox{0.86}{
\bfig
 \putmorphism(-170,500)(1,0)[``\hat w]{980}1a
 \putmorphism(-170,120)(1,0)[``\hat u']{500}1a
 \putmorphism(250,120)(1,0)[\phantom{G(B)}` `\hat v']{550}1a
\putmorphism(-180,500)(0,-1)[\phantom{Y_2}``l]{380}1l
\putmorphism(820,500)(0,-1)[\phantom{Y_2}``r]{380}1r
\put(260,310){\fbox{$\omega_2^*$}}
\efig}
\qquad
\scalebox{0.86}{
\bfig
 \putmorphism(-170,500)(1,0)[``\hat u]{500}1a
 \putmorphism(250,500)(1,0)[\phantom{F(f)}``\hat v]{550}1a
 \putmorphism(-170,120)(1,0)[``\hat u']{500}1a
 \putmorphism(250,120)(1,0)[\phantom{G(B)}`. `\hat v']{550}1a
\putmorphism(-180,500)(0,-1)[\phantom{Y_2}``l]{380}1l
\putmorphism(820,500)(0,-1)[\phantom{Y_2}``r]{380}1r
\put(260,310){\fbox{$\omega_3^*$}}
\efig}
$$
Let $F:\Dd\to\Dd$ be a pseudodouble functor and consider the composite 2-cells
$$\xi_1^*=\frac{F_{\hat u,\hat v}}{F(\omega_1^*)}, \qquad \xi_2^*=\frac{F(\omega_2^*)}{F_{\hat u',\hat v'}^{-1}}, \qquad 
\xi_3^*=\threefrac{F_{\hat u,\hat v}}{F(\omega_3^*)}{F_{\hat u',\hat v'}^{-1}}$$
where $F_{\hat u,\hat v}$ and $F_{\hat u,\hat v}^{-1}$ present the laxity and the colaxity compositor, respectively. 
Then the horizontally globular 2-cells $\hat\xi_1, \hat\xi_2, \hat\xi_3$ (as in \prref{essence}) are $\theta$ 2-cell isomorphisms. 
\end{lma}

\begin{proof}
First of all, observe that the 1-cells on the edges of $\xi_1^*,\xi_2^*,\xi_3^*$ are compositions of the same type 
as in $\omega_1^*,\omega_2^*,\omega_3^*$. The domain and codomain 1h-cells of $\hat\xi_1$ are as indicated: 
$[F(\hat u)\,\,\vert\,\, F(\hat v)\,\,\vert\,\, \widehat{F(r)}]\Rightarrow [\widehat{F(l)}\,\,\vert\,\,F(\hat w)]$.  
By \prref{alg}, 6a) and 2) the domain and codomain of $\hat\xi_1$ are companions of 
$\threefrac{F(u)}{F(v)}{F(r)}=F\left(\threefrac{u}{v}{r}\right)$ and $F(\frac{l}{w})$, respectively, and these are equal because of the identity vertically globular 2-cell $\omega_1$. Thus there is a $\theta$ 2-cell between the domain and codomain of $\hat\xi_1$. Now as in the proof of 
\leref{horiz is teta}, substituting $\theta$ in \equref{teta-property} by $\hat\xi_1$, applying \axiomref{lx.f.c-nat} and 
the assignment $\hat\omega\mapsto\omega$ from \prref{essence}, one sees that 
\equref{teta-property} holds, hence $\hat\xi_1$ is a $\theta$ isomorphism. The proof for $\hat\xi_2, \hat\xi_3$ is similar. 
\qed\end{proof}

Coming back to the proof of \thref{Kl-prem}: the composites of (the 2-cell components of) $(y')^*\ot D, y^*\ot D, A\ot q^*$ with the 
corresponding compositors are of the form of $\xi^*_1, \xi^*_2, \xi^*_3$ from the lemma, and so they are in 1-1 correspondence with $\theta$ isomorphisms (and the same holds for their variations). 

Now the same reasoning as in the proof of part 3. of \prref{omega*} applies here 
to finalize the proof. Those four types of 2-cells 
$(y')^*_{A',B,C}\ot D, \,\, \hat s^{\alpha_{A',B,C},1_D}, \,\,  (y')^*_{A',BC,D}, \,\, \alpha_{\hat s_{A',B},C,D}$ 
we express - using the respective 1-1 correspondences - in terms of their corresponding $\theta$ isomorphisms. Any of the four pentagon 
equations $E_i^*, i=1,2,3,4$ gets the form $\threefrac{[\Id\vert\iota]}{\hat E_i}{[\Epsilon\vert\Id]}$, {\em i.e.} $[\Id\vert\iota]$ and $[\Epsilon\vert\Id]$ are composed to both sides of $\hat E_i$, where $\hat E_i$ is an equation as in \prref{essence}, part 3.. 
Since by \prref{essence} the equations $\hat E_i,i=1,2,3,4$ hold true, the four pentagon-equations are also fulfilled. 

The proof for the six triangles is similar. 
\qed\end{proof}

The above proof offers a great evidence of the advantage of working in a double categorical setting in the cases that the assumption 
on the existence of the appropriate companions is fulfilled. We were spared of the tedious checking that the four pentagons and six triangles commute. Instead, we had to take an account on the type of the 2-cells involved in the pasted diagrams constituting the equation/axiom in question, and make sure that they all come as companion-lifts of vertically globular 2-cells that are identities. The rest flows by \leref{teta}, \prref{essence} and \prref{omega*}. Having carried out this easier and shorter task, we are then also able to draw the desired consequences for (underlying) bicategories (recall \thref{strenghts thm} and \prref{premon isofib new}):

\begin{cor} \colabel{Kl-prem}
Under the assumptions of \thref{Kl-prem}, $\HH(\hat T)$ 
is a bistrong pseudomonad on $\HH(\Dd)$, and the Kleisli bicategory $\HH(\Kk l(\hat T))=\Kl(\HH(\hat T))$ is premonoidal.
\end{cor}

This is a double categorical version of \cite[Theorem 2 on page 21]{HF}, \cite[Proposition 5.10]{HF2} and simplifies its proof. The relation between the two 
is as follows. If we are given a bistrong pseudomonad $S$ on a monoidal bicategory $\B$ that comes from a 
vertically bistrong vertical double monad $T$ on a monoidal double category $\Bb$, then by \coref{Kl-prem} the Kleisli bicategory 
$\Kl(S)$ is premonoidal. We explain in more details what it means for $S$ and $\B$ to come from $T$ and $\Bb$. Namely, that 
given a bistrong pseudomonad $(S, \mu_1, \eta_1)$ with a left and a right strength $s_l, s_r$ on a monoidal bicategory 
$(\B, \ot_1, \alpha_1, \lambda_1, \rho_1)$, we can find a monoidal double category $(\Bb, \ot, \alpha, \lambda, \rho)$ and a 
vertically bistrong vertical double monad $(T, \mu, \eta)$ on $\Bb$, with a left vertical strength $t_l$ 
and a right vertical strength $t_r$, so that $\alpha, \lambda, \rho$ and $\mu, \eta, t_l, t_r$ are invertible liftable vertical transformations, 
and so that $\HH\left((\Bb, \ot, \hat\alpha, \hat\lambda, \hat\rho)\right)=(\B, \ot_1, \alpha_1, \lambda_1, \rho_1)$ and 
$\HH\left((\hat T, \hat\mu, \hat\eta, \hat t_l, \hat t_r)\right)=(S, \mu_1, \eta_1, s_l, s_r)$. 

A fully horizontal version of \thref{Kl-prem} is possible. Namely, that for a {\em bistrong horizontal} double monad $S$ on a 
horizontally monoidal double category $\Dd$ the Kleisli pseudodouble category $\Kk l(S)$ is premonoidal in a sense not studied 
in this paper: where the premonoidality constraints are {\em horizontal} pseudonatural tarnsformations. Though, this perspective is out of the interest here. 

\medskip

In conclusion, by working in a setting of double categories, by means of lifting the vertical structures, which are simpler and possess less data, to the horizontal ones, 
one may obtain results concerning horizontal or mixed structures (expressed by non-globular 2-cells),  
analogous to those known in bicategories, with the benefit of passing through significantly less laborous proofs. This is possible 
under the assumption that the suitable vertical transformations are invertible and liftable, meaning that their vertical 1-cell components are have companions.

\bigskip

{\bf Acknowledgements.} The author was supported by the Science Fund of the Republic of Serbia, Grant No.~7749891, Graphical Languages - GWORDS.
My thanks to Hugo Paquet and Philip Saville for nice conversations and comments regarding their construction.

\vspace{1cm}

{\bf {\Large Appendix A.0}} 

\vspace{0,6cm}

The horizontal composition of modifications is induced on components by the horizontal composition of the corresponding 2-cells: 
$$[\Theta\vert\Theta'](A)=
\scalebox{0.86}{
\bfig
\putmorphism(-100,250)(1,0)[F(A)`G(A)`\alpha(A)]{550}1a
 \putmorphism(430,250)(1,0)[\phantom{F(A)}`H(A) `\alpha'(A)]{580}1a

 \putmorphism(-100,-200)(1,0)[F(A)`G(A)`\beta(A)]{550}1a
 \putmorphism(450,-200)(1,0)[\phantom{F(A)}`H(A). `\beta'(A)]{580}1a

\putmorphism(-100,250)(0,-1)[\phantom{Y_2}``\alpha_0(A)]{450}1l
\putmorphism(420,250)(0,-1)[\phantom{Y_2}``]{450}1r
 \putmorphism(400,250)(0,-1)[\phantom{Y_2}``\beta_0(A)]{450}0r
\putmorphism(1020,250)(0,-1)[\phantom{Y_2}``\beta' _0(A)]{450}1r
\put(40,20){\fbox{$\Theta_A$}}
\put(700,20){\fbox{$\Theta_A' $}}
\efig}
$$

The vertical composition of modifications is induced on components by the vertical composition of the corresponding 2-cells:
$$(\frac{\Theta}{\Sigma })(A)=
\scalebox{0.86}{
\bfig
 \putmorphism(-150,470)(1,0)[F(A)`G(A)  `\alpha(A)]{440}1a
\putmorphism(-160,470)(0,-1)[\phantom{Y_2}`F(A) `\alpha_0(A)]{400}1l
\putmorphism(-160,80)(0,-1)[\phantom{Y_2}`F(A)`\alpha_0'(A)]{400}1l
\putmorphism(300,80)(0,-1)[\phantom{Y_2}`G(A).`\beta_0' (A)]{400}1r
\putmorphism(300,470)(0,-1)[\phantom{Y_2}`G(A) `\beta_0 (A)]{400}1r
\put(-20,280){\fbox{$\Theta_A$}}

\putmorphism(-180,70)(1,0)[\phantom{F(A)}``\beta(A)]{380}1a
\putmorphism(-100,-300)(1,0)[\phantom{Y_2}` `\gamma(A)]{300}1b
\put(-20,-130){\fbox{$\Sigma_A$}}
\efig}
$$ 
It is clear that the associativity and unitality of modifictaions in both horizontal and vertical direction hold strictly.

\pagebreak

{\bf {\Large Appendix A.}} 

\vspace{0,6cm}

{\bf The 24 axioms determining interrelations of 12 horizontal, 12 vertical pseudodouble transformations and 6 modifications from 
\seref{24 axioms}} 
(to simplify the annotation in the diagrams, we will write $(-,-)$ both for $-\ltimes-$ and $-\rtimes-$, which one is meant will be clear from the context): 

\medskip

\noindent \axiom{$(f\ltimes,g,C)$} \vspace{-0,3cm}
$$\scalebox{0.86}{
\bfig
\putmorphism(-150,500)(1,0)[``((A,g),C)]{600}1a
 \putmorphism(480,500)(1,0)[` `((f,B'),C)]{640}1a
 \putmorphism(-150,50)(1,0)[``((f,B),C)]{600}1a
 \putmorphism(470,50)(1,0)[` `((A',g),C)]{660}1a

\putmorphism(-180,500)(0,-1)[\phantom{Y_2}``=]{450}1l
\putmorphism(1100,500)(0,-1)[\phantom{Y_2}``=]{450}1r
\put(230,280){\fbox{$(f\ltimes-\vert_g,C)$}}

\putmorphism(-170,-400)(1,0)[` `(f,(B,C))]{640}1b
 \putmorphism(470,-400)(1,0)[` `(A',(g,C))]{640}1b

\putmorphism(-180,50)(0,-1)[\phantom{Y_2}``\alpha_{A,B,C}]{450}1l %
\putmorphism(450,50)(0,-1)[\phantom{Y_2}``]{450}1l
\putmorphism(660,50)(0,-1)[\phantom{Y_2}``\alpha_{A',B,C}]{450}0l 
\putmorphism(1100,50)(0,-1)[\phantom{Y_3}``\alpha_{A',B',C}]{450}1r
\put(-40,-180){\fbox{$\alpha_{f,B,C}$}} 
\put(670,-180){\fbox{$\alpha_{A',g,C}$}}
\efig}
\quad
=
\quad
\scalebox{0.86}{
\bfig
\putmorphism(-150,500)(1,0)[``((A,g),C)]{600}1a
 \putmorphism(480,500)(1,0)[` `((f,B'),C)]{640}1a

 \putmorphism(-150,50)(1,0)[``(A,(g,C))]{600}1a
 \putmorphism(450,50)(1,0)[``(f,(B',C))]{640}1a

\putmorphism(-180,500)(0,-1)[\phantom{Y_2}``\alpha_{A,B,C}]{450}1l
\putmorphism(450,500)(0,-1)[\phantom{Y_2}``]{450}1r
\putmorphism(250,500)(0,-1)[\phantom{Y_2}``\alpha_{A,B',C}]{450}0r
\putmorphism(1100,500)(0,-1)[\phantom{Y_2}``\alpha_{A',B',C}]{450}1r
\put(-90,280){\fbox{$\alpha_{A,g,C}$}}
\put(680,280){\fbox{$\alpha_{f,B',C}$}}

\putmorphism(-150,-400)(1,0)[` `(f,(B,C))]{640}1b
 \putmorphism(490,-400)(1,0)[` `(A',(g,C))]{640}1b

\putmorphism(-180,50)(0,-1)[\phantom{Y_2}``=]{450}1l
\putmorphism(1100,50)(0,-1)[\phantom{Y_3}``=]{450}1r
\put(270,-200){\fbox{$f\ltimes-\vert_{(g,C)}$}}
\efig}
$$
for every left central 1h-cell $f\colon A\to A'$ and any 1h-cell $g\colon B\to B'$ (observe that in the rectangular diagram for 
the 2-cell $(f\ltimes-\vert_g,C)$ we omitted the compositor 2-cells of the pseudodouble functor $-\rtimes C$ on top and bottom);

\medskip
\noindent \axiom{$(f\ltimes,v,C)$} \vspace{-0,6cm}

$$\scalebox{0.86}{
\bfig
 \putmorphism(-120,500)(1,0)[` `=]{550}1a
 \putmorphism(450,500)(1,0)[` `((f,B),C)]{550}1a
\putmorphism(-140,520)(0,-1)[` `\alpha_{A,B,C}]{480}1l
\put(0,50){\fbox{$\alpha_{A,v,C}$}}
\putmorphism(-150,-380)(1,0)[` `=]{540}1a
\putmorphism(-140,80)(0,-1)[``(A,(v ,C))]{450}1l
\putmorphism(430,50)(0,-1)[` `\alpha_{A,\tilde B,C}]{450}1l
\putmorphism(430,520)(0,-1)[` `((A,v),C)]{480}1l
\put(470,300){\fbox{$(f\w\ltimes\w-\vert_v,C)$}}
\putmorphism(430,50)(1,0)[``((f,\tilde B),C)]{540}1a
\putmorphism(1000,80)(0,-1)[``\alpha_{A',\tilde B,C}]{450}1r
\putmorphism(1000,520)(0,-1)[``((A',v),C)]{480}1r
\putmorphism(450,-380)(1,0)[``(f,(\tilde B,C))]{540}1b
\put(580,-170){\fbox{$\alpha_{f,\tilde B,C}$}}
\efig}
\quad=\quad
\scalebox{0.86}{
\bfig
 \putmorphism(-150,500)(1,0)[` `((f,B),C)]{600}1a
 \putmorphism(450,500)(1,0)[` `=]{540}1a
\putmorphism(-180,520)(0,-1)[` `\alpha_{A,B,C}]{450}1l
\put(0,280){\fbox{$\alpha_{f,B,C}$}}
\putmorphism(-150,-380)(1,0)[` `(f,(\tilde B,C))]{500}1b
\putmorphism(-180,80)(0,-1)[``(A,(v,C))]{450}1l
\putmorphism(450,80)(0,-1)[``(A',(v,C))]{450}1r
\putmorphism(450,520)(0,-1)[` `\alpha_{A',B,C}]{450}1r
\put(580,50){\fbox{$\alpha_{A',v,C}$}}
\putmorphism(-150,50)(1,0)[``(f,(B,C))]{500}1a
\putmorphism(1000,80)(0,-1)[``\alpha_{A',\tilde B,C}]{450}1r
\putmorphism(1000,520)(0,-1)[``((A',v),C)]{450}1r
\putmorphism(450,-380)(1,0)[``=]{520}1b
\put(-100,-170){\fbox{$f\w\ltimes\w-\vert_{(B,C)}$}}
\efig}
$$
for every left central 1h-cell $f\colon A\to A'$ and any 1v-cell $v\colon B\to \tilde B$;

\noindent \axiom{$(u\ltimes,g,C)$} \vspace{-0,3cm}
$$\scalebox{0.86}{
\bfig
 \putmorphism(-120,500)(1,0)[` `=]{550}1a
 \putmorphism(450,500)(1,0)[` `((A,g),C)]{550}1a
\putmorphism(-140,520)(0,-1)[` `\alpha_{A,B,C}]{480}1l
\put(0,50){\fbox{$\alpha_{u,B,C}$}}
\putmorphism(-150,-380)(1,0)[` `=]{540}1a
\putmorphism(-140,80)(0,-1)[``(u,(B,C))]{450}1l
\putmorphism(430,50)(0,-1)[` `\alpha_{\tilde A,B,C}]{450}1l
\putmorphism(430,520)(0,-1)[` `((u,B),C)]{480}1l
\put(470,300){\fbox{$(u\w\ltimes\w-\vert_g,C)$}}
\putmorphism(430,50)(1,0)[``((\tilde A,g),C)]{540}1a
\putmorphism(1000,80)(0,-1)[``\alpha_{\tilde A,B',C}]{450}1r
\putmorphism(1000,520)(0,-1)[``((u,B'),C)]{480}1r
\putmorphism(450,-380)(1,0)[``(\tilde A,(g,C))]{540}1b
\put(580,-170){\fbox{$\alpha_{\tilde A,g,C}$}}
\efig}
\quad=\quad
\scalebox{0.86}{
\bfig
 \putmorphism(-150,500)(1,0)[` `((A,g),C)]{600}1a
 \putmorphism(450,500)(1,0)[` `=]{540}1a
\putmorphism(-180,520)(0,-1)[` `\alpha_{A,B,C}]{450}1l
\put(0,280){\fbox{$\alpha_{A,g,C}$}}
\putmorphism(-150,-380)(1,0)[` `(\tilde A,(g,C))]{500}1b
\putmorphism(-180,50)(0,-1)[``(u,(B,C))]{450}1l
\putmorphism(450,50)(0,-1)[``((u,B'),C)]{450}1r
\putmorphism(450,520)(0,-1)[` `\alpha_{A,B',C}]{450}1r
\put(600,50){\fbox{$\alpha_{u,B',C}$}}
\putmorphism(-150,50)(1,0)[``(A,(g,C))]{500}1a
\putmorphism(1000,50)(0,-1)[``\alpha_{\tilde A,B',C}]{450}1r
\putmorphism(1000,520)(0,-1)[``(u,(B',C))]{450}1r
\putmorphism(450,-380)(1,0)[``=]{520}1b
\put(-100,-170){\fbox{$u\w\ltimes\w-\vert_{(g,C)}$}}
\efig}
$$
for every left central 1v-cell $u\colon A\to \tilde A$ and any 1h-cell $g\colon B\to B'$; 

\pagebreak

\noindent \axiom{$(u\ltimes,v,C)$} \vspace{-0,8cm}

$$
\scalebox{0.8}{
\bfig
 \putmorphism(-170,500)(1,0)[` `=]{480}1a
\putmorphism(-180,500)(0,-1)[` `((u,B),C)]{450}1l
\put(-80,250){\fbox{$\alpha_{u,B,C}$}}
\putmorphism(-170,-380)(1,0)[` `=]{500}1a
\putmorphism(-180,50)(0,-1)[``\alpha_{\tilde A,B,C}]{450}1l
\putmorphism(300,50)(0,-1)[``(u,(B,C))]{450}1l
\putmorphism(300,500)(0,-1)[` `\alpha_{A,B,C}]{450}1r
\putmorphism(-700,50)(1,0)[``=]{480}1a
\putmorphism(-700,50)(0,-1)[``((\tilde A,v),C)]{450}1l
\putmorphism(-700,-380)(0,-1)[``\alpha_{\tilde A,\tilde B,C}]{450}1l
\putmorphism(-700,-800)(1,0)[``=]{500}1a 
\putmorphism(-180,-380)(0,-1)[` `(\tilde A,(v,C))]{450}1r
\put(-610,-630){\fbox{$\alpha_{\tilde A, v,C}$}}

\putmorphism(360,50)(1,0)[``=]{480}1a
\putmorphism(870,50)(0,-1)[``(A,(v,C))]{450}1r
\putmorphism(300,-380)(0,-1)[` `]{450}1r 
\putmorphism(870,-380)(0,-1)[` `(u,(\tilde B,C))]{450}1r
\putmorphism(350,-800)(1,0)[` `=]{530}1a 
\put(360,-200){\fbox{$u\ltimes-\vert_{(v,C)}$}}
\efig}
=
\scalebox{0.8}{
\bfig
 \putmorphism(20,500)(1,0)[` `=]{520}1a
\putmorphism(0,500)(0,-1)[` `((u,B),C)]{450}1l
\put(30,0){\fbox{$(u\w\ltimes\w-\vert_v,C)$}}
\putmorphism(0,-400)(1,0)[` `=]{520}1a %
\putmorphism(0,50)(0,-1)[``((\tilde A,u),C)]{450}1l
\putmorphism(560,50)(0,-1)[``((u,\tilde B),C)]{450}1r
\putmorphism(560,500)(0,-1)[` `]{450}1l 
\putmorphism(550,50)(1,0)[``=]{480}1a %
\putmorphism(1060,50)(0,-1)[``\alpha_{A,\tilde B,C}]{450}1r
\putmorphism(1060,-400)(0,-1)[``(u,(\tilde B,C))]{450}1r
\putmorphism(530,-800)(1,0)[``=]{450}1a 
\putmorphism(560,-400)(0,-1)[` `\alpha_{\tilde A,\tilde B,C}]{450}1l
\put(640,-460){\fbox{$\alpha_{u,\tilde B,C}$}}
\putmorphism(1060,500)(1,0)[``=]{480}1a
\putmorphism(1060,500)(0,-1)[` `((A,v),C)]{450}1l
\putmorphism(1500,500)(0,-1)[` `\alpha_{A,B,C}]{450}1r
\putmorphism(1500,50)(0,-1)[` `(A,(v,C))]{450}1r
\putmorphism(1040,-400)(1,0)[` `=]{470}1a
\put(1110,200){\fbox{$\alpha_{A,v,C}$}}
\efig}
$$
for every left central 1v-cell $u\colon A\to \tilde A$ and any 1v-cell $v\colon B\to\tilde B$;

\noindent to the above four axioms the analogous four axioms 
$$\axiom{$(f,\rtimes g,C)$}, \,\, \axiom{$(f,\rtimes v,C)$}, \,\, \axiom{$(u,\rtimes g,C)$}, \,\,\axiom{$(u,\rtimes v,C)$}$$ 
hold for right central 1h-cell $g$, right central 1v-cell $v$ and any 1h-cell $f$ and 1v-cell $u$, 
and additionally to the above eight axioms the analogous $8\cdot 2=16$ axioms, which we denote symbolically like this 
$$\axiom{$(f\ltimes,B,h)$}, \,\,\axiom{$(f\ltimes,B,z)$}, \,\,\axiom{$(u\ltimes, B,h)$}, \,\,\axiom{$(u\ltimes, B,z)$}$$ 
$$\axiom{$(f,B,\rtimes h)$}, \,\,\axiom{$(f,B,\rtimes z)$}, \,\,\axiom{$(u,B,\rtimes h)$}, \,\,\axiom{$(u,B,\rtimes z)$}$$ 
$$\axiom{$(A,g\ltimes,h)$}, \,\,\axiom{$(A,g\ltimes,z)$}, \,\,\axiom{$(A, v\ltimes,h)$}, \,\,\axiom{$(A, v\ltimes,z)$}$$ 
$$\axiom{$(A,g,\rtimes h)$}, \,\,\axiom{$(A,g,\rtimes z)$}, \,\,\axiom{$(A,v,\rtimes h)$}, \,\,\axiom{$(A,v,\rtimes z)$},$$
where in $p\ltimes-$ the corresponding 1-cell $p$ is left central and in $-\rtimes p$  the corresponding 1-cell 
$p$ is right central.

\vspace{1cm}

{\bf {\Large Appendix B.}} 

\vspace{1cm}



\begin{prop} \prlabel{char df} \cite[Proposition 3.3]{Fem:Bif}
A lax double functor $\F\colon\Aa\to\Lax_{hop}(\Bb, \Cc)$ of double categories consists of the following: \\
1. lax double functors 
$$(-,A)\colon\Bb\to\Cc\quad\text{ and}\quad (B,-)\colon\Aa\to\Cc$$ 
such that $(-,A)\vert_B=(B,-)\vert_A=(B,A)$, 
for objects $A\in\Aa, B\in\Bb$, \\
2. 2-cells
$$
\scalebox{0.86}{
\bfig
 \putmorphism(-150,50)(1,0)[(B,A)`(B', A)`(k, A)]{600}1a
 \putmorphism(450,50)(1,0)[\phantom{A\ot B}`(B', A') `(B', K)]{680}1a
\putmorphism(-180,50)(0,-1)[\phantom{Y_2}``=]{450}1r
\putmorphism(1100,50)(0,-1)[\phantom{Y_2}``=]{450}1r
\put(350,-190){\fbox{$(k,K)$}}
 \putmorphism(-150,-400)(1,0)[(B,A)`(B,A')`(B,K)]{600}1a
 \putmorphism(450,-400)(1,0)[\phantom{A\ot B}`(B', A') `(k, A')]{680}1a
\efig}
$$

$$
\scalebox{0.86}{
\bfig
\putmorphism(-150,50)(1,0)[(B,A)`(B,A')`(B,K)]{600}1a
\putmorphism(-150,-400)(1,0)[(\tilde B, A)`(\tilde B,A') `(\tilde B,K)]{640}1a
\putmorphism(-180,50)(0,-1)[\phantom{Y_2}``(u,A)]{450}1l
\putmorphism(450,50)(0,-1)[\phantom{Y_2}``(u,A')]{450}1r
\put(0,-180){\fbox{$(u, K)$}}
\efig}
\quad
\scalebox{0.86}{
\bfig
\putmorphism(-150,50)(1,0)[(B,A)`(B',A)`(k,A)]{600}1a
\putmorphism(-150,-400)(1,0)[(B, \tilde A)`(B', \tilde A) `(k,\tilde A)]{640}1a
\putmorphism(-180,50)(0,-1)[\phantom{Y_2}``(B,U)]{450}1l
\putmorphism(450,50)(0,-1)[\phantom{Y_2}``(B',U)]{450}1r
\put(0,-180){\fbox{$(k,U)$}}
\efig}
$$

$$
\scalebox{0.86}{
\bfig
 \putmorphism(-150,500)(1,0)[(B,A)`(B,A) `=]{600}1a
\putmorphism(-180,500)(0,-1)[\phantom{Y_2}`(B, \tilde A) `(B,U)]{450}1l
\put(-20,50){\fbox{$(u,U)$}}
\putmorphism(-150,-400)(1,0)[(\tilde B, \tilde A)`(\tilde B, \tilde A) `=]{640}1a
\putmorphism(-180,50)(0,-1)[\phantom{Y_2}``(u,\tilde A)]{450}1l
\putmorphism(450,50)(0,-1)[\phantom{Y_2}``(\tilde B, U)]{450}1r
\putmorphism(450,500)(0,-1)[\phantom{Y_2}`(\tilde B, A) `(u,A)]{450}1r
\efig}
$$ 
in $\Cc$ for every 1h-cells $A\stackrel{K}{\to} A'$ and $B\stackrel{k}{\to} B'$ and 1v-cells $A\stackrel{U}{\to} \tilde A$ and 
$B\stackrel{u}{\to} \tilde B$ which satisfy: 

\noindent $\bullet$ \quad \axiom{($1_B,K$)}  
$$
\scalebox{0.86}{
\bfig
 \putmorphism(-210,420)(1,0)[(B,A)`(B,A)`=]{550}1a
\putmorphism(-210,50)(1,0)[(B,A)`(B,A) `(1_B,A)]{560}1a
\putmorphism(350,400)(1,0)[\phantom{F(A)}` (B,A') `(B,K)]{600}1a
\putmorphism(360,50)(1,0)[\phantom{F(A)}`(B,A') `(B,K)]{600}1a

\putmorphism(-170,420)(0,-1)[\phantom{Y_2}``=]{350}1l
\putmorphism(320,420)(0,-1)[\phantom{Y_2}``]{370}1r
\putmorphism(300,420)(0,-1)[\phantom{Y_2}``=]{370}0r
\put(-120,250){\fbox{$(-,A)_B$}}

\putmorphism(860,420)(0,-1)[\phantom{Y_2}``=]{350}1r
\put(450,240){\fbox{$\Id_{(B,K)}$}}
\putmorphism(-180,50)(0,-1)[\phantom{Y_2}``=]{350}1r
\putmorphism(860,50)(0,-1)[\phantom{Y_2}``=]{350}1r
 \putmorphism(-150,-300)(1,0)[(B,A)`(B,A')`(B,K)]{500}1a
 \putmorphism(350,-300)(1,0)[\phantom{A\ot B}`(B, A') `(1_B, A')]{600}1a
\put(170,-110){\fbox{$(1_B,K)$}}
\efig}
\quad=\quad
\scalebox{0.86}{
\bfig
 \putmorphism(-260,200)(1,0)[(B,A)`\phantom{F(A)} `(B,K)]{550}1a
\putmorphism(330,200)(1,0)[(B,A')`(B,A')`=]{600}1a
\putmorphism(-210,200)(0,-1)[\phantom{Y_2}``=]{370}1l
\putmorphism(320,200)(0,-1)[\phantom{Y_2}``]{370}1l
\putmorphism(340,200)(0,-1)[\phantom{Y_2}``=]{370}0l
\putmorphism(960,200)(0,-1)[\phantom{Y_2}``=]{370}1r
 \putmorphism(-260,-170)(1,0)[(B,A)`\phantom{Y_2}`(B,K)]{520}1a
 \put(420,30){\fbox{$(-,A')_{B}$}} 
\putmorphism(360,-170)(1,0)[(B,A')`(B,A') `(1_B,A')]{620}1a
\put(-160,30){\fbox{$\Id_{(B,K)}$}}
\efig}
$$


\noindent $\bullet$ \quad \axiom{($k,1_A$)} \vspace{-0,9cm}
$$ 
\scalebox{0.86}{
\bfig
 \putmorphism(-150,200)(1,0)[(B,A)`(B,A)`=]{500}1a
\putmorphism(360,200)(1,0)[\phantom{F(A)}`(B,A') `(k,A)]{500}1a
\putmorphism(830,200)(0,-1)[\phantom{Y_2}``=]{350}1r
\put(460,30){\fbox{$\Id_{(k,A)}$}}

\putmorphism(-180,200)(0,-1)[\phantom{Y_2}``=]{370}1l
\putmorphism(320,200)(0,-1)[\phantom{Y_2}``]{370}1r
\putmorphism(300,200)(0,-1)[\phantom{Y_2}``=]{370}0r
 \putmorphism(-150,-170)(1,0)[(B,A)`(B,A)`(B,1_A)]{500}1a
 \put(-140,30){\fbox{$(B,-)_A$}} 
\putmorphism(350,-170)(1,0)[\phantom{F(A)}`(B',A) `(k,A)]{560}1a
\efig}
\quad
=
\quad
\scalebox{0.86}{
\bfig
\putmorphism(-150,420)(1,0)[(B,A)` \phantom{Y_2}`(k,A)]{450}1a
\putmorphism(370,420)(1,0)[(B',A)` (B',A) `=]{500}1a
\putmorphism(-170,420)(0,-1)[\phantom{Y_2}``=]{350}1l
\putmorphism(350,420)(0,-1)[\phantom{Y_2}``=]{350}1l
\putmorphism(860,420)(0,-1)[\phantom{Y_2}``=]{350}1r
 \putmorphism(420,50)(1,0)[\phantom{Y_2}`(B',A)`(B',1_A)]{520}1a
\putmorphism(-150,50)(1,0)[(B,A)` (B',A) `(k,A)]{500}1a
 \put(410,250){\fbox{$(B',-)_A$}} 
\put(-120,250){\fbox{$\Id_{(k,A)}$}}

\putmorphism(-180,50)(0,-1)[\phantom{Y_2}``=]{350}1r
\putmorphism(860,50)(0,-1)[\phantom{Y_2}``=]{350}1r
\putmorphism(-150,-300)(1,0)[(B,A)`(B,A)`(B,1_A)]{500}1a
 \putmorphism(350,-300)(1,0)[\phantom{A\ot B}`(B', A) `(k, A)]{560}1a
\put(200,-120){\fbox{$(k,1_A)$}}
\efig}
$$

\noindent where the 2-cells $(-,A)_{B}$ and $(B,-)_A$ come from laxity of the lax double functors $(-,A)$ and $(B,-)$ 

\noindent $\bullet$ \quad  \axiom{($u,1_A$)} 
$$
\scalebox{0.86}{
\bfig
\putmorphism(-250,500)(1,0)[(B,A)`(B,A)` =]{550}1a
 \putmorphism(-250,50)(1,0)[(\tilde B,A)`(\tilde B,A)` =]{550}1a
 \putmorphism(-250,-400)(1,0)[(\tilde B,A)`(\tilde B,A)` (\tilde B,1_A)]{550}1a

\putmorphism(-280,500)(0,-1)[\phantom{Y_2}``(u,A)]{450}1l
 \putmorphism(-280,70)(0,-1)[\phantom{F(A)}` `=]{450}1l

\putmorphism(300,500)(0,-1)[\phantom{Y_2}``(u,A)]{450}1r
\putmorphism(300,70)(0,-1)[\phantom{Y_2}``=]{450}1r
\put(-150,270){\fbox{$Id^{(u,A)}$}}
\put(-190,-150){\fbox{$(\tilde B,-)_A$}}
\efig}
=
\scalebox{0.86}{
\bfig
\putmorphism(-250,500)(1,0)[(B,A)`(B,A)` =]{550}1a
 \putmorphism(-250,50)(1,0)[(B,A)`(B,A)` (B,1_A)]{550}1a
 \putmorphism(-250,-400)(1,0)[(\tilde B,A)`(\tilde B,A)` (\tilde B,1_A)]{550}1a

\putmorphism(-280,500)(0,-1)[\phantom{Y_2}``= ]{450}1l
 \putmorphism(-280,70)(0,-1)[\phantom{F(A)}` `(u,A)]{450}1l

\putmorphism(300,500)(0,-1)[\phantom{Y_2}``=]{450}1r
\putmorphism(300,70)(0,-1)[\phantom{Y_2}``(u,A)]{450}1r
\put(-160,290){\fbox{$(B,-)_A$}}
\put(-140,-150){\fbox{$(u,1_A)$}}
\efig}
$$

\noindent $\bullet$ \quad \axiom{($1_B,U$)} 
$$\scalebox{0.86}{
\bfig
\putmorphism(-250,500)(1,0)[(B,A)`(B,A)` =]{550}1a
 \putmorphism(-250,50)(1,0)[(B,A)`(B,A)` (1_B,A)]{550}1a
 \putmorphism(-250,-400)(1,0)[(B,\tilde A)`(B,\tilde A)` (1_B,\tilde A)]{550}1a

\putmorphism(-280,500)(0,-1)[\phantom{Y_2}``= ]{450}1l
 \putmorphism(-280,70)(0,-1)[\phantom{F(A)}` `(B,U)]{450}1l

\putmorphism(300,500)(0,-1)[\phantom{Y_2}``=]{450}1r
\putmorphism(300,70)(0,-1)[\phantom{Y_2}``(B,U)]{450}1r
\put(-170,290){\fbox{$(-,A)_B$}}
\put(-150,-160){\fbox{$(1_B,U)$}}
\efig}\quad
=
\scalebox{0.86}{
\bfig
\putmorphism(-250,500)(1,0)[(B,A)`(B,A)` =]{550}1a
 \putmorphism(-250,50)(1,0)[(B,\tilde A)`(B,\tilde A)` =]{550}1a
 \putmorphism(-250,-400)(1,0)[(B,\tilde A)`(B,\tilde A)` (\tilde B,1_{\tilde A})]{550}1a

\putmorphism(-280,500)(0,-1)[\phantom{Y_2}``(B,U)]{450}1l
 \putmorphism(-280,70)(0,-1)[\phantom{F(A)}` `=]{450}1l

\putmorphism(300,500)(0,-1)[\phantom{Y_2}``(B,U)]{450}1r
\putmorphism(300,70)(0,-1)[\phantom{Y_2}``=]{450}1r
\put(-170,270){\fbox{$Id^{(B,U)}$}}
\put(-180,-160){\fbox{$(-, \tilde A)_B$}}
\efig}
$$

\noindent $\bullet$ \qquad \axiom{($1^B,K$)}  \quad $(1^B,K)=Id_{(B,K)}$ \qquad\hspace{-0,16cm}\text{and}\qquad 
$\bullet$ \quad  \axiom{($k,1^A$)} \qquad $(k,1^A)=Id_{(k,A)}$ 

\noindent $\bullet$ \qquad  \axiom{($1^B,U$)} \quad $(1^B,U)=Id^{(B,U)}$ \qquad\text{and}\qquad \hspace{-0,22cm} 
$\bullet$ \quad  \axiom{($u,1^A$)} \qquad $(u,1^A)=Id^{(u,A)}$

\noindent $\bullet$ \quad  \axiom{($k'k,K$)} 
$$
\scalebox{0.78}{
\bfig
 \putmorphism(450,650)(1,0)[(B', A) `(B'', A) `(k', A)]{680}1a
 \putmorphism(1140,650)(1,0)[\phantom{A\ot B}`(B'', A') ` (B'', K)]{680}1a

 \putmorphism(-150,200)(1,0)[(B, A) `(B', A)`(k, A)]{600}1a
 \putmorphism(450,200)(1,0)[\phantom{A\ot B}`(B', A') `(B', K)]{680}1a
 \putmorphism(1130,200)(1,0)[\phantom{A\ot B}`(B'', A') ` (k', A')]{680}1a

\putmorphism(450,650)(0,-1)[\phantom{Y_2}``=]{450}1r
\putmorphism(1750,650)(0,-1)[\phantom{Y_2}``=]{450}1r
\put(1000,420){\fbox{$ (k',K)$}}

 \putmorphism(-150,-250)(1,0)[(B, A)`(B, A') `(B,K)]{640}1a
 \putmorphism(480,-250)(1,0)[\phantom{A'\ot B'}`(B', A') `(k, A')]{680}1a

\putmorphism(-180,200)(0,-1)[\phantom{Y_2}``=]{450}1l
\putmorphism(1120,200)(0,-1)[\phantom{Y_3}``=]{450}1r
\put(310,-50){\fbox{$ (k,K)$}}

 \putmorphism(1170,-250)(1,0)[\phantom{A\ot B}`(B'', A') ` (k', A')]{650}1a
\putmorphism(450,-250)(0,-1)[\phantom{Y_2}``=]{450}1r
\putmorphism(1750,-250)(0,-1)[\phantom{Y_2}``=]{450}1r

 \putmorphism(480,-700)(1,0)[(B, A') `(B'', A'') `(k'k, A')]{1320}1a
\put(920,-470){\fbox{$ (-,A')_{k'k}$}}
\efig}
=\quad
\scalebox{0.78}{
\bfig
 \putmorphism(-150,500)(1,0)[(B, A) `(B', A)`(k,A)]{580}1a
 \putmorphism(450,500)(1,0)[\phantom{(B, A)} `(B'', A) `(k',A)]{660}1a
\putmorphism(-180,500)(0,-1)[\phantom{Y_2}``=]{450}1r
\put(240,270){\fbox{$ (-,A)_{k'k}$}}

 \putmorphism(-150,50)(1,0)[(B,A)` `(k'k,A)]{1080}1a
 \putmorphism(1080,50)(1,0)[(B'',A)`(B'', A') ` (B'',K)]{680}1a

\putmorphism(1050,500)(0,-1)[\phantom{Y_2}``=]{450}1r

 \putmorphism(-150,-400)(1,0)[(B, A)`(B, A') `(B,K)]{640}1a
 \putmorphism(530,-400)(1,0)[\phantom{Y_2X}`(B'', A') `(k'k,A')]{1220}1a
\put(570,-200){\fbox{$ (k'k,K)$}}

\putmorphism(-180,50)(0,-1)[\phantom{Y_2}``=]{450}1l
\putmorphism(1700,50)(0,-1)[\phantom{Y_3}``=]{450}1r
\efig}
$$ 
where $(-,A)_{k'k}$ is the 2-cell from the laxity of $(-,A)$ 

\noindent $\bullet$ \quad  \axiom{($k,K'K$)}  
$$
\scalebox{0.78}{
\bfig
 \putmorphism(450,500)(1,0)[(B', A) `(B', A') `(B', K)]{680}1a
 \putmorphism(1140,500)(1,0)[\phantom{A\ot B}`(B', A'') ` (B', K')]{680}1a

 \putmorphism(-150,50)(1,0)[(B, A) `(B', A)`(k,A)]{600}1a
 \putmorphism(450,50)(1,0)[\phantom{A\ot B}`(B', A'') `(B', K'K)]{1350}1a

\putmorphism(450,500)(0,-1)[\phantom{Y_2}``=]{450}1r
\putmorphism(1750,500)(0,-1)[\phantom{Y_2}``=]{450}1r
\put(880,290){\fbox{$ (B',-)_{K'K}$}}

 \putmorphism(-150,-400)(1,0)[(B, A)`(B, A'') `(B, K'K)]{980}1a
 \putmorphism(780,-400)(1,0)[\phantom{A'\ot B'}`(B', A'') `(k,A'')]{980}1a

\putmorphism(-180,50)(0,-1)[\phantom{Y_2}``=]{450}1l
\putmorphism(1750,50)(0,-1)[\phantom{Y_3}``=]{450}1r
\put(560,-200){\fbox{$ (k,K'K)$}}

\efig}
=
\scalebox{0.78}{
\bfig
 \putmorphism(-150,450)(1,0)[(B,A)`(B',A)`(k,A)]{600}1a
 \putmorphism(450,450)(1,0)[\phantom{A\ot B}`(B', A') `(B',K)]{680}1a

 \putmorphism(-150,0)(1,0)[(B,A)`(B,A')`(B,K)]{600}1a
 \putmorphism(450,0)(1,0)[\phantom{A\ot B}`(B', A') `(k,A')]{680}1a
 \putmorphism(1120,0)(1,0)[\phantom{A'\ot B'}`(B', A'') `(B', K')]{660}1a

\putmorphism(-180,450)(0,-1)[\phantom{Y_2}``=]{450}1r
\putmorphism(1100,450)(0,-1)[\phantom{Y_2}``=]{450}1r
\put(350,210){\fbox{$(k,K)$}}
\put(1000,-250){\fbox{$(k,K')$}}

 \putmorphism(450,-450)(1,0)[\phantom{A''\ot B'}` (B, A'') `(B, K')]{680}1a
 \putmorphism(1100,-450)(1,0)[\phantom{A''\ot B'}`(B', A'') ` (k, A'')]{660}1a

\putmorphism(450,0)(0,-1)[\phantom{Y_2}``=]{450}1l
\putmorphism(1750,0)(0,-1)[\phantom{Y_2}``=]{450}1r
 \putmorphism(-150,-450)(1,0)[(B,A)`(B,A')`(B,K)]{600}1a
\putmorphism(-180,-450)(0,-1)[\phantom{Y_2}``=]{450}1r
\putmorphism(1100,-450)(0,-1)[\phantom{Y_2}``=]{450}1r
 \putmorphism(-150,-900)(1,0)[(B,A)`(B,A'')`(B, K'K)]{1280}1a
\put(260,-670){\fbox{$(B,-)_{K'K}$}}
\efig}
$$
where $(B,-)_{K'K}$ is the 2-cell from the laxity of $(B,-)$ 

\noindent $\bullet$ \quad   \axiom{($u, K'K$)}  
$$
\scalebox{0.86}{
\bfig
\putmorphism(-150,500)(1,0)[(B,A)`(B,A')`(B,K)]{600}1a
 \putmorphism(470,500)(1,0)[\phantom{F(A)}`(B,A'') `(B, K')]{600}1a
 \putmorphism(-150,50)(1,0)[(B,A)`(B,A'')`(B,K'K)]{1220}1a

\putmorphism(-180,500)(0,-1)[\phantom{Y_2}``=]{450}1r
\putmorphism(1080,500)(0,-1)[\phantom{Y_2}``=]{450}1r
\put(250,290){\fbox{$(B,-)_{K'K}$}}

\putmorphism(-150,-400)(1,0)[(\tilde B,A)`(\tilde B,A'') `(\tilde B,K'K)]{1200}1a

\putmorphism(-180,50)(0,-1)[\phantom{Y_2}``(u,A)]{450}1l
\putmorphism(1080,50)(0,-1)[\phantom{Y_3}``(u,A'')]{450}1r
\put(250,-160){\fbox{$(u,K'K)$}} 
\efig}
=
\scalebox{0.86}{
\bfig
\putmorphism(-150,500)(1,0)[(B,A)`(B,A')`(B,K)]{600}1a
 \putmorphism(470,500)(1,0)[\phantom{F(A)}`(B,A'') `(B, K')]{600}1a

 \putmorphism(-150,50)(1,0)[(\tilde B, A)`(\tilde B,A')`(\tilde B,K)]{600}1a
 \putmorphism(470,50)(1,0)[\phantom{F(A)}`(B'',\tilde A) `(\tilde B,K')]{620}1a

\putmorphism(-180,500)(0,-1)[\phantom{Y_2}``(u,A)]{450}1l
\putmorphism(450,500)(0,-1)[\phantom{Y_2}``]{450}1r
\putmorphism(300,500)(0,-1)[\phantom{Y_2}``(u,A')]{450}0r
\putmorphism(1080,500)(0,-1)[\phantom{Y_2}``(u,A'')]{450}1r
\put(-40,280){\fbox{$(u,K)$}}
\put(620,280){\fbox{$(u,K')$}}

\putmorphism(-150,-400)(1,0)[(\tilde B, A)`(\tilde B,A'') `(\tilde B,K'K)]{1200}1a

\putmorphism(-180,50)(0,-1)[\phantom{Y_2}``=]{450}1l
\putmorphism(1080,50)(0,-1)[\phantom{Y_3}``=]{450}1r
\put(260,-160){\fbox{$(\tilde B,-)_{K'K}$}}

\efig}
$$

\noindent $\bullet$ \quad  \axiom{($k'k, U$)} 
$$
\scalebox{0.86}{
\bfig
\putmorphism(-150,500)(1,0)[(B,A)`(B',A)`(k,A)]{600}1a
 \putmorphism(450,500)(1,0)[\phantom{F(A)}`(B'',A) `(k',A)]{620}1a

 \putmorphism(-150,50)(1,0)[(B,\tilde A)`(B',\tilde A)`(k,\tilde A)]{600}1a
 \putmorphism(470,50)(1,0)[\phantom{F(A)}`(B'',\tilde A) `(k',\tilde A)]{620}1a

\putmorphism(-180,500)(0,-1)[\phantom{Y_2}``(B,U)]{450}1l
\putmorphism(450,500)(0,-1)[\phantom{Y_2}``]{450}1r
\putmorphism(240,500)(0,-1)[\phantom{Y_2}``(B',U)]{450}0r
\putmorphism(1080,500)(0,-1)[\phantom{Y_2}``(B'',U)]{450}1r
\put(-40,280){\fbox{$(k,U)$}}
\put(620,280){\fbox{$(k',U)$}}

\putmorphism(-150,-400)(1,0)[(B,\tilde A)`(B'',\tilde A) `(k'k,\tilde A)]{1200}1a

\putmorphism(-180,50)(0,-1)[\phantom{Y_2}``=]{450}1l
\putmorphism(1080,50)(0,-1)[\phantom{Y_3}``=]{450}1r
\put(260,-160){\fbox{$(-,\tilde A)_{k'k}$}}
\efig}
=
\scalebox{0.86}{
\bfig
\putmorphism(-150,500)(1,0)[(B,A)`(B',A)`(k,A)]{600}1a
 \putmorphism(450,500)(1,0)[\phantom{F(A)}`(B'',A) `(k',A)]{620}1a
 \putmorphism(-150,50)(1,0)[(B,A)`(B'',A)`(k'k, A)]{1220}1a

\putmorphism(-180,500)(0,-1)[\phantom{Y_2}``=]{450}1r
\putmorphism(1080,500)(0,-1)[\phantom{Y_2}``=]{450}1r
\put(260,290){\fbox{$(-,A)_{k'k}$}}

\putmorphism(-150,-400)(1,0)[(B,\tilde A)`(B'',\tilde A) `(k'k,\tilde A)]{1200}1a

\putmorphism(-180,50)(0,-1)[\phantom{Y_2}``(B,U)]{450}1l
\putmorphism(1080,50)(0,-1)[\phantom{Y_3}``(B'',U)]{450}1r
\put(300,-180){\fbox{$(k'k,U)$}} 
\efig}
$$

\noindent $\bullet$ \quad   \axiom{($\frac{u}{u'}, K$)} \qquad $(\frac{u}{u'}, K)=\frac{(u,K)}{(u', K)}$  \qquad\text{and}\qquad 
$\bullet$ \quad  \axiom{($k,\frac{U}{U'}$)}  \qquad $(k,\frac{U}{U'})=\frac{(k,U)}{(k,U')}$ 

\noindent $\bullet$ \quad  \axiom{($u,\frac{U}{U'}$)}
$$(u,\frac{U}{U'})=
\scalebox{0.86}{
\bfig
 \putmorphism(-150,500)(1,0)[(B,A)`(B,A) `=]{600}1a
\putmorphism(-180,500)(0,-1)[\phantom{Y_2}`(B, \tilde A) `(B,U)]{450}1l
\put(0,50){\fbox{$(u,U)$}}
\putmorphism(-150,-400)(1,0)[(\tilde B, \tilde A)`(\tilde B, \tilde A) `=]{640}1a
\putmorphism(-180,50)(0,-1)[\phantom{Y_2}``(u,\tilde A)]{450}1l
\putmorphism(450,50)(0,-1)[\phantom{Y_2}``(\tilde B, U)]{450}1r
\putmorphism(450,500)(0,-1)[\phantom{Y_2}`(\tilde B, A) `(u,A)]{450}1r
\putmorphism(-820,50)(1,0)[(B, \tilde A)``=]{520}1a
\putmorphism(-820,50)(0,-1)[\phantom{(B, \tilde A')}``(B,U')]{450}1l
\putmorphism(-820,-400)(0,-1)[(B, \tilde A')`(\tilde B, \tilde A')`(u,\tilde A')]{450}1l
\putmorphism(-820,-850)(1,0)[\phantom{(B, \tilde A)}``=]{520}1a
\putmorphism(-150,-400)(0,-1)[(\tilde B, \tilde A)`(\tilde B, \tilde A') `(\tilde B, U')]{450}1r
\put(-650,-630){\fbox{$(u,U')$}}
\efig}
$$

\noindent $\bullet$ \quad  \axiom{($\frac{u}{u'},U$)} 
$$(\frac{u}{u'},U)=
\scalebox{0.86}{
\bfig
 \putmorphism(-150,500)(1,0)[(B,A)`(B,A) `=]{600}1a
\putmorphism(-180,500)(0,-1)[\phantom{Y_2}`(B, \tilde A) `(B,U)]{450}1l
\put(0,50){\fbox{$(u,U)$}}
\putmorphism(-150,-400)(1,0)[(\tilde B, \tilde A)` `=]{500}1a
\putmorphism(-180,50)(0,-1)[\phantom{Y_2}``(u,\tilde A)]{450}1l
\putmorphism(450,50)(0,-1)[\phantom{Y_2}`(\tilde B, \tilde A)`(\tilde B, U)]{450}1r
\putmorphism(450,500)(0,-1)[\phantom{Y_2}`(\tilde B, A) `(u,A)]{450}1r
\putmorphism(450,50)(1,0)[\phantom{(B, \tilde A)}`(\tilde B, A)`=]{620}1a
\putmorphism(1070,50)(0,-1)[\phantom{(B, \tilde A')}``(u',A)]{450}1r
\putmorphism(1070,-400)(0,-1)[(\tilde B', A)`(\tilde B', \tilde A)`(\tilde B', U)]{450}1r
\putmorphism(450,-850)(1,0)[\phantom{(B, \tilde A)}``=]{520}1a
\putmorphism(450,-400)(0,-1)[\phantom{(B, \tilde A)}`(\tilde B', \tilde A) `(U',\tilde A)]{450}1l
\put(600,-630){\fbox{$(u',U)$}}
\efig}
$$

\noindent $\bullet$ \quad  \axiom{$(k,K)$-l-nat}  
$$
\scalebox{0.86}{
\bfig
 \putmorphism(-150,500)(1,0)[(B,A)`(B', A)`(k, A)]{600}1a
 \putmorphism(450,500)(1,0)[\phantom{A\ot B}`(B', A') `(B', K)]{680}1a
 \putmorphism(-150,50)(1,0)[(B,A)`(B,A')`(B,K)]{600}1a
 \putmorphism(450,50)(1,0)[\phantom{A\ot B}`(B', A') `(k, A')]{680}1a

\putmorphism(-180,500)(0,-1)[\phantom{Y_2}``=]{450}1r
\putmorphism(1100,500)(0,-1)[\phantom{Y_2}``=]{450}1r
\put(350,260){\fbox{$(k,K)$}}

\putmorphism(-150,-400)(1,0)[(\tilde B,A)`(\tilde B,A') `(\tilde B, K)]{640}1a
 \putmorphism(450,-400)(1,0)[\phantom{A'\ot B'}` (\tilde B',A') `(l,A')]{680}1a

\putmorphism(-180,50)(0,-1)[\phantom{Y_2}``(u,A)]{450}1l
\putmorphism(450,50)(0,-1)[\phantom{Y_2}``]{450}1r
\putmorphism(300,50)(0,-1)[\phantom{Y_2}``(u,A')]{450}0r
\putmorphism(1100,50)(0,-1)[\phantom{Y_2}``]{450}1r
\putmorphism(1080,50)(0,-1)[\phantom{Y_2}``(v,A')]{450}0r
\put(-20,-180){\fbox{$(u,K)$}}
\put(660,-180){\fbox{$(\omega,A')$}}

\efig}
\quad=\quad
\scalebox{0.86}{
\bfig
 \putmorphism(-150,500)(1,0)[(B,A)`(B', A)`(k, A)]{600}1a
 \putmorphism(450,500)(1,0)[\phantom{A\ot B}`(B', A') `(B', K)]{680}1a

 \putmorphism(-150,50)(1,0)[(\tilde B,A)`(\tilde B',A)`(l, A)]{600}1a
 \putmorphism(450,50)(1,0)[\phantom{A\ot B}`(\tilde B',A') `(\tilde B',K)]{680}1a
\putmorphism(-180,500)(0,-1)[\phantom{Y_2}``]{450}1l
\putmorphism(-160,500)(0,-1)[\phantom{Y_2}``(u,A)]{450}0l
\putmorphism(450,500)(0,-1)[\phantom{Y_2}``]{450}1l
\putmorphism(610,500)(0,-1)[\phantom{Y_2}``(v,A)]{450}0l 
\putmorphism(1120,500)(0,-1)[\phantom{Y_3}``(v,A')]{450}1r
\put(-40,270){\fbox{$(\omega,A)$}} 
\put(650,270){\fbox{$(v,K)$}}
\putmorphism(-150,-400)(1,0)[(\tilde B,A)`(\tilde B,A') `(\tilde B, K)]{640}1a
 \putmorphism(450,-400)(1,0)[\phantom{A'\ot B'}` (\tilde B',A') `(l,A')]{680}1a

\putmorphism(-180,50)(0,-1)[\phantom{Y_2}``=]{450}1l
\putmorphism(1120,50)(0,-1)[\phantom{Y_3}``=]{450}1r
\put(300,-200){\fbox{$(l,K)$}}

\efig}
$$

\noindent $\bullet$ \quad  \axiom{$(k,K)$-r-nat}  
$$
\scalebox{0.86}{
\bfig
 \putmorphism(-150,500)(1,0)[(B,A)`(B', A)`(k, A)]{600}1a
 \putmorphism(450,500)(1,0)[\phantom{A\ot B}`(B', A') `(B', K)]{680}1a
 \putmorphism(-150,50)(1,0)[(B,A)`(B,A')`(B,K)]{600}1a
 \putmorphism(450,50)(1,0)[\phantom{A\ot B}`(B', A') `(k, A')]{680}1a

\putmorphism(-180,500)(0,-1)[\phantom{Y_2}``=]{450}1r
\putmorphism(1100,500)(0,-1)[\phantom{Y_2}``=]{450}1r
\put(350,260){\fbox{$(k,K)$}}

\putmorphism(-180,50)(0,-1)[\phantom{Y_2}``(B,U)]{450}1l
\putmorphism(450,50)(0,-1)[\phantom{Y_2}``]{450}1r
\putmorphism(300,50)(0,-1)[\phantom{Y_2}``(B,V)]{450}0r
\putmorphism(1100,50)(0,-1)[\phantom{Y_2}``]{450}1r
\putmorphism(1080,50)(0,-1)[\phantom{Y_2}``(B',V)]{450}0r
\put(-20,-180){\fbox{$(B,\zeta)$}}
\put(660,-180){\fbox{$(k,V)$}}

\putmorphism(-150,-400)(1,0)[(B,\tilde A)`(B,\tilde A') `(B,L)]{640}1a
 \putmorphism(450,-400)(1,0)[\phantom{A'\ot B'}` (B',\tilde A') `(k,\tilde A')]{680}1a
\efig}
\quad=\quad
\scalebox{0.86}{
\bfig
 \putmorphism(-150,500)(1,0)[(B,A)`(B', A)`(k, A)]{600}1a
 \putmorphism(450,500)(1,0)[\phantom{A\ot B}`(B', A') `(B', K)]{680}1a

 \putmorphism(-150,50)(1,0)[(B,\tilde A)`(B',\tilde A)`(k,\tilde A)]{600}1a
 \putmorphism(450,50)(1,0)[\phantom{A\ot B}`(B',\tilde A') `(B',L)]{680}1a

\putmorphism(-180,500)(0,-1)[\phantom{Y_2}``]{450}1l
\putmorphism(-160,500)(0,-1)[\phantom{Y_2}``(B,U)]{450}0l
\putmorphism(450,500)(0,-1)[\phantom{Y_2}``]{450}1l
\putmorphism(610,500)(0,-1)[\phantom{Y_2}``(B',U)]{450}0l 
\putmorphism(1120,500)(0,-1)[\phantom{Y_3}``(B',V)]{450}1r
\put(-40,270){\fbox{$(k,U)$}} 
\put(650,270){\fbox{$(B', \zeta)$}}

\putmorphism(-180,50)(0,-1)[\phantom{Y_2}``=]{450}1l
\putmorphism(1120,50)(0,-1)[\phantom{Y_3}``=]{450}1r
\put(300,-200){\fbox{$(k,L)$}}

\putmorphism(-150,-400)(1,0)[(B,\tilde A)`(B,\tilde A') `(B,L)]{640}1a
 \putmorphism(450,-400)(1,0)[\phantom{A'\ot B'}` (B',\tilde A') `(k,\tilde A')]{680}1a
\efig}
$$

\noindent $\bullet$ \quad  \axiom{$(u,U)$-l-nat} 
$$
\scalebox{0.86}{
\bfig
 \putmorphism(-150,500)(1,0)[(B,A)`(B,A) `=]{600}1a
 \putmorphism(450,500)(1,0)[(B,A)` `(k,A)]{450}1a
\putmorphism(-180,500)(0,-1)[\phantom{Y_2}`(B, \tilde A) `(B,U)]{450}1l
\put(-20,50){\fbox{$(u,U)$}}
\putmorphism(-170,-400)(1,0)[(\tilde B, \tilde A)` `=]{480}1a
\putmorphism(-180,50)(0,-1)[\phantom{Y_2}``(u,\tilde A)]{450}1l
\putmorphism(450,50)(0,-1)[\phantom{Y_2}`(\tilde B, \tilde A)`(\tilde B, U)]{450}1l
\putmorphism(450,500)(0,-1)[\phantom{Y_2}`(\tilde B, A) `(u,A)]{450}1l
\put(600,260){\fbox{$(\omega,A)$}}
\putmorphism(430,50)(1,0)[\phantom{(B, \tilde A)}``(l, A)]{500}1a
\putmorphism(1070,50)(0,-1)[\phantom{(B, A')}`(\tilde B', \tilde A)`]{450}1r
\putmorphism(1050,50)(0,-1)[``(\tilde B',U)]{450}0r
\putmorphism(1070,500)(0,-1)[(B', A)`(\tilde B', A)`]{450}1r
\putmorphism(1050,500)(0,-1)[``(v,A)]{450}0r
\putmorphism(450,-400)(1,0)[\phantom{(B, \tilde A)}``(l, \tilde A)]{500}1a
\put(600,-170){\fbox{$(l,U)$}}
\efig}
\quad=\quad
\scalebox{0.86}{
\bfig
 \putmorphism(-150,500)(1,0)[(B,A)`(B',A) `(k,A)]{600}1a
 \putmorphism(450,500)(1,0)[\phantom{(B,A)}` `=]{450}1a
\putmorphism(-180,500)(0,-1)[\phantom{Y_2}`(B, \tilde A) `]{450}1l
\putmorphism(-160,500)(0,-1)[` `(B,U)]{450}0l
\put(620,50){\fbox{$(v,U)$}}
\putmorphism(-170,-400)(1,0)[(\tilde B, \tilde A)` `(l, \tilde A)]{470}1a
\putmorphism(-180,50)(0,-1)[\phantom{Y_2}``]{450}1l
\putmorphism(-160,50)(0,-1)[\phantom{Y_2}``(u,\tilde A)]{450}0l
\putmorphism(450,50)(0,-1)[\phantom{Y_2}`(\tilde B', \tilde A)`(v,\tilde A)]{450}1r
\putmorphism(450,500)(0,-1)[\phantom{Y_2}`(B', \tilde A) `(B',U)]{450}1r
\put(0,260){\fbox{$(k,U)$}}
\putmorphism(-190,50)(1,0)[\phantom{(B, \tilde A)}``(k, \tilde A)]{500}1a
\putmorphism(1070,50)(0,-1)[\phantom{(B, A')}`(\tilde B', \tilde A)`(\tilde B',U)]{450}1r
\putmorphism(1070,500)(0,-1)[(B', A)``(v,A)]{450}1r
\putmorphism(1100,500)(0,-1)[`(\tilde B', A)`]{450}0r
\putmorphism(450,-400)(1,0)[\phantom{(B, \tilde A)}``=]{500}1b
\put(0,-170){\fbox{$(\omega, \tilde A)$}}
\efig}
$$

\noindent $\bullet$ \quad \axiom{$(u,U)$-r-nat} 
$$
\scalebox{0.86}{
\bfig
 \putmorphism(-150,500)(1,0)[(B,A)`(B,A) `=]{600}1a
 \putmorphism(550,500)(1,0)[` `(B,K)]{380}1a
\putmorphism(-180,500)(0,-1)[\phantom{Y_2}`(B, \tilde A) `(B,U)]{450}1l
\put(-20,50){\fbox{$(u,U)$}}
\putmorphism(-170,-400)(1,0)[(\tilde B, \tilde A)` `=]{480}1a
\putmorphism(-180,50)(0,-1)[\phantom{Y_2}``(u,\tilde A)]{450}1l
\putmorphism(450,50)(0,-1)[\phantom{Y_2}`(\tilde B, \tilde A)`(\tilde B, U)]{450}1l
\putmorphism(450,500)(0,-1)[\phantom{Y_2}`(\tilde B, A) `(u,A)]{450}1l
\put(620,280){\fbox{$(u,K)$}}
\putmorphism(430,50)(1,0)[\phantom{(B, \tilde A)}``(\tilde B,K)]{500}1a
\putmorphism(1070,50)(0,-1)[\phantom{(B, A')}`(\tilde B, \tilde A')`]{450}1r
\putmorphism(1090,50)(0,-1)[\phantom{(B, A')}``(\tilde B,V)]{450}0r
\putmorphism(1070,500)(0,-1)[(B, A')`(\tilde B, A')`]{450}1r
\putmorphism(1090,500)(0,-1)[``(u,A')]{450}0r
\putmorphism(450,-400)(1,0)[\phantom{(B, \tilde A)}``(\tilde B, L)]{500}1a
\put(620,-170){\fbox{$ (\tilde{B},\zeta)$ } } 
\efig}
\quad=\quad
\scalebox{0.86}{
\bfig
 \putmorphism(-150,500)(1,0)[(B,A)`(B,A') `(B,K)]{600}1a
 \putmorphism(450,500)(1,0)[\phantom{(B,A)}` `=]{450}1a
\putmorphism(-180,500)(0,-1)[\phantom{Y_2}`(B, \tilde A) `]{450}1l
\putmorphism(-160,500)(0,-1)[\phantom{Y_2}` `(B,U)]{450}0l
\put(620,50){\fbox{$(u,V)$}}
\putmorphism(-170,-400)(1,0)[(\tilde B, \tilde A)` `(\tilde B, L)]{500}1a
\putmorphism(-180,50)(0,-1)[\phantom{Y_2}``]{450}1l
\putmorphism(-160,50)(0,-1)[\phantom{Y_2}``(u,\tilde A)]{450}0l
\putmorphism(450,50)(0,-1)[\phantom{Y_2}`(\tilde B, \tilde A')`(u,\tilde A')]{450}1r
\putmorphism(450,500)(0,-1)[\phantom{Y_2}`(B, \tilde A') `(B,V)]{450}1r
\put(0,260){\fbox{$(B,\zeta)$}}
\putmorphism(-190,50)(1,0)[\phantom{(B, \tilde A)}``(B,L)]{500}1a
\putmorphism(1070,50)(0,-1)[\phantom{(B, A')}`(\tilde B, \tilde A')`(\tilde B,V)]{450}1r
\putmorphism(1070,500)(0,-1)[(B, A')``(u,A')]{450}1r
\putmorphism(1110,500)(0,-1)[`(\tilde B, A')`]{450}0r
\putmorphism(450,-400)(1,0)[\phantom{(B, \tilde A)}``=]{500}1b
\put(0,-170){\fbox{$(u,L)$}}
\efig}
$$
for any 2-cells 
\begin{equation} \eqlabel{omega-zeta}
\scalebox{0.86}{
\bfig
\putmorphism(-150,175)(1,0)[B` B'`k]{450}1a
\putmorphism(-150,-175)(1,0)[\tilde B`\tilde B' `l]{440}1b
\putmorphism(-170,175)(0,-1)[\phantom{Y_2}``u]{350}1l
\putmorphism(280,175)(0,-1)[\phantom{Y_2}``v]{350}1r
\put(0,-15){\fbox{$\omega$}}
\efig}
\quad\text{and}\quad
\scalebox{0.86}{
\bfig
\putmorphism(-150,175)(1,0)[A` A'`K]{450}1a
\putmorphism(-150,-175)(1,0)[\tilde A`\tilde A' `L]{440}1b
\putmorphism(-170,175)(0,-1)[\phantom{Y_2}``U]{350}1l
\putmorphism(280,175)(0,-1)[\phantom{Y_2}``V]{350}1r
\put(0,-15){\fbox{$\zeta$}}
\efig}
\end{equation}
in $\Bb$, respectively $\Aa$. 

The data of the points 1. and 2. of this proposition comprise the definition of a {\em lax double quasi-functor} $H:\Aa\times\Bb\to\Cc$. 
\end{prop}

\vspace{1,3cm}

{\bf {\Large Appendix C.}} 

\vspace{0,8cm}

{\bf Table 1 of \cite[Proposition 3.3]{Fem:Bif} enriched by additional interpretations: }

\begin{table}[H]
\begin{center}
\begin{tabular}{ c c c } 
New axiom & \hspace{0,3cm} Origin from $\F\colon\Aa\to\llbracket\Bb,\Cc\rrbracket$ & \hspace{0,3cm}  additional interpretation\\ [0.5ex]
\hline
2-cell $(k,K)$ & part 3 of $(-,K)$ being a h.o.t. 
\\ [1ex]   
2-cell $(u,K)$ & part 2 of $(-,K)$ being a h.o.t.  \\ [1ex]
2-cell $(k,U)$ & part 2 of $(-,U)$ being a v.l.t.  \\ [1ex]
2-cell $(u,U)$ & part 3 of $(-,U)$ being a v.l.t.  \\ [1ex]
\hline
\axiomref{($1_B,K$)}  & \axiomref{h.o.t.-2} of $(-,K)$ \\ [1ex]    \hdashline[0.5pt/5pt]
\axiomref{($k,1_A$)}  & \axiomref{m.ho-vl.-1} of unitor \\ [1ex] 
  & $\F_A\colon \Id_{(-,A)}\Rrightarrow (-,1_A)$ & \axiomref{h.l.t.-2} of $(k,-)$ \\ [1ex]   \hdashline[0.5pt/5pt]
\axiomref{($1^B,K$)}  & \axiomref{h.o.t.-4} of $(-,K)$ \\ [1ex]    \hdashline[0.5pt/5pt]
\axiomref{($u,1_A$)}  & \axiomref{m.ho-vl.-2} of unitor \\ [1ex] 
  & $\F_A\colon \Id_{(-,A)}\Rrightarrow (-,1_A)$ & \axiomref{v.o.t.\x 2} of $(u,-)$ \\ [1ex]   \hdashline[0.5pt/5pt]
\axiomref{($1_B,U$)}  & \axiomref{v.l.t.\x 2} of $(-,U)$ \\ [1ex]   \hdashline[0.5pt/5pt]
\axiomref{($k,1^A$)}  & \axiomref{lx.f.v2} of $\F$ (is an equality of v.l.t.) \\ [1ex] 
  & evaluated at $k$ & \axiomref{h.l.t.-4} of $(k,-)$ \\ [1ex]    \hdashline[0.5pt/5pt]
\axiomref{($1^B,U$)}  & \axiomref{v.l.t.\x 4} of $(-,U)$ \\ [1ex]    \hdashline[0.5pt/5pt]
\axiomref{($u,1^A$)}  & \axiomref{lx.f.v2} of $\F$ (is an equality of v.l.t.) \\ [1ex] 
  & evaluated at $u$ & \axiomref{v.o.t.\x 4} of $(u,-)$ \\ [1ex]  \hdashline[0.5pt/5pt]
\hline
\axiomref{($k'k,K$)}  & \axiomref{h.o.t.-1} of $(-,K)$ \\ [1ex]    \hdashline[0.5pt/5pt]
\axiomref{($k,K'K$)}  & \axiomref{m.ho-vl.-1} of compositor \\ [1ex] 
  & $\F_{LK}\colon (-,L)(-,K)\Rrightarrow (-,LK)$ & \axiomref{h.l.t.-1} of $(k,-)$ \\ [1ex]   \hdashline[0.5pt/5pt]
\axiomref{($\frac{u}{u'},K$)}  & \axiomref{h.o.t.-3} of $(-,K)$ \\ [1ex]    \hdashline[0.5pt/5pt]
\axiomref{($u,K'K$)}  & \axiomref{m.ho-vl.-2} of compositor \\ [1ex] 
  & $\F_{LK}\colon (-,L)(-,K)\Rrightarrow (-,LK)$ & \axiomref{v.o.t.\x 1} of $(u,-)$ \\ [1ex]   \hdashline[0.5pt/5pt]
\axiomref{($k'k,U$)}  & \axiomref{v.l.t.\x 1} of $(-,U)$ \\ [1ex]   \hdashline[0.5pt/5pt]
\axiomref{($k,\frac{U}{U'}$)}  & \axiomref{lx.f.v1} of $\F$ (is an equality of v.l.t.) & \axiomref{h.l.t.-3} of $(k,-)$ \\ [1ex] 
  & evaluated at $k$ {h.o.t.-3} of $(k,-)$ \\ [1ex]  \hdashline[0.5pt/5pt]
\axiomref{($u,\frac{U}{U'}$)}  & \axiomref{lx.f.v1} of $\F$ (is an equality of v.l.t.) \\ [1ex] 
  & evaluated at $u$ & \axiomref{v.o.t.\x 3} of $(u,-)$ \\ [1ex]
\axiomref{($\frac{u}{u'},U$)}  & \axiomref{v.l.t.\x 3} of \axiomref{$-,U$} \\ [1ex]   \hdashline[0.5pt/5pt]
\hline
\axiomref{$(k,K)$-l-nat}  & \axiomref{h.o.t.-5} of \axiomref{$-,K$} & \axiomref{m.hl-vo.-1} of \axiomref{$\omega,-$} \\ [1ex]   
\axiomref{$(k,K)$-r-nat}  & \axiomref{m.ho-vl.-1} of \axiomref{$-,\zeta$} & \axiomref{h.l.t.-5} of $(k,-)$ \\ [1ex]   
\axiomref{$(u,U)$-l-nat}  & \axiomref{v.l.t.\x 5} of \axiomref{$-,U$} & \axiomref{m.hl-vo.-2} of \axiomref{$\omega,-$} \\ [1ex]   
\axiomref{$(u,U)$-r-nat}  & \axiomref{m.ho-vl.-2} of \axiomref{$-,\zeta$} & \axiomref{v.o.t.\x 5} of $(u,-)$ \\ [1ex]   
\end{tabular}
\caption{Generation of a lax double quasi-functor $\Aa\times\Bb\protect\to\Cc$}
\label{table:12}
\end{center}
\end{table}

\pagebreak

{\bf {\Large Appendix D.}}

\vspace{0,5cm}

\begin{defn} \cite[Definition 4.3]{Fem:Bif} \delabel{lax v tr cubical}
A vertical lax transformation $\theta_0\colon (-,-)_1\Rightarrow (-,-)_2$ between lax double quasi-functors $(-,-)_1,(-,-)_2\colon 
\Aa\times\Bb\to\Cc$ is given by: for each $A\in\Aa$ a vertical lax transformation $\theta_0^A\colon (-,A)_1\Rightarrow(-,A)_2$ and 
for each $B\in\Bb$ a vertical lax transformation $\theta_0^B\colon (B,-)_1\Rightarrow(B,-)_2$, both of lax double functors, such that 
$(\theta_0^A)_B=(\theta_0^B)_A$ and such that 

\medskip
\noindent \axiom{$VLT^q_1$} \vspace{-0,8cm}

$$
\scalebox{0.8}{
\bfig
 \putmorphism(-170,500)(1,0)[(B,A)_1`(B,A)_1 `=]{500}1a
\putmorphism(-180,500)(0,-1)[\phantom{Y_2}`(B, A)_2 `(\theta_0^A)_B]{450}1l
\put(-50,250){\fbox{$(\theta_0^A)^u$}}
\putmorphism(-150,-400)(1,0)[(\tilde B, A)_2`(\tilde B, A)_2 `=]{500}1a
\putmorphism(-180,50)(0,-1)[\phantom{Y_2}``(u,A)_2]{450}1l
\putmorphism(300,50)(0,-1)[\phantom{Y_2}``(\theta_0^A)_{\tilde B}]{450}1l
\putmorphism(300,500)(0,-1)[\phantom{Y_2}`(\tilde B, A)_1 `(u,A)_1]{450}1r
\putmorphism(-720,50)(1,0)[(B, A)_2``=]{400}1a
\putmorphism(-720,50)(0,-1)[\phantom{(B, \tilde A')}``(B,U)_2]{450}1l
\putmorphism(-720,-400)(0,-1)[(B, \tilde A)_2`(\tilde B, \tilde A)_2`(u,\tilde A)_2]{450}1l
\putmorphism(-720,-850)(1,0)[\phantom{(B, \tilde A)}``=]{370}1a 
\putmorphism(-180,-400)(0,-1)[`(\tilde B, \tilde A)_2 `(\tilde B, U)_2]{450}1r
\put(-610,-630){\fbox{$(u,U)_2$}}

\putmorphism(430,50)(1,0)[`(\tilde B, A)_1`=]{450}1a
\putmorphism(870,50)(0,-1)[\phantom{Y_2}`(\tilde B, \tilde A)_1`(\tilde B, U)_1]{450}1r
\putmorphism(300,-400)(0,-1)[`(\tilde B, \tilde A)_2 `(\tilde B,U)_2]{450}1r
\putmorphism(870,-400)(0,-1)[`(\tilde B, \tilde A)_2 `(\theta_0^{\tilde B})_{\tilde A}]{450}1r
\putmorphism(450,-850)(1,0)[` `=]{260}1a 
\put(430,-200){\fbox{$(\theta_0^{\tilde B})^U$}}
\efig}
=
\scalebox{0.8}{
\bfig
 \putmorphism(0,500)(1,0)[(B,A)_1`(B,A)_1 `=]{500}1a
\putmorphism(0,500)(0,-1)[\phantom{Y_2}`(B, A)_2 `(\theta^A_0)_B]{450}1l
\put(70,-200){\fbox{$(\theta^B_0)^U$}}
\putmorphism(-40,-400)(1,0)[(B, \tilde A)_2` `=]{330}1a %
\putmorphism(0,50)(0,-1)[\phantom{Y_2}``(B,U)_2]{450}1l
\putmorphism(450,50)(0,-1)[\phantom{Y_2}`(B, \tilde A)_2`(\theta^B_0)_{\tilde A}]{450}1r
\putmorphism(450,500)(0,-1)[\phantom{Y_2}`(B, \tilde A)_1 `(B,U)_1]{450}1l
\putmorphism(450,50)(1,0)[\phantom{(B, \tilde A)}`(B, \tilde A)_1`=]{620}1a %
\putmorphism(970,50)(0,-1)[\phantom{(B, \tilde A')}``(u,\tilde A)_1]{450}1r
\putmorphism(970,-400)(0,-1)[(\tilde B, \tilde A)_1`(\tilde B, \tilde A)_2`(\theta^{\tilde A}_0)_{\tilde B}]{450}1r
\putmorphism(450,-850)(1,0)[\phantom{(B, \tilde A)}``=]{350}1a 
\putmorphism(450,-400)(0,-1)[\phantom{(B, \tilde A)}`(\tilde B, \tilde A) `(u,\tilde A)_2]{450}1l
\put(600,-630){\fbox{$(\theta^{\tilde A}_0)^u$}}

 \putmorphism(970,500)(1,0)[(B,A)_1`(B,A)_1 `=]{500}1a
\putmorphism(970,500)(0,-1)[\phantom{Y_2}` `(B,U)_1]{450}1l
\putmorphism(1500,500)(0,-1)[\phantom{Y_2}`(\tilde B, A)_1 `(u,A)_1]{450}1r
\putmorphism(1500,50)(0,-1)[\phantom{Y_2}`(\tilde B, \tilde A)_1 `(\tilde B,U)_1]{450}1r
\putmorphism(1110,-400)(1,0)[` `=]{260}1a
\put(1050,200){\fbox{$(u,U)_1$}}
\efig}
$$
for every 1v-cells $U\colon A\to \tilde A$ and $u\colon B\to\tilde B$; 



\noindent \axiom{$VLT^q_2$} \vspace{-0,3cm}
$$\scalebox{0.86}{
\bfig
 \putmorphism(-150,500)(1,0)[(B,A)_1`(B,A)_1 `=]{600}1a
 \putmorphism(450,500)(1,0)[\phantom{(B,A)_1}` `(B,K)_1]{450}1a
\putmorphism(-200,500)(0,-1)[\phantom{Y_2}`(B, A)_2 `(\theta_0^A)_B]{450}1l
\put(-20,50){\fbox{$(\theta^A_0)^u$}}
\putmorphism(-170,-400)(1,0)[(\tilde B, A)_2` `=]{500}1a
\putmorphism(-200,50)(0,-1)[\phantom{Y_2}``(u,A)_2]{450}1l
\putmorphism(450,50)(0,-1)[\phantom{Y_2}`(\tilde B, A)_2 `(\theta_0^A)_{\tilde B}]{450}1l
\putmorphism(450,500)(0,-1)[\phantom{Y_2}`(\tilde B, A)_1 `(u,A)_1]{450}1l
\put(600,260){\fbox{$(u,K)_1$}}
\putmorphism(450,50)(1,0)[\phantom{(B, \tilde A)}``(\tilde B,K)_1]{500}1a
\putmorphism(1070,50)(0,-1)[\phantom{(B, A')}`(\tilde B, A')_2`(\theta^{\tilde B}_0)_{A'}]{450}1r
\putmorphism(1070,500)(0,-1)[(B, A')_1`(\tilde B, A')_1`(u,A')_1]{450}1r
\putmorphism(450,-400)(1,0)[\phantom{(B, \tilde A)}``(\tilde B,K)_2]{500}1a
\put(600,-170){\fbox{$(\theta_0^{\tilde B})_K$}}
\efig}
\quad=\quad
\scalebox{0.86}{
\bfig
 \putmorphism(-150,500)(1,0)[(B,A)_1`(B,A')_1 `(B,K)_1]{600}1a
 \putmorphism(450,500)(1,0)[\phantom{(B,A)}` `=]{540}1a
\putmorphism(-180,500)(0,-1)[\phantom{Y_2}`(B, A)_2 `(\theta^A_0)_B]{450}1l
\put(0,280){\fbox{$(\theta_0^B)_K$}}
\putmorphism(-180,-400)(1,0)[(\tilde B,A)_2` `(\tilde B,K)_2]{500}1a
\putmorphism(-180,50)(0,-1)[\phantom{Y_2}``(u,A)_2]{450}1l
\putmorphism(450,50)(0,-1)[\phantom{Y_2}`(\tilde B, A')_2`(u, A')_2]{450}1r
\putmorphism(450,500)(0,-1)[\phantom{Y_2}`(B, A')_2 `(\theta^{A'}_0)_B]{450}1r
\put(620,50){\fbox{$(\theta_0^{A'})^u$}}
\putmorphism(-180,50)(1,0)[\phantom{(B, \tilde A)}``(B,K)_2]{500}1a
\putmorphism(1130,50)(0,-1)[\phantom{(B, A')}`(\tilde B, A')_2`(\theta_0^{A'})_{\tilde B}]{450}1r
\putmorphism(1130,500)(0,-1)[(B, A')_1`(\tilde B, A')_1`(u,A')_1]{450}1r
\putmorphism(450,-400)(1,0)[\phantom{(B, \tilde A)}``=]{520}1b
\put(0,-170){\fbox{$(u,K)_2$}}
\efig}
$$
for every 1h-cell $K\colon A\to A'$ and 1v-cell $u\colon B\to \tilde B$, 

\medskip
\noindent \axiom{$VLT^q_3$} \vspace{-0,6cm}

$$\scalebox{0.86}{
\bfig
 \putmorphism(-150,500)(1,0)[(B,A)_1`(B,A)_1 `=]{600}1a
 \putmorphism(450,500)(1,0)[\phantom{(B',A)_1}` `(k,A)_1]{450}1a
\putmorphism(-200,500)(0,-1)[\phantom{Y_2}`(B, A)_2 `(\theta_0^B)_A]{450}1l
\put(-20,50){\fbox{$(\theta^B_0)^U$}}
\putmorphism(-170,-400)(1,0)[(B, \tilde A)_2` `=]{500}1a
\putmorphism(-200,50)(0,-1)[\phantom{Y_2}``(B,U)_2]{450}1l
\putmorphism(450,50)(0,-1)[\phantom{Y_2}`(B, \tilde A)_2 `(\theta_0^B)_{\tilde A}]{450}1l 
\putmorphism(450,500)(0,-1)[\phantom{Y_2}`(B, \tilde A)_1 `(B,U)_1]{450}1l
\put(600,260){\fbox{$(k,U)_1$}}
\putmorphism(450,50)(1,0)[\phantom{(B, \tilde A)}``(k,\tilde A)_1]{500}1a %
\putmorphism(1070,50)(0,-1)[\phantom{(B, A')}`(B', \tilde A)_2`(\theta^{\tilde A}_0)_{B'}]{450}1r
\putmorphism(1070,500)(0,-1)[(B, A')_1`(B', \tilde A)_1`(B',U)_1]{450}1r
\putmorphism(450,-400)(1,0)[\phantom{(B, \tilde A)}``(k,\tilde A)_2]{500}1a
\put(600,-170){\fbox{$(\theta_0^{\tilde A})_k$}}
\efig}
\quad=\quad
\scalebox{0.86}{
\bfig
 \putmorphism(-150,500)(1,0)[(B,A)_1`(B',A)_1 `(k,A)_1]{600}1a
 \putmorphism(450,500)(1,0)[\phantom{(B,A)}` `=]{540}1a
\putmorphism(-180,500)(0,-1)[\phantom{Y_2}`(B, A)_2 `(\theta^A_0)_B]{450}1l
\put(0,280){\fbox{$(\theta_0^A)_k$}}
\putmorphism(-180,-400)(1,0)[(\tilde B,A)_2` `(k,\tilde A)_2]{500}1a
\putmorphism(-180,50)(0,-1)[\phantom{Y_2}``(B,U)_2]{450}1l
\putmorphism(450,50)(0,-1)[\phantom{Y_2}`(B', \tilde A)_2`(B', U)_2]{450}1r
\putmorphism(450,500)(0,-1)[\phantom{Y_2}`(B', A)_2 `(\theta^{B'}_0)_A]{450}1r
\put(640,50){\fbox{$(\theta_0^{B'})^U$}}
\putmorphism(-180,50)(1,0)[\phantom{(B, \tilde A)}``(k,A)_2]{500}1a
\putmorphism(1150,50)(0,-1)[\phantom{(B, A')}`(B',\tilde A)_2`(\theta_0^{\tilde A})_{B'}]{450}1r
\putmorphism(1150,500)(0,-1)[(B', A)_1`(B',\tilde A)_1`(B', U)_1]{450}1r
\putmorphism(450,-400)(1,0)[\phantom{(B, \tilde A)}``=]{520}1b
\put(0,-170){\fbox{$(k,U)_2$}}
\efig}
$$
for every 1v-cell $U\colon A\to \tilde A$ and 1h-cell $k\colon B\to B'$, and 

 \pagebreak

\medskip
\noindent \axiom{$VLT^q_4$} \vspace{-0,3cm}
$$\scalebox{0.86}{
\bfig
\putmorphism(-150,500)(1,0)[(B,A)_1`(B',A)_1`(k,A)_1]{600}1a
 \putmorphism(480,500)(1,0)[\phantom{F(A)}`(B',A')_1 `(B',K)_1]{640}1a
 \putmorphism(-150,50)(1,0)[(B,A)_1`(B,A')_1`(B,K)_1]{600}1a
 \putmorphism(470,50)(1,0)[\phantom{F(A)}`(B',A')_1 `(k,A')_1]{660}1a

\putmorphism(-180,500)(0,-1)[\phantom{Y_2}``=]{450}1r
\putmorphism(1100,500)(0,-1)[\phantom{Y_2}``=]{450}1r
\put(320,280){\fbox{$(k,K)_1$}}

\putmorphism(-150,-400)(1,0)[(B,A)_2`(B,A')_2 `(B,K)_2]{640}1a
 \putmorphism(490,-400)(1,0)[\phantom{F(A')}` (B',A')_2 `(k,A')_2]{640}1a

\putmorphism(-180,50)(0,-1)[\phantom{Y_2}``(\theta_0^A)_B]{450}1l %
\putmorphism(450,50)(0,-1)[\phantom{Y_2}``]{450}1l
\putmorphism(610,50)(0,-1)[\phantom{Y_2}``(\theta_0^B)_{A'}]{450}0l 
\putmorphism(1120,50)(0,-1)[\phantom{Y_3}``(\theta_0^{A'})_{B'}]{450}1r
\put(-40,-180){\fbox{$(\theta^B_0)_K$}} 
\put(620,-180){\fbox{$(\theta_0^{A'})_k$}}
\efig}
\quad
=
\quad
\scalebox{0.86}{
\bfig
\putmorphism(-150,500)(1,0)[(B,A)_1`(B',A)_1`(k,A)_1]{600}1a
 \putmorphism(480,500)(1,0)[\phantom{F(A)}`(B',A')_1 `(B',K)_1]{640}1a

 \putmorphism(-150,50)(1,0)[(B,A)_2`(B',A)_2`(k,A)_2]{600}1a
 \putmorphism(450,50)(1,0)[\phantom{F(A)}`(B',K)_2 `(B', K)_2]{640}1a

\putmorphism(-180,500)(0,-1)[\phantom{Y_2}``(\theta_0^A)_B]{450}1l
\putmorphism(450,500)(0,-1)[\phantom{Y_2}``]{450}1r
\putmorphism(300,500)(0,-1)[\phantom{Y_2}``(\theta_0^{B'})_A]{450}0r
\putmorphism(1100,500)(0,-1)[\phantom{Y_2}``(\theta_0^{B'})_{A'}]{450}1r
\put(-40,280){\fbox{$(\theta_0^{A})_k$}}
\put(700,280){\fbox{$(\theta_0^{B'})_K$}}

\putmorphism(-150,-400)(1,0)[(B,A)_2`(B,A')_2 `(B,K)_2]{640}1a
 \putmorphism(490,-400)(1,0)[\phantom{F(A')}` (B',A')_2 `(k,A')_2]{640}1a

\putmorphism(-180,50)(0,-1)[\phantom{Y_2}``=]{450}1l
\putmorphism(1120,50)(0,-1)[\phantom{Y_3}``=]{450}1r
\put(320,-200){\fbox{$(k,K)_2$}}
\efig}
$$
for every 1h-cells $K\colon A\to A'$ and $k\colon B\to B'$. 
\end{defn}

\end{document}